\documentclass[12pt]{amsart}

\usepackage{quiver}

\usepackage{ourstyle}

\usepackage{mathtools}

\newcommand{\ph}{\varphi}
\newcommand{\pf}{\mathfrak{p}}
\DeclareMathOperator{\Coker}{Coker}

\usepackage{hyperref}

\title
{
	$p$-Adic periods and Selmer scheme images
}

\author{David Corwin and Ishai Dan-Cohen}

\thanks{The second author's work on this article was supported by ISF grants 726/17 and 621/21.}

\date{\today}

\begin{document}

	\begin{abstract} The Chabauty--Kim method was developed with the aim of approaching effective Faltings', the problem of explicitly determining the finite set of rational points on a hyperbolic curve. This method has seen success with the more particular Quadratic Chabauty method, but this method still applies only to certain curves. Previous applications of Chabauty--Kim beyond the quadratic level, as pursued by the authors, by S. Wewers, and by others, use mixed Tate motives and the $p$-adic period map of Chatzistamatiou-\"Unver to approach the particular hyperbolic curve $\Poneminusthreepoints$.
   
%    Previous applications of Chabauty--Kim beyond Quadratic Chabauty, as pioneered by Dan-Cohen--Wewers, used mixed-Tate motives of Deligne--Goncharov and the $p$-adic period map of Chatzistamatiou--\"Unver.
		
		The main purpose of this article is to lay foundations for extending the above approach to more general hyperbolic curves, in particular by defining an analogous $p$-adic period map for more general categories of motives and their non-conjectural cousins such as systems of realizations and $p$-adic Galois representations. We use this to describe a general setup for non-abelian Chabauty for an arbitrary hyperbolic curve.
		
		Our period map also connects the study of $p$-adic iterated integrals with Goncharov's theory of \textit{motivic} iterated integrals, and allows us to investigate Goncharov's conjectures from a $p$-adic point of view. In particular, it suggests the possibility of evaluating syntomic regulators by writing elements of $K$-theory in terms of motivic iterated integrals. Lastly, it forms the basis for a certain generalization of the $p$-adic period conjecture of Yamashita for mixed Tate motives well-suited to applications in Chabauty--Kim theory.
	\end{abstract}
	
	\maketitle

	\setcounter{tocdepth}{2}
	\tableofcontents
	
	\raggedbottom
	\SelectTips{cm}{11}

	\setcounter{tocdepth}{1}
	%\tableofcontents

	%%%%%%%%%%%%%%
	\section{Introduction}
	%%%%%%%%%%%%%
	
	\subsection{$p$-adic periods and multiple polylogarithms}\label{sec:intro_intro}
	
	The $p$-adic period map of \cite{ChatUnv13} from the theory of mixed Tate motives (\cite{DelGon05}) plays a central role in drawing a fruitful connection between two areas of investigation. One is the study of the structure of mixed Tate motivic Galois groups as pursued by Goncharov in his ``beyond standard'' motivic conjectures (\cite{GonGal}, \cite{GonMot}, \cite{GonICM}, \cite{GonMPMTM01}, etc.) as well as by Brown \cite{BrownMTMZ,BrownDecomp} and others. A one-sentence summary of this multi-faceted subject might go as follows: one defines certain special functions $f \in \Oo(G)$ on various motivic Galois groups $G$ (\textit{motivic multiple polylogarithms}, \textit{unipotent multiple polylogarithms}) and attempts to understand which portions of $\Oo(G)$ are spanned by which families of special functions. 
	
	The other is the study of $p$-adic multiple polylogarithms and $p$-adic multiple zeta values (see, for instance, \cite{FurushoI, FurushoII, BesserFurusho, FurushoJafari, FurushoPenagonHexagon, UnverCyclotomic, JarossayDynamical, FresanGil}), where one of the main interests is in understanding relations between rational values of multiple polylogarithms with a view towards Grothendieck's philosophy, according to which all relations should \emph{come from geometry}. Indeed, Yamashita has conjectured a $p$-adic mixed-Tate analog of the Grothendieck period conjecture (\cite[Conjecture 4.1]{AndreGaloisMotTrans09}), as follows. There's a ring homomorphism
	\[
	\m{per}_p : \Oo(G) 
	\twoheadrightarrow \Oo(U) \to \Qp,
	\]
	dual to the period point $\eta$ of \cite{ChatUnv13} defined on the quotient Hopf algebra associated to the prounipotent radical $U$ of $G$, which sends motivic multiple polylogarithms to $p$-adic multiple polylogarithms. Yamashita conjectures (\cite[Conjecture 4]{YamashitaBounds10}\footnote{Technically, Yamashita uses the period point $\varphi_{\pf} = F_{\pf}^{-1} \tau(q)^{-1}$ due to \cite[5.28]{DelGon05}, which produces Deligne $p$-adic polylogarithms (\cite[19.6]{Deligne89}, c.f. \cite[\S 3.2]{YamashitaBounds10}) rather than those of Furusho (\cite{FurushoI}), which are output by $\eta$ and therefore also our period point. However, \cite[Propositions 3.10, 3.12]{YamashitaBounds10} imply that these two versions of the conjecture are equivalent, at least for $N \le 4$.}) that this is injective for mixed Tate motives over $\Qb(\mu_N)$.
	
	\subsection{Mixed Tate Chabauty--Kim}\label{sec:CK_intro}
	
	These two areas meet with particular force in the approach to Chabauty--Kim theory in the mixed Tate setting developed by the authors and others. Chabauty--Kim theory (described in detail in \ref{sec:CK_setup_intro}), also known as non-abelian Chabauty's Method, is a programme started by Minhyong Kim (\cite{kim09}) to $p$-adically bound the number of rational points on an arbitrary hyperbolic curve $X$ over a number field.
	
	Kim's programme has so far found significant applications via the Quadratic Chabauty method of \cite{BalakrishnanDograI,BalakrishnanDograII}, in particular in determining the set of rational points on a previously-inaccessible curve (\cite{CursedCurve}). Nonetheless, Quadratic Chabauty, like the original Chabauty--Coleman method that came before it, is limited to hyperbolic curves whose arithmetic and geometric data satisfy certain inequalities. Faltings' Theorem, on the other hand, applies to all hyperbolic curves. One would therefore like to compute the Chabauty--Kim method for all hyperbolic curves, which conjecturally implies an effective approach to Faltings' Theorem by \cite[Conjecture 3.1]{nabsd} (c.f. also \cite{KimEffective}).

    Work by Dan-Cohen and Wewers in the setting of mixed Tate motives \cite{MTMUE,MTMUEII}, while limited in scope, was particularly successful at going beyond the quadratic level. Subsequent work in this direction includes \cite{PolGonI,PolGonII, CKTwo, BrownUnit, BettsEtAl,JarLilSaeWeiZeh24}. The goal of the present article is to lay the foundations for applying this approach to general hyperbolic curves, which necessarily goes beyond mixed Tate motives. We first review previous work in the latter setting.
	
	\subsubsection{}\label{sec:MT_CK}
	
	Let $X$ be the projective line with at least 3 punctures, which are required to be distinct over every closed point of an open subscheme $Z$ of $\Spec \ZZ$. Such a curve can be viewed as a \emph{non-proper hyperbolic} smooth curve. Siegel's Theorem says, in a parallel way to Faltings', that $X(Z)$ is finite, and this finiteness was the original testing ground for Kim's method (\cite{kim05}).
	
	The work of Deligne-Goncharov (\cite{DelGon05}) provides a lifting of the rational prounipotent completion of the fundamental group of $X(\CC)$ (at an appropriate base point) to a prounipotent group object $\pi_1^{\MT}(X)$ of the Tannakian category $\MT(Z,\Qb)$ of mixed Tate motives over $Z$. We consider finite type quotients $\ppi$ of $\pi_1^{\MT}(X)$. One may regard such a $\ppi$ as a unipotent $\QQ$-group $\ppi_{\dR}$ equipped with an action of the mixed Tate motivic Galois group $G=G^{\MT}(Z)$ of $Z$ at the de Rham fiber functor. Letting $U=U^{\MT}(Z)$ denote the unipotent radical of $G(Z)$, we have $G/U \simeq \Gm$. Thus there's a $\Gm$-action on the moduli space of cocycles $Z^1(U; \ppi_{\dR})$, and it is proved in \cite[Proposition 5.2]{MTMUE} that $H^1(G;\ppi_{\dR}) \simeq Z^1(U;\ppi_{\dR})^{\Gm}$.
	
	By work of \cite{ChatUnv13}, we have for $p \in Z$ a special point $\uU_p \in U(\Qp)$ corresponding to the homomorphism $\m{per}_p$ of \S \ref{sec:intro_intro}. The properties of this map show (\cite[\S 4.11]{MTMUE},\cite[\S 2.4]{PolGonI}) that the map from the global to local Selmer variety in Kim's theory described in \S \ref{sec:CK_setup_intro} may be identified with evaluation at $\uU_p$:
	\[
	Z^1(U; \ppi_{\dR})^\Gm_\Qp \xto{\ev_{\uU_p}}
	(\ppi_{\dR})_\Qp.
	\]
	
	In our previous approach to the mixed Tate case, the problem of bounding the image of $X(Z)$ in $X(\Zp)$ (for some $p \in Z$) is thus reduced to the problem of understanding the (scheme-theoretic) image of $\ev_{\uU_p}$.
	%\commentDC{DC: The trivial $U$-action and the freeness of $U$ are NOT the reason why we are led to consider the image of the universal evaluation map. This universal evaluation map is relevant even when neither of these is true, and the real reason it's relevant is because of the relationship between the $p$-adic period map and Kim's local unipotent Albanese map. ID: If the action is nontrivial, then the "geometric step" involves arithmetic work.}
	Moreover, the abstract structure of $U$ as a prounipotent group with $\Gm$-action is known (it's free on the rational algebraic $K$-theory of $Z$), and when $\ppi_{\dR}$ has trivial $U$-action, $Z^1(U; \ppi)^\Gm_\Qp$ is readily computable.\footnote{Work-in-progress of the two authors and M. L\"udtke computes even when $\ppi_{\dR}$ has nontrivial $U$-action.} But $\ev_{\uU_p}$ is just the pullback along $\uU_p \colon \Spec{\Qp} \to U$ of the universal cocycle evaluation map (\cite[Definition 2.20]{PolGonI}):
	\[
	Z^1(U; \ppi_{\dR})^\Gm \times U
	\xto{\ev_{\ppi}}
	\ppi_{\dR} \times U.
	\]
	
	The period conjecture (\cite[Conjecture 4]{YamashitaBounds10}) says that $\uU_p$ is generic, which implies that the image of the former is the pullback of the image of the latter. In particular, letting $I^\m{mot}$ denote the ideal vanishing on the image of $\ev_{\ppi}$ and $I$ the ideal vanishing on the image of $\ev_{\uU_p}$, we have
	\[
	I^\m{mot} \otimes_{\Oc(U)} \Qp \subseteq I;
	\]
	if the period conjecture holds, then they are equal.
	
	The problem of computing the image of $\ev_{\ppi}$ in terms of free generators for (suitable portions of) $U$ is purely geometric and is known as the \emph{geometric step}. The more subtle part is computing coordinates on $U$ on which $\per_p$ may be computed, which ultimately involves some intricate manipulations with motivic iterated integrals, such as in \cite[\S 4.3]{PolGonI}.
	
	Our goal here is to extend the above picture beyond the limited setting of mixed Tate motives, so our initial focus is on aspects that are largely formal.\footnote{As the reader will see, non-formal Tannakian considerations do arise.}

	\subsection{Beyond the mixed Tate setting}
	
	For a non-rational hyperbolic curve, one must pass beyond the (Artin-)Tate setting. Here, abelian categories of motives are only conjectural (see \S \ref{sec:mixed_motives} and \ref{app:mixed_motives}). We are therefore forced to work with systems of realizations if we want our setup to be non-conjectural. Given the central role of $p$-adic Hodge theory in $p$-adic realization, we prefer to use \emph{motivic structures} of Fontaine--Perrin-Riou (\cite[\S III]{FPR91}) in \S \ref{sec:PSPM}-\ref{sec:motivic_structures}, which include the $p$-adic Hodge comparison as part of the data. Nonetheless, our formal results (\S \ref{sec:arith_hodge_path}-\ref{sec:geom_loc_real}) apply equally to other substitutes for the abelian category of mixed motives, such as appropriate categories of Galois representations.
	
	%\comment{When you say ``for concreteness'', are you saying you're just considering the $GL_2$ case instead of a general reductive group? Maybe at least make that clear that you're considering a special case - though I see no reason not to mention $\Gb$, and also I think it would be good to explain how it connects to $X$, since that's the ultimate application. But also, we are considering a general reductive group in this paper, right? My writeup on my website is just about extending from $GL_2$ to the general case, and I thought you felt this comes before that? ID: Yes, we're considering the $GL_2$ case here only for concreteness. I've added the work ``temporarily."}
	
	\subsubsection{}
	For concreteness, let us temporarily focus attention on a special case, namely the category $\mathbf{SM}_{Z}(\Qb,E)$ (\ref{sssec:SM_A,S},\ref{category_of_lisse_structures}) of \textit{lisse mixed-$E$-motivic structures over $Z$} associated to a non-CM elliptic curve $E/\basefield$ with good reduction over an open subscheme $Z \subseteq \Spec{\Oc_{\basefield}}$. This applies whenever $X$ is of \emph{mixed-elliptic type}, for example when $X$ is an affine open subscheme of an elliptic curve or a genus $2$ curve whose Jacobian is isogeneous to the square of an elliptic curve, such as the curves with LMFDB labels 38416.a.614656.1 and 614656.a.614656.1.
	
	This is, roughly speaking, the Serre subcategory category of \textit{lisse motivic structures over $Z$} (\ref{category_of_structures}) generated by $h_1(E)$. Its Tannakian Galois group at the fiber functor induced by de Rham realization sits in an exact sequence 
	\[
	1 \to U(Z,E) \to G(Z,E) \to \GL_2 \to 1. 
	\]
	The group $U=U(Z,E)$ consists of those Tannakian loops which induce the identity on associated graded for the weight filtration. Let $F^0U \subset U$ (\S \ref{WH_groups_diagram}) be the subgroup of automorphisms that preserve the Hodge filtration. For $\pf \in Z$, our main results (\S \ref{sec:hodge_frob_together}) lead to the construction of a special $\Qp$-point $\uU_\pf^R$ (resp. $\uU_\pf^L$) (\S \ref{sec:period_loops_mot_gal}) of $U$ (a ``unipotent $p$-adic period loop") which is well-suited to our applications to Chabauty--Kim theory. This point depends on a certain choice, and is well-defined only modulo the left- (resp. right-) action of $F^0U$.
	
	Let us focus on $\uU_{\pf}^L$. The homomorphism
	\[
	\m{per}^L_\pf  \colon \Oo \big(
	 U / F^0 U\big) \to \basefield_{\pf}
	\]
	associated to $\uU_\pf^L$ should probably not be regarded as \textit{the} $p$-adic period map.\footnote{Instead, we suggest the terminology \textit{unipotent $p$-adic period map}.} In one of several drawbacks, by its very construction, $\m{per}_\pf^L$ provides pure motives with no nonzero $p$-adic periods. While in the mixed Tate setting, the only period of a pure motive ($2 \pi i$) may reasonably be considered to be ``$p$-adically zero'', this appears not to be the case for more general motives. A construction of a candidate for the $p$-adic period map which is better in this regard (although perhaps less suited to Chabauty--Kim) has been carried out concurrently with this work by Ancona--Fr\u{a}\c{t}il\u{a} \cite{AnFr2025algebraic}. See also Dan-Cohen \cite{dancohen2025andreperiodsmixedtate}.
	
	%\commentDC{DC: Don't we want to discuss the general case of a hyperbolic curve over an integer ring? I thought this paper is about the general case. In fact, the only difference between the paper on my website that you had objected to and what's on the arxiv is that the former discusses the general case, so then I'm confused why the former paper caused concern if we're not focusing on the general case here.}
	
	\subsubsection{}
	
	More generally, in order to apply Chabauty--Kim to an arbitrary curve $X$, we must consider a category
	\[
	\mathbf{SM}_{Z}(\Qb,J)
	\]
	of \textit{lisse mixed-$J$-motivic structures over $Z$}, where $J$ is the Jacobian of $X$, which is the same as above but generated by $h_1(J)$. We describe this category in detail in \S \ref{sec:motivic_structures}-\ref{sec:SM_A}. We have a de Rham fiber functor over $\basefield$ giving rise to a motivic Galois group $G(Z,J)$ (\ref{notation:GZ_UZ}) over $\basefield$, sitting in a similar exact sequence
	\[
	1 \to U(Z,J) \to G(Z,J) \to \Gb(J) \to 1,
	\]
	where $\Gb(J)$ is a reductive group closely related to both the Mumford--Tate group and the $p$-adic monodromy group of $J$, and $U=U(Z,J)$ is generated (conjecturally freely) by certain $\Ext^1$ groups in $\mathbf{SM}_{Z}(\Qb,J)$, which are in turn conjecturally isomorphic to pieces of certain rational algebraic $K$-theory or higher Chow groups of powers of $E$.
	
	Similarly to the case of $E$, for $\pf \in Z$ above $p$, we have a homomorphism
	\[
	\m{per}_\pf^L  \colon \Oo \big(
	U(Z,J)/F^0U(Z,J) \big) \to \basefield_{\pf}.
	\]
	
	\begin{ssRemark}\label{remark:computing_period_map}
	The reader may be left wondering how to compute the $p$-adic period maps on explicit elements of $\Oc(U)$. For this, we note that the commutativity of \ref{ssec:conclusion}(\ref{eqn:CK_diagram2}) below allows us to compute the maps on `motivic iterated integrals'. Specifically, given $\om \in \Oc(\pi_{\dR})$ and $z \in \Xc(Z)$, we may define a `motivic iterated integral'
	\[
	I^{\mot}(\om;b,z) \coloneqq \kappa_{\ppi}(z)^{\#}(\omega) \in \Oc(U),
	\]
	depending on a choice of $\tannloop^\m{H}_{\pf}$. If $\omega$ is left-invariant under $F^0$, then
	\[
	\per_{\pf}^L\left(I^{\mot}(\om;b,z)\right) = \int_b^z \omega
	\]
	in the sense of Coleman integration.

    \end{ssRemark}

	\subsection{Chabauty--Kim setup}\label{sec:CK_setup_intro} Let $\basefield$ be a number field, $Z=  \Spec{\Oc_{\basefield,S}}$ for a finite subset $S \subseteq S_f(\basefield)$, and $\Xc \to Z$ a smooth morphism whose generic fiber $X$ is a hyperbolic curve (i.e., with non-abelian geometric fundamental group) along with smooth proper compactification $\Xc \hookrightarrow \Xc' \to Z$ such that the complement $\Yc = \Xc' \setminus \Xc$ is \'etale over $Z$. We denote by $Y,X',X$ their generic fibers, respectively.
	
	The theorems of Faltings and Siegel show that $\Xc(Z)$ is always finite. Note that if $X=X'$, then $\Xc(Z)=X(\basefield)$, and this is the case of Faltings. The open problem is to compute $\Xc(Z)$ given $\Xc$.
	
	Fix a place $\pf \in S_p(\basefield) \cap Z$ and a finite-dimensional Galois-equivariant quotient $\ppi_p$ of $\pi_1^{\et,\un}(\overline{X})_{\Qp}$, the $\Qp$-prounipotent completion of the profinite group $\pi_1^{\et}(\overline{X})$. Kim's approach to this problem (\cite{kim09},\cite{nabsd}) uses a diagram
	\[
	\xymatrix{
		X(Z) \ar[r] \ar[d]^{\kappa_{\ppi_p}} & X(Z_{\pf}) \ar[d]_{\kappa_{\pf,\ppi_p}} \ar[dr]^{\int_{\pf}} \\
		H^1_{f,S}(\basefield;\ppi_p)_{\Qp} \ar[r]_-{\mathrm{loc}_{\ppi_p}} & H^1_f(\basefield_{\pf};\ppi_p) \ar[r]_-{\log_{\mathrm{BK}}}^-{\sim} & \operatorname{Res}_{\basefield_{\pf}/\Qp} (F^0 \ppi_{\dR} \backslash \ppi_{\dR})
	}
	\tag{Kim's Cutter}
	\]
	known as Kim's cutter, where $\ppi_{\dR} \coloneqq \operatorname{D}_{\dR}(\ppi_p)$.
	
	We refer to the composition
	\[
	\rm{LocReal} = \rm{LocReal}_{\ppi_p} \coloneqq \log_{\mathrm{BK}} \circ \mathrm{loc}_{\ppi_p}\]
	as the \emph{Localization-Realization} map of Chabauty--Kim theory. Computing the \emph{image} of this map is the most important and nontrivial step in using this theory to bound $\Xc(Z)$. \cite[\S 2]{KimEffective} refers to the computability of this map as ``Hypothesis [H]'' and shows in the case $\basefield=\Qb$ that assuming [H] and the Bloch--Kato conjectures, one may effectively compute a function on $\Xc(\Oc_{\pf})$ vanishing on $\Xc(Z)$.
	
	When $X'$ has genus $0$, the approach mentioned in \S \ref{sec:MT_CK} replaces the $\Qp$-scheme $H^1_{f,S}(\basefield;\pi_1^{\et,\un}(\overline{X})_{\Qp})$ by a $\Qb$-form $H^1(\mathbf{MT}(Z,\Qb);\pi_1^{\MT}(X))$, where $\mathbf{MT}(Z,\Qb)$ is the category of mixed Tate motives of \cite[1.6]{DelGon05}, and $\pi_1^{\MT}(X)$ is the fundamental group in this category due to \cite[5.2]{DelGon05}. The Localization-Realization map is then quite simply, in the notation of \S \ref{sec:MT_CK}, the evaluation map
	\[
	Z^1(U; \ppi_{\dR})^\Gm_\Qp \xto{\ev_{\uU_p}}
	(\ppi_{\dR})_\Qp.
	\]
	
	In \S \ref{sec:path_torsor_\basefield,S}, we define a similar group scheme object $\pi_1^{\mathbf{SM}}(X)$ in $\mathbf{SM}_{\basefield,S}(\Qb)$ for an \emph{arbitrary} curve $X$ and in \ref{sec:CK_setup}(\ref{eqn:CK_diagram}) define an alternative version of Kim's Cutter relative to a finite-dimensional motivic quotient $\ppi$ of $\pi_1^{\mathbf{SM}}(X)$, namely:
	\begin{equation}\label{eqn:CK_diagram_intro}
		\xymatrix{
			\Xc(Z) \ar[r] \ar[d]^{\kappa_{\ppi}} & \Xc(\Oc_{\pf}) \ar[d]_{\kappa_{\pf,\ppi}} \ar[dr]^{\int_{\pf}} \\
			H^1(\mathbf{SM}_{\basefield,S}(\Qb);\ppi) \ar[r]_-{\mathrm{loc}_{\ppi}} & H^1_f(\basefield_{\pf};\ppi_p) \ar[r]^-{\sim}_-{\log_{\mathrm{BK}}} & \operatorname{Res}_{\basefield_{\pf}/\Qp} F^0 (\ppi_{\dR})_{\basefield_{\pf}} \backslash (\ppi_{\dR})_{\basefield_{\pf}}.
		}
	\end{equation}

	In \S \ref{sec:abstract_CK_maps} (resp. \S \ref{left_variant}), we define a more abstract version of $\rm{LocReal}$ called `$\beta^R$' (resp. `$\beta^L$') for any effective group object in a ``weight-filtered Tannakian category'' (Definition \ref{defn:weight_filtr_cat}) with Hodge-filtered Frobenius-equivariant fiber functor (Definition \ref{defn:hodge_frob}). The description of $\log_{\mathrm{BK}}$ in \cite[3.4.3]{nabsd} (c.f. \ref{sec:BK_log}) implies that this is the same as $\rm{LocReal}$ of the previous paragraph.
	
	We thus set out to study the (scheme-theoretic) image $\opnm{\overline{Im}} \rm{LocReal}$ of the associated map of schemes over $\Qp$
	\[
	\m{LocReal}  \colon
	H^1(\mathbf{SM}_{\basefield,S}(\Qb);\ppi) \times_{\Qb} \Spec \Qp
	\to
	\operatorname{Res}_{\basefield_{\pf}/\Qp} F^0 (\ppi_{\dR})_{\basefield_{\pf}} \backslash (\ppi_{\dR})_{\basefield_{\pf}}
	\]
	(\ref{left_variant}).
	
	Our main result is as follows:
	
	\begin{thm*}\label{thm:intro_locreal}
		The (Weil restriction of the) localization-realization map of Chabauty--Kim theory may be identified with the evaluation map:
		\[
		Z^1(U(Z,J); \ppi_{\dR})^{\Gb(J)}_\Qp
		\xto{\ev_{\uU_\pf}}
		\operatorname{Res}_{\basefield_{\pf}/\Qp} F^0 (\ppi_{\dR})_{\basefield_{\pf}} \backslash (\ppi_{\dR})_{\basefield_{\pf}}.
		\]
	\end{thm*}
	\noindent

	See the Theorem of \S \ref{thm:final_theorem} for a precise statement, as well as Theorems \ref{thm:right_commutes}, \ref{thm:left_commutes} for the abstract formulation.
	
	\subsection{Remarks on Chabauty--Kim computations}\label{sec:remarks_on_computation}
	
	In order to compute the image of the localization-realization map, one needs to know the structure of $U(Z,J)$. This group is generated (conjecturally freely) by the $\Ext$ groups $\Ext^1_{\mathbf{SM}_{\basefield,S}(\Qb)}(\Qb(0),V)$, as $V$ ranges over irreducible representations of $\Gb(J)$ (equivalently, irreducible objects of $\mathbf{SM}_{Z}(\Qb,J)$).
	
	Such groups receive maps from certain rationalized algebraic $K$-theory groups (\S \ref{sec:sm_mot_coh}-(\ref{eqn:k_ch_mc}),(\ref{eqn:map_from_k_theory})), which are generally believed to be isomorphisms. The Beilinson conjectures imply that the dimensions of most of these groups are equal to dimensions of certain absolute Hodge cohomology groups, whose dimensions may be computed in terms of the Hodge structure on $V_{\dR}$.
	
	But there is also a natural map
	\[
	\Ext^1_{\mathbf{SM}_{\basefield,S}(\Qb)}(\Qb(0),V) \otimes_{\Qb} \Qp
	\to
	H^1_{f,S}(G_{\basefield};V_p),
	\]
	which \cite[III.4.1.5]{FPR91} conjectures to be an isomorphism (c.f. Conjecture \ref{conj:FPR_adm}). Arithmetic duality theorems allow one to compute their dimensions conditionally on the vanishing of Bloch--Kato Selmer groups in non-negative weight, as carried out in \cite[\S 5]{CorwinMECK} and described in more generality in \cite[\S 2.7]{BettsCorwinLeonhardt}.
	
	In fact, one may circumvent $\Ext^1_{\mathbf{SM}_{\basefield,S}(\Qb)}(\Qb(0),V)$ entirely and work directly with the groups $H^1_{f,S}(G_{\basefield};V_p)$ by using $\Qp$-Tannakian categories of $p$-adic Galois representations (\ref{sec:sg_galois_reps}). Our axiomatic setup applies equally in this setting. This has the disadvantage that one loses the potentially significant Zariski closure of the period loop (Conjecture \ref{conj:period_conjecture}). On the other hand, it has the advantage that one may sometimes prove the desired vanishing of Bloch--Kato Selmer groups using Iwasawa theory. This direction is pursued in the parallel work \cite{CorwinTSV,CorwinMECK} by the first author.
	
	In particular, \cite[\S 2.4]{CorwinMECK}, explains how one may use Kato's Iwasawa Main Conjecture to prove $H^1_f(G_{\Qb};V_p(E)(1))=0$ for a non-CM elliptic curve $E$ with ordinary reduction at $p$. This in turn allows one to unconditionally apply the Chabauty--Kim method to the case of $\Zb[1/\ell]$-points on a punctured elliptic curve of rank $1$, in which Quadratic Chabauty does not apply. M. L\"udtke and the first author have already computed Chabauty--Kim functions in this case (\cite{CLinitialreport2026}). Note that all of that work depends crucially on the main theorem of the present article, which already demonstrates its significance.
	
	%The data of a mixed Tate filtered $\phi$-module over $\Qp$ is equivalent to the data of a ``weight'' filtered vector space together with two splittings. The underlying vector space has a unique automorphism which is unipotent with respect to the weight filtration and interchanges the two splittings. Varying the filtered $\phi$ module, these automorphisms assemble to a $p$-adic point of a Tannakian Galois group which maps to the Tannakian Galois group of the category of mixed Tate motives over any open integer scheme $Z$ containing $p$. The associated ring homomorphism\[\m{per}_p \colon A(Z) \to \Qp \] is the $p$-adic period map. 

	\subsection{Proof sketch}
	
	While \cite{ChatUnv13} works only with filtered $\varphi$-modules, we work with global objects (Galois representations, motives, or systems of realizations), for the following reason.
	%% DC082: I added the above sentence to introduce the paragraph 
	Even if $Y \to \Spec \Zp$ is a smooth and proper morphism, the admissible filtered $\phi$-modules $H^i_{\rig}(Y_\Qp)$ may fail to be semisimple. Thus, beyond the mixed Tate setting, we don't have a useful weight-filtered category (see the next paragraph) of admissible filtered $\phi$-modules. For this reason, constructions that center on the category of mixed Tate filtered $\phi$-modules don't generalize readily. Instead, we construct our period loop $\uU_p$ while working directly with an appropriate global Tannakian category (systems of realizations, $p$-adic Galois representations, or motives, as described in \S \ref{sec:examples_cats}).
	
	For example, suppose we work with a subcategory $\Tc$ (such as $\mathbf{SM}_{K,S}(\Qb,J)$ in \ref{sssec:SM_A,S}) of the category of \textit{motivic structures} (\S \ref{sec:motivic_structures}). In our axiomatic setup, an important feature of $\Tc$ is that it is \emph{weight-filtered} (\ref{defn:weight_filtr_cat}), which implies that the induced action on associated graded is always semisimple. De Rham realization
	% For a discussion of sentences beginning with a name that begins with a lower-case letter, see https://blog.apastyle.org/apastyle/2013/11/is-it-sometimes-okay-to-begin-a-sentence-with-a-lowercase-letter.html
	$\om^{\dR}$ provides a $\basefield$-rational fiber functor which is \emph{Hodge filtered} (\ref{defn:hodge_filtr}); moreover, it is isomorphic over $\basefield_{\pf}$ to a fiber functor $\om^{\cri}_0$ over $\basefield_{\pf,0}$ that is \emph{Frobenius-equivariant} (\ref{defn:frob_equiv}), assuming $\pf \in Z$. In this setting, our main theorems are as follows:
	
	\begin{Theorem*}[\ref{thm:hodge-filtered_weight-splitting} below]
		Tannakian paths $\tannloop^\m{H} \colon (\om^{\dR} \circ \gr^W) \to \om^{\dR}$ that preserve Hodge filtrations form a trivial (in particular, nonempty!) torsor under $F^0U$.
	\end{Theorem*}
	
	\noindent
	We refer to such paths as \emph{arithmetic Hodge paths}.
	
	\begin{Theorem*}[\ref{thm:arith_frob_path} below]
		There exists a unique Tannakian path 
		\[
		\tannloop^\m{cr}_0 \colon (\om^{\cri}_0 \circ \gr^W) \to \om^{\cri}_0
		\]
		which respects Frobenii.
	\end{Theorem*}
	
	\noindent
	We refer to the base-change $\tannloop^\m{cr}$ of $\tannloop^\m{cr}_0$ to $\basefield_{\pf}$ as the \emph{arithmetic crystalline path}. We define our (left) period loop by (\S \ref{sec:hodge_frob_together})
	\[
	\uU_p^L \coloneqq \tannloop^\m{cr} \circ (\tannloop^\m{H})\inv
	\]
	for any arithmetic Hodge path $\tannloop^\m{H}$. Our proof of Theorem \ref{thm:hodge-filtered_weight-splitting} builds on P. Ziegler's study of filtered fiber functors in \cite{ZieglerGrFil15},\footnote{Similar results have also appeared in work of Hain-Matsumoto (\cite[Lemma 3.1]{HainMatsumoto}) and Hain (\cite{Hain_HodgeDeRham_Chapter2016,HainGoldman20}). See the paragraph following Theorem \ref{thm:hodge-filtered_weight-splitting}.}
while our proof of Theorem \ref{thm:arith_frob_path} is based on Besser's proof in \cite{BesserColeman} of the uniqueness and existence of \textit{geometric} crystalline paths.

    In \S \ref{sec:geom_loc_real}, we prove the abstract form of Theorem \ref{thm:intro_locreal}. One subtle point in defining the diagram in this abstract setting is to define a Hodge filtration on torsors over an arbitrary $\basefield_{\pf}$-algebra $R$. For this, we prove a form of Tannaka duality for the base-change of a Tannakian category to an arbitrary base ring (Theorem \ref{theorem:tannakian_reconstr_groupoids}), which we then use in \S \ref{base_change_fil_gr} to base-change an arbitrary filtered fiber functor, generalizing a result of \cite{ZieglerFiltGeneralBase}.

    Most of \S \ref{sec:p-adic_periods_for_motivics_galoisr} spells out the concrete consequences of the abstract results of \S \ref{sec:geom_loc_real} in arithmetic settings. Then \S \ref{sec:ck_application} sets up the non-abstract Chabauty--Kim diagram. A key technical point is the construction in \S \ref{sec:fund_grp_and_path_torsors}--\ref{sec:tangential_basepoints} of the motivic fundamental group and path torsors in the category of motivic structures defined in \S \ref{sec:motivic_structures} for a curve over a number field.
	
	%% ####

	%\commentDC{DC: If you want, I'd be happy to write a version of what I'd like as an introduction...but at least for now, I left everything as comments as you asked. ID: Thanks.}

	\subsection{Acknowledgements}
	We wish to thank Giuseppe Ancona, Richard Hain, Owen Patashnick, Paul Ziegler, Joseph Ayoub, and Cl\'ement Dupont for helpful conversations and email exchanges. We thank Amnon Besser, Minhyong Kim, Martin L\"udtke, Benjamin Moore for their support and encouragement.
	
	The first author thanks Ben Moonen for mentioning Jannsen's mixed motives for absolute Hodge cycles.
	
	\subsection{Notations and Conventions}\label{ssec:notations_and_conventions}
	
	When we refer to the \emph{intersection} of two full subcategories of a larger category, we mean the full subcategory on the class of objects that belong to each subcategory. For a union of a filtered collection of subcategories, we similarly mean the full subcategory on all objects contained in at least one of these subcategories. In a given monoidal category $\Cc$, $\mathbbm{1}=\mathbbm{1}_{\Cc}$ denotes a unit object. If $\Cc$ is a symmetric monoidal category, we let $\Cc_{\rig}$ denote the full subcategory of dualizable objects of $\Cc$. For an object of a category consisting of modules with extra data (filtration, grading, Frobenius, etc.), we let $\forg$ denote the functor to the underlying category of modules (usually vector spaces). All functors between additive categories are assumed additive unless otherwise stated. As explained in \ref{sec:tann_prelim}, the term `tensor category' is always specified as `$X$ tensor category' for $\emptyset \neq X \subseteq \{\mathrm{additive}, \mathrm{abelian}, k-\mathrm{linear}\}$; a `tensor functor' refers to an additive (and linear if relevant) symmetric monoidal functor.
    
    %\comment{Let's either add more here, or move this to where it's needed for the first time. ID: Now it's ok.}
	
	We use the concept of Tannakian category throughout, which is reviewed in Appendix \ref{appendix:tannakian_categories}.
	
	An \emph{algebraic group} refers to a finite-type group scheme over a field, while a \emph{pro-algebraic group} refers to a group scheme that is an inverse limit of algebraic groups.
	
	%We freely use the concepts of Tannakian category and fiber functor over a general base scheme. We refer to, rings, schemes, group( scheme)s, and torsors in a Tannakian category in the sense of \cite[5.3-5.4]{Deligne89}.
	
	For a field extension $\coeff'/\coeff$, we let $\Ex_{\coeff}^{\coeff'}$ the extension of scalars functor from $\coeff$ to $\coeff'$ (for vector spaces, algebras, schemes, Tannakian categories, etc.).
	
	If $V$ is a vector space, we denote the associated free pronilpotent Lie algebra by $\nN(V)$, the associated prounipotent group by $U(V)$, the completed universal enveloping algebra by $\Uu(V)$ and the coordinate ring by $A(V)$. If $V$ is a representation of an algebraic group $\GG$, then there's an induced action on all of the above.
	
	For a field $\basefield$, we let $G_{\basefield}$ denote the absolute Galois group of $\basefield$ (for some implicit choice of separable closure). If $\basefield$ is a number field, we let $S_{\infty}(\basefield)$ denote the set of infinite places of $\basefield$,  $S_f(\basefield)$ the set of finite places of $\basefield$, and $S(\basefield) \coloneqq S_{\infty}(\basefield) \cup S_f(\basefield)$. For a subset $S \subseteq S(\basefield)$, let $\Oc_{\basefield,S}$ denote the ring of $S$-integers in $\basefield$, i.e. the set of elements of the form $\frac{\alpha}{\beta}$ for $\alpha,\beta \in \Oc_{\basefield}$ and $(\beta)$ supported in $S$, and $G_{\basefield,S}$ the maximal quotient of $G_{\basefield}$ unramified outside $S$ (equivalently, $\pi_1^{\et}(\Oc_{\basefield,S})$). For a prime number $\ell$, we let $S_\ell(\basefield)$ denote the subset of $S_f(\basefield)$ above $\ell$. For $v \in S(\basefield)$, we let $\basefield_v$ denote the completion of $\basefield$ at $v$, $G_v$ the absolute Galois group of $\basefield_v$ (note that the embedding $G_v \hookrightarrow G_{\basefield}$ is defined only up to conjugation by $G_{\basefield}$), and (if $v \in S_f(\basefield)$) $\Oc_{K,v}$ the Zariski local ring, $\Ov$ the (completed) integer ring, and $\mf_v$ its maximal ideal. If $\pf \in S_p(\basefield)$, we let $K_{\pf,0}$ denote the maximal unramified extension of $\Qp$ in $K_{\pf}$. For a scheme $X$ over $\basefield$, we let $X_v \coloneqq \Ex_{\basefield}^{\basefield_v} X$ and $\overline{X} \coloneqq \Ex_{K}^{\overline{\basefield}} X$.

	%%%%%%%%%%%%%%%%%%%%%%%%%%%%%%%%%%%%%%
	\section{Weight-Filtered Tannakian Categories}\label{sec:wt_fil_tann}
	%%%%%%%%%%%%%%%%%%%%%%%%%%%%%%%%%%%%%%%%
	
	%\commentDC{ID: I've recently had trouble with a journal whose formatting made untitled subsections and subsubsections look strange. I've created commands `spar' (alias `paragraph') and `sspar' to generate numbered paragraphs (`spar' for 2-digit numbers, `sspar' for three digit numbers) without sectioning. DC: Okay - is there something I should do about that? ID: No.}
	
	\subsection{Filtered objects in Abelian categories}
	
	We recall some basic notions about filtered objects in additive and (pre-)abelian categories. Some of this can be found (in the abelian case) in \cite[\href{https://stacks.math.columbia.edu/tag/0120}{Tag 0120}]{stacks-project}.

	\subsubsection{Notation} If $\{M_i\}_I$ is a collection of subobjects of an object $M$ of a pre-abelian\footnote{is additive with all finite (co)limits, hence has (co)kernels, hence (co)images} category $\Ac$ for which $\oplus_{i \in I} M_i$ is defined (e.g., if $I$ is finite or $\Ac$ has all coproducts of size $|I|$), we denote by $\sum_{i \in I} M_i$ the image of $\oplus_{i \in I} M_i \to M$; it is the minimal subobject of $M$ containing $M_i$ for all $i$.
	
	\subsubsection{Definition}\label{defn:filtered_objects} For an additive\footnote{The only non-abelian additive $\Ac$ we consider are categories of vector bundles ($\FilLF$ in \ref{filtrations_in_families}) and categories of filtered objects ($\Fil^W_F\LF$ in \ref{notation:fil_X^X}). The former can be embedded in the abelian category of quasi-coherent sheaves, and the latter is generally pre-abelian.} category $\Ac$, we let $\mathrm{Fil}(\Ac)$ denote the additive category of tuples $\Mm$ consisting of an object $M$ of $\Ac$ and a decreasing\footnote{but see \ref{rem:increasing_vs_decreasing}} filtration $\{\mathrm{fil}^n \Mm\}_{n \in \Zb}$ of $M$.
	
	For $M,N \in \Ac$ underlying $\Mm,\Nc \in \mathrm{Fil}(\Ac)$, $\homo_{\mathrm{Fil}(\Ac)}(\Mm,\Nc)$ is the set of $f \in \homo_{\Ac}(M,N)$ for which $f(\mathrm{fil}^n \Mm) \subseteq \mathrm{fil}^n \Nc$ for all $n \in \Zb$.

    \subsubsection{}\label{sssec:filt_is_quasi-abelian}

    If $\Ac$ is abelian, then $\Fil(\Ac)$ is quasi-abelian (\cite[Definition 1.1.3]{SchneidersQA99}) by \cite[Theorem 3.9]{schapira2013derivedcategoryfilteredobjects} and is therefore pre-abelian.
	
	\subsubsection{Definition}\label{defn:forg_fil_gr} For $\Mm$ consisting of $M \in \Ac$ and $\{\mathrm{fil}^n \Mm\}_{n \in \Zb}$, we have functors
	\[
	\forg = \forg_{\m{fil}} \colon \mathrm{Fil}(\Ac) \to \Ac,
	\]
	\[
	\mathrm{fil}^n \colon \mathrm{Fil}(\Ac) \to \Ac,
	\]
	and if $\Ac$ is pre-abelian,
	\[
	\gr^n \colon \mathrm{Fil}(\Ac) \to \Ac
	\]
	for all $n \in \Zb$ defined by $\forg(\Mm) \coloneqq M$, $\mathrm{fil}^n(\Mm) \coloneqq \mathrm{fil}^n \Mm$, and $\gr^n(\Mm) = \gr^n \Mm \coloneqq \mathrm{fil}^n \Mm/\mathrm{fil}^{n+1} \Mm$.
	
	\subsubsection{Definition}\label{defn:fin_sep_ex} A filtration $\{\mathrm{fil}^n M\}_{n \in \Zb}$ on an object $M \in \Ac$ is said to be \emph{separated} if the only subobject contained in all $\fil^n M$ is $0$ and \emph{exhaustive} if the only subobject containing all $\fil^n M$ is $M$. It is \emph{finite} if $\exists \, n_1,n_2 \in \Zb$ such that $\mathrm{fil}^{n_1} M = 0$ and $\mathrm{fil}^{n_2} M = M$.
	
	We denote by $\mathrm{Fil}^{\fin}\Ac$ (resp. $\mathrm{Fil}^{\sex}\Ac$) the subcategory of $\Mm \in \mathrm{Fil}(\Ac)$ with finite (resp. separated and exhaustive) filtration. If $\Ac$ is Artinian and Noetherian (e.g., the category of finitely-generated modules over an Artinian ring such as a field), then $\mathrm{Fil}^{\fin}(\Ac) =\mathrm{Fil}^{\sex}(\Ac)$.

    \subsubsection{}\label{sssec:filtered_tensor}

    	If $\Ac$ is a (symmetric) pre-abelian tensor category, we define a (symmetric) additive tensor structure on $\mathrm{Fil}^{\fin}(\Ac)$ in general and on $\Fil(\Ac)$ if $\Ac$ has countable coproducts, for which
	\[
	\mathrm{fil}^n (\Mm \otimes \Nc) = \sum_{a+b=n} \mathrm{fil}^a \Mm \otimes \mathrm{fil}^b \Nc
	\]
	
	\subsubsection{}\label{Vect_FilVect} For a field $\coeff$, we let $\Vect_\coeff$ denote the category of finite-dimensional vector spaces over $\coeff$. We let $\FilVect_{\coeff}$ denote $\mathrm{Fil}^{\fin}(\Vect_{\coeff}) = \mathrm{Fil}^{\sex}(\Vect_{\coeff})$.

    \subsubsection{Convention} \emph{Whenever discussing exactness or strictness in $\Fil(\Ac)$ or $\Gr(\Ac)$, we assume that $\Ac$ has notions of image and exact sequence.}
	
	\subsubsection{Definition}\label{defn:strict_exact} A morphism $f \colon \Mm \to \Nc$ of filtered objects is \emph{strict} or \emph{admissible} if $\Im(\restr{f}{\mathrm{fil}^n \Mm}) = \mathrm{fil}^n \Nc \cap \operatorname{Im}{f} \subseteq N$ for all $n \in \Zb$. A sequence $\Lc \xrightarrow{f} \Mm \xrightarrow{g} \Nc$ is \emph{exact} if it is exact as a sequence in $\Ac$ and $f$ is strict. In particular, $\forg$ (\ref{defn:forg_fil_gr}) is exact.

    \subsubsection{\somewhatrem} That $\Fil(\Ac)$ is quasi-abelian when $\Ac$ is abelian (\ref{sssec:filt_is_quasi-abelian}) implies that $\Fil(\Ac)$ forms an \emph{exact category} (\cite[Definition 2.1]{BuehlerExactCats10}).
	
	\subsubsection{\notreallyrem}\label{rem:SES_filt} A sequence $0 \to \Lc \xrightarrow{f} \Mm \xrightarrow{g} \Nc \to 0$ in $\Fil(\Ac)$ is short exact iff it is short exact in $\Ac$, and $f,g$ are both strict.
	
	\subsubsection{\notreallyrem}\label{rem:increasing_vs_decreasing} Definition \ref{defn:filtered_objects} was stated in terms of decreasing filtrations, such as the Hodge filtrations. All of the above results and definitions apply equally to increasing filtrations by the rule that an increasing filtration is a decreasing filtration with negated indices.

	\subsection{Categories of graded objects}\label{sec:graded}
	
	For an additive category $\Ac$, we let $\mathrm{Gr}(\Ac)$ denote the category of functors $\Mm$ from the discrete category $\Zb$ into $\Ac$, where $\Mm_i \coloneqq \Mm(i)$. 
    
 \subsubsection{}\label{sssec:gr_abelian}   If $\Ac$ is abelian, then so is $\Gr(\Ac)$.
	
	\subsubsection{}\label{defn:gr_fin} A graded object $\Mm$ is \emph{finite} if $\Mm_i = 0$ for all but finitely many $i$, and the category of such $\Mm$ is denoted $\mathrm{Gr}^{\fin}(\Ac)$.

	\subsubsection{}\label{gr_forg} If $\Ac$ admits countable coproducts (e.g., quasicoherent sheaves, countably-generated modules), we have a functor 
	\[
	\forg = \forg_{\gr} \colon \mathrm{Gr}(\Ac) \to \Ac
	\]
	given by $\forg(\Mm) \coloneqq \bigoplus_{i \in \ZZ} \Mm_i$. More generally, for arbitrary additive $\Ac$ (e.g., coherent sheaves or vector bundles), we still have
	\[
	\forg = \forg_{\gr} \colon \mathrm{Gr}^{\fin}(\Ac) \to \Ac,
	\]
	defined in the same way. We sometimes think of $\Mm \in \Gr(\Ac)$ as $\forg(\Mm)$ along with the collection of subobjects $\Mm_i$.

\subsubsection{}\label{sssec:gr_tensor} If $\Ac$ is an additive tensor category that admits countable coproducts, we have an additive tensor structure on $\mathrm{Gr}(\Ac)$ for which
	\[
	(\Mm \otimes \Nc)_n = \bigoplus_{a+b=n} \Mm_a \otimes \Nc_b.
	\]

Without the assumption of countable coproducts, the same formula defines an additive tensor structure on $\mathrm{Gr}^{\fin}(\Ac)$.

In either case, $\forg_{\gr}$ is monoidal.

	\subsubsection{}\label{defn:functors_fil_gr}
	When $\Ac$ is pre-abelian, we have a tensor functor
	\[
	\gr \colon \Fil(\Ac) \to \Gr(\Ac)
	\]
	defined by $\gr(\Mm)_i = \gr^i \Mm$. When $\Ac$ admits countable coproducts (resp. when restricting to finite gradings) we have a tensor functor
	\[
	\fil \colon \Gr(\Ac) \to \Fil(\Ac)
	\]
    (resp.
    \[
	\fil \colon \Gr^{\fin}(\Ac) \to \Fil^{\fin}(\Ac))
    \]
	sending $\Mm \in \Gr(\Ac)$ to the filtered object defined by $M \coloneqq \forg_{\gr}(\Mm)$ and $\fil^i M \coloneqq \bigoplus_{j \ge i} \Mm_j$.\footnote{Or $j \le i$ when considering increasing filtrations, c.f. \ref{rem:increasing_vs_decreasing}.}

\subsubsection{}\label{sssec:grfil_id}

	There are canonical isomorphisms of tensor functors
    \[
    \forg_{\fil} \circ \fil \xrightarrow{\sim} \forg_{\gr}
    \]
	\[
	\gr \circ \fil \xrightarrow{\sim} \id.
	\]
	
	\subsubsection{}\label{exact_gr} A sequence $\Lc \to \Mm \to \Nc$ is exact if $\Lc_i \to \Mm_i \to \Nc_i$ is exact for all $i \in \ZZ$. The functors $\forg_{\gr}$, $\fil$, and $\gr$ are all exact, the latter by the following lemma:

	\begin{sLemma}\label{lemma:exactness_conditions} Consider a sequence $\Lc \xrightarrow{f} \Mm \xrightarrow{g} \Nc$ in $\Fil(\Ac)$ for $\Ac$ abelian. Assuming either the underlying sequence $L \to M \to N$ is exact in $\Ac$ or that $\Mm,\Nc$ are exhaustive, the following are equivalent:
	
	\begin{enumerate}
		\item\label{item:exact} $\Lc \xrightarrow{f} \Mm \xrightarrow{g} \Nc$ is exact
		\item\label{item:exact_for_n} $\mathrm{fil}^n \Lc \xrightarrow{\fil^n f} \mathrm{fil}^n \Mm \xrightarrow{\fil^n g} \mathrm{fil}^n \Nc$ is exact for all $n \in \Zb$
		%    \item\label{item:exact_for_graded} $\gr \Lc \xrightarrow{\gr f} \gr \Mm \xrightarrow{\gr g} \gr \Nc$ is exact
	\end{enumerate}
	
	Either of the above implies
	\begin{enumerate}[resume]
		\item\label{item:exact_for_graded} $\gr \Lc \to \gr \Mm \to \gr \Nc$ is exact,
	\end{enumerate}
	and the converse holds if $\Lc,\Mm,\Nc$ are finite and $g \circ f = 0$.
	%% decided to assume separated and exhaustive
	
	\subsubsection{Proof of \ref{lemma:exactness_conditions}}

	Let $n \in \Zb$. Then we have a diagram
	\[
	\xymatrix{
		\mathrm{fil}^n \Lc \ar[r]_-{\mathrm{fil}^n f} \ar[d]_-{\iota_{L,n}}
		& \mathrm{fil}^n \Mm \ar[r]_-{\mathrm{fil}^n g} \ar[d]_-{\iota_{M,n}}
		& \mathrm{fil}^n \Nc \ar[d]^-{\iota_{N,n}} \\
		L \ar[r]^-{f}
		& M \ar[r]^-{g}
		& N
	}
	\]
	
	Note that
	\[\tag{*}
	\Ker(g) \cap \mathrm{fil}^n \Mm = \Ker(g \circ \iota_{M,n}) = \Ker(\iota_{L,n} \circ \mathrm{fil}^n g) = \Ker(\mathrm{fil}^n g)\]
	because $\iota_{N,n}$ is injective.
	
	Assume (\ref{item:exact_for_n}) and exhaustiveness for $\Mm,\Nc$. If $\alpha \in \Ker(g)$,\footnote{We phrase the proof as if $\Ac$ is a category of modules, noting that this suffices by the Freyd--Mitchell Emebedding Theorem.} then $\alpha \in \fil^n \Mm$ for sufficiently small $n$ by exhaustiveness, so $\alpha \in \Ker(\mathrm{fil}^n g)$ by (*), and thus $\alpha \in \operatorname{Im}(\fil^n f)$ by exactness of the top, hence $\alpha \in \operatorname{Im}(f)$. If $\alpha \in \operatorname{Im}(f)$, say $\alpha = f(\beta)$, then $\beta \in \fil^n \Nc$ for sufficiently small $n$, so $\alpha \in \operatorname{Im}(\fil^n f) = \Ker(\mathrm{fil}^n g) = \Ker(g) \cap \mathrm{fil}^n \Mm$, i.e., $\alpha \in \Ker(g)$. It follows that $L \xrightarrow{f} M \xrightarrow{g} N$ is exact.
	
	Assuming $L \xrightarrow{f} M \xrightarrow{g} N$ is exact, we have by (*) that $\operatorname{Im}(f) \cap \mathrm{fil}^n \Mm = \Ker(g) \cap \mathrm{fil}^n \Mm = \Ker(\mathrm{fil}^n g)$. Strictness of $f$ (at $n$) is the condition that the left-hand side is $\operatorname{Im}(\mathrm{fil}^n f)$, and exactness of the top is the condition that the right-hand side is $\operatorname{Im}(\mathrm{fil}^n f)$, so we have shown (\ref{item:exact}) $\iff$ (\ref{item:exact_for_n}).
	
	For either statement about (\ref{item:exact_for_graded}), we may assume $g \circ f = 0$. In this case, $\fil^n g \circ \fil^n f = 0$ for all $n$, so $\gr^n g \circ \gr^n f = 0$.
	
	We thus get a short exact sequence of chain complexes
	\[
	\xymatrix{
		0 \ar[d] & 0\ar[d] & 0\ar[d]\\
		\fil^{n+1} \Lc \ar[r]^-{\fil^{n+1} f} \ar[d] &
		\fil^{n+1} \Mm \ar[r]^-{\fil^{n+1} g} \ar[d] &
		\fil^{n+1} \Nc \ar[d]\\
		\fil^{n} \Lc \ar[r]_-{\fil^{n} f} \ar[d] &
		\fil^{n} \Mm \ar[r]_-{\fil^{n} g} \ar[d] &
		\fil^{n} \Nc \ar[d]\\
		\gr_n \Lc \ar[r]_-{\gr^n f} \ar[d] &
		\gr_n \Mm \ar[r]_-{\gr^n g} \ar[d] &
		\gr_n \Nc \ar[d]\\
		0 & 0 & 0
	}
	\]
	and thus a long exact sequence
	%\comment{ID: Please fix this overfull H-box}
	\begin{align*}
		\Ker(\fil^{n+1} g)/\operatorname{Im}(\fil^{n+1} f) \to \Ker(\fil^{n} g)/\operatorname{Im}(\fil^{n} f)\\
		\to \Ker(\gr^n g)/\operatorname{Im}(\gr^n f) \to \Coker(\fil^{n+1} g)
	\end{align*}
	
	If we assume (\ref{item:exact_for_n}), then $\Ker(\gr_n g)/\operatorname{Im}(\gr_n f) = 0$ for all such $f$ and $g$, which implies (\ref{item:exact_for_graded}).
	
	Conversely, assume (\ref{item:exact_for_graded}) and that $\Lc,\Mm,\Nc$ are finite. Choose $m>>0$ so that $\fil^{m} \Lc =\fil^{m} \Mm = \fil^{m} \Nc = 0$. Then clearly $0 \to \fil^{n} \Lc \to \fil^{n} \Mm \to \fil^{n} \Nc \to 0$ is exact for $n \ge m$.
	
	Now we apply downward induction on $n$ to show that $\Ker(\fil^{n} g)/\operatorname{Im}(\fil^{n} f) =0$ for all $n$. We have shown this for $n+1$. Then $\Ker(\fil^{n+1} g)/\operatorname{Im}(\fil^{n+1} f) = \Ker(\gr^n g)/\operatorname{Im}(\gr^n f) = 0$ by induction, so $\Ker(\fil^{n} g)/\operatorname{Im}(\fil^{n} f) =0$.

	\subsubsection{\notreallyrem}\label{rem:exact_lemma_increasing} As in \ref{rem:increasing_vs_decreasing}, Lemma \ref{lemma:exactness_conditions} applies equally to exact sequences of objects with increasing instead of decreasing filtration.
	
	\subsubsection{\remish}\label{rem:filn_is_exact} Lemma \ref{lemma:exactness_conditions} implies that $\mathrm{fil}^n \colon \mathrm{Fil}(\Ac) \to \Ac$ is exact in the sense that it sends exact sequences (in the sense of \ref{defn:strict_exact}) to exact sequences.

    \end{sLemma}

	\subsection{Functors into categories of filtered objects}
	
	\subsubsection{}\label{functors_into_fil} Let $\Ac,\Bc$ be additive categories. A functor $\phi \colon \Ac \to \mathrm{Fil}(\Bc)$ is equivalent to giving a functor $\om \colon \Ac \to \Bc$ and a decreasing sequence of subfunctors $\{F^n \om\}_{n \in \Zb}$ of $\om$. Indeed, given $\phi$, we set $\om \coloneqq \forg \circ \phi$ and $F^n \om \coloneqq \mathrm{fil}^n \circ \phi$. Given $\om$, $\{F^n \om\}$, and $M \in \Ac$ we set $\phi(M) \in \mathrm{Fil}(\Bc)$ to be $\om(M) \in \Bc$ with filtration given by $\mathrm{fil}^n \phi(M) \coloneqq F^n \om(M)$.
	
	\subsubsection{}\label{exactness_into_fil} If $\phi$ is exact and $f \in \Hom_{\Ac}(\Lc,\Mm)$, then $\phi(\Lc \xrightarrow{f} \Mm \to \Coker(f))$ is exact in $\Fil(\Bc)$, so $\phi(f)$ is strict. It follows easily that $\phi$ is exact iff $\om$ is exact and $\phi(f)$ is strict for all $f \in \operatorname{Mor}(\Ac)$. 
	
	By Lemma \ref{lemma:exactness_conditions} (for $\Ac$ abelian), $\phi$ is exact (\ref{rem:filn_is_exact}) iff $\om$ and each $F^n \om$ is exact. If $\phi$ lands in $\Fil^{\fin}(\Bc)$, then it is exact iff $\gr \circ \phi$ is exact.
	
	\subsubsection{} A symmetric monoidal structure on $\phi$ is equivalent to a symmetric monoidal structure on $\om$ for which
	\[
	F^n \omega (M \otimes N) = \sum_{a+b=n} F^a \om(M) \otimes F^b \om(N)
	\]
	for all $n$ under the isomorphism $\om(M \otimes N) \xrightarrow{\sim} \om(M) \otimes \om(N)$.
	
	\subsubsection{}\label{sssec:filtered_fiber_functor} If $\Tc$ is a Tannakian category (\ref{defn:tannakian_category}) over $\coeff$, a \emph{filtered fiber functor} (\cite[Definition 4.4(i)]{ZieglerGrFil15}; functor with exact $\otimes$-filtration in \cite[2.1.1]{SaavedraRivanoBook}) over $\coeff'$ is an exact tensor functor $\phi \colon \Tc \to \FilVect_{\coeff'}$. We will return to this in \S\ref{defn:graded_filtered_fib_func}.

\subsubsection{}\label{sssec:Gr_Fil_functors}

If $\om \colon \Ac \to \Bc$ is an additive functor between abelian (resp. additive) categories, then we have induced functors $\Fil(\om) \colon \Fil(\Ac) \to \Fil(\Bc)$ (resp. $\Gr(\om) \colon \Gr(\Ac) \to \Gr(\Bc)$), the former by setting $\fil^n(\om(A)) \coloneqq \Im(\om(\fil^n(A)) \to \om(A))$. In this way, $\Fil$ and $\Gr$ can be viewed as functors from the category of abelian (resp. additive) categories to the category of quasi-abelian (resp. additive) categories. The same is true for $\Fil^{\fin}$ and $\Gr^{\fin}$, as well as $\Fil^{\sex}$ if $\om$ respects countable filtered limits and colimits.

We always have
\[\forg_{\fil} \circ \Fil(\om) = \om \circ \forg_{\fil}\]
\[\forg_{\gr} \circ \Gr^{\fin}(\om) = \om \circ \forg_{\gr}\]
\[
\fil \circ \Gr^{\fin}(\om) = \Fil^{\fin}(\om) \circ \fil.
\]
If countable coproducts exist and are respected by $\om$, then
\[\forg_{\gr} \circ \Gr(\om) = \om \circ \forg_{\gr}\]
\[
\fil \circ \Gr(\om) = \Fil(\om) \circ \fil.
\]
If $\om$ is left exact, then \[\fil^n \circ \Fil(\om) = \om \circ \fil^n,\]
and if $\om$ is exact, then
\[
\gr^n \circ \Fil(\om) = \om \circ \gr^n
\]
\[
\gr \circ \Fil(\om) = \Gr(\om) \circ \gr.
\]

	\subsection{Filtrations on Tannakian categories}
	
	For the rest of \S \ref{sec:wt_fil_tann}, we fix a field $\coeff$ of characteristic $0$ and work exclusively with a Tannakian category $\Tc$ over $\coeff$.

	%\subsubsection{Definition}
	%A \emph{filtration} on a Tannakian category $\Tc$ is an increasing sequence $\{W_n\}_{n \in \Zb}$ of subfunctors of the identity functor such that for any $M \in \Ob(\Tc)$, the canonical inclusion
	%\[
	%W_n(W_n M) \hookrightarrow W_n M
	%\] is an isomorphism for every $n$, the induced filtration $\{W_n M\}_{n \in \Zb}$ %of $M$ is finite, and for any $M,N \in \Ob(\Tc)$ and $n \in \Zb$, we have
	%\[
	%W_n(M \otimes N) = \sum_{p+q=n} W_p M \otimes W_q N.
	%\]

	\subsubsection{Definition}
	A \emph{filtration} $W$ on $\Tc$ is a tensor functor $\xi_W \colon \Tc \to \mathrm{Fil}^{\mathrm{fin}}\Tc$ such that $\mathrm{forg}_{\fil} \circ \xi_W = \id_{\Tc}$.
	
	We set $W_n \coloneqq \mathrm{fil}^{-n} \circ \phi$, so that $\{W_n M\}_{n \in \Zb}$ is an increasing filtration of any $M \in \Tc$ (c.f. \ref{rem:increasing_vs_decreasing}).
	
	A Tannakian category $\Tc$ with a choice of filtration is called a \emph{filtered Tannakian category}.

	%% Not true if we think about the image of nilpotent morphism
	%\begin{rem}
	%An equivalent definition is a choice of filtration on each object that is functorial, i.e., such that every morphism respects the filtration.
	%\end{rem}
	
	\subsubsection{}\label{defn:gr_in_filtered_tannakian}
	For a filtered Tannakian category $\Tc$, $M \in \Ob(\Tc)$, and $n \in \Zb$, we set
	\[
	\gr_n M \coloneqq W_n M/W_{n-1}M.
	\]
	
	Note that $\gr_n$ is a functor from $\Tc$ to itself. In the notation of \ref{defn:forg_fil_gr}, we have $\gr_n = \gr^{-n} \circ \xi_W$.

    We also set \begin{equation}\label{eqn:gr^W}\gr^W(M) \coloneqq \bigoplus_{n \in \Zb} \gr_n M.\end{equation} In the notation of \ref{defn:functors_fil_gr}, $\gr^W = \forg_{\gr} \circ \gr \circ \xi_W$.

	\subsubsection{Definition/Theorem}\label{dfthm:strict}
	A filtration $\{W_n\}$ on a Tannakian category $\Tc$ is \emph{strict} if any of the following equivalent conditions holds:
	\begin{enumerate}
		\item\label{item:0} The functor $\xi_W$ is exact.
		\item\label{item:i} Every morphism in $\Tc$ is strict for the filtration
		\item\label{item:ii} Each functor $W_n$ is exact
		\item\label{item:iii} The associated graded functor $\gr^W \colon \Tc \to \Tc$ is exact
		%    \item\label{item:iv} For $M$ and $N$ pure of weights $m$ and $n$ respectively, $\homo_{\Tc}(M,N)=0$ if $m \neq n$.
	\end{enumerate}
	
	A Tannakian category $\Tc$ equipped with a strict filtration is called a \emph{strictly filtered Tannakian category}.
	
	\subsubsection{Proof of \ref{dfthm:strict}}
	
	Since $\mathrm{forg} \circ \xi_W = \id_{\Tc}$ is exact, \ref{exactness_into_fil} says (\ref{item:0}) is equivalent to $\xi_W(f)$ being strict for all $f \in \operatorname{Mor}(\Tc)$, which is (\ref{item:i}).
	
	The equivalence between (\ref{item:0}) and (\ref{item:ii}) (resp. \ref{item:iii}) is immediate from the first (resp. second) part of Lemma \ref{lemma:exactness_conditions}, noting again that $\mathrm{forg} \circ \xi_W = \id_{\Tc}$ is exact, and $\xi_W$ lands in $\mathrm{Fil}^{\mathrm{fin}}\Tc$.
	
	\subsubsection{\remish}
	If $\xi_W$ is a strict filtration on $\Tc$ and $\om$ a choice of fiber functor, then $W_n$ turns $\om$ into a filtered fiber functor (Definition \ref{defn:weight_filtered_fib}). Conversely, if $\om$ is a filtered fiber functor over $\coeff$ whose Tannakian Galois group preserves the filtration on every object, then this filtration comes from a strict filtration on $\Tc$ (\ref{PG_weight-filtered}).

\subsubsection{Definition}\label{defn:concentrated_in_deg_n}
We say that an object $M$ in a filtered Tannakian category $\Tc$ is \emph{concentrated in degree $n \in \Zb$} if $0 = W_{n-1} M \subseteq W_n M = M$.

	\begin{ssLemma}
		
		\label{lemma:strictness_crit}
		Suppose $\Tc$ has a filtration. Then the filtration is strict iff for any $M,N \in \Ob(\Tc)$ and $m \neq n \in \Zb$, no nonzero subquotient of $\gr_m M$ is isomorphic to a nonzero subquotient of $\gr_n N$.

		\begin{proof}

Suppose the filtration is strict. Then every morphism in $\Tc$ is strict for the filtration. In particular, if $0 \to A \to B \xto{f} C \to 0$ is an exact sequence, then $W_i A = W_i B \cap A$, and $W_i C = f(W_i B)$ for all $i \in \Zb$. From this, it is easy to see that $\gr_m M$ (resp. $\gr_n N$) and any subquotient thereof is concentrated in degree $m$ (resp. $n$). Thus if $L$ is isomorphic to a subquotient of both $\gr_m M$ and $\gr_n N$ for $m<n$, then $L = W_m L \subseteq W_{n-1} L = 0$, so $L=0$.

Conversely, fix a morphism $f \colon A \to B$ in $\Tc$ and $n \in \Zb$. Let $C$ denote the preimage in $A$ of $W_n B$, so that $C/\Ker{f} \simeq f(C) = W_n B \cap f(A)$. Modding out on both sides by $W_n A$, we get $C/(W_n A + \Ker{f}) \simeq (W_n B \cap f(A))/f(W_n A)$. It suffices to show that this is always zero. 
			
			%Note for $m \ge n$ that $(W_n A + \Ker{f}) \cap W_m A = W_n A + (\Ker{f} \cap W_m A)$.
			
			Suppose $D \coloneqq C/(W_n A + \Ker{f})$ is nonzero. Then for $i \ge n$, $V_i \coloneqq (C \cap W_i A)/(W_n A + (\Ker{f} \cap W_i A))$ is an increasing sequence of subobjects of $D$ with $V_n = 0$ and $V_i = D$ for $i$ sufficiently large (e.g., such that $W_i A = A$). Let $m > n$ such that $V_m \neq 0$ but $V_{m-1} = 0$. Then $V_m = V_m/V_{m-1} = (C \cap W_m A)/((\Ker{f} \cap W_m A) + C \cap W_{m-1} A)$ is a quotient of $(C \cap W_m A)/(C \cap W_{m-1} A)$, which is the image of $C$ in $\gr_m A$. In particular, it is a subquotient of $\gr_m A$.
			
			On the other hand, for $i \le n$, let $U_i \coloneqq (W_i B \cap f(A))/(W_i B \cap f(W_n A))$. We clearly have $U_n = D$, while $U_i = 0$ for $i$ sufficiently small. Let $j \le n$ so that $U_j \cap V_m \neq 0$ while $U_{j-1} \cap V_m = 0$. Then $U_j \cap V_m$ is a nonzero subobject of $V_m$ and hence a subquotient of $\gr_m A$.
			
			Now $U_j/U_{j-1}$ is a quotient of $(W_j B \cap f(A))/(W_{j-1} B \cap f(A))$, which is the image of $W_j B \cap f(A)$ in $\gr_j B$, hence a subquotient of $\gr_j B$. On the other hand, $U_j \cap V_m$ maps isomorphically onto its image in $U_j/U_{j-1}$, so $U_j \cap V_m$ is a subquotient of $\gr_j B$. Since $j \le n < m$, we find that $U_j \cap V_m = 0$, a contradiction. Therefore, $D$ is zero, and the filtration is strict.

			%On the other hand, for $i \ge n$, let $U_i \coloneqq (W_n B \cap f(A))/f(W_i A \cap C)$. Note that $W_n A \subseteq C$, so $U_n = D$. For $i$ sufficiently large, we have $W_i A \cap C = C$, so that $f(W_i A \cap C) = f(C) = W_n B \cap f(A)$, so $U_i = 0$. Let $j$ be such that $U_j \neq 0$ while $U_{j+1}=0$. Then $W_n B \cap f(A) = f(W_{j+1} A \cap C)$, so $U_j = f(W_{j+1} A \cap C)/f(W_j A \cap C)

			%First, we verify that the filtration is strict. Fix a morphism $f \colon A \to B$ in $\mathrm{MMot}(k,\Qb)$ and $n \in \Zb$. Let $C$ denote the preimage in $A$ of $f(A) \cap W_n B$; clearly $C$ contains $W_n A$. We can restrict $f$ to a map $C \to W_n B$ and then get an induced map $g \colon C/W_n A \to (f(A) \cap W_n B)/f(W_n A)$. Note that this map is surjective because $C$ surjects onto $f(A) \cap W_n B$.
			
			%Have $A \to B$ and $D \subseteq B$ and $C$ is kernel of map $A \to B \to B/D$. Then cokernel of $C \to A$ should be the same as the kernel of $B/D \to B/A$
			
			%Suppose otherwise, i.e. that there is a morphism $f \colon A \to B$ and $n \in \mathbb{Z}$ such that $f(W_n A) \subsetneq f(A) \cap W_n B$. We consider the induced map $A/W_n A \to B/f(W_n A)$. 
			
			%Clearly $C$ contains $W_n A$. Furthermore, if $C=W_n A$, then because $f(A) \cap W_n B = f(C)$ (why???**), we would get $C=
			
		\end{proof}

    \subsubsection{}\label{sssec:double_gr^W}

If the filtration is strict, note $\gr_n M$ is concentrated in degree $n$. It follows that \begin{equation}\label{eqn:W_n_of_gr^W}W_n \gr^W M = \sum_{m \le n} \gr_m M\end{equation} and hence that
\[
\gr_n^W \gr^W M = \sum_{m \le n} \gr_m M \big / \sum_{m \le n-1} \gr_m M \simeq \gr_n M,
\]
and thus
\begin{equation}\label{eqn:double_gr^W_M}
\gr^W \gr^W M \coloneqq \bigoplus_i \gr^W_i \gr^W M \simeq \bigoplus_i \gr^W_i M = \gr^W M.
\end{equation}
Furthermore, any $f \in \Hom_{\Tc}(\gr^W M,\gr^W N)$ respects the decomposition (\ref{eqn:gr^W}) by Lemma \ref{lemma:strictness_crit} and thus equals its own associated graded. Thus the isomorphism (\ref{eqn:double_gr^W_M}) upgrades to an isomorphism of functors
\begin{equation}\label{eqn:double_gr^W_functor}
\gr^W \circ \gr^W \simeq \gr^W \colon \Tc \to \Tc.
\end{equation}

	\end{ssLemma}

	\subsubsection{Definition}
	\label{defn:weight_filtr_cat}
	We say that a strict filtration $W_n$ on a Tannakian category $\Tc$ is a \emph{weight filtration} if for any $M \in \Ob(\Tc)$ and $n \in \Zb$,
	\[
	\gr_n^W M
	\]
	is semisimple in $\Tc$. In this case, $\gr^W M$ is semisimple.
	
	A Tannakian category $\Tc$ with a choice of weight filtration is called a \emph{weight-filtered Tannakian category}.
	%\commentDC{ID: Let's add a reference to Goncharov's similar notion of a \emph{mixed category over $\mathcal{P}$} in \S3 of \cite{GonMEM}. Curiously, Goncharov doesn't assume $\Mm_\Pp$ is Tannakian. DC: Added, but I think this is more relevant to the section you put in an appendix (which tbh I don't think it should be there). ID: Thanks.}
	
	\subsubsection{\yesrem}
	\label{rem:goncharov_reference}
	Goncharov lists some of the defining properties of a weight-filtered Tannakian category when describing categories of the form $\Mm_{\Pp}$ in \cite[\S 3]{GonMEM}, especially at the bottom of p.15 of loc.cit. We mention these briefly in Remark \ref{sec:motives_subcategory}.
	
	\subsubsection{\notreallyrem}
	\label{rem:filtr_sub}
	If $\Tc$ is a filtered Tannakian category, and $\Tc' \subseteq \Tc$ is a Tannakian subcategory, then $\Tc'$ acquires a natural filtration. If the filtration on $\Tc$ is strict (resp. a weight filtration), then so is the filtration on $\Tc'$. For the weight filtration, note that if $M \in \Tc'$ is semisimple in $\Tc$, then it is semisimple in $\Tc'$, so $\Tc'$ is weight-filtered if $\Tc$ is.
	
	% Is the converse true?? ** If $M$ is semisimple in $\Tc'$, suppose it has a subobject in $\Tc$.

	\subsection{Hodge and Frobenius fiber functors}

	\subsubsection{Definition}
	\label{defn:hodge_filtr}
	If $\Tc$\footnote{Technically a strictly filtered Tannakian category is a pair $(\Tc,\xi_W)$, but we suppress this in the notation.} is a strictly filtered Tannakian category over a field $\coeff$ and $\coeff'/\coeff$ a field extension, a \emph{Hodge-filtered fiber functor over $\coeff'$} on $\Tc$ is a filtered fiber functor (\ref{sssec:filtered_fiber_functor})
	\[
	\om_H \colon \Tc \to \mathrm{FilVect}_{\coeff'}.
	\]
	
	By abuse of language, we sometimes refer to $\om \coloneqq \forg \circ \om_H$ as a Hodge-filtered fiber functor, viewing the filtration as a sequence of subfunctors $F^p\om$ as in \ref{functors_into_fil}.
	
	%%% We *could* phrase this more easily as an exact tensor-functor to filtered vector space. Well, maybe - I'm not sure if the condition of the filtration being finite is automatic or not. It might be automatic because of rigidity.

	%\subsubsection{Remark} By \cite[Lemma 4.1]{ZieglerGrFil15}, the exactness of the $F^p \om$ is equivalent to $\om_H$ being a filtered fiber functor in the sense of \cite[Definition 4.4(i)]{ZieglerGrFil15}, which is equivalent to $\om(f)$ being strict for the Hodge filtration for every morphism $f$ of $\Tc$. We will return to this in \ref{sec:graded_filtered_fib_func}.

	\subsubsection{Definition}\label{defn:theory_of_weights}
	A \emph{theory of weights} on a field $\coeff'_0$ is a collection of disjoint subsets $w_n=w_n(\coeff'_0)$ of $\overline{\coeff'_0}$ indexed by $n \in \Zb$ such that $1_{\coeff'_0} \in w_0$, if $\alpha \in w_n$ and $\beta \in w_m$, then $\alpha \beta \in w_{n+m}$, and if $\alpha^m \in w_{nm}$, then $\alpha \in w_n$.

	\subsubsection{Example}
	\label{ex:q_weil_weight}
	If $\coeff$ is of characteristic $0$ and $q$ is a positive real number, we define the \emph{$q$-Weil weight structure} by setting $W_n$ to be the set of $\alpha \in \overline{\coeff}$ such that for any embedding $\iota \colon \Qb(\alpha) \hookrightarrow \Cb$, we have $|\iota(\alpha)| = q^{\frac{n}{2}}$.

	\subsubsection{Definition}\label{defn:ph-mod}
	Let $R$ a $\coeff$-algebra equipped with a finite-order $\coeff$-algebra automorphism $\phi$. A \emph{$\ph$-module over $R$} is a an $R$-module equipped with a $\phi$-linear automorphism $\ph$. $\ph$-Modules over $R$ form an abelian tensor category with 
	\[
	\ph_{X \otimes Y} = \ph_X \otimes \ph_Y.
	\]
	We denote the tensor category of filtered $\ph$-modules over $R$ by $\ph\m{Mod}(R)$.
	
	\subsubsection{\notreallyrem}\label{rem:ph^d_is_linear}
	Let $d$ be the order of $\phi$ acting on $R$. Then $\ph^d$ is an $R$-linear automorphism of $M$.

	\begin{ssDefinition}
		\label{defn:frob_equiv}
		Let $\Tc$ be a strictly filtered Tannakian category over a field $\coeff$, $\coeff_0'/\coeff$ a field extension, $\{w_n\}$ a theory of weights on $\coeff'_0$, and $\phi$ an automorphism of $\coeff'_0$ of order $d$ fixing $\coeff$. A \emph{Frobenius-equivariant fiber functor over $(\coeff'_0,\phi)$} is an exact tensor functor
		\[
		\om_{\phi} \colon \Tc \to \ph\m{Mod}(\coeff'_0)
		\]
		such that for $M=\gr_n^W(M)$, the eigenvalues of $\ph^d$ acting on $\om_{\phi}(M)$ are in $w_{dn}$. We set $\om_0 \coloneqq \forg_{\phi} \circ \om_{\phi}$.
	\end{ssDefinition}
	
	\begin{ssDefinition}
		\label{defn:hodge_frob}
		If $\coeff'/\coeff'_0/\coeff$ are field extensions where $\coeff'_0$ has a theory of weights and a finite-order automorphism $\phi$, a \emph{Hodge-filtered Frobenius-equivariant fiber functor} $\om_{H,\phi} = (\om_H, \om_{\phi}, \eta)$ over $(\coeff',\coeff'_0,\phi)$ is a triple of a Hodge-filtered fiber functor $\om_H$ valued in $\FilVect_{\coeff'}$, a Frobenius-equivariant fiber functor $\om_{\phi}$ valued in $\ph\m{Mod}(\coeff'_0)$, and a natural $\otimes$-isomorphism \[\om \xrightarrow[\raisebox{0.25 em}{\smash{\ensuremath{\sim}}}]{\eta} \Ex_{\coeff'_0}^{\coeff'}  \circ \om_0.\] We use $\eta$ to implicitly identify the two functors over $\coeff'$.
	\end{ssDefinition}
	
	\subsubsection{\remish}
	\label{rem:adm}
	A Hodge-filtered Frobenius-equivariant fiber functor over $(\coeff',\coeff'_0,\phi)$ naturally outputs filtered $\ph$-modules over $\coeff'$. One may define such a functor to be \emph{admissible} if its image lands in admissible filtered $\ph$-modules, where the Newton polygon is defined in an obvious way using the theory of weights on $\coeff'_0$.
	
	\subsubsection{\notreallyrem}
	\label{rem:hodge_frob_sub}
	If $\Tc$ is a strictly filtered Tannakian category with a Hodge-filtered Frobenius-equivariant fiber functor $\om_{H,\phi}$, and $\Tc' \subseteq \Tc$ is a Tannakian subcategory, then $\om_{H,\phi}$ restricts to a Hodge-filtered Frobenius-equivariant fiber functor on $\Tc'$.

	%let $\om$ be a fiber functor valued in an extension $F/k$, and let $\{w_n\}$ be a theory of weights on $F$. Then a \emph{Frobenius morphism on $\om$} is a tensor automorphism $\phi$ of $\om$ such that for all $M \in \Tc$, the eigenvalues of $\om(\phi)$ acting on $\om(\gr M)$ are all in $w_n$.
	%\end{defn}
	
	%\begin{defn}
	%A \emph{Frobenius-equivariant fiber functor over $K$} on a weight-filtered Tannakian category $\Tc$ is a fiber functor $\om$ valued in $\Vect_K$ along with a choice of Frobenius morhpisms on $\om$.
	%\end{defn}
	
	%Definition: Frobenius-filtered
	%[Ishai: why "filtered"?  ] fiber functor (an automorphism of the fiber functor [ID: only over $\Qp$] that acts by eigenvalues of weight $n$ on the $n$th weight-graded piece)

	%%A \emph{$\phi$-equivariant fiber functor over $R$} is a monoidal functor
	%%\[
	%%F: T \to \phi\m{Mod}(R)
	%%\]
	%%whose composite with the forgetful functor to $\m{Mod}(R)$ is a fiber functor, such that if $E \in T$ is semisimple of negative weight, then the action of $\phi$ on $F(E)$ is invertible.
	%%\footnote{ID: This is just a suggestion.}
	%%% DC: probably don't want to phrase it this way,

	\section{Examples of Weight-Filtered Tannakian Categories}\label{sec:examples_cats}
	
	\subsection{Strictly geometric Galois representations}\label{sec:sg_galois_reps}

    We recall some constructions that originally appeared in the parallel work \cite[\S 2]{CorwinMECK} and explain how they fit under the general definitions above.
	
	\subsubsection{Definition}
	Let $\basefield$ be a number field and $p$ a prime number, and let 
	\[
	\Rep_{\Qp}^{\mathrm{g}}(G_{\basefield})
	\]
	denote the category of $p$-adic representations of $G_{\basefield}$ that are unramified almost everywhere and de Rham at every place of $\basefield$ above $p$. We let
	\[
	\Rep_{\Qp}^{\mathrm{sg}}(G_{\basefield})
	\]
	denote the subcategory of representations $V$ that are \emph{strongly geometric}, meaning $V$ has a finite increasing filtration $W^{\bullet} V$, known as the \emph{motivic weight filtration}, for which
	\[
	\gr^W_n V
	\]
	is pure of weight $n$ at almost all (unramified) places in the sense of \cite{WeilII}. This category is strictly filtered by Lemma \ref{lemma:strictness_crit}.
	
	\subsubsection{}
	We let $\Rep_{\Qp}^{\rm sg}(G_{\basefield})^{\wss}$ be the subcategory of $V \in \Rep_{\Qp}^{\rm sg}$ for which $\gr_n^W V$ is semisimple for each $n$. Then $\Rep_{\Qp}^{\mathrm{sg}}(G_{\basefield})^{\wss}$ is a neutral weight-filtered Tannakian category over $\coeff=\Qp$ in the sense of Definition \ref{defn:weight_filtr_cat}.
	
	\subsubsection{}\label{hodge_rep}
	Fixing a place $\pf$ of $\basefield$ above $p$, the functor
	\[
	\om^{\dR}_H \colon V \mapsto \om^{\dR}(V) \coloneqq V^{\dR} \coloneqq \operatorname{D}_{\dR}(V) \in \mathrm{FilVect}_{\basefield_{\pf}}
	\]
	is a Hodge-filtered fiber functor on $\Rep_{\Qp}^{\mathrm{sg}}(G_{\basefield})$ over $\coeff'=\basefield_{\pf}$ by \cite[Theorem 5.29]{fontaineouyang}, and $\om^{\dR} \coloneqq \forg \circ \om^{\dR}_H$.
	%% See also \cite[Theorem 8.2.11]{BrinonConrad}
	%% \cite[I.2.2.2(iii)]{FPR}, which proves Fil^0 is exact. Can get the others by Tate twist.

	%Note that if $v \notin S$, then $\Rep_{\Qp}^{\mathrm{sf,S}}(G_{\basefield}) \subseteq \Cc_v$
	
	%%% ** don't see any reason for this hereWe let \[\Rep_{\Qp}^{\mathrm{ss}}(G_{\basefield})\] denote the full subcategory of $\Rep_{\Qp}^{\mathrm{sg}}(G_{\basefield})$ consisting of semisimple objects, or equivalently, those objects for which the weight filtration splits.
	
	\subsection{\texorpdfstring{$S$}{S}-integral variants}
	\label{sec:S-integral}
	For a subset $S \subseteq S_f(\basefield)$, let $\Rep_{\Qp}^{\rm f, S}(G_{\basefield})$ denote the full subcategory of $V \in \Rep_{\Qp}^{\rm g}(G_{\basefield})$ with good reduction\footnote{A $p$-adic representation of $G_{\basefield}$ has \emph{good reduction} at $v$ if it is unramified (resp. crystalline) for $v \notin S_p(\basefield)$ (resp. $v \in S_p(\basefield)$).} at all $v \notin S$. Let $\Rep_{\Qp}^{\mathrm{sf,S}}(G_{\basefield})$ denote the full subcategory of $V \in \Rep_{\Qp}^{\mathrm{f,S}}(G_{\basefield}) \cap \Rep_{\Qp}^{\mathrm{sg}}(G_{\basefield})$ such that:
	\begin{itemize}
		\item For $v \notin S$, the Frobenius at $v$ acts by weight $n$ on $\gr^W_n V$ (resp. $\operatorname{D}_{\cris}(\gr^W_n V)$) for $v \notin S_p(\basefield)$ (resp. $v \in S_p(\basefield)$).
	\end{itemize}
	We similarly define $\Rep_{\Qp}^{\mathrm{sf,S}}(G_{\basefield})^{\wss} \coloneqq \Rep_{\Qp}^{\mathrm{sf,S}}(G_{\basefield}) \cap \Rep_{\Qp}^{\mathrm{sg}}(G_{\basefield})^{\wss}$
	
	\subsubsection{Definition}\label{defn:CCv}
	
	The category $\Rep_{\Qp}^{\rm sf, S}(G_{\basefield})$ has a Frobenius-equivariant fiber functor as long as there is $\pf \in S_p(\basefield) \setminus S$. If we fix $\pf \in S_p(\basefield)$, we may simply work with
	\[
	\Rep_{\Qp}^{\mathrm{sf},\cri(\pf)}(G_{\basefield}) \coloneqq \bigcup_{\pf \notin S} \Rep_{\Qp}^{\rm sf, S}(G_{\basefield}).
	\]
	
	Then the functor
	\[
	V \mapsto \om^{\cri}(V) \coloneqq V_{\cri} \coloneqq \operatorname{D}_{\cris}(V)
	\]
	is a Frobenius-equivariant fiber functor on $\Rep_{\Qp}^{\mathrm{sf},\cri(\pf)}(G_{\basefield})$ over $(\coeff'_0,\phi)=(\basefield_{\pf,0},\Frob_{\pf})$ with respect to the standard Frobenius and the $q$-Weil weight structure, where $q$ is the size of the residue field of $\basefield_{\pf}$.
	
	\subsubsection{}\label{hodge_frob_rep}
	
	The inclusion $B_{\cris} \subseteq B_{\dR}$ induces a functorial tensor isomorphism
	\[
	V^{\dR} \xrightarrow[\raisebox{0.25 em}{\smash{\ensuremath{\sim}}}]{\eta^{\et}} V^{\cris} \otimes_{\basefield_{\pf,0}} \basefield_{\pf}
	\]
	making $\boxed{(\om_H,\om_{\phi},\eta)=(\om^{\dR}_H,\om^{\cri},\eta^{\et})}$ into a Hodge-filtered Frobenius-equivariant fiber functor over $\boxed{(\coeff',\coeff'_0,\phi)=(\basefield_{\pf},\basefield_{\pf,0},\Frob_{\pf})}$ on the strictly filtered Tannakian category $\boxed{\Tc=\Rep_{\Qp}^{\mathrm{sf},\cri(\pf)}(G_{\basefield})}$ over $\boxed{k=\Qp}$ and thus also on its weight-filtered subcategory $\Rep_{\Qp}^{\mathrm{sf},\cri(\pf)}(G_{\basefield})^{\wss}$. Note that it is admissible in the sense of Remark \ref{rem:adm}.

\subsection{Systems of realizations and pre-motivic structures}\label{sec:PSPM}

While categories of Galois representations suffice for some purposes, it is sometimes better to have $\Qb$-linear Tannakian categories of motives. Such a category is known to exist with all the desired properties only in the category of mixed (Artin-)Tate motives, which apply only when $X$ is rational. Therefore, in general, we must consider the unconditional categories of systems of realizations. For a treatment of the conjectural abelian category of mixed motives, see Appendix \ref{app:mixed_motives}.

%The appendix discusses the conjectural abelian category of mixed motives and shows that it is a weight-filtered Tannakian category conditional upon a set of axioms found in the literature.

There are numerous versions of the category of systems of realizations in the literature, for example due to Deligne, Huber, and Jannsen. We focus on Fontaine--Perrin-Riou's category $\mathbf{SPM}_\basefield(\Qb)$ of pre-motivic structures over $\basefield$ with $\Qb$-coefficients, as it has a built-in $p$-adic Hodge comparison isomorphism, which is particularly relevant to $p$-adic periods and the Chabauty--Kim method. We show how to define the subcategory $\mathbf{SM}_\basefield(\Qb)$ of realizations of motives (which \cite[III.4.1.2]{FPR91} assumes but does not define). We make remarks about comparisons between the different categories (\ref{sec:lit_comparison}) and note in particular that $\mathbf{SM}_\basefield(\Qb)$ should be equivalent to the category of mixed motives for absolute Hodge cycles (Remark \ref{rem:absolute_hodge_motives}).

%Jannsen's category $\underline{\mathrm{MM}}_k$ (\cite[\S 4.1]{JannsenBook90}) of mixed motives for absolute Hodge cycles.

\subsubsection{Brief Overview of \S \ref{sec:PSPM}-\ref{sec:SM_A}}

For the convenience of the reader who wishes to take the rest of \S \ref{sec:examples_cats} as a black box, we give a quick guide to the definition of our replacement $\mathbf{SM}_\basefield(\Qb)$ (Definition \ref{defn:motivic_structures}) for the category of mixed motives over $\mathbb{Q}$.

In \ref{PSMP}, we define the category $\mathbf{PSPM}_\basefield(\Qb)$ abstract systems of realizations, and in Definition \ref{sec:SPM}, we define the category $\mathbf{SPM}_\basefield(\Qb)$ of systems with a weight filtration.

In Theorem \ref{thm:real}, we define a realization functor from Voevodsky's triangulated category of geometric motives to the bounded derived category of $\mathbf{SPM}_{\basefield}(\Qb)$. Finally, in Definition \ref{defn:motivic_structures}, we define \[\mathbf{SM}_\basefield(\Qb)\] as the Tannakian subcategory of $\mathbf{SPM}_{\basefield}(\Qb)$ generated by all cohomology objects of complexes in the image of the realization functor.

In \S \ref{sec:variants_SPM}-\ref{sec:pol_wt_fil}, we discuss various subcategories of $\mathbf{SPM}_\basefield(\Qb)$. The most important one is $\mathbf{SPM}_{\basefield}(\Qb)^{\Del,\pol}$, as it contains $\mathbf{SM}_\basefield(\Qb)$; it consists of those that have good reduction and appropriately-weighted Frobenii outside $S$ for some finite set $S$ of places of $\basefield$ and have a polarization. The latter property ensures that the category is weight-filtered. These sections are not strictly needed for the definition of $\mathbf{SM}_{\basefield}(\Qb)$ but imply useful properties of it.

In Definition \ref{defn:SM_beta} and \S \ref{sec:SM_A}, we define some important subcategories of $\mathbf{SM}_\basefield(\Qb)$, including those with good reduction outside a fixed set $S$ of places and those with socle generated by the cohomology of a single fixed variety. In \S \ref{sec:sm_mot_coh}, we discuss the relation between $\Ext$ in this category and motivic cohomology and $K$-theory.

\subsubsection{}
\label{PSMP}

We recall, for a number field $\basefield$, Fontaine--Perrin-Riou's category $\mathbf{PSPM}_\basefield(\Qb)$ (\cite[III.2.1.1]{FPR91}). An object $M$ of $\mathbf{PSPM}_\basefield(\Qb)$ consists of the modules:
\begin{itemize}
	\item A vector space $M_{\dR}$ over $\basefield$ with a finite decreasing Hodge filtration $F^i M_{\dR}$
	\item For each $v \in S_{\infty}(\basefield)$, a $\Qb$-vector space $M_{B,v}$ with an action of $G_{v}$
	\item For each prime number $\ell$, an object $M_{\ell} \in \Rep_{\Ql}^{\mathrm{g}}(G_{\basefield})$
	%% The fact that it is in ^{\mathrm{g}} is because they say it is "pseudo-g\'eometrique", which is what I call geometric
\end{itemize}
and the comparison isomorphisms:
\begin{itemize}
	\item For each pair of a prime number $\ell$ and $v \in S_{\infty}(\basefield)$, an isomorphism
	\[
	\iota_{\ell,v} \colon M_{B,v} \otimes_{\Qb} \Ql \xrightarrow{\raisebox{-1 em}{\smash{\ensuremath{\sim}}}} M_{\ell}
	\]
	compatible with the action of $G_{v}$
	\item For each $v \in S_{\infty}(\basefield)$, a $G_v$-equivariant isomorphism
	\[
	\iota_v \colon M_{B,v} \otimes_{\Qb} \overline{\basefield_v} \xrightarrow{\raisebox{-1 em}{\smash{\ensuremath{\sim}}}} M_{\dR} \otimes_{\basefield} \overline{\basefield_v},
	\]
	where $G_v$ acts on $M_{B,v}$ and both copies of $\overline{\basefield_v}$
	\item For each prime number $p$ and $v \in S_p(\basefield)$, an isomorphism
	\[
	\iota_v \colon M_p \otimes_{\Qp} B_{\dR,v} \xrightarrow{\raisebox{-1 em}{\smash{\ensuremath{\sim}}}} M_{\dR} \otimes_{\basefield} B_{\dR,v}
	\]
	of $B_{\dR,v}$-vector spaces compatible with action of $G_v$ and the Hodge filtration.
\end{itemize}

A morphism $f \colon M \to N$ is a collection of morphisms of all associated modules (Betti, de Rham, $\ell$-adic), that preserve the associated structures (Hodge filtration, conjugation action, $G_{\basefield}$-action) and are compatible with all comparison isomorphisms. The tensor product is defined in an obvious way, making this category Tannakian.

\subsubsection{Fiber Functors}\label{sec:fiber_func_PSPM}

For each $v \in S_{\infty}(\basefield)$, we have a $\Qb$-linear fiber functor $\om^v = \om^{B,v}$ given by $M \mapsto M_{B,v}$. We also have a $\basefield$-linear fiber functor $\om^{\dR}$ sending $M$ to $M_{\dR}$, and for each prime number $\ell$, a $\Ql$-linear fiber functor $\om^{\et,\ell} = \om^{\ell}$ sending $M$ to $M_{\ell}$. The comparison isomorphisms give natural isomorphisms between all these fiber functors over the appropriate base extensions.

When we remember the Galois representation on $M_{\ell}$, we get a functor $\mathrm{real}_{\et,\ell} \colon \mathbf{PSPM}_\basefield(\Qb) \to \Rep_{\Qb_{\ell}}^{\mathrm{g}}(G_{\basefield})$.

%\subsubsection{}
%Given a finite subset $S \subseteq S_f(\basefield)$, we let $\mathbf{PSPM}_{\basefield,S}(\Qb)$ denote the full subcategory of $M \in \mathbf{PSPM}_\basefield(\Qb)$ for which $M_{\ell} \in \Rep_{\Ql}^{\mathrm{f,S}}(G_{\basefield})$ (c.f. \S \ref{sec:S-integral}) for all $\ell$.

\subsubsection{Definition}
\label{sec:SPM}

Following \cite[III.2.1.4]{FPR91}, we define a \emph{pre-motivic structure} $M$ to be an object $M \in \mathbf{PSPM}_\basefield(\Qb)$ along with an increasing finite filtration $\{W_n M\}_{n \in \Zb}$ in $\mathbf{PSPM}_\basefield(\Qb)$ such that
\begin{itemize}
	\item For each $v \in S_{\infty}(\basefield)$, the vector space $M_{B,v} \otimes_{\Qb} \Rb$ is a mixed Hodge structure over $\basefield_v$ when equipped with the weight filtration coming from $\{W_n M\}$, the Hodge filtration coming from $\iota_v$, and if $F_v \simeq \Rb$, the Frobenius action from $G_v$.
\end{itemize}

As in \cite[III.2.1.4]{FPR91}, we denote the category of such objects and morphisms in $\mathbf{PSPM}_\basefield(\Qb)$ that respect the filtration by $\mathbf{SPM}_\basefield(\Qb)$. As mentioned in loc.cit., it is a neutral $\Qb$-linear Tannakian category, with unit object the cohomology of $\Spec{\basefield}$.

By III.2.1.2 of loc.cit., all morphisms are strictly compatible with the weight filtration, so $\mathbf{SPM}_\basefield(\Qb)$ is a strictly filtered Tannakian category over $\coeff=\Qb$ with Hodge-filtered fiber functor $\om^{\dR}_H$.

\subsubsection{}\label{sssec:de_rham_filtered}

When we remember the filtration on $M_{\dR}$, we get a filtered\footnote{The subfunctors $F^p \om^{\dR}$ are exact by \cite[Th\'eor\`eme 2.3.5(iii)]{DeligneHodgeII}.} fiber functor (\ref{sssec:filtered_fiber_functor}) $\om^{\dR}_H \colon \mathbf{PSPM}_\basefield(\Qb) \to \FilVect_{\basefield}$.

\subsection{Variants of pre-motivic structures}\label{sec:variants_SPM}

We describe various Tannakian subcategories of $\mathbf{SPM}_\basefield(\Qb)$, most of which are known to (and all of which are expected to) contain realizations of motives.

\subsubsection{}\label{variant_hub}

We let $\mathbf{SPM}_\basefield(\Qb)^{\Hub}$ denote the subcategory of $M \in \mathbf{SPM}_\basefield(\Qb)$ such that
\begin{itemize}
	\item For each $\ell$, the $G_{\basefield}$-representation $\gr_n^W M_{\ell}$ is pure of weight $n$ at almost all unramified $v \in S_f(\basefield)$ in the sense of \cite{WeilII}.
\end{itemize}

%%% Huber says on p.87 that the Galois modules are all constructible, which means unramified almost everywhere

\subsubsection{}\label{variant_del}

\emph{Convention: Any appearance of `$\basefield,S$' as a subscript after $\mathbf{M}$ is interchangeable with `$\Oc_{\basefield,S}$' or `$Z$' where $Z=\Spec{\Oc_{\basefield,S}}$.}

For a subset $S \subseteq S_f(\coeff)$, we let $\mathbf{SPM}_{\basefield,S}(\Qb) = \mathbf{SPM}_{Z}(\Qb)$ denote the subcategory of $M \in \mathbf{SPM}_\basefield(\Qb)$ for which $M_{\ell} \in \Rep^{\rm{sf,S}}_{\Ql}(G_{\basefield})$ for all $\ell$, and
\[
\mathbf{SPM}_{\basefield}(\Qb)^{\Del} \coloneqq \bigcup_{\underset{\m{finite}}{S \subseteq S_f(\basefield)}}
\mathbf{SPM}_{\basefield,S}(\Qb).
\]

Note that $\mathbf{SPM}_{\basefield}(\Qb)^{\Del} \subseteq \mathbf{SPM}_{\basefield}(\Qb)^{\Hub}$.

\subsubsection{}\label{frob_cris_spm}

For any $\pf \in S_p(\basefield) \setminus S$, we have a Frobenius-equivariant fiber functor over $(\coeff'_0,\phi)=(\basefield_{\pf,0},\mathrm{Frob}_{\pf})$ on the category $\mathbf{SPM}_{\basefield,S}(\Qb)$. In fact, the union
\[
\mathbf{SPM}_{\basefield,\cri(\pf)}(\Qb) \coloneqq \mathbf{SPM}_{\basefield,S_f(\basefield)\setminus \{\pf\}}(\Qb) = \mathbf{SPM}_{\Oc_{K,\pf}}(\Qb)
\]
over all such $S$ has a Frobenius-equivariant fiber functor $\om^{\cri}$ over $(\basefield_{\pf,0},\mathrm{Frob}_{\pf})$ given by
\[
M \mapsto \om^{\cri}(M) \coloneqq M_{\cri} \coloneqq \operatorname{D}_{\cris}(M_p).\]
%\cite[1.4(AM4)]{Deligne89}

\subsubsection{}\label{variant_pol}

We let $\mathbf{SPM}_\basefield(\Qb)^{\pol}$ denote the full subcategory of $M \in \mathbf{SPM}_\basefield(\Qb)$ for which each $\gr^W_n M$ is \emph{polarizable of weight $n$} in the following sense:
\begin{itemize}
	\item There exists an isomorphism $\gr^W_n M \to (\gr^W_n M)^{\vee}(-n)$ in $\mathbf{SPM}_\basefield(\Qb)$ whose restriction to the Betti and de Rham realizations is a polarization of Hodge structures of weight $n$.
\end{itemize}

\subsubsection{}\label{tann_subcats}

We claim the categories $\mathbf{SPM}_{\basefield,S}(\Qb)^{\Del} \subseteq \mathbf{SPM}_{\basefield}(\Qb)^{\Del} \subseteq \mathbf{SPM}_{\basefield}(\Qb)^{\Hub}$ and $\mathbf{SPM}_\basefield(\Qb)^{\pol}$ are Tannakian subcategories of $\mathbf{SPM}_{\basefield}(\Qb)$.

For the first claim, note that the property of being unramified/crystalline and of weight $n$ at a particular place $v$ is stable under taking subobjects, quotients, tensor products, and duals, so this follows.

For the second, see \ref{proof:semisimplicity} for kernels and quotients. 
%note that a polarization restricts to a polarization of any subobject by the positive-definiteness. For quotients, note that the quotient is isomorphic to the orthogonal complement of the kernel under the polarization.
Finally, noting that the weight-graded pieces of a tensor product (resp. dual) are the tensor products (resp. duals) of the weight-graded pieces and that one can take the tensor product (resp. dual) of polarization(s), $\mathbf{SPM}_\basefield(\Qb)^{\pol}$ is stable under tensor product and dual.

\subsubsection{\yesrem}
One might wish to define the full subcategory $\mathbf{SPM}_\basefield(\Qb)^{\comp}$ of $M \in \mathbf{SPM}_\basefield(\Qb)$ for which $\{M_{\ell}\}_{\ell}$ form a (weakly) compatible system of $\ell$-adic representation in the sense of \cite[p.81]{TaylorGaloisReps}. Similarly, we could define $\mathbf{SPM}_\basefield(\Qb)^{\scomp}$ denote those for which $\{M_{\ell}\}_{\ell}$ is strongly compatible in the sense of loc.cit. Note that \cite[Conjecture 1.3(2)]{TaylorGaloisReps} says that $\comp=\scomp$.

Unfortunately, these subcategories are not abelian, since the kernel and cokernel of a map between compatible systems need not be compatible.\footnote{For example, if we take $M$ to be the motive of a degree $2$ extension of $\basefield$ but then twist the Betti-\'etale comparison isomorphisms for a single $\ell$ by the involution that switches the trivial and nontrivial characters, then there is a map from $M$ to itself whose kernel is not in $\mathbf{SPM}_\basefield(\Qb)^{\comp}$.} In particular, $\mathbf{SPM}_\basefield(\Qb)^{\comp}$ is not a Tannakian category, so we cannot apply the Nori formalism to it as in \S \ref{sec:proof_of_real}. Nonetheless, the realization of the motive of a single smooth proper variety is contained in $\mathbf{SPM}_\basefield(\Qb)^{\comp}$ (c.f. Remark \ref{rem:scomp}).

If we restrict to objects in $\mathbf{SPM}_\basefield(\Qb)^{\comp}$, all of whose subobjects in $\mathbf{SPM}_\basefield(\Qb)$ are also in $\mathbf{SPM}_\basefield(\Qb)^{\comp}$, then we get an abelian category. However, we don't know that the cohomology of a variety satisfies this additional property. Such a property would follow from the absolute Hodge conjecture and an argument as in \cite{SawinMO_Pure_motives_compatible_systems}.

\subsubsection{}

For $\alpha_1, \cdots,\alpha_m \in \{\Hub,\Del,\comp,\scomp,\pol\}$ and \newline $\beta \in \{\{\basefield\},\{\basefield,S\},\{\basefield,\cri(\pf)\}\}$, we let
\[
\mathbf{SPM}_{\beta}(\Qb)^{\alpha_1,\cdots,\alpha_m}  \coloneqq 
\mathbf{SPM}_{\beta}(\Qb) \cap
\bigcap_{j=1}^m \mathbf{SPM}_{\basefield}(\Qb)^{\alpha_j}.
\]

For example, we will show that motives have realizations in $\mathbf{SPM}_\basefield(\Qb)^{\Del,\pol}$.

\subsubsection{Comparison with the literature}
\label{sec:lit_comparison}

There is an obvious functor from $\mathbf{SPM}_\basefield(\Qb)$ to Jannsen's category $\underline{\mathrm{MR}}_k$ (\cite[Definition 2.1]{JannsenBook90}).

There is also a functor from $\mathbf{SPM}_\basefield(\Qb)^{\Hub}$ to Huber's $\mathcal{MR}$ (\cite[11.1]{HuberMMRealizDerivCat}), where $A_{\sigma,\ell} \coloneqq A_{\sigma} \otimes_{\Qb} \Ql$ and $A_{\sigma,\Cb} \coloneqq A_{\sigma} \otimes_{\Qb} \Cb$.
%%% Not expected to be an equivalence because Huber doesn't require any p-adic Hodge theory conditions.

Finally, there is are obvious functors from $\mathbf{SPM}_{\Qb}(\Qb)^{\Del}$ to Deligne's \emph{syst\`emes de r\'ealisations} (\cite[1.9]{Deligne89}) and from $\mathbf{SPM}_{\Qb,S}(\Qb)^{\Del}$ to those that are lisse over $\Spec{\Zb} \setminus S$ (\cite[1.15]{Deligne89}) by taking $M_{\cris,p} = \operatorname{D}_{\cris}(M_p)$ for $p$ such that $M_p$ is crystalline (resp. $p \notin S$).

None of these functors is expected to be an equivalence, as none of the target categories require the Galois representations to be de Rham.

%%% refer to Huber and Deligne, mention Jannsen (Are Huber and Deligne different in that there's absence of Frobenius at infinity?** Unless this is hiding in Huber and it has to do with $\sigma$ vs $\overline{\sigma}$??)

\subsection{Polarization and weight-filtered structures}\label{sec:pol_wt_fil}

\begin{ssProposition}
	\label{prop:semisimplicity}
	
	For $M \in \mathbf{SPM}_\basefield(\Qb)^{\pol}$, $\gr_n^W M$ is semisimple in $\mathbf{SPM}_\basefield(\Qb)$. In particular, $\mathbf{SPM}_\basefield(\Qb)^{\pol}$ is a weight-filtered Tannakian category over $\coeff=\Qb$.
	
	By \ref{rem:filtr_sub} and \ref{tann_subcats}, so are $\mathbf{SPM}_\basefield(\Qb)^{\pol,\Hub}$, $\mathbf{SPM}_\basefield(\Qb)^{\pol,\Del}$, and $\mathbf{SPM}_{\basefield,S}(\Qb)^{\pol}$.
\end{ssProposition}

%%% This is a weight-filtered Tannakian category (with Hodge-filtered fiber functor?). Key point is showing that weight-graded pieces are semisimple. (Also see how Ishai does Propositions - is the proof have another numbering? or the same?)

%% Mention Hodge-filtered fiber functor

\subsubsection{Proof}\label{proof:semisimplicity}

The proof is a modification of that of \cite[Lemma 1.1]{JannsenBook90}. Given $M \in \mathbf{SPM}_\basefield(\Qb)^{\pol}$ and $A \subseteq \gr_n^W M$ in $\mathbf{SPM}_\basefield(\Qb)$, we define $B \subseteq \gr_n^W M$ to be the orthogonal complement of $A$ under the polarization. Note that $B$ can also be defined by taking the cokernel $B'$ of
\[
A \to \gr^W_n M \to (\gr^W_n M)^{\vee}(-n)
\]
and then dualizing and twisting the map $(\gr^W_n M)^{\vee}(-n) \to B'$ to get an inclusion $B \to \gr^W_n M$, which shows that $B$ is an object of $\mathbf{SPM}_\basefield(\Qb)$. The polarization property implies by a standard argument that $A$ and $B$ are complementary, so that $\gr^W_n M \simeq A \oplus B$ in $\mathbf{SPM}_\basefield(\Qb)$.

Note that the polarization restricts to a pairing on any subobject in $\mathbf{SPM}_\basefield(\Qb)$, and this restriction is in fact non-degenerate (and also positive-definite) by the positive-definiteness of a polarization. In particular, $\gr^W_n M \simeq A \oplus B$ is in $\mathbf{SPM}_\basefield(\Qb)^{\pol}$ as well. This implies $\gr^W_n M$ is also semisimple in $\mathbf{SPM}_\basefield(\Qb)^{\pol}$, and therefore $\mathbf{SPM}_\basefield(\Qb)^{\pol}$ is a weight-filtered Tannakian category.

\subsubsection{}\label{sec:hodge_frob_spm}

Now that we have a weight-filtered Tannakian category, let us recall the extra structure on its fiber functors.

The category $\mathbf{SPM}_\basefield(\Qb)$ (and therefore any Tannakian subcategory) has a $\basefield$-linear Hodge-filtered fiber functor $\om^{\dR}_H$ (\ref{sssec:de_rham_filtered}). Let $\om^{\dR,\pf}_H$ denote its base-change along $\basefield \hookrightarrow \basefield_{\pf}$.

As noted in \ref{frob_cris_spm}, $\mathbf{SPM}_{\Oc_{\basefield,\pf}}(\Qb)$ (and therefore all its Tannakian subcategories) has a Frobenius-equivariant fiber functor $\om^{\cri}$ over $(\basefield_{\pf,0},\Frob_{\pf})$.

We set $\eta^{\mot} \coloneqq \eta^{\et} \circ (\iota_{\pf}^{-1})^{G_{\pf}}$, where
\[
(\iota_{\pf}^{-1})^{G_{\pf}}
\colon
M_{\dR} = (M_{\dR} \otimes_{\basefield} B_{\dR,\pf})^{G_{\pf}} \to (M_p \otimes_{\Qp} B_{\dR,v})^{G_{\pf}} = V^{\dR}.
\]
Then $\boxed{(\om_H,\om_{\phi},\eta)=(\om^{\dR,\pf}_H,\om^{\cri},\eta^{\mot})}$ is a \emph{Hodge-filtered Frobenius-equivariant fiber functor} over $\boxed{(\coeff',\coeff'_0,\phi) =(\basefield_{\pf},\basefield_{\pf,0},\Frob_{\pf})}$ on the \emph{strictly filtered Tannakian category} $\boxed{\Tc = \mathbf{SPM}_{\Oc_{\basefield,\pf}}(\Qb)}$ over $\boxed{\coeff = \Qb}$ and thus also on its \emph{weight-filtered Tannakian} subcategory $\boxed{\mathbf{SPM}_{\Oc_{\basefield,\pf}}(\Qb)^{\pol}}$.

\subsubsection{\somewhatrem}\label{rem:compat_hodge_frob}

The Hodge-filtered Frobenius-equivariant fiber functor of \ref{hodge_frob_rep} pulls back to that of \ref{sec:hodge_frob_spm} via the functor $\mathrm{real}_{\et,\ell} \colon \mathbf{SPM}_{\Oc_{\basefield,\pf}}(\Qb) \to \Rep_{\Qb_{\ell}}^{\mathrm{sf,\cri(\pf)}}(G_{\basefield})$. Note in particular that $\iota_{\pf}$ provides the isomorphism between $\om^{\dR} \circ \mathrm{real}_{\et,\ell}$ (\ref{hodge_rep},\ref{sec:fiber_func_PSPM}) and $\om^{\dR,\pf}$ (\ref{sec:hodge_frob_spm}).

\subsection{Motivic structures}\label{sec:motivic_structures}

While $\mathbf{SPM}_{\basefield,S}(\Qb)^{\Del,\pol}$ might seem like a good category to work with, the $\Ext$ groups are generally not finite-dimensional $\Qb$-vector spaces. This is in particular because the compatibility between $M_{\ell}$ for different $\ell$ sees only the characteristic polynomial, which does not see extensions. One could therefore cook up independent extension classes at many different $\ell$, creating an extension space that is at best finite rank over the adeles. This issue is discussed briefly in \cite[p.194]{HuberMMRealizDerivCat}.

To regulate this chaos, it seems necessary to restrict to systems of realizations that come from geometry (equivalently, from motives). In terms of $\Ext$ groups, this corresponds roughly to looking at the image
\[
K_m(Y)^{(n)} \to \prod_{v \in S_{\infty}(\basefield)} H^{2n-m}_{\mathcal{D}}(Y_v;\basefield_v(n)) \times \prod_{\ell} H^1_{\rm{g}}(G_{\basefield};H^{2n-m-1}_{\et}(\overline{Y};\Ql(m)))
\]
for a variety $Y/\basefield$ and $m>0$, which is a (conjecturally) finite-dimensional $\Qb$-vector space inside a space most naturally considered as a module over the adeles.

This is what \cite{FPR91} refers to as ``motivic systems'' and denotes $\mathbf{SM}_\basefield(\Qb)$. However, \cite[III.4.1.2]{FPR91} does not specify what this means precisely, and in particular, what a ``mixed motive'' is. They conjecture furthermore (\cite[III.4.1.5]{FPR91}) that such a category should induce isomorphisms (\ref{conj:FPR_adm}) of $\Ql$-vector spaces
\[
\Ext^1_{\mathbf{SM}_\basefield(\Qb)}(\Qb,M) \otimes_{\Qb} \Ql \xrightarrow{\sim} H^1_{\rm{g}}(G_{\basefield};M).
\]
This property is especially important for us because the dimensions of the latter (and thus by the conjecture, of the former) follow from the Bloch--Kato conjectures, as mentioned in \S \ref{sec:remarks_on_computation}.

We claim that ``mixed motive'' should be taken in the sense of Voevodsky's triangulated category $\DM_{\gm}(\basefield,\Qb)$ (\cite[\S 2.1]{VoevodskyTriCa00},\cite[11.1.10]{CisDegTriCat}) of \emph{geometric motives} and construct a realization functor from this category to the bounded derived category $D^b(\mathbf{SPM}_\basefield(\Qb)^{\Del,\pol})$ and define $\mathbf{SM}_\basefield(\Qb)$ as the subcategory generated by the image of $H^0$.

\subsubsection{Motivic Structures Associated to Pairs of Varieties}\label{sssec:motivic_associated_to_pairs}

As a preliminary step, we describe the object $h^i_{\mathbf{SPM}}(X,Y)$ of $\mathbf{SPM}_\basefield(\Qb)^{\Del,\pol}$ associated to the $i$th cohomology of a pair $X,Y$ of smooth varieties over $\basefield$.

We let $M(X,Y)$ denote the geometric motive given by the complex $0 \to Y \hookrightarrow X \to 0$ with $X$ in degree $0$. We get all the objects and most of the compatibilities and conditions defining $\mathbf{SPM}_{\basefield}(\Qb)^{\Hub}$ by applying the contravariant functor $R_{\mathcal{MR}}$ of \cite[Theorem 2.3.3]{HuberMMRealizDerivCat} to the Voevodsky motive $M(X,Y)$ and taking $H^{i}$.
%%% Note that Huber doesn't even assume the representation is de Rham at all places, but the next sentence deals with this. Consider making this a footnote. Or putting it in the section on comparison with the literature.
The one missing piece is the $p$-adic Hodge comparison and its compatibility with maps; but this follows by applying \cite[Theorem 1.2]{NizUnique20} to $M(X,Y)$.

Note that any pair $(X,Y)$ has a nice model over $\Oc_{\basefield,S}$ for sufficiently large $S$, which shows that $A^*(X,Y)$ in fact lies in $\mathbf{SPM}_{\basefield}(\Qb)^{\Del}$.

For the polarization, note that most of the data of the polarization comes from \cite[Proposition 21.2.3]{HuberMMRealizDerivCat}. The only missing part is the compatibility of the polarizations with $p$-adic comparison isomorphisms. The polarization on a smooth projective variety is constructed using the Lefschetz class of a hyperplane section and Poincar\'e duality. The former is compatible with the $p$-adic comparison isomorphism by \cite[(A1.3)]{TsujiSurvey02} and the latter by \cite[(A1.2),(A1.5)]{TsujiSurvey02}.

We remark that for $M=h^i_{\mathbf{SPM}}(X,Y)$, $M_{\dR}=H^i_{\dR}(X,Y)$, $M_{B,v} = H^i_{B,v}(X,Y;\Qb) \coloneqq H^i_{\mathrm{sing}}(X(\basefield_v),Y(\basefield_v);\Qb)$, and $M_{\ell} = H^i_{\et}(X_{\overline{\basefield}},Y_{\overline{\basefield}};\Ql)$.

\begin{ssTheorem}\label{thm:real}

There is a realization functor 
\[\mathrm{real}_{\mathbf{SPM}} \colon \DM_{\gm}(\basefield,\Qb) \to D^b(\mathbf{SPM}_\basefield(\Qb)^{\Del,\pol})
\]
such that
\[
H^{-i}(\mathrm{real}_{\mathbf{SPM}}(M(X,Y))) = h_i^{\mathbf{SPM}}(X,Y) \coloneqq h^i_{\mathbf{SPM}}(X,Y)^{\vee}.
\]

\end{ssTheorem}

\begin{ssDefinition}\label{defn:motivic_structures}
	We define the category $\mathbf{SM}_\basefield(\Qb)$ of \emph{motivic structures} as the Tannakian subcategory of $\mathbf{SPM}_\basefield(\Qb)^{\Del,\pol}$ generated by the essential image of \[H^0 \circ \mathrm{real}_{\mathbf{SPM}} \colon \DM_{\gm}(\basefield,\Qb) \to \mathbf{SPM}_\basefield(\Qb)^{\Del,\pol}.\] By \ref{rem:filtr_sub} and Proposition \ref{prop:semisimplicity}, it is a weight-filtered Tannakian category.
\end{ssDefinition}

We believe this is the \textit{meaning} of \cite[III.4.1.2]{FPR91} and for $\basefield=\Qb$ of \cite[Definition 1.11]{Deligne89}, updated to reflect the availability of $\DM_{\gm}(\basefield,\Qb)$. (See also the \textit{hope} expressed in the paragraph immediately preceding Definition 1.11 of loc.cit.)

\subsubsection{\yesrem}\label{rem:absolute_hodge_motives}

The category $\mathbf{SM}_\basefield(\Qb)$ is very close to various categories of mixed motives for absolute Hodge cycles: the category $\mathcal{MM}$ of \cite[22.1.3]{HuberMMRealizDerivCat}, the category $\underline{MM}_\basefield$ of \cite[Definition 4.1]{JannsenBook90}, and $\mathcal{MM}_{\mathrm{AH}}$ of \cite[Definition 6.3.11]{HuberMSBook17}. The main difference for us is that our category of mixed realizations includes $p$-adic Hodge comparison; we felt this was important both because our application (Chabauty--Kim) heavily uses $p$-adic Hodge theory and because it allows us to invoke Fontaine--Perrin-Riou's conjectures (\ref{conj:FPR_adm}) comparing with Bloch--Kato Selmer groups.

It's not clear that Jannsen's $\underline{MM}_\basefield$ includes realizations of all complexes of varieties; see the discussion in the middle of p.194 of \cite{HuberMMRealizDerivCat}. We want all realizations of complexes of varieties so as to be able to realize path torsors as in \ref{sec:path_torsor_SM}.

\subsubsection{Definition}\label{defn:SM_beta}

For $\beta \in \{\{\basefield\},\{\basefield,S\},\{\basefield,\cri(\pf)\}\}$, we let (c.f. \S \ref{ssec:notations_and_conventions}):
\[
\mathbf{SM}_{\beta}(\Qb) \coloneqq \mathbf{SM}_\basefield(\Qb) \cap \mathbf{SPM}_{\beta}(\Qb).
\]

Note then that by \ref{sec:hodge_frob_spm} and \ref{rem:hodge_frob_sub},
\[
\mathbf{SM}_{\basefield,\cri(\pf)}(\Qb) = \mathbf{SM}_{\Oc_{\basefield,\pf}}(\Qb)
\]
and its subcategories $\mathbf{SM}_{\basefield,S}(\Qb)$ for $\pf \notin S$ are weight-filtered Tannakian categories with Hodge-filtered Frobenius-equivariant fiber functor $(\om^{\dR,\pf}_H,\om^{\cri},\eta^{\mot})$ over $(\basefield_{\pf},\basefield_{\pf,0},\Frob_{\pf})$.

\subsubsection{Proof of Theorem \ref{thm:real}}\label{sec:proof_of_real}

We use \cite[Theorem 4.10]{DegNiz18} (described also in \cite[\S 3.2]{ChiarBruLaz19}), noting that the arguments work the same if we take $F$ to be any field of characteristic $0$. The theorem associates a canonical monoidal realization functor $R_A \colon \DM_{\gm}(\basefield,\Qb) \to D^b(\mathscr{A})$ to:
\begin{itemize}
	\item An $F$-linear Tannakian category $\mathscr{A}$ and fiber functor $\om \colon \mathscr{A} \to \Vect_F$
	\item A representation $A^* \colon \Delta_g \to \mathscr{A}$ for which $\om \circ A^*$ is the singular representation (see below).
\end{itemize}

Here, $\Delta_g$ is Nori's directed graph of good pairs, which is a localization of $\Delta_g^{\eff}$ at multiplication by $(\Gm,\{1\},1)$. The directed graph $\Delta_g^{\eff}$ has vertices $3$-tuples $(X,Y,i)$ with $X,Y$ quasi-projective (reduced but not necessarily integral) varieties over $\basefield$ such that $H^j_{\et}(\overline{X},\overline{Y};\Qp)=0$ for $j \neq i$. The edges are everything of the form either $(X,Y,i) \to (X',Y',i)$ for a morphism $f \colon (X,Y) \to (X',Y')$ of pairs or a boundary map $(Y,Z,i) \to (X,Y,i+1)$ for a chain $Z \subseteq Y \subseteq X$. It has an obvious product structure given by Cartesian product of pairs of varieties.

We take $F=\Qb$, $\mathscr{A}=\mathbf{SPM}_\basefield(\Qb)^{\Del,\pol}$ and wish to define $A^*$. As in \cite[Lemma 8.2.4]{HuberMSBook17}, it suffices to define $A^*$ on $\Delta_g^{\eff}$ in such a way that it is monoidal and such that $(\Gm,\{1\},1)$ maps to an invertible object.

We can define $A^*(X,Y,i) \coloneqq h^i_{\mathbf{SPM}}(X,Y)$ as in \ref{sssec:motivic_associated_to_pairs}. Monoidality follows from that of \cite[Theorem 2.3.3]{HuberMMRealizDerivCat} and by functoriality and compatible with cup products of $p$-adic period isomorphisms \cite[A.1.1-2]{TsujiSurvey02}. Finally, $h^1_{\mathbf{SPM}}(\Gm,\{1\})$ is invertible because $1 = \dim_{\Qb} H^1_{\dR}(\Gm,\{1\}) = \dim_{\Qb} H^1_{B,\sigma}(\Gm,\{1\};\Qb) = \dim_{\Ql} H^1_{\et}((\Gm)_{\overline{\basefield}},\{1\};\Ql)$, so we are done.

\subsubsection{\remish}\label{cor:real_sm}

Note that for any $(X,Y,i)$
\[
h^i_{\mathbf{SM}}(X,Y) \coloneqq h^i_{\mathbf{SPM}}(X,Y) = H^{-i}(\mathrm{real}_{\mathbf{SPM}}(M(X,Y)))^{\vee} \in \mathbf{SM}_\basefield(\Qb).
\]
The proof of Theorem \ref{thm:real} applied to $\mathbf{SM}_\basefield(\Qb)$ in place of $\mathbf{SPM}_\basefield(\Qb)$ shows that $\mathrm{real}_{\mathbf{SPM}}$ factors through a functor
\[
\mathrm{real}_{\mathbf{SM}} \colon \DM_{\gm}(\basefield,\Qb) \to D^b(\mathbf{SM}_\basefield(\Qb)).
\]

\subsubsection{Notation}\label{notation_tate_object}

We set
\[\Qb(1) = \Qb(1)_{\mathbf{SM}} \coloneqq h^1_{\mathbf{SM}}(\Gm).
\]
It satisfies $\Qb(1)_{\dR} = \gr_F^1 \Qb(1)_{\dR}$, its $\ell$-adic realization is the $\ell$-adic cyclotomic character for each $\ell$, and $\Qb(1)_{\dR} = 2 \pi i K \iota_v(\Qb(1)_{B,v}) \subseteq \Qb(1)_{\dR} \otimes_F \overline{F_v}$ for each $v \in S_{\infty}(\basefield)$.

We set $\Qb(n) \coloneqq \Qb(1)^{\otimes n}$, so that in particular $\Qb(0) = \mathbbm{1}_{\mathbf{SM}_{\basefield}(\Qb)}$ and $\Qb(-1) = \Qb(1)^{\vee}$.

For $M \in \mathbf{SM}_{\basefield}(\Qb)$, we set
\[
M(n) \coloneqq M \otimes \Qb(n).
\]

\subsubsection{Conjecture}\label{conj:FPR_adm}

In \cite[III.4.1.5]{FPR91}, Fontaine--Perrin-Riou conjecture that $\mathbf{SM}_{\basefield,S}(\Qb)$ is $S$-admissible. As described in \cite[III.3.3.1]{FPR91}, this says in particular that for $M \in \mathbf{SM}_{\basefield,S}(\Qb)$ and any prime number $\ell$, the natural map
\[
\Ql \otimes_{\Qb} \Ext^1_{\mathbf{SM}_{\basefield,S}(\Qb)}(\Qb(0),M) \to H^1_{f,S}(G_{\basefield};M_{\ell})
\]
is an isomorphism.

This would imply that $\dim_{\Qb} \Ext^1_{\mathbf{SM}_{\basefield,S}(\Qb)}(\Qb(0),M) = \dim_{\Ql} H^1_{f,S}(G_{\basefield};M_{\ell})$. As mentioned in \S \ref{sec:remarks_on_computation} and described in more detail in \cite[\S 2.7]{BettsCorwinLeonhardt}, it is possible to give conjectural formulas for these dimensions.

\subsubsection{\yesrem}\label{rem:scomp}

\cite[III.4.1.5]{FPR91} posits also that every object $M$ of $\mathbf{SM}_{\basefield,S}(\Qb)$ is ``L-admissible''. This means that in fact $\mathbf{SM}_{\basefield,S}(\Qb) \subseteq \mathbf{SPM}_\basefield(\Qb)^{\Del,\pol,\scomp}$.

What's known is that for a smooth proper variety $X/\basefield$, we have $H^i(\mathrm{real}_{\mathbf{SM}}(M(X))) \in \mathbf{SPM}_\basefield(\Qb)^{\Del,\pol,\comp}$. Replacing $\comp$ with $\scomp$ is Conjecture $C_{WD}(X,m)$ of \cite[2.4.3]{Fontaine94}.

%%% Motivation: FPR define this but don't give a real definition (see also Deligne's issue). We define this precisely, as we will need it. We also note conjecture on the Ext-groups.

\subsection{Maps from motivic cohomology and $K$-theory}\label{sec:sm_mot_coh}

Let $X$ be a smooth quasi-projective variety over $\basefield$. By \cite[Theorem 3.4]{LevineBlochRevisited} and \cite[Theorem 19.1]{LNMC2006}, we have isomorphisms
\begin{equation}\label{eqn:k_ch_mc}
	K_{2n-k}(X)^{(n)}_{\Qb} \xrightarrow[\raisebox{0.25 em}{\smash{\ensuremath{\sim}}}]{\che} \CH^{n}(X,2n-k)_{\Qb} \xrightarrow[\raisebox{0.25 em}{\smash{\ensuremath{\sim}}}]{\cyc_{\mot}} H^{k}_{\mot}(X;\Qb(n)).
\end{equation}

By \cite[p.197]{VoevodskyTriCa00}, we have an isomorphism
\begin{equation}\label{eqn:Hmot_ext_isom}
	H^{k}_{\mot}(X;\Qb(n)) \simeq \homo_{\DM_{\gm}(\basefield,\Qb)}(M(X),\Qb(n)[k]),
\end{equation}
hence an isomorphism with all of (\ref{eqn:k_ch_mc}).

Via the functor $\mathrm{real}_{\mathbf{SM}}$ of Remark \ref{cor:real_sm}, we have a map
\begin{equation}\label{eqn:DM_gm_to_SM_K}
	\homo_{\DM_{\gm}(\basefield,\Qb)}(M(X),\Qb(n)[k]) \to \homo_{D^b(\mathbf{SM}_\basefield(\Qb))}(\mathrm{real}_{\mathbf{SM}}(M(X)),\Qb(n)[k]).
\end{equation}

Taking $H^{-k}$ gives us a map
\[
\beta_k \colon \homo_{D^b(\mathbf{SM}_\basefield(\Qb))}(\mathrm{real}_{\mathbf{SM}}(M(X)),\Qb(n)[k]) \to \homo_{\mathbf{SM}_\basefield(\Qb))}(H^{-k}(\mathrm{real}_{\mathbf{SPM}}(M(X,Y))),\Qb(n)).
\]

Since $H^{-k}(\mathrm{real}_{\mathbf{SPM}}(M(X,Y))) \simeq h_{\mathbf{SM}}^k(X)^{\vee}$, this latter group is
\[
\homo_{\mathbf{SM}_\basefield(\Qb))}(h_{\mathbf{SM}}^k(X)^{\vee},\Qb(n))
\simeq
\homo_{\mathbf{SM}_\basefield(\Qb))}(\Qb,h_{\mathbf{SM}}^k(X)(n)).
\]

The latter is trivial unless $k=2n$, in which case it is the space of absolute Hodge cycles in degree $\coeff$. Either way, it is zero on $\CH^{n}_{\mathrm{hom}}(X,2n-k)_{\Qb}$, the (higher\footnote{Of course $\CH^n(X,m) = \CH^n_{\mathrm{hom}}(X,m)$ if $m>0$.}) Chow group of cycles homologically equivalent to zero. The fact that $\mathrm{real}_{\mathbf{SM}}(M(X))$ is supported in negative degrees means the spectral sequence computing $\homo_{D^b(\mathbf{SM}_\basefield(\Qb))}$ gives us a map
\[
\Ker{\beta_k} \to \Ext^1_{\mathbf{SM}_\basefield(\Qb))}(h_{\mathbf{SM}}^{k-1}(X)^{\vee},\Qb(n))
\simeq \Ext^1_{{\mathbf{SM}_\basefield(\Qb))}}(\Qb,h_{\mathbf{SM}}^{k-1}(X)(n))
\]
and thus a motivic Abel--Jacobi map
\begin{equation}\label{eqn:map_from_k_theory}
	\AJ_{\mot} \colon \CH^{n}_{\mathrm{hom}}(X,2n-k)_{\Qb} \to \Ext^1_{{\mathbf{SM}_\basefield(\Qb))}}(\Qb,h_{\mathbf{SM}}^{k-1}(X)(n)).
\end{equation}

The composition
\[\CH^{n}_{\mathrm{hom}}(X,2n-k)_{\Qb} \xrightarrow{\AJ_{\mot}} \Ext^1_{{\mathbf{SM}_\basefield(\Qb))}}(\Qb,h_{\mathbf{SM}}^{k-1}(X)(n)) \xrightarrow{\mathrm{real}_{\et}} H^1(G_{\basefield};H_{\et}^{k-1}(\overline{X};\Qp(n)))\]
is the \'etale Abel--Jacobi map $\AJ_{\et}$ induced from the \'etale cycle class map $\cyc_{\et} \colon \CH^{n}(X,2n-k)_{\Qb} \to H^k_{\et}(X;\Qp(n))$ and the Hochschild--Serre spectral sequence.

%below probably will be deleted
%The spectral sequence computing $\homo_{D^b(\mathbf{SM}_\basefield(\Qb))}$ gives us a map
%\[
%\homo_{D^b(\mathbf{SM}_\basefield(\Qb))}(\Qb,\mathrm{real}_{\mathbf{SM}}(M(X))[k](n)) \to \homo_{\mathbf{SM}_\basefield(\Qb))}(\Qb,h_{\mathbf{SM}}^k(X)(n)).
%\]

%Assuming $m \coloneqq 2n-k \neq 0$, this last group vanishes by weights. We therefore get a map $\homo_{D^b(\mathbf{SM}_\basefield(\Qb))}(\Qb,\mathrm{real}_{\mathbf{SM}}(M(X))[k](n)) \to \Ext^1_{\mathbf{SM}_\basefield(\Qb))}(\Qb,H^{k-1}(M(X))(n))$

% composed map
%\[
%K_{m}(X)^{(n)}_{\Qb} \simeq H^{k}_{\mot}(X;\Qb(n)) \simeq \homo_{\DM_{\gm}(\basefield,\Qb)}(\Qb,M(X)[k](n)) \to 
%\]

\subsection{Subcategories generated by a single variety}\label{sec:SM_A}

%###

The categories $\mathbf{SM}_{\basefield,S}(\Qb)$ and $\Rep_{\Qp}^{\mathrm{sf,S}}(G_{\basefield})^{\wss}$, although they have (conjecturally in the first case) finite-dimensional Ext-groups for $S$ finite, are nonetheless too big to describe explicitly. We need to restrict to objects whose weight-graded pieces are in the subcategory generated by the cohomology of a fixed variety. For Galois representations, this appears in \cite[\S 4]{CorwinMECK} for $A=E$ an elliptic curve and \cite[\S 4]{CorwinTSV} in general.

\subsubsection{}\label{sssec:SM_A}

Given a smooth projective variety $A/\basefield$, let $\mathbf{SM}_{\basefield}(\Qb,A)^{\sesi}$ (resp. $\Rep_{\Qp}^{\mathrm{sg}}(G_{\basefield},A)^{\sesi}$) denote the Tannakian subcategory of $\mathbf{SM}_{\basefield}(\Qb)$ (resp. $\Rep_{\Qp}^{\mathrm{sg}}(G_{\basefield})$) generated by $\{h^n_{\mathbf{SM}}(A)\}_{n \in \Zb}$ (resp. $\{H^n_{\et}(\overline{A};\Qp)\}_{n \in \Zb}$). Our main focus is the case where $A$ is an abelian variety, in which case $h^1_{\mathbf{SM}}(A)$ (resp. $H^1_{\et}(\overline{A};\Qp)$) generates $\mathbf{SM}_{\basefield}(\Qb,A)^{\sesi}$ (resp. $\Rep_{\Qp}^{\mathrm{sg}}(G_{\basefield},A)^{\sesi}$). We denote by
\[
\mathbf{SM}_{\basefield}(\Qb,A)
\]
(resp.
\[
\Rep_{\Qp}^{\mathrm{sg}}(G_{\basefield},A)
)\]
the subcategory of $M$ in $\mathbf{SM}_{\basefield}(\Qb)$ (resp. $\Rep_{\Qp}^{\mathrm{sg}}(G_{\basefield})$) whose weight-graded pieces lie in $\mathbf{SM}_{\basefield}(\Qb,A)^{\sesi}$ (resp. $\Rep_{\Qp}^{\mathrm{sg}}(G_{\basefield},A)^{\sesi}$).

\subsubsection{}\label{sssec:SM_A,S}

If $A$ has good reduction outside $S$, then $\mathbf{SM}_{\basefield,S}(\Qb,A)^{\sesi} \coloneqq \mathbf{SM}_{\basefield}(\Qb,A)^{\sesi}$ (resp. $\Rep_{\Qp}^{\mathrm{sf,S}}(G_{\basefield},A)^{\sesi} \coloneqq \Rep_{\Qp}^{\mathrm{sg}}(G_{\basefield},A)^{\sesi}$) is contained in $\mathbf{SM}_{\basefield,S}(\Qb)$ (resp. $\Rep_{\Qp}^{\mathrm{sf,S}}(G_{\basefield})$), and we let
\[
\mathbf{SM}_{\basefield,S}(\Qb,A)
\]
(resp.
\[
\Rep_{\Qp}^{\mathrm{sf,S}}(G_{\basefield},A)
)\]
denote the full subcategory on the intersection of $\mathbf{SM}_{\basefield}(\Qb,A)$ and $\mathbf{SM}_{\basefield,S}(\Qb)$ (resp. $\Rep_{\Qp}^{\mathrm{sg}}(G_{\basefield},A)$ and $\Rep_{\Qp}^{\mathrm{sf,S}}(G_{\basefield})$).

\subsubsection{}\label{sssec:SM_A_semisimplicity}

In general, every object of $\mathbf{SM}_{\basefield}(\Qb,A)^{\sesi}$ is semisimple by Proposition \ref{prop:semisimplicity}. The latter is true for $\Rep_{\Qp}^{\mathrm{sg}}(G_{\basefield},A)^{\sesi}$ if $A$ satisfies the Grothendieck--Serre Semisimplicity Conjecture (GSSC), as discussed in \cite[\S 4.1]{CorwinTSV}. This is true if $A$ is an abelian variety, by Faltings' work. Only under this assumption do we have $\Rep_{\Qp}^{\mathrm{sg}}(G_{\basefield},A)^{\sesi} \subseteq \Rep_{\Qp}^{\mathrm{sg}}(G_{\basefield})^{\wss}$ and thus $\Rep_{\Qp}^{\mathrm{sg}}(G_{\basefield},A) \subseteq \Rep_{\Qp}^{\mathrm{sg}}(G_{\basefield})^{\wss}$.

%Note that $\mathbf{SM}_{\basefield,S}(\Qb,A)$ is the thick subcategory of $\mathbf{SM}_{\basefield,S}(\Qb)$ generated by $\mathbf{SM}_{\basefield,S}(\Qb,A)^{\sesi}$.

When working with motivic structures (resp. if $A$ satisfies the GSSC), $\mathbf{SM}_{\basefield,S}(\Qb,A)$ (resp. $\Rep_{\Qp}^{\mathrm{sf,S}}(G_{\basefield},A)$) is the thick subcategory of $\mathbf{SM}_{\basefield,S}(\Qb)$ (resp. $\Rep_{\Qp}^{\mathrm{sf,S}}(G_{\basefield})^{\wss}$) generated by $\mathbf{SM}_{\basefield,S}(\Qb,A)^{\sesi}$ (resp. $\Rep_{\Qp}^{\mathrm{sf,S}}(G_{\basefield},A)^{\sesi}$).

The categories $\mathbf{SM}_{\basefield}(\Qb,A)$ and $\Rep_{\Qp}^{\mathrm{sg}}(G_{\basefield},A)$ have the important property that their semisimple subcategory has a tensor generator (\ref{sec:algebraic_tannakian}), which implies that the maximal reductive quotient of their fundamental groups are algebraic groups (i.e., finite-type as schemes over $\Qb$ or $\Qp$, respectively).

\subsubsection{}\label{sec:THE_CATS_FOR_CK}

If $A$ has good reduction at $\pf \in S_p(\basefield)$, we let
\[
\mathbf{SM}_{\Oc_{\basefield,\pf}}(\Qb,A)
\]
(resp.
\[
\Rep_{\Qp}^{\mathrm{sf,\cri(\pf)}}(G_{\basefield},A)
)\]
denote the full subcategory on the intersection of $\mathbf{SM}_{\basefield}(\Qb,A)$ and $\mathbf{SM}_{\Oc_{\basefield,\pf}}(\Qb)$ (resp. $\Rep_{\Qp}^{\mathrm{sg}}(G_{\basefield},A)$ and $\Rep_{\Qp}^{\mathrm{sf,\cri(\pf)}}(G_{\basefield})$). As in \ref{defn:SM_beta}, \ref{sec:hodge_frob_spm} and \ref{rem:hodge_frob_sub} imply these categories have Hodge-filtered Frobenius-equivariant fiber functor $(\om^{\dR,\pf}_H,\om^{\cri},\mathrm{id})$ over $(\basefield_{\pf},\basefield_{\pf,0},\Frob_{\pf})$. By \ref{rem:hodge_frob_sub}, we get the same on its weight-filtered Tannakian subcategory $\mathbf{SM}_{\basefield,S}(\Qb,A)$ (resp. $\Rep_{\Qp}^{\mathrm{sf,S}}(G_{\basefield},A)$).

\subsection{Mixed motives}\label{sec:mixed_motives}

Another example of a weight-filtered Tannakian category is expected to be provided by the conjectural abelian category of mixed motives. We discuss this in more detail in Appendix \ref{app:mixed_motives} and explain how conjectured properties of this category imply that it is a weight-filtered Tannakian category in our sense.

%%%%%%%%%%%%%%
\section{Construction of Arithmetic Hodge Path}\label{sec:arith_hodge_path}
%%%%%%%%%%

\subsection{Introduction}

The goal of this section is to show (Theorem \ref{thm:hodge-filtered_weight-splitting}) that a weight-filtered Tannakian category $\Tc$ over $\coeff$ with Hodge-filtered fiber functor $\om_H$ over an extension $\coeff'/\coeff$ has a functorial splitting of the weight filtration over $\coeff'$ that respects the Hodge filtration. Such a splitting is given by a $\otimes$-isomorphism
\[
\Psi \colon  \om^{\gr} \coloneqq \om \circ \gr^W \xrightarrow{\sim} \om
\]
which we call an \emph{arithmetic Hodge path}. The condition of respecting the Hodge filtration says that for $M \in \Tc$, $\Psi_M^{-1}$ sends
\[
F^i \om(M)
\]
to
\[
\bigoplus_{n \in \Zb} \Im(F^i \om(M) \cap \om(W_n M) \to \om(\gr_n^W M)) = F^i\om(\gr_n^W M),
\]
the latter equality following from \ref{lemma:WH}.

\subsubsection{Overview of \S \ref{sec:arith_hodge_path}}\label{sssec:overview_of_hodge_path_section}

Most of \S \ref{sec:arith_hodge_path} is devoted to explaining what is meant by such a splitting; this is ultimately defined as a \emph{Hodge-filtered weight-splitting} in Definition \ref{hodge-filtered_weight_splitting}.

More precisely, in \S \ref{ssec:some_linear_algebra} have a mostly informal discussion of the concept of splitting of a filtration in the setting of vector spaces. In \S \ref{grfil_in_families}-\ref{sec:graded_galois_groups_diagram}, we discuss technicalities on graded and filtered fiber functors. In \S \ref{defn:weight_filtered_fib}-\ref{example:E}, we spell out that formalism more explicitly in the still-abstract setup from \S \ref{sec:wt_fil_tann} of a strictly-filtered Tannakian category with a Hodge-filtered fiber functor.

Finally, in \S \ref{prop:U^W_H_is_prounip}-\ref{thm:hodge-filtered_weight-splitting}, we prove the main theorem, which is a multi-step reduction to a form of Levi's Theorem.

\subsection{Some linear algebra}\label{ssec:some_linear_algebra}

We fix a field $\coeff$ of characteristic $0$. We begin with a discussion of some elementary linear algebra. Then, starting in \ref{grfil_in_families}, we will upgrade this discussion in two ways: (1) Instead of working with vector spaces (and their subspaces) over $\coeff$ we will work with finite locally free $\Oo_S$-modules (and locally split submodules) for $S$ an arbitrary $\coeff$-scheme. (2) Instead of working with individual $\Oo_S$-modules we will work with a neutral Tannakian category over $\coeff$ and fiber functors with various extra structures (filtrations, gradings) over $S$.

\subsubsection{}
By a \emph{filtered vector space} we mean an object of $\FilVect_{\coeff} $ as in \ref{Vect_FilVect}, i.e., a finite-dimensional vector space with a finite (equivalently, separated and exhaustive) filtration. We may use an increasing filtration (usually denoted $W$) or a decreasing filtration (usually denoted $F$). As in \ref{rem:increasing_vs_decreasing}, we may convert an increasing filtration $W$ to a decreasing filtration $F$ by setting $F^i = W_{-i}$. This operation defines an isomorphism of symmetric monoidal categories. In \ref{30522A}-\ref{defn:weight_filtered_fib}, we will for the most part work agnostically with \emph{filtered vector spaces} without specifying if our filtrations are increasing or decreasing. 

\subsubsection{}\label{30522A}
We distinguish two notions of \emph{splitting} of a filtration. If $V$ is a filtered vector space, a \emph{parabolic splitting} is an isomorphism of filtered vector spaces 
\[
\tag
{P}
\fil \gr V \overset {\Psi} \simeq V.
\]
A \emph{unipotent splitting} is a linearly independent collection $\{V_i\}_i $ of subspaces $V_i \subset V$ (i.e. for vectors $v_i \in V_i$, $\sum_i v_i = 0$ implies $v_i = 0$ for all $i$) such that the filtration is equal to the filtration induced by $\{V_i\}_i $. 

Every unipotent splitting $\{V_i\}_i$ gives rise to a parabolic splitting given by
\[
\tag{U}
V = \sum_i V_i \simeq \bigoplus _i V_i.
\]
This map from unipotent splittings to parabolic splittings is evidently injective with image given by those filtered isomorphisms $\Psi$ as in (P) for which the composition
\[
\gr V
\xto[\ref{defn:functors_fil_gr}(*)]{\sim}
\gr \fil \gr V 
\xto[{\gr \Psi}]{\sim}
\gr V 
\]
is $\id_{\gr V}$.

Conversely, a parabolic splitting $\Psi$ gives rise to a unipotent splitting given by setting
\[
V_i \coloneqq \Psi(\gr_i V).
\]

For example, if $V$ is concentrated in filtered degree $0$, then every automorphism of $V$ is a parabolic splitting, but there's only one unipotent splitting, given by $V_0\coloneqq V$.

\subsubsection{}\label{11}
We let $\GrVect_{\coeff} \coloneqq \Gr^{\fin}(\Vect_{\coeff})$ in the sense of \ref{defn:gr_fin}. When the field $\coeff$ stays fixed, we allow ourselves to drop the subscript `$k$'.

By a \emph{graded and filtered vector space} we mean a vector space equipped with a grading and a filtration subject to \textit{no} compatibility conditions. A graded and filtered vector space is a filtered object of the category of graded vector spaces if and only if it is a graded object of the category of filtered vector spaces, as we now explain.

Suppose $V = \bigcup_{i} F^iV$ is a filtered object of $\GrVect$. This means we're given gradings
\[
F^iV = \bigoplus_j \gr^W_j F^iV
\]
preserved by the inclusions
\[
F^i V \subset F^{i+1}V.
\]

Setting $\gr^W_j V \coloneqq \bigcup_i \gr^W_j F^i V$ gives a graded and filtered vector space structure on $V$.

We may retrieve the filtered object of $\GrVect$ from the underlying graded and filtered vector space via
\[
\gr^W_j F^iV = \gr^W_jV \cap F^iV.
\]
Furthermore, a graded and filtered vector space comes from a filtered object of $\GrVect$ if and only if for each $i$, the subspaces
\[
\gr^W_j V \cap F^iV \subset F^iV
\]
span $F^iV$.

Suppose now that
\[
\tag{*}
V = \bigoplus_j \gr_j^W V  = \bigoplus_j V_j
\]
is a graded object of $\FilVect$. This means we're given (separated, exhaustive) filtrations
\[
V_j = \bigcup_i F^i  (V_j)
\]
which are subject to \textit{no} conditions. We may form the \emph{underlying} vector space with filtration and grading by endowing $V$ with the direct sum filtration
\[
F^i V = \bigoplus_{j} F^i ( V_j).
\]
We then have
\[
F^iV \cap V_j = 
\left(
\bigoplus_j F^iV_j
\right) 
\cap V_j
= F^i(V_j)
\]
which shows that the graded object (*) is uniquely determined by the underlying filtered and graded vector space. Furthermore, a graded and filtered vector space 
\[
\bigcup_i F^i V = V = \bigoplus_j V_j
\]
comes from a graded object of $\FilVect$ if and only if the filtration \textit{is} the direct sum filtration associated to the induced filtrations 
\[
V_j = \bigcup_i ( F^iV \cap V_j)
\]
of the graded pieces. To conclude:

\begin{ssLemma}\label{lemma:graded_and_filtered}
	
	A graded and filtered vector space
	%\comment{ID: I've added a definition at the beginning of the paragraph. DC: Good! I still think it might be good to separate that out as a definition and then make the proof after the statement of the lemma. ID: I prefer not to, since any way we eventually need locally free modules over rings.} 
	determines a graded object of $\FilVect_{\coeff}$ if and only it it determines a filtered object of $\GrVect_{\coeff}$. 
	
\end{ssLemma}

\subsubsection{Proof of Lemma \ref{lemma:graded_and_filtered}}
The conditions described in \ref{11} are tautologically the same.

\subsubsection{Definition}
\label{defn:compat_gr_fil}
We will say the filtration and the grading are \emph{compatible} if they satisfy these conditions.

\subsubsection{Corollary}\label{cor:graded_filtered_equivalent}
Let $\bigcup F^iV = V = \bigoplus V_j$ be a graded and filtered vector space. The following conditions are equivalent:
\begin{enumerate}
	\item The filtration and grading are compatible.
	\item
	For each $i$, the intersections $V_j \cap F^iV$ span $F^iV$. 
	\item
	The filtration $F$ of $V$ is the direct sum of the filtrations
	\[
	V_j = \bigcup_i (V_j \cap F^iV)
	\]
	of the subspaces $V_j$.
	\item
	The $\Gm$-action associated to the grading preserves the filtration. 
\end{enumerate}
Indeed, as noted above, (2) and (3) are tautologically the same. On the other hand, (4) is clearly equivalent to the condition that $V$ determine a filtered object of $\Rep \Gm \simeq \GrVect$ which, as noted above, is equivalent to the other conditions.

\subsubsection{}
\label{hodge_on_grW}
Let $V$ be a Hodge and weight filtered vector space (i.e. equipped with an increasing filtration $W$ and a decreasing filtration $F$ both assumed finite). The Hodge filtration on $V$ endows $\gr^W(V)$ with a Hodge filtration given by defining $F^i \gr^W_j V$ to be the image of 
\[
\tag{*}
F^iV \cap W_jV 
\subset
W_jV \surj \gr^W_j V
\]
and then taking the direct-sum filtration. By construction, the weight grading and Hodge filtration on $\gr^W V$ are compatible.

\begin{ssLemma}\label{lemma:preserve_hodge}
	A unipotent splitting $\Psi \colon \gr^W V \simeq V$ of the weight filtration preserves Hodge filtrations if and only if the associated representation
	\[
	\chi \colon \Gm \to \GL V
	\]
	lands in $P(F) \coloneqq \{g \in \GL V \, \mid \, g(F^i V) = F^i V \, \forall \, i \in \Zb\}$.
	
\end{ssLemma}

\subsubsection{Proof of Lemma \ref{lemma:preserve_hodge}}
%Let $V_i$ be the image of $\gr^W_iV$ in $V$. Assume first that $\Psi$ preserves Hodge filtrations. For each $i$, the ismorphism $V_i \simeq \gr^W_iV$ induces a filtration on $V_i$ and the Hodge filtration on $V$ is the direct sum filtration associated to these. Let $v \in F^iV$ and $z$ a point of $\Gm$.

By \ref{cor:graded_filtered_equivalent}, $\chi$ lands in $P(F)$ iff the grading on $V$ determined by $\Psi$ is compatible with the Hodge filtration on $V$.

Since the Hodge filtration on $\gr^W V$ is compatible with the grading, this compatibility holds for $V$ if $\Psi$ preserves the Hodge filtration.

Conversely, suppose the grading associated to $\Psi$ and the Hodge filtration on $V$ are compatible. Let $V_j \coloneqq \Psi(\gr^W_j V)$, and for $\alpha \in V$, let $\alpha_j$ denote the projection of $\alpha$ onto $V_j$, so that $\alpha_j \in V_j$ for all $j$, and $\alpha = \sum_j \alpha_j$.

Suppose $\alpha \in F^i V$. The compatibility implies that $\alpha_j \in F^i V$ for all $j$. Thus $\alpha_j \in F^i V \cap V_j \subseteq F^i V \cap W_j V$.

The subspace $F^i \gr^W_j V$ is defined as the image of $F^i V \cap W_j V$ under the projection $W_j V \twoheadrightarrow W_j V/W_{j-1} V = \gr^W_j V$, so the image of $\alpha_j \in W_j V$ is in $F^i \gr^W_j V$.

But the composition $V_j \subseteq W_j V \twoheadrightarrow W_j V/W_{j-1} V = \gr^W_j V$ is the same as $\restr{\Psi^{-1}}{V_j}$. Thus $\Psi^{-1}(\alpha_j) \in F^i \gr^W_j V \subseteq F^i \gr^W V$.

It follows that
\[
\Psi^{-1}(\alpha) = \Psi^{-1}\left(\sum_j \alpha_j\right) = \sum_j \Psi^{-1}(\alpha_j) \in F^i \gr^W V.
\]

In other words, $\Psi^{-1}$ sends $F^i V$ into $F^i \gr^W V$. Since
\[
\dim F^i V = \sum_j \dim (F^i V \cap W_j V)/(F^i V \cap W_{j-1} V) = \sum_j \dim F^i \gr^W_j V = \dim F^i \gr^W V,
\]
$\Psi^{-1}$ induces an isomorphism $F^i V \xto{\sim} F^i \gr^W V$, so that $\Psi(F^i \gr^W V) = F^i V$, and we are done.

\subsection{Gradings and filtrations in families}\label{grfil_in_families}
%Our construction of arithmetic Hodge paths takes place in the axiomatic setting of \textit{weight-filtered Tannakian categories} and \textit{Hodge-filtered fiber functors}. Our treatment is based heavily on the work of Ziegler \cite{ZieglerGrFil15}, and we refer to the latter for a through discussion of foundational issues. 

We recall some notation from \cite[\S 3]{ZieglerGrFil15}. If $S$ is a scheme, we define $\LF(S)$ to be the category of finite locally free $\Oo_S$-modules (\cite[\href{https://stacks.math.columbia.edu/tag/01C6}{Tag 01C6}(2)]{stacks-project}), also known as `vector bundles'.

\subsubsection{Notation}

We let $\GrQCoh(S)$ and $\FilQCoh(S)$ denote $\Gr(\QCoh(S))$ and $\Fil^{\sex}(\QCoh(S))$, respectively.

\subsubsection{Definition}\label{gradings_in_families} We define $\GrLF(S)$ to be the subcategory of $\Mm \in \Gr(\LF(S))$ for which $\bigoplus_{i \in \ZZ} \Mm_i \in \LF(S)$.\footnote{This can be larger than $\Gr^{\fin}(\LF(S))$, for example if $S$ has infinitely many connected components.} By \cite[Remark 3.2]{ZieglerGrFil15}, $\GrLF(S) = \GrQCoh(S)_{\rig}$ (c.f. \S \ref{ssec:notations_and_conventions}).

%define $\GrLF(S)$ to be the category of locally free finite rank $\Oo_S$-modules $\Mm$ together with a decomposition 
%\[
%\Mm = \bigoplus_{n \in \ZZ} \Mm_n
%\]
%into locally free finite rank submodules (see \S 3.1 and Remark 3.2 of loc.cit.). We equip $\GrLF(S)$ with the evident symmetric monoidal structure for which $(\Mm \otimes \Nc)_n = \bigoplus_{a+b=n} \Mm_a \otimes \Nc_b$.

\subsubsection{Definition}
\label{filtrations_in_families}
As in \cite[\S 4.1]{ZieglerGrFil15}, we define $\FilLF(S)$ to be the full subcategory of $\Mm \in \Fil^{\sex}(\LF(S))$ for which the inclusion $\mathrm{fil}^n(\Mm) \hookrightarrow \Mm$ is locally split\footnote{For $S$ affine, it is in fact split, so $\fil(\gr(M)) \simeq M$ for all $M \in \FilLF(S)$.} for all $n$. By \cite[Lemma 4.2]{ZieglerGrFil15}, $\FilLF(S) = \FilQCoh(S)_{\rig}$. We set $\FilLF^{\fin}(S) \coloneqq \Fil^{\fin}(\LF(S)) \cap \FilLF(S)$.

\subsubsection{}\label{grfil_functors_S}

We have the exact tensor functors $\fil$ and $\gr$\footnote{Since $\LF(S)$ is not even pre-abelian, $\gr$ lands a priori in $\GrQCoh(S)$; however, the condition of being locally split ensures the graded piece are in $\LF(S)$.} as in \ref{defn:functors_fil_gr} between the categories $\GrLF(S)$ and $\FilLF(S)$. We also have the usual functors $\forg_{\fil}, \fil^n, \gr^n$ (\ref{defn:forg_fil_gr}) and $\forg_{\gr}$ (\ref{gr_forg}) landing in $\LF(S)$.

\subsubsection{}\label{LF_Vect}

Note that for a field $F$, $\LF(F) \coloneqq \LF(\Spec{F}) = \Vect_F$. Similarly, $\FilLF(F) = \FilVect_F$, and $\GrLF(F) = \GrVect_F$.

%There are evident exact tensor functors
%\[
%\FilLF(S)
%{\underset{\gr}
	%{\overset{\m{fil}}
		%{\leftrightarrows}}
	%\GrLF(S).
	%}
%\]
%For instance, if our filtrations are decreasing and $\Mm = \bigoplus_{n \in \ZZ} \Mm_i$ is an object of $\GrLF(S)$, then the associated filtration is given by
%\[
%F^n \Mm = \bigoplus_{i \ge n }\Mm_i.
%\]

\begin{sDefinition}\label{defn:graded_filtered_fib_func} \emph{For the rest of \S \ref{sec:arith_hodge_path}, let $\Tc$ be a Tannakian category over $\coeff$.}
	
	A \emph{graded (resp. filtered) fiber functor over $S$} is an exact $\coeff$-linear tensor functor 
	\[
	\phi \colon \Tc \to \GrQCoh(S) \hspace{5mm} (\mathrm{resp.  } \, \phi \colon \Tc \to \FilQCoh(S))
	\]

\subsubsection{}\label{sssec:grfil_values_in_rig} As in \cite[Remarks 3.2, 4.5]{ZieglerGrFil15},  any graded (resp. filtered) fiber functor takes values in $\GrLF(S)$ (resp. $\FilLF(S)$) by \cite[2.7]{DelTann90}. C.f. the footnote in \ref{sec:fiber_functors}.
	
	%\subsubsection{Filtered Fiber Functors}\label{sec:filtered_fiber_functors}
	
	%\subsubsection{Definition} As in \ref{defn:strict_exact}, a sequence $0 \to A \to B \to C \to 0$ in $\FilLF(S)$ is \emph{exact} if it is exact as a sequence of $\Oo_S$-modules and if all morphisms are strict.
	
	%\subsubsection{Remark} By \cite[Lemma 4.1]{ZieglerGrFil15}, a sequence $0 \to A \to B \to C \to 0$ is exact iff $0 \to F^n A \to F^n B \to F^n C \to 0$ is exact for all $n \in \Zb$. The proof, while not given in loc.cit., is similar to that of Definition/Theorem \ref{dfthm:strict}.
	
	%\subsubsection{Definition}\label{defn:filtered_fiber_functor}
	%Let $\Tc$ be a Tannakian category over $\coeff$ and $S$ a $\coeff$-scheme. A \emph{filtered fiber functor} $\phi$ over $S$ is a $k$-linear exact $\otimes$-functor
	%\[
	%\phi \colon \Tc \to \FilLF(S).
	%\]
	
	%\subsubsection{Remark}

\subsubsection{}\label{sssec:fil_gr_of_fiber_functors}

If $\phi$ is a graded (resp. filtered) fiber functor, then $\fil \circ \phi$ (resp. $\gr \circ \phi$) is filtered (resp. graded).

If $\phi$ is graded, then by \ref{sssec:grfil_id},
\[
\gr \circ \fil \circ \phi \simeq \phi
\]
canonically.

\end{sDefinition}

\subsection{Subgroups associated to a filtration}
\label{sec:subgroups_filtration}

In \S \ref{sec:subgroups_filtration}-\ref{prop:ranks_equal}, we discuss some general results and definitions regarding filtered fiber functors. While our notation mostly reflects decreasing filtrations, we recall as in \ref{rem:increasing_vs_decreasing} that results and definitions apply equally in either case, and we will apply this material both to (decreasing) Hodge-filtered fiber functors $\om_H$ (Definition \ref{defn:hodge_filtr}) and (increasing) weight-filtered fiber functors $\om^W$ (Definition \ref{defn:weight_filtered_fib}).

\subsubsection{}\label{base_change_filt_fib} Given $f \colon S' \to S$ over $\coeff$, a graded (resp. filtered) fiber functor $\phi$ over $S$ may be base-extended to a graded (resp. filtered) fiber functor over $S'$ in an evident way:
\[
(\phi)_{S'}(E) = f^*\phi(E)
\quad
\text{for all }
E \in \Tc.
\]

\subsubsection{}\label{isom_filt_grad_fib} For any two graded (resp. filtered) fiber functors $\phi_1$, $\phi_2$ on $\Tc$ over $S$, we let $\underline{\Isom}^\otimes(\phi_1, \phi_2)$ be the presheaf
\[
\m{Sch}_{/S} \to \m{Sets}
\]
which sends $S' \to S$ to the set of tensor isomorphisms 
\[
(\phi_1)_{S'} \xrightarrow{\sim} (\phi_2)_{S'}.
\]
By \cite[Theorems 3.11, 4.7]{ZieglerGrFil15}, $\underline{\Isom}^\otimes(\phi_1, \phi_2)$ is representable by an affine $S$-scheme, which is flat in the graded case.

\subsubsection{}\label{defn:subgroups_filtration} Let $\phi \colon \Tc \to \FilVect_{S}$ be a filtered fiber functor, and let $\phi^{\gr} \coloneqq \gr \circ \phi \colon \Tc \to \GrVect_{S}$. Following Ziegler \cite[Definition 4.7]{ZieglerGrFil15}, we define the associated \emph{parabolic Galois group} $P(\phi)$ (in terms of points with values in an arbitrary $S$-scheme $f \colon S' \to S$) by 
\[
P(\phi)(S') \coloneqq \Aut^\otimes(f^* \circ \phi),
\]
and the \emph{Levi Galois group} $L(\phi)$ by
\[
L(\phi)(S') \coloneqq \Aut^\otimes(f^* \circ \phi^{\gr}) = \Aut^\otimes(\gr \circ f^* \circ \phi).
\]

\subsubsection{}\label{defn:subgroups_filtration2}

Taking the associated graded of an isomorphism of filtered objects provides a morphism
\[
\m{ss} \colon P(\phi) \to L(\phi),
\]
and we define the \emph{unipotent Galois group} as
\[
U(\phi) = \Ker(\m{ss});
\]
By \ref{isom_filt_grad_fib}, these groups are all representable by affine $S$-groups. 

\subsection{Splittings of filtered fiber functors}
\label{sec:splittings_fib}
Let
\[
\phi \colon \Tc \to \FilLF(S)
\]
be a filtered fiber functor over a $\coeff$-scheme $S$. We use the $\otimes$-functors $\m{fil}$ and $\gr$ to define two notions of splitting. %Note that $\gr \circ \m{fil}$   is canonically isomorphic to the identity functor, while $\m{fil} \circ \gr$ is not.

A \emph{parabolic splitting} of $\phi$ is a $\otimes$-isomorphism of filtered fiber functors 
\[
\tag{P}
\m{fil} \circ \phi^{\gr}
\overset {\Psi} {\simeq} 
\phi
\]
A \emph{unipotent splitting} (also simply a \emph{splitting} as in \cite[Definition 4.9(i)]{ZieglerGrFil15}) is a graded fiber functor
\[
\ga \colon \Tc \to \GrLF(S)
\]
such that
\[
\phi = \m{fil} \circ \ga.
\]

\subsubsection{}\label{sssec:unip_sec_gives_par}

Equivalently, a unipotent splitting is an equivalence class of parabolic splittings modulo the left action of $L(\phi)$ acting by automorphisms of filtered fiber functors on $\phi^{\gr}$, as we now explain.

A parabolic splitting $\Psi$ gives rise to a unipotent splitting $\ga$ given by $\ga(E) \coloneqq \phi(E)$ with grading given by
\[
\gr_i \ga(E) \coloneqq  \Psi_E ( \gr_i \phi^{\gr}(E) ).
\]

If we precompose $\Psi$ with an automorphism of $\phi^{\gr}$, the resulting graded fiber functor $\gamma$ is the same. Thus we get a map from parabolic splittings modulo $L(\phi)$ to unipotent splittings.

If parabolic splittings $\Psi,\Psi'$ give rise to the same $\ga$, then $\Psi^{-1} \circ \Psi'$, a priori only an automorphism of $\m{fil} \circ \phi^{\gr}$, in fact preserves the grading and is thus an automorphsim of $\phi^{\gr}$; i.e., an element of $L(\phi)$. Thus this map is injective.

To see that the map is surjective, note that a unipotent splitting gives rise to a parabolic splitting given by the isomorphism of filtered fiber functors
\[
\tag{U}
\m{fil} \circ \phi^{\gr}
=
\m{fil} \circ \gr \circ \phi
=
\m{fil} \circ \gr \circ \m{fil} \circ \ga
\overset {\ref{sssec:fil_gr_of_fiber_functors}} \simeq
\m{fil} \circ \ga
=
\phi.
\]

Thus the concept of a \emph{splittable} filtered fiber functor (\cite[Definition 4.9(ii)]{ZieglerGrFil15}) is independent of which notion of splitting we choose.

\subsubsection{}\label{sssec:section_par_unip_split}

The proof of surjectivity in \ref{sssec:unip_sec_gives_par} in fact provides a canonical section of the projection from parabolic splittings to unipotent splittings. The image is given by those $\otimes$-isomorphisms of filtered fiber functors whose composition with $\gr$ is the identity automorphism of $\phi^{\gr}$. 

\subsubsection{\notreallyrem}
\label{splittings_torsor}
Parabolic splittings forms a pseudo-torsor\footnote{torsor or empty set with trivial action} under the parabolic Galois group $P(\phi)$; unipotent splittings form a pseudo-torsor under the unipotent Galois group $U(\phi)$.

\subsection{Graded Galois groups}
\label{sec:graded_galois_groups}
Let
\[
\phi \colon \Tc\to \FilLF(S)
\]
be a filtered fiber functor. In addition to the parabolic and Levi Galois groups, as well as the usual Galois (or \textit{fundamental}) group
\[
G(\phi) = \underline{\Aut}^\otimes (\forg_{\m{fil}} \circ \phi), 
\]
we define the \emph{(usual and parabolic) graded Galois groups of $\Tc$ at $\phi$} by
\[
G^\m{gr}(\phi) = \underline{\Aut}^\otimes(\forg_{\m{gr}} \circ \phi^{\gr})
\]
\[
P^\m{gr}(\phi) = \underline{\Aut}^\otimes(\fil \circ \phi^{\gr})
\]
(functorial in $S' \to S$ in the usual way). Note that $L(\phi) \subseteq G^{\gr}(\phi)$.
% We also define the \emph{cocharacter $\chi(\phi)$ associated to $\phi$} to be the cocharacter associated to the graded fiber functor $\m{gr} \circ \phi$ as in paragraph \ref{associated_cocharacter}.

\subsection{Associated cocharacter}
\label{associated_cocharacter}
Let
\[
\phi \colon \Tc\to \FilVect_{S}
\]
be a filtered fiber functor. We define the \emph{associated cocharacter}
\[
\chi(\phi) \colon \Gm \to L(\phi)
\]
as follows. We work with points with values in an arbitrary affine $S$-scheme $\Spec{R}$. Fix $z \in R^{\times} = \Gm(R)$ and $E \in \Tc$. Then
\[
\phi_R^{\gr}(E) = \phi^{\gr}(E) \otimes_{S} R
= \bigoplus_i \phi^{\gr}(E)_i \otimes_{S} R
\]
is a graded $R$-module, and we have $z$ act on it via scalar multiplication by $z^i$ on the $i$th summand. Evidently, $\chi(\phi)(z)$ preserves gradings, and hence belongs to $L(\phi)(R)$. 

\begin{ssLemma}
	\label{lemma:cochar_central}
	The cocharacter $\chi(\phi)$ associated to a filtered fiber functor $\phi$ is \emph{central} in the Levi Galois group $L(\phi)$; moreover, $L(\phi)$ is its centralizer in $\Gr^{\gr}(\phi)$.
\end{ssLemma}

\begin{proof}

The structure of a grading is the same as a $\Gm$-action. Thus an automorphism preserves the grading iff it commutes with the $\Gm$-action, so $L(\phi)$ is the centralizer of $\chi(\Gm)$.

%It follows from the construction that for any $E,R,z$ as above, $\phi_R(E)$ is central in the set of graded automorphisms of $\phi(E)$. Conversely, an isomorphism preserves the grading iff it commutes with $\chi(\phi)$; so the centralizer is the set of automorphisms of $\phi^{\gr}$.
\end{proof}

\subsection{Diagram of groups associated to a filtered fiber functor}\label{sec:graded_galois_groups_diagram}

In addition to the map $\m{ss}$ defined in \ref{defn:subgroups_filtration2}, there is a map $\m{ss}^{\gr} \colon P^{\gr}(\phi) \to L(\phi) = \underline{\Aut}^{\otimes}(\phi^{\gr})$ defined by applying the functor $\gr$ to an automorphism of $\fil \circ \phi^{\gr}$ and using the isomorphism $\gr \circ \m{fil} \simeq \m{id}_{\GrLF(S)}$. The functor $\fil$ similarly defines a section $r \colon L(\phi) \to P^{\gr}(\phi)$, fitting into a diagram of groups over $S$:
% https://q.uiver.app/#q=WzAsNixbMSwwLCJQKFxccGhpKSJdLFsyLDAsIkcoXFxwaGkpIl0sWzEsMSwiTChcXHBoaSkiXSxbMiwxLCJHXlxcbWF0aHJte2dyfShcXHBoaSkiXSxbMCwwXSxbMCwxLCJcXG1hdGhiYntHfV9tIl0sWzAsMV0sWzAsMl0sWzIsMywiIiwyLHsib2Zmc2V0IjoxfV0sWzMsMiwiIiwyLHsib2Zmc2V0IjoxfV0sWzUsMiwiXFxjaGkoXFxwaGkpIiwyXV0=
\[\begin{tikzcd}
	{} & {P(\phi)} & {G(\phi)} & {}\\
	{\mathbb{G}_m} & {L(\phi)} & {P^\mathrm{gr}(\phi)} & {G^\mathrm{gr}(\phi)}
	\arrow[hookrightarrow, from=1-2, to=1-3]
	\arrow["\m{ss}", from=1-2, to=2-2]
	\arrow["r"', shift right, from=2-2, to=2-3]
	\arrow["\m{ss}^{\gr}"', shift right, from=2-3, to=2-2]
	\arrow["{\chi(\phi)}"', from=2-1, to=2-2]
	\arrow[hookrightarrow, from=2-3, to=2-4]
\end{tikzcd},\]
with $\chi(\phi)$ central in $L(\phi)$.

\begin{sDefinition}
	\label{defn:weight_filtered_fib}
	Any fiber functor $\om \colon \Tc \to \QCoh(S)$ on a strictly filtered Tannakian category $(\Tc,\xi_W)$ over $\coeff$ upgrades canonically to a(n increasing) filtered fiber functor $\om^W$ defined by
\[\om^W \coloneqq \Fil^{\fin}(\om) \circ \xi_W \colon \Tc \to \Fil^{\fin}(\Tc) \to \Fil^{\fin} \QCoh(S) \subseteq \FilQCoh(S)\] (c.f. \ref{sssec:Gr_Fil_functors}), or on objects by $W_i \om^W(E) \coloneqq \om(W_i E)$; we define a \emph{weight-filtered fiber functor} to be a filtered fiber functor of this sort.
	
	We also set \[\om^{W,\gr} \coloneqq \gr \circ \om^W \colon \Tc \to \GrLF(S)\] \[\om^{\gr} \coloneqq \forg_{\m{gr}} \circ \om^{W,\gr} \colon \Tc \to \LF(S).\] If $\phi = \om^W$, then $ \om^{W,\gr} = \phi^{\gr}$, and $\om^{\gr} = \forg_{\gr} \circ \phi^{\gr}$ in the notation of \ref{defn:subgroups_filtration}.

	\subsubsection{}\label{PG_weight-filtered} For such $\om^W$, $P(\om^W) = G(\om^W) = \pi_1(\Tc,\om)$, and $P^{\gr}(\om^W) = G^{\gr}(\om^W) = \pi_1(\Tc,\om^{\gr})$.

\subsubsection{\yesrem}\label{rem:is_strictly_filtered_if_P=G}
	
	Conversely, suppose $\Tc$ is a Tannakian category over $\coeff$, and $\om^W \colon \Tc \to \FilVect_{\coeff}$ is any (increasing) filtered fiber functor for which $P(\om^W)=G(\om^W)$. Then there is a unique strict filtration on $\Tc$ for which $\om^W$ is the associated weight-filtered fiber functor. %a fact we will use in the proof of Theorem \ref{thm:hodge-filtered_weight-splitting}. (no longer use it because connectedness isn't an issue)

	\subsubsection{\notreallyrem}\label{rem:gr_commute}
	Let
	\[
	\om^W \colon \Tc \to \FilLF(S)
	\]
	be a weight-filtered fiber functor. Since $\om^W$ is exact, for each $E \in \Tc$ and $i\in \ZZ$ we have an exact sequence 
	\[
	0 \to \om^W(W_{i-1}E)
	\to \om^W(W_iE) \to \om^W(\gr_i^W E) \to 0
	\]
	in $\FilLF(S)$, from which a canonical isomorphism
	\[
	\tag{*}
	\gr^W_i \om^W E \simeq \om(\gr^W_i E).
	\]
    and thus
    \begin{eqnarray*}
    \om^{\gr}(E) \coloneqq \forg_{\gr} \circ \gr \circ \om^W(E)
    &=&
    \bigoplus_i \gr_i \om^W(E)\\
    &\overset{(*)}{\simeq}&
    \bigoplus_i \om(\gr^W_i E)\\
    &\simeq&
    \om\left(\bigoplus_i \gr^W_i E\right)\\
    &=&
    \om(\gr^W E),\end{eqnarray*}
    
which are in $E$ and thus provide a canonical isomorphism of tensor functors
\begin{equation}\label{eqn:gr_commute}	\om^{\gr} \simeq \om \circ \gr^W
\footnote{We may also see this last point by noting that $\om^{W,\gr} \coloneqq \gr \circ \om^W = \gr \circ \Fil^{\fin}(\om) \circ \xi_W
    \overset{\ref{sssec:Gr_Fil_functors}}{=} \Gr^{\fin}(\om) \circ \gr \circ \xi_W$ and thus that
$
\om^{\gr} \coloneqq \forg_{\gr} \circ \om^{W,\gr} = \forg_{\gr} \circ \Gr^{\fin}(\om) \circ \gr \circ \xi_W
    \overset{\ref{sssec:Gr_Fil_functors}}{=}
    \om \circ \forg_{\gr} \circ \gr \circ \xi_W
    \overset{\ref{defn:gr_in_filtered_tannakian}}{=}
    \om \circ \gr^W$.
}\end{equation}

%    \begin{eqnarray*}
%    \om^{\gr} \coloneqq \forg_{\gr} \circ \om^{W,\gr} &=& \forg_{\gr} \circ \Gr^{\fin}(\om) \circ \gr \circ \xi_W\\
%    &\overset{\ref{sssec:Gr_Fil_functors}}{=}&
%    \om \circ \forg_{\gr} \circ \gr \circ \xi_W\\
%    &\overset{\ref{defn:gr_in_filtered_tannakian}}{=}&
%    \om \circ \gr^W.
%    \end{eqnarray*}
    
%  &=& \forg_{\gr} \circ \Gr^{\fin}(\om) \circ \gr \circ \xi_W\\
 %   &\overset{\ref{sssec:Gr_Fil_functors}}{=}&
  %  \om \circ \forg_{\gr} \circ \gr \circ \xi_W\\
   % &\overset{\ref{defn:gr_in_filtered_tannakian}}{=}&
%    \om \circ \gr^W.
%    \end{eqnarray*}

%    \[
 %   \gr \om^W E \coloneqq \{\gr^W_i \om^W E\}_{\i \in \Zb}
%    \simeq
%    \{\om(\gr^W_i E)\}_{\i \in \Zb}
%    \]

In fact, we can replace the RHS by $\om^W \circ \gr^W$ to get a functor valued in $\FilLF(S)$ and then note that this upgrades to an isomorphism
\begin{equation}\label{eqn:gr_commute2}
\fil \circ \om^{W,\gr} \simeq \om^W \circ \gr^W
\end{equation}
of filtered fiber functors since $W_n$ of either side is $\bigoplus_{i \le n} \gr_i \om^W(E) \simeq \bigoplus_{i \le n} \om(\gr_i^W E)$ by \ref{sssec:double_gr^W}(\ref{eqn:W_n_of_gr^W}).
    
We will often use (\ref{eqn:gr_commute}) and (\ref{eqn:gr_commute2}) to commute $\om^W$ and $\gr^W$ without further mention.

\subsubsection{\remish}\label{remark:gr^W_in_terms}

For $E \in \Tc$, there are two natural Galois actions on $\om^{\gr}(E) \simeq \om(\gr^W E)$. The first is via the natural action of $G^{\gr}$ on $\om^{\gr}(E)$, while the second is via the natural action of $G$ on $\om(\gr^W E)$.

The actions are related in the following way. The action of $G$ on $\om(\gr^W E)$ is trivial when restricted to $U^W$ and thus factors through an action of $L^W$ on $\om(\gr^W E)$. There is also an action of $L^W$ on $\om^{\gr}(E)$ via $L^W \hookrightarrow G^{\gr}$. Then the isomorphism (\ref{eqn:gr_commute2}) intertwines these $L^W$-actions.

In this way, when $\coeff'=\coeff$, one may actually view $\gr^W$ as the operation of taking the fiber at $\om^{\gr}$, letting $G$ act via $G = P^W \xto{\sesi} L^W \hookrightarrow G^{\gr}$, then viewing this vector space with $G$-action as an object of $\Tc$ via $\om$. C.f. \ref{if_S=k}-\ref{groupoid_action_on_bitorsor}.

\subsubsection{\remish}\label{rem:gr_on_T^ss}

%% ***** ultimately seems unnecessary even if maybe a little nice
    
If $\Tc$ is weight-filtered, then $\gr^W \colon \Tc \to \Tc$ lands in the full subcategory $\Tc^{\sesi}$ of semisimple objects. Conversely, if $M \in \Tc^{\sesi}$, then $M \simeq \gr^W M$, so $\Tc^{\sesi}$ is the essential image of $\gr^W$. The isomorphism \ref{sssec:double_gr^W}(\ref{eqn:double_gr^W_functor}) thus implies that
\[
\restr{\gr^W}{\Tc^\m{ss}} \simeq \id_{\Tc^\m{ss}}
\]

% In the reference, probably need to also show that $\gr^W$ of a morphism is the same morphism. It would follow from the isomorphism of functors $\gr^W \circ \gr^W \simeq \gr^W$ and a statement that every morphism in $\Tc^{\sesi}$ comes from $\gr^W$. Then again, this last fact might be easiest to prove by proving that statement...

%If $M \in \Tc^{\sesi}$, then  implies that $\gr^W M \simeq M$.

%We thus have a square
%	\[
%	\tag{*}
%	\xymatrix{
%		\Tc \ar[d]_-{\gr^W} \ar[r]^-{\om^W} 
%		&
%		\FilLF(S) \ar[d]^-{\gr}
%		\\
%		\Tc^\m{ss} \ar[r]_-{\restr{\om^W}{\Tc^\m{ss}}} 
%		&
%		\GrLF(S),
%	}
%	\]

	\subsubsection{\yesrem}\label{rem:can_split}
	Every weight-filtered fiber functor over $\coeff'$ on a strictly filtered Tannakian category over $\coeff$ admits a ($\coeff'$-rational) splitting. Indeed, this is a special case of Ziegler \cite[Theorem 4.15]{ZieglerGrFil15}. However, under our extra assumption that the filtration exists on the level of the category itself, this is much easier to prove, as we now explain. Classical Tannaka duality guarantees the existence of an isomorphism of fiber functors 
	\[
	\forg_{\gr} \circ \om^{W,\gr}
	=
	\om^{\gr}
	\overset{\Psi^\circ} \simeq
	\om
	=
	\forg_{\m{fil}} \circ \om^W
	\]
	fpqc locally. Since the filtration exists on the level of the category itself, $\Psi^\circ$ respects the filtration and hence upgrades canonically to an isomorphism $\Psi$ of filtered fiber functors, i.e. to a parabolic splitting. As we've seen, this gives us a unipotent splitting fpqc locally. Since unipotent splittings form a torsor under a prounipotent group, this guarantees the existence of a $\coeff'$-rational unipotent splitting.

\end{sDefinition}

\subsection{Notation}
\label{notation:fil_X^X}
It will be convenient for us to decorate `$\FilVect$' and '$\FilLF$' so as to specify a symbol to be used to denote the associated filtrations and also to indicate whether filtrations are increasing or decreasing: a superscript $X$ as in `$\m{Fil}^X$' or `$\FilLF^X$' indicates increasing filtrations denoted by $X$; a subscript $X$ as in `$\m{Fil}_X$' or `$\FilLF_X$' indicates decreasing filtrations denoted by $X$. If $\Tc$ is a Tannakian category, by a \emph{Hodge-filtered fiber functor} we mean a filtered fiber functor which takes values in $\FilLF_F$. If $\Tc$ is a strictly filtered Tannakian category, by a \emph{weight-filtered fiber functor} we mean the filtered fiber functor with values in $\FilLF^W$ induced by an ordinary fiber functor $\Tc \to \LF$. We decorate symbols such as $\fil$, $\gr$, $\forg$ (\ref{defn:forg_fil_gr}, \ref{gr_forg}, \ref{defn:functors_fil_gr}) and others with $^W$ and $_F$.

We set $\FilLF^W_F(S)$ to be the full subcategory of $\Fil^W(\FilLF_F(S))$ for which forgetting the Hodge filtration lands in $\FilLF^W(S)$ and for which the inclusions $W_i \subseteq W_{i+1}$ are strict for the Hodge filtration. We let $\operatorname{GrFilLF}_F^W(S)$ be the full subcategory of $\Gr^W(\FilLF_F(S))$ for which $\forg_{\gr}^W$ is also in $\FilLF_F(S)$. Note that both are additive tensor categories but not abelian.

%$ \coloneqq \m{Fil}^W \m{Fil}_F \LF$.% When working with (filtered or graded) \textit{vector spaces}, we allow ourselves to drop the field from the notation.

\subsection{Weight and Hodge together}\label{WH}

Given a Hodge-filtered fiber functor $\om_H \colon \Tc \to \FilLF_F(S)$, we may set $W_i \om^W_H(E) \coloneqq \om_H(W_i E)$ as in Definition \ref{defn:weight_filtered_fib} to get an object of $\FilLF^W_F(S)$.

We can then apply $\gr^W \colon \FilLF^W_H(S) \to \operatorname{GrFilLF}_F^W(S)$, with Hodge filtration determined by \ref{hodge_on_grW} (in other words, we form the associated graded by taking quotients in the category $\m{Fil}_F\LF(S)$), to $\om_H^W$ to get a functor $\om_H^{W,\gr}$. We set $\om_H^{\gr} \coloneqq \forg^W \circ \om_H^{W,\gr}$. % = \forg^W \circ \om_H^{W,\gr}$.

\begin{ssLemma}\label{lemma:WH} The isomorphism $\fil^W(\om^{W,\gr}(E)) \simeq \om^W(\gr^W E)$ of \ref{rem:gr_commute}(\ref{eqn:gr_commute2}) preserves the Hodge filtration, i.e., upgrades to an isomorphism $\fil^W \om_H^{W,\gr}(E) \simeq \om_H^W(\gr^W E)$ in $\FilLF^W_F(S)$.
	
\end{ssLemma}

\subsubsection{Proof}\label{proof:WH} Let $i,j \in \ZZ$. By exactness of $\om_H$, the map $\om_H(W_j E) \to \om_H(E)$ is strict, so $F^i \om(W_j E) = \om(W_j) \cap F^i \om(W)$. Since $0 \to W_{j-1} E \to W_j E \to \gr_j^W E \to 0$ is exact in $\Tc$, the sequence
\[
0 \to \om_H(W_{j-1} E) \to \om_H(W_j E) \to \om_H(\gr_j^W E) \to 0
\]
is exact in $\FilLF_F(S)$. Thus $F^i \om(\gr_j^W E)$ is the image of $F^i \om(W_j E) = \om(W_j) \cap F^i \om(W)$, which agrees with \ref{hodge_on_grW} and thus \ref{WH}.

\begin{sDefinition}
	\label{hodge-filtered_weight_splitting}
	Let $\Tc$ be a strictly filtered Tannakian category over a field $\coeff$ and let
	\[
	\om_H \colon \Tc \to \FilLF_F(S)
	\]
	be a Hodge-filtered fiber functor over a $\coeff$-scheme $S$. 
	
	Recall we have $\om$ valued in $\LF(S)$ as in \ref{defn:hodge_filtr}, $\om^W, \om^{W,\gr}, \om^{\gr}$ valued in $\FilLF^W(S)$, $\GrLF^W(S)$, $\LF(S)$ as in \ref{defn:weight_filtered_fib}, and $\om_H^W$, $\om_H^{W,\gr}$ and $\om_H^{\gr}$ valued in $\Fil^W_F \LF(S)$, $\Gr^W  \Fil_F \LF(S)$, $\Fil_F \LF(S)$ as in \ref{WH}.
	
	Any $\otimes$-isomorphism as in the following diagram
	% https://q.uiver.app/#q=WzAsMixbMCwwLCJUIl0sWzEsMCwiXFxtYXRocm17RmlsfV9GXFxtYXRocm17VmVjdH0oUykiXSxbMCwxLCJcXG9tZWdhX0ggXFxjaXJjIFxcbWF0aHJte2dyfV5XIiwwLHsiY3VydmUiOi00fV0sWzAsMSwiXFxvbWVnYV9IIiwyLHsiY3VydmUiOjR9XSxbMiwzLCJwIiwyLHsic2hvcnRlbiI6eyJzb3VyY2UiOjIwLCJ0YXJnZXQiOjIwfX1dXQ==
	\[\begin{tikzcd}
		\Tc & {\mathrm{Fil}_F\LF(S)}
		\arrow[""{name=0, anchor=center, inner sep=0}, "{\om_H^{\gr}}", curve={height=-24pt}, from=1-1, to=1-2]
		\arrow[""{name=1, anchor=center, inner sep=0}, "{\om_H}"', curve={height=24pt}, from=1-1, to=1-2]
		\arrow["\Psi"', shorten <=6pt, shorten >=6pt, Rightarrow, from=0, to=1]
	\end{tikzcd}\]
	upgrades canonically to a $\otimes$-isomorphism as in the following diagram
	% https://q.uiver.app/#q=WzAsNSxbMCwxLCJUIl0sWzEsMSwiXFxtYXRocm17RmlsfV5XX0ZcXG1hdGhybXtWZWN0fShTKSJdLFsyLDAsIlxcbWF0aHJte0ZpbH1fRlxcbWF0aHJte1ZlY3R9KFMpIl0sWzIsMiwiXFxtYXRocm17RmlsfV5XXFxtYXRocm17VmVjdH0oUykiXSxbMywyLCJcXG1hdGhybXtHcn1eV1xcbWF0aHJte1ZlY3R9KFMpIl0sWzAsMSwiXFxvbWVnYV5XX0ggXFxjaXJjIFxcbWF0aHJte2dyfV5XIiwwLHsiY3VydmUiOi00fV0sWzAsMSwiXFxvbWVnYV5XX0giLDIseyJjdXJ2ZSI6NH1dLFsxLDMsIlxcbWF0aHJte2Zvcmd9X0YiLDJdLFsxLDJdLFswLDMsIlxcb21lZ2FeVyIsMix7Im9mZnNldCI6MiwiY3VydmUiOjV9XSxbMyw0LCJcXG1hdGhybXtncn1eVyIsMl0sWzUsNiwicCIsMix7InNob3J0ZW4iOnsic291cmNlIjoyMCwidGFyZ2V0IjoyMH19XV0=
	\[\begin{tikzcd}
		&& {\FilLF_F(S)} \\
		\Tc & {\FilLF^W_F(S)} \\
		&& {\FilLF^W(S)} & {\GrLF^W(S)}
		\arrow[""{name=0, anchor=center, inner sep=0}, "{\om^W_H \circ \mathrm{gr}^W \simeq \fil^W \circ \om^{W,\gr}_H}", curve={height=-24pt}, from=2-1, to=2-2]
		\arrow[""{name=1, anchor=center, inner sep=0}, "{\om^W_H}"', curve={height=24pt}, from=2-1, to=2-2]
		\arrow["{\forg_F}"', from=2-2, to=3-3]
		\arrow["{\forg^W}"', from=2-2, to=1-3]
		\arrow["{\om^W}"', shift right=2, curve={height=30pt}, from=2-1, to=3-3]
		\arrow["{\mathrm{gr}^W}"', from=3-3, to=3-4]
		\arrow["\Psi^W"', shorten <=6pt, shorten >=6pt, Rightarrow, from=0, to=1]
	\end{tikzcd},\]
	where $\forg^W \circ \Psi^W = \Psi$.
	
	Now $\forg_F \circ \Psi^W$ is an isomorphism from $\fil^W \circ \om^{W,\gr}$ to $\om^W$. By a \emph{Hodge-filtered weight-splitting}, we mean a $\otimes$-isomorphism $\Psi$ as above, for which $\forg_F \circ \Psi^W$ is a unipotent splitting of $\om^W$; recall this means that its further composition with $\gr^W$ is equal to the identity natural automorphism of $\om^{W,\gr}$.
\end{sDefinition}

For the remainder of \S \ref{sec:arith_hodge_path}, we fix a weight-filtered Tannakian category $\Tc$ and a Hodge-filtered fiber functor $\om_H$ over $\coeff'$.% As above, we denote the underlying fiber functor by $\om$ and the associated weight-filtered fiber functor by $\om^W$.

\subsection{Diagram of groups associated to a Hodge-filtered fiber functor}
\label{WH_groups_diagram}
Associated to the filtered fiber functors $\om_H$ and $\om^W$ we have \textit{Levi}, \textit{Parabolic}, \textit{unipotent}, and \textit{graded} Galois groups; we decorate accordingly. Specifically, we set
\begin{itemize}
	\item
	$G = \underline{\Aut}^\otimes(\om)$
	\item
	$G^{\gr} = \underline{\Aut}^\otimes(\om \circ \gr^W) = \underline{\Aut}^{\otimes}(\om^{\gr}) = P^{W,\gr}$
	\item
	$P^W = \underline{\Aut}^\otimes(\om^W) = G$
	;
	$P_H = \underline{\Aut}^\otimes(\om_H) = \underline{\Aut}^\otimes(\om_H^W)$
	\item
	$L^W = \underline{\Aut}^\otimes(\om^{W,\gr}) = \underline{\Aut}^{\otimes}( \restr{\om^W}{\Tc^\m{ss}})$
	\item
	$U^W = \Ker (P^W \to L^W)$
	
\end{itemize}
as in \ref{defn:subgroups_filtration}. We also have the central cocharacter $\chi^W \colon \Gb_m \to L^W$ associated to $\om^W$ (\ref{associated_cocharacter}). We let $U^W_H$ be the kernel of the induced map
\[
P_H \hookrightarrow G = P^W \twoheadrightarrow L^W.
\]
These fit into a commutative diagram of groups over $\coeff'$ as follows:
% https://q.uiver.app/#q=WzAsOSxbMSwxLCJQXlciXSxbMiwxLCJHIl0sWzEsMiwiTF5XIl0sWzIsMiwiR15XIl0sWzAsMV0sWzMsMSwiUF9IIl0sWzAsMiwiXFxtYXRoYmJ7R31fbSJdLFs0LDAsIlVfSF5XIl0sWzEsMCwiVV5XIl0sWzUsMSwiIiwxLHsic3R5bGUiOnsidGFpbCI6eyJuYW1lIjoiaG9vayIsInNpZGUiOiJib3R0b20ifX19XSxbMCwxLCJcXHNpbSJdLFswLDJdLFsyLDMsIiIsMix7Im9mZnNldCI6MX1dLFszLDIsIiIsMix7Im9mZnNldCI6MX1dLFs2LDIsIlxcY2hpXlciLDJdLFs3LDUsIiIsMSx7InN0eWxlIjp7InRhaWwiOnsibmFtZSI6Imhvb2siLCJzaWRlIjoidG9wIn19fV0sWzgsMF0sWzcsOCwiIiwwLHsic3R5bGUiOnsidGFpbCI6eyJuYW1lIjoiaG9vayIsInNpZGUiOiJib3R0b20ifX19XV0=
\[\begin{tikzcd}\tag{*}
	& {U^W} &&& {U_H^W} \\
	{} & {P^W} & G & {P_H} \\
	{\mathbb{G}_m} & {L^W} & {G^{\gr}}
	\arrow[hook', from=2-4, to=2-3]
	\arrow["\sim", from=2-2, to=2-3]
	\arrow["\m{ss}", from=2-2, to=3-2]
	\arrow["r"', shift right, from=3-2, to=3-3]
	\arrow["\m{ss}^{\gr}"', shift right, from=3-3, to=3-2]
	\arrow["{\chi^W}"', from=3-1, to=3-2]
	\arrow[hook, from=1-5, to=2-4]
	\arrow[from=1-2, to=2-2]
	\arrow[hook', from=1-5, to=1-2]
\end{tikzcd}\]

\subsubsection{\yesrem} The condition of $\Tc$ being \emph{weight-filtered} (\ref{defn:weight_filtr_cat}) rather than just \emph{strictly filtered} (\ref{dfthm:strict}) ensures that $L^W$ is reductive, although we will not need this until \S \ref{sec:arith_frob_path}.

\subsection{Example: non-CM elliptic curve}
\label{example:E}
%Suppose $G$ is reductive. Then any parabolic subgroup $P \subset G$ gives rise in an evident way to a filtered fiber functor
%\[
%\phi \colon \Rep G \to \FilVect_{\coeff}
%\]
%such that $P(\phi) = P$. Saavedra-Rivano \cite[IV.2.4.3.2]{SaavedraRivanoBook} proves the converse: if $\phi$ is any filtered fiber functor, then $P(\phi)$ is a parabolic subgroup of $G$. In particular, the rank of $P$ is equal to the rank of $G$. 

Consider $\mathbf{SM}_{\QQ}(\Qb,E)$ (\ref{sec:SM_A}) for $E$ an elliptic curve over $\QQ$ without complex multiplication, and let
\[
\om_H^\dR \colon\mathbf{SM}_{\QQ}(\Qb,E) \to \FilVect_{\QQ}
\]
and $\om^{\dR}$ as in \ref{sec:fiber_func_PSPM}.

We have
\[
L^W = \GL(\om^{\dR}(h^1_{\mathbf{SM}}(E))) = \GL(H^1_{\dR}(E)) \simeq \GL_{2},
\]
which follows either by the fact that the Mumford--Tate group of $h_1(E)$ is $\GL_2$ (\cite[5.4]{MoonenIntroMT}) or by Serre's Open Image Theorem (\cite{SerreOpenImageI}).

Then the image of $P_H$ in $L^W$, denoted $F^0 \GL_2$ is the Borel (hence parabolic) subgroup of $\GL(H^1_{\dR}(E))$ fixing the subspace generated by a nonzero invariant differential form on $E$.

\begin{sProposition}\label{prop:U^W_H_is_prounip}
	The affine $\coeff'$-group $U^W_H$ is prounipotent.
\end{sProposition}

\begin{proof}
	Evidently, it's enough to show that $U^W$ is prounipotent. This has been noted before; 
	see, for instance, \cite[\S 3]{GonMEM}. One way to see this is to note that the left-regular representation of $G$ provides a faithful Ind-unipotent representation of $U^W$.
\end{proof}

%\begin{ssLemma}\label{lemma:reduction_to_K''}
%	For a given $\om_H$ and algebraic $\Tc$, it suffices to show that $(\om_H)_{\coeff''}$ has a Hodge-filtered weight-splitting for some overfield $\coeff''$ of $\coeff'$.
%\end{ssLemma}

\begin{ssLemma}\label{lemma:Spl^W_H_pseudotorsor}
For a $\coeff'$-scheme $S$, let $\underline{\Spl^W_H}(\om_H)(S)$ denote the set of Hodge-filtered weight-splittings of $(\om_H)_S$. Then $\underline{\Spl^W_H}(\om_H)$ forms an effective\footnote{representable by a scheme} pseudo-torsor under $U^W_H$ over $\Spec{\coeff'}$ for the fpqc topology.
\end{ssLemma}

\begin{proof}[Proof of Lemma \ref{lemma:Spl^W_H_pseudotorsor}]
	
%For a $\coeff'$-scheme $S$, let $\underline{\Spl^W_H}(\om_H)(S)$ denote the set of Hodge-filtered weight-splittings of $(\om_H)_S$.

A priori this is a presheaf on $\coeff'$-schemes, but we claim it is an affine $\coeff'$-scheme. Indeed, $\underline{\Isom}^{\otimes}(\om_H^{\gr},\om_H)$ is an affine $\coeff'$-scheme by \ref{isom_filt_grad_fib}, and $\underline{\Spl^W_H}(\om_H)$ is the fiber of $\id_{\om^{\gr}}$ under the map $\underline{\Isom}^{\otimes}(\om_H^{\gr},\om_H) \xto{\forg_F} \underline{\Isom}^{\otimes}(\om^{\gr},\om) \xto{\gr^W} \underline{\Aut}^{\otimes}(\om^{\gr})$, thus a closed subscheme.
	
	Next, note that $U^W_H=U^W_H(\om)$ is the kernel of the map
	\[
	P_H = \underline{\Aut}^\otimes(\om^W_H) \xrightarrow{\gr^W \circ \forg_F} 
	\underline{\Aut}^\otimes(\om^{W,\gr}) = L^W.
	\]
	From this, it is clear that if $a \in \underline{\Spl^W_H}(\om_H)(S)$ and $b \in U^W_H(S)$, then $b \circ a \in \underline{\Spl^W_H}(\om_H)(S)$. Similarly, if $a,a' \in \underline{\Spl^W_H}(\om_H)(S)$, then $a' \circ a^{-1} \in U^W_H(S)$. Thus $\underline{\Spl^W_H}(\om_H)$ is a pseudo-torsor (as a sheaf) under $U^W_H$.
	
%	Suppose $\underline{\Spl^W_H}(\om_H)(\coeff'')$ is nonempty for some overfield $\coeff''/\coeff'$. Then $(\underline{\Spl^W_H}(\om_H))_{\coeff''} \simeq (U^W_H)_{\coeff''}$. Since $\Tc$ is algebraic, $G(\om)$ and therefore $U^W_H$ is finite type and thus fppf. Thus $\underline{\Spl^W_H}(\om_H)$ is fpqc-locally fppf and hence fppf.
	
%	Finally, since $H^1_{\mathrm{fppf}}(\coeff';U^W_H) = H^1_{\et}(\coeff';U^W_H)$ is trivial because $U^W_H$ is unipotent, the set $\underline{\Spl^W_H}(\om_H)(\coeff')$ is nonempty, and we are done.
\end{proof}

\begin{sProposition}
	\label{prop:ranks_equal}Suppose $\Tc$ is algebraic (\ref{sec:algebraic_tannakian}). Then the ranks of $P_H$, $G$, $P^W$, $L^W$, and $G^{\gr}$ are all equal.
\end{sProposition}

\begin{proof}

The assumption on $\Tc$ implies that all the groups in \ref{WH_groups_diagram}(*) are algebraic (i.e., of finite type), so they have a well-defined notion of rank.

Let $M$ be a generator of $\Tc$. Then $P_H$ is the subgroup of $G$ fixing the Hodge filtration on $\om_H(M)$, so $P_H$ is parabolic, and therefore $r \coloneqq \rk{P_H}=\rk{G}$. Since $G=P^W$, $\rk{P^W}=\rk{G}=r$. Since $G \simeq G^{\gr}$, we have $\rk{G^{\gr}} = \rk{G} = r$.

Since the kernel of $\sesi$ is unipotent, $\rk{\sesi(P^W)} = \rk{P^W} =r$. Since a subgroup has smaller or equal rank, we have $r=\rk{\sesi(P^W)} \le \rk{L^W} \le \rk{G^{\gr}}=r$, so in fact $\rk{L^W}=r$.
\end{proof}

\begin{sTheorem}
	\label{thm:hodge-filtered_weight-splitting}
	Let $\coeff$ be a field of characteristic zero, $\Tc$ a strictly filtered Tannakian category over $\coeff$, $\coeff'$ an overfield of $\coeff$, and
	\[
	\om_H \colon \Tc \to \m{Fil}_F\Vect_{\coeff'}
	\]
	a Hodge-filtered fiber functor. Then $\om_H$ admits a $\coeff'$-rational Hodge-filtered weight-splitting; in other words, $\underline{\Spl^W_H}(\om_H)(\coeff')$ is nonempty.
\end{sTheorem}

The proof spans \ref{proof:finite_type}--\ref{general_G}. In the case $\Tc$ is algebraic, the result reduces to a form of Levi's Theorem (\cite[Theorem 5.5]{ZieglerGrFil15}), as described in \ref{proof:finite_type}. A similar idea appears in \cite[Lemma 3.1]{HainWeighted03} and \cite[\S 10.2]{HainGoldman20}, although without the condition of respecting the Hodge filtration; c.f. also \cite[Proposition 3.1]{Hain_HodgeDeRham_Chapter2016}. We reduce to this case in \ref{general_G}.

\subsubsection{}\label{proof:finite_type}

Assume that $\Tc$ and hence $G$ is algebraic. Let $T$ be a maximal torus of $P_H$. Since $\rk{P_H}=\rk{G}$, $T$ is also a maximal torus of $G=P^W$. Furthermore, $T \cap U^W$ is trivial, so $\sesi(T) \simeq T$ is a subtorus of $L^W$ of rank $\dim T = \rk{P_H} \overset{\ref{prop:ranks_equal}}{=} \rk{L^W}$. Therefore $\sesi(T)$ is a maximal torus of $L^W$. Since $\chi^W$ is central (Lemma \ref{lemma:cochar_central}), its image lands in $\sesi(T)$.

We may therefore lift $\chi^W$ to a cocharacter of $T$ and therefore of $P_H$:
\[
\tilde \chi^W  \colon \Gm \to P_H
\]
Since $\tilde \chi^W$, regarded as a cocharacter of $P^W$, lifts $\chi^W$, it splits the weight filtration by \cite[Lemma 4.11]{ZieglerGrFil15}. Since it lands in $P_H$, the splitting preserves Hodge filtrations by Lemma \ref{lemma:preserve_hodge}.

\subsubsection{}\label{general_G}

Next, following \cite[\S 5.5]{ZieglerGrFil15} we argue indirectly that the same holds for general $\Tc$. (Meanwhile, we avoid the main issues dealt with there by restricting attention to $\coeff$ of characteristic zero.) The collection $I$ of algebraic full Tannakian subcategories is a \textit{set} since $\Tc$, by definition, is essentially small. As a matter of notation, we treat $I$ as an index set, and we denote the Tannakian category associated to $i\in I$ by $\Tc_i$. Let $G_i = \underline{\Aut}^\otimes(\om\restr{}{\Tc_i})$ be the Tannakian Galois group of $\Tc_i$. Then $G_i$ is a finite type quotient of $G$. In particular, the case considered above applies to show that $\underline{\Spl^W_H}(\restr{\om_H}{\Tc_i})$ possesses a $\coeff'$-point.

Whenever $\Tc_i \subset \Tc_j$, there's an associated map 
\[
\underline{\Spl^W_H}(\restr{\om_H}{\Tc_i}) 
\from
\underline{\Spl^W_H}(\restr{\om_H}{\Tc_j}).
\]
These assemble to a diagram of affine $\coeff'$-schemes, and a straightforward verification shows that 
\[
\underline{\Spl^W_H}(\om_H) = \varprojlim_{i \in I} 
\underline{\Spl^W_H}(\restr{\om_H}{\Tc_i}).
\]
As explained in \cite[Lemma 2.9]{ZieglerGrFil15}, $I$ is a lattice and is in particular filtered. Additionally, the base $\Spec{\coeff'}$ occurs within the diagram as
\[
\Spec{\coeff'} = \underline{\Spl^W_H}(\restr{\om_H}{\langle \mathbbm{1}_{\Tc} \rangle}).
\]
%where $\langle \mathbbm{1} \rangle$ denotes the Tannakian subcategory generated by $\mathbbm{1}$.
Hence, Proposition 8.3.8 of EGA IV \cite{EGAIV3} applies to show that $\underline{\Spl^W_H}(\om_H)$ is faithfully flat over $\coeff'$ (equivalently, nonempty). Consequently, $\underline{\Spl^W_H}(\om_H)$ is an fpqc torsor under the prounipotent group $U^W_H$.

Let $U \coloneqq U^W$ and $U_i \coloneqq U(\om^W\restr{}{\Tc_i})$. In \cite[Theorem 4.17]{ZieglerGrFil15}, it is shown that every fpqc torsor under $U$ is trivial (apply the theorem to $\varphi = \om^W$) using an intricate argument involving $\varprojlim^1$. This is done by putting compatible decreasing filtrations $\{U^{\alpha}\}_{\alpha \in \Nb},\{U^{\alpha}_i\}_{\alpha \in \Nb}$ on $U$ and $U_i$ and showing (in Lemma 4.29 of loc.cit.) that
\[
{\varprojlim_i}^1 U^{\alpha}_i/U^{\alpha+1}_i
\]
is trivial for each $\alpha \ge 1$. The proof uses Theorem 2.17 of loc.cit., for which there are four conditions.

We may put the same filtration on the subgroup $U^W_H$ of $U=U^W$
%and on the finite-type versions $(U^W_H)_i \coloneqq U^W_H(\om^W\restr{}{\Tc_i})$
, where $(U^W_H)^\alpha$ is the subgroup acting trivially on $W_{j+\alpha} \om(M)/W_j \om(M)$ for all $M \in \Tc$ and $j \in \Zb$. We let $(U^W_H)^\alpha_i$ denote the image of $(U^W_H)^\alpha$ in $G_i$; note that it respects the Hodge filtration and acts trivially on the associated graded for the weight filtration. We wish to show that 
\[
{\varprojlim_i}^1 (U^W_H)^\alpha_i/(U^W_H)^{\alpha+1}_i
\]
is trivial for each $\alpha \ge 1$. In fact, all the conditions in Theorem 2.17 follow from the corresponding conditions in the proof of Lemma 4.29 of loc.cit.: (i) because it is the same lattice, (ii) because in characteristic $0$, an algebraic subgroup of an abelian unipotent group is again abelian unipotent, and (iii),(iv) because we have defined $(U^W_H)^\alpha_i$ as the image of $(U^W_H)^\alpha$.

Then the proof of Lemma 4.30 and thus of Theorem 4.17 of loc.cit. works in the same way (in particular, for fixed $i$, we also have $(U^W_H)^\alpha_i$ trivial for sufficiently large $\alpha$) to show that every fpqc torsor under $U^W_H$ is trivial. Therefore, $\underline{\Spl^W_H}(\om_H)(\coeff')$ is nonempty, and we are done.

%One could try to apply the article ``A Remark on Torsors under Affine Group Schemes'' by Michael Wibmer. But in trying to prove his Lemma 2.4 in our case (not algebraically closed), I'm stuck showing that the map $U(\coeff') \to U_i(\coeff')$ is surjective, which seems to boil down to showing that a torsor under a prounipotent group is always trivial...

\qed

%\subsection{Remark}\label{rem:hodge_thm_k_k'}
%If $\Tc$ a strictly filtered Tannakian category over $\coeff$, and
%\[
%\om_H \colon \Tc \to \m{Fil}_F\Vect_{\coeff'}
%\]
%a Hodge-filtered fiber functor for an extension $\coeff'/\coeff$, then Theorem \ref{thm:hodge-filtered_weight-splitting} shows that $\om_H$ admits a $\coeff'$-rational Hodge-filtered weight-splitting by applying it to $\Tc_{\coeff'}$ (\ref{arbitrary_S}) and then restricting the splitting to $\Tc$.

%%%%%%%%%%%
\section{Construction of Unipotent \texorpdfstring{$p$}{p}-adic Period Map}\label{sec:arith_frob_path}
%%%%%%%%%%%%

\begin{sTheorem}
\label{thm:arith_frob_path}
Let $\Tc$ be a weight-filtered Tannakian category over $\coeff$, $\coeff'_0/\coeff$ an extension with a theory of weights and finite order automorphism $\phi$ fixing $\coeff$, and
\[
\om_{\phi} \colon
\Tc \to \phi\m{Mod}(\coeff'_0)
\]
a Frobenius-equivariant fiber functor over $(\coeff'_0,\phi)$ (Definition \ref{defn:frob_equiv}).

%% below is not needed because already defined earlier
%We denote the associated fiber functor 
%\[
%\Tc \to \phi\m{Mod}(\coeff'_0) \to 
%\Vect_{\coeff'_0}
%\]
%by $\om$. We indicate precomposition with $\gr^W$ by adding a superscript `$\gr$'.

Then there exists a unique $\otimes$-isomorphism of Frobenius-equivariant fiber functors
	\[
	\tannloop^{\cri}_{\phi} \colon \om_{\phi}^{\gr}
	\to
	\om_{\phi},
	\]
	where $\om_{\phi}^{\gr}$ has Frobenius action via $\om_{\phi}^{\gr} \simeq \om_{\phi} \circ \gr^W$ (\ref{rem:gr_commute}). In other words, in the notation of \ref{defn:frob_equiv}, there exists a unique Tannakian path 
	\[
	\tannloop^{\cri}_0 \colon \om_{0}^{\gr}
	\to
	\om_{0}
	\]
	such that for each $E\in \Tc$, the square
	\[
	\tag{*}
	\xymatrix{
		\omgr_0(E) 
		\ar[r]^-{\tannloop^\m{cr}_0} \ar[d]_{\gr^W(\phi)}
		& 
		\om_0(E)
		\ar[d]^\phi
		\\
		\omgr_0(E)
		\ar[r]_-{{\tannloop^\m{cr}_0}} 
		& 
		\om_0(E)
	}
	\]
	commutes.
	
\end{sTheorem}

The proof spans \S \ref{sec:bitorsor}-\ref{sec:phi_on_spl}.

\subsection{Bitorsor of splittings}\label{sec:bitorsor}

Since $\phi$ does not in general commute with arbitrary $\coeff'_0$-linear maps, we cannot upgrade $\om_{\phi}$ to a functor from $\Tc_{\coeff'_0}$ (\ref{arbitrary_S}) to $\phi$-modules. Since our proof uses Frobenius actions on groups such as $U^W$ and on the torsor of splittings, we must realize these as objects of $\Tc$ itself. Since our only given fiber functor $\om_0$ is over $\coeff'_0$ rather than $\coeff$, we will need to use the theory described in \S \ref{sec:reconstruction_non_neutral}, beginning in \ref{groupoid_action_on_bitorsor}.

Suppose $\Tc$ is strictly filtered. Let $\om \coloneqq \om_0 = \forg_{\phi} \circ \om_{\phi}$, and in the notation of \S \ref{sec:graded_galois_groups_diagram}, consider $G(\om)=G(\om^W)$, $G^{\gr}(\om)=G^{\gr}(\om^W)$, $L(\om)=L(\om^W)$, $\m{ss}$, $\m{ss}^{\gr}$, and $r$. The fiber functor $\omega$ gives rise to a $\coeff$-groupoid scheme $\Gc(\om)$ with scheme of objects $S_0 \coloneqq \Spec{\coeff'_0}$ as in \S \ref{sec:reconstruction_non_neutral}. Note that $G(\om)$ is the fiber of $\Gc(\om)$ over the diagonal $S_0 \xto{\Delta} S_0 \times_{\coeff} S_0$. We recall $U(\om)=U(\om^W) = \Ker(\m{ss})$ and define $U^{\gr}(\om) = U^{\gr}(\om^W) \coloneqq \Ker(\m{ss}^{\gr})$. Finally, we let
\[
\underline{\Spl^W}(\om)
\]
denote the scheme of unipotent splittings of $\om^W$ as in \cite[Definition 4.10(iv)]{ZieglerGrFil15}; for an $S_0$-scheme $S$, $\Spl(\om)(S)$ is the set of unipotent splittings (\S \ref{sec:splittings_fib}) of $\om_{S}^W$. By \cite[Theorem 4.15]{ZieglerGrFil15}, $\underline{\Spl^W}(\om)$ is a $U(\om) \text{-} U^{\gr}(\om)$-bitorsor over $S_0$.

\subsubsection{}\label{if_S=k}

If $S_0=\Spec{\coeff}$, we can upgrade this to a bitorsor in $\Tc$ as follows:

In the notation of \S \ref{sec:reconstruction_non_neutral}, we define a left action $\star$ of $g \in G(\om)$ on $p \in \underline{\Spl^W}(\om)$ by
\[
g\star p  = g \circ p \circ r(\m{ss}(g))\inv.
\]

We also denote the conjugation action of $G(\om)$ on $U(\om)$ and via $r \circ \m{ss}$ on $U^{\gr}(\om)$ by $\star$. This gives rise to group schemes $\underline{U}$ and $\underline{U^{\gr}}$ in $\Tc$.\footnote{We remark that $\underline{U^{\gr}}$ is isomorphic in $\Tc$ to $\gr^W \underline{U}$ by Remark \ref{remark:gr^W_in_terms}.}

With these definitions, the composition action
\[\tag{*}
\circ \colon
U(\om) \times \underline{\Spl^W}(\om) \times U^{\gr}(\om) \to \underline{\Spl^W}(\om)
\]
is $\star$-equivariant and gives $\underline{\Spl^W}(\om)$ the structure of a $G(\om)$-equivariant $U(\om) \text{-} U^{\gr}(\om)$-bitorsor. In particular, it underlies via $\om_0$ a $\underline{U} \text{-} \underline{U^{\gr}}$-bitorsor $\underline{\Spl^W}$ in the category $\Tc$.

We now wish to apply \ref{rep_of_groupoid} to do this even when $S_0 \neq \Spec{\coeff}$.

\subsubsection{}\label{groupoid_action_on_bitorsor}
We define an action of $\Gc=\Gc(\om)$ on $U(\om), \underline{\Spl^W}(\om), U^{\gr}(\om)$ as follows. We work with points in an arbitrary $\coeff$-scheme $T$.

Given $g \in \Gc(T)$ corresponding to $(b(g),s(g)) \in (S_0 \times_{\coeff} S_0)(T)$ and $g \in \Isom^{\otimes}(s(g)^* \om,b(g)^* \om)$, we define
\[
g^{\m{ss}} \in \Isom^{\otimes}(s(g)^* \om^{\gr},b(g)^* \om^{\gr})
\]
by composing $g$ with the associated graded functor $\gr^W \colon \Tc \to \Tc$.

Given such $g,s(g),b(g)$ and $u \in U(\om)_{s(g)}(T)$, we define $g \star u \in U(\om)_{b(g)}(T)$ by the commutative square
\[\tag{*}
\xymatrix{
	s(g)^* \om \ar[r]^-u \ar[d]_-{g} & s(g)^* \om \ar[d]^-{g}\\
	b(g)^* \om \ar[r]^-{g \star u} & b(g)^* \om
}
\]

Note that the group scheme $\underline{U}$ in $\Tc$ thus defined is a $\Tc$-subgroup scheme of $\pi(\Tc)$, the latter in the sense of \cite[D\'efinition 6.1]{Deligne89} or \cite[D\'efinition 8.13]{DelTann90}. In the notation of \cite[6.6]{Deligne89}, $U$ is the subgroup $H$ associated to $\Tc_1 = \Tc^{\m{ss}}$.

Given $p \in \underline{\Spl^W}(\om)_{s(g)}(T)$, we define $g \star p \in \underline{\Spl^W}(\om)_{b(g)}(T)$ by the commutative square
\[
\xymatrix{
	s(g)^* \omgr \ar[r]^-p \ar[d]_-{g^{\m{ss}}} & s(g)^* \om \ar[d]^-{g}\\
	b(g)^* \omgr \ar[r]^-{g \star p} & b(g)^* \om
}
\]

Given $u \in U^{\gr}(\om)_{s(g)}(T)$, we define $g \star u \in U^{\gr}(\om)_{b(g)}(T)$ by the commutative square
\[
\xymatrix{
	s(g)^* \omgr \ar[r]^-u \ar[d]_-{g^{\m{ss}}} & s(g)^* \omgr \ar[d]^-{g^{\m{ss}}}\\
	b(g)^* \omgr \ar[r]^-{g \star u} & b(g)^* \omgr
}
\]

Note that when $s(g)=b(g)$, this is precisely the action $\star$ of $G(T)$ as defined in \ref{if_S=k}.

It is a simple diagram chase that the actions of $U(\om)$ and $U^{\gr}(\om)$ on $\underline{\Spl^W}(\om)$ commute with the action of $\Gc$. By \ref{rep_of_groupoid}, we get a $\underline{U} \text{-} \underline{U^{\gr}}$-bitorsor $\underline{\Spl^W}$ in the category $\Tc$.

\begin{sLemma}\label{lemma:Uab_V_Ext}
	In the situation of \S \ref{sec:bitorsor}, suppose furthermore that $\Tc$ is weight-filtered, so that $\GG \coloneqq L(\om_0)$ is reductive and $\Ext^1_{\Tc}(\mathbbm{1}_{\Tc},\mathbbm{1}_{\Tc}) = 0$. Then there's an isomorphism in $\Pro-\Tc$
	\[
	a \colon
	\underline{U}^\m{ab} \xrightarrow{\sim} 
	\prod_{V \, \mathrm{ simple}} 
	\Ext^1_{\Tc}(\mathbbm{1}_{\Tc}, V)^\vee \otimes V.
	\]
\end{sLemma}

\begin{proof}
	This is well known when $S_0=\Spec{k}$; see for example \cite[Proposition 6.1]{BrownNotes17} or \cite[Proposition 3.1.1]{eskandari2024tannakianfundamentalgroupsblended} for the dual statement, or \cite[(3.4)]{Hain_HodgeDeRham_Chapter2016} for a related non-dual statement. Let us define the map and check equivariance in general; the map will be an isomorphism since it is after base-change to $S_0$.
	
	Given $V$, we define the $V$-component 
	\[
	a_V \colon \underline{U}^\m{ab} \to
	M_V \coloneqq \Ext^1_{\Tc}(\mathbbm{1}_{\Tc}, V)^\vee \otimes_{\coeff} V
	\]
	of $a$ as follows. If $V$ is the trivial representation, then the Ext group is trivial so the target is terminal. Suppose $V$ is not trivial, and let 
	\[
	0 \to M_V \to E_V \to \mathbbm{1}_{\Tc} \to 0
	\]
	be the universal extension.
	
	As usual, we work with values in an arbitrary affine $S_0$-scheme $\Spec{R}$. Let $1_E$ be an element of $\om(E_V)(R)$ mapping to $1 \in R$. Then for $u \in U(\om)(R)$,
	\[
	a_V(u) \coloneqq -1_E + u(1_E)\footnote{This strange way of writing it is to stress the similarity with $(\de P)(u) = (\gamma^\m{gr})\inv u(\gamma^\m{gr})$ in \S \ref{sec:abstract_CK_maps}.} \in \om(M_V)(R).
	\]
	
	We must verify that it is well-defined (independent of $1_E$). For this, note that since $V$ is simple, we have $V = \gr^W V$. Therefore, $U(\om)$ acts trivially on $\om(V)$ and hence on $\om(M_V) = \Ext^1(\mathbbm{1}_{\Tc}, V)^{\vee} \otimes_{\coeff} \om(V)$. It follows that if $1_E'$ is another lift of $1$ to $\om(E_V)(R)$, we have
	\begin{eqnarray*}
		\left[-1_E + u(1_E)\right] - \left[-1_E' + u(1_E')\right]
		&=&
		\left[-1_E + 1_E'\right] - \left[-u(1_E) + u(1_E')\right]\\
		&=&
		\left[-1_E + 1_E'\right] - u\left[-1_E + 1_E'\right]\\
		&=& 0
	\end{eqnarray*}
	because $-1_E + 1_E' \in \om(M_V)(R)$.
	
	We can use well-definedness to show that the map is a group homomorphism. We have
	\begin{eqnarray*}
		a_V(u_1 u_2) = -1_E + u_1(u_2(1_E))
		&=&
		-1_E + u_1(1_E) + u_1(-1_E + u_2(1_E))\\
		&=&
		a_V(u_1) + u_1(a_V(u_2))\\
		&=&a_V(u_1)+a_V(u_2),
	\end{eqnarray*}
	the last by the fact that $U(\om)$ acts trivially on $\om(M_V)$. Since the target is abelian, the map factors through $U(\om)^{\ab} = \om(\underline{U}^{\ab})$.
	
	We have thus defined a homomorphism of group schemes over $S_0$ from $\om(\underline{U}^{\ab})$ to $\om(M_V)$. We now must show that it is equivariant for the action of $\Gc(\om)$ to get a map $\underline{U}^{\ab} \to M_V$.
	
	Let $T$ be a $\coeff$-scheme and $g \in \Gc(T)$. We denote (by abuse of notation) $s(g),b(g)$ by $s,b$ so that we have morphisms of $\coeff$-schemes $s,b \colon T \to S_0$. We write $T_s$ and $T_b$ for these two different ways of viewing $T$ as an $S_0$-scheme. We have a map $g \colon s^* \om \to b^* \om$. We must show that the following diagram is commutative:
	\[
	\xymatrix{
		s^* \om(U^{\ab})(T) \ar[r]^-{g_{U^{\ab}}} \ar[d]_-{s^* a_V} & b^* \om(U^{\ab})(T) \ar[d]^-{b^* a_V}\\
		s^* \om(M_V)(T) \ar[r]^-{g_{M_V}} & b^* \om(M_V)(T)
	}.
	\]
	
	Fix $1_E \in \om(E_V)(S_0)$. Note that we have an isomorphism
	\[
	g_{E_v} \colon s^* \om(E_V) \to b^* \om(E_V)
	\]
	restricting to $g_{M_V}$ on $s^* \om(M_V)$.
	Since $g$ is a tensor-isomorphism, it induces the trivial action on $\om(\mathbbm{1}_{\Tc})$. In particular, it sends $s^* 1_E \in \om(E_V)(T_s)$ to an element of $\om(E_V)(T_b)$ projecting to $1 \in \om(\mathbbm{1}_{\Tc})(T_b) = \Oc(T_b)$. In particular, we find that for $u_b \in b^*\om(\underline{U}^{\ab})(T)$,
	\[
	b^* a_V(u_b) = -g_{E_V}(s^* 1_E) + u_b(g_{E_V}(s^* 1_E))
	=
	- b^* 1_E + u_b(b^* 1_E).
	\]
	
	Letting $u_s \in s^* \om(\underline{U}^{\ab})(T)$, we find that
	\begin{eqnarray*}
		b^* a_V(g_{U^{\ab}}(u_s))
		&=&
		b^* a_V(g \star u_s)\\
		&=&
		-g_{E_V}(s^* 1_E) + (g \star u_s)(g_{E_V}(s^* 1_E))\\
		&=& - g_{E_V}(s^* 1_E) + g_{E_V}(u_s(s^* 1_E)),
	\end{eqnarray*}
	where the last step follows by commutativity of \ref{groupoid_action_on_bitorsor}(*). But
	\begin{eqnarray*}
		- g_{E_V}(s^* 1_E) + g_{E_V}(u_s(s^* 1_E))
		&=&
		g_{E_V}(- s^* 1_E + u_s(s^* 1_E))\\
		&=&
		g_{M_V}(- s^* 1_E + u_s(s^* 1_E))\\
		&=&
		g_{M_V}(s^* a_V(u_s)),
	\end{eqnarray*}
	so we are done.
	
\end{proof}

\subsection{Proof of Theorem \ref{thm:arith_frob_path}}\label{sec:phi_on_spl}

We apply Lemma \ref{lemma:Uab_V_Ext} to the $\Tc$ of Theorem \ref{thm:arith_frob_path} and $\om=\om_0$. We get an object $\underline{U}$ of $\Tc$ and an isomorphism $a \colon
\underline{U}^\ab \xrightarrow{\sim} 
\prod_{V \, \mathrm{ simple}} 
\Ext^1_{\Tc}(\mathbbm{1}_{\Tc}, V)^\vee \otimes V.$ Note that for $V$ simple, there is a unique $n \in \Zb$ for which $V = \gr^W_n V$, which we call its \emph{weight}. The weight-filtered property of $\Tc$ implies that $\Ext^1_{\Tc}(\mathbbm{1}_{\Tc}, V)=0$ unless $n<0$. This implies that $\underline{U}^{\ab}$ is concentrated in negative weights. Applying $\om_{\phi}$, we find that the action of $\phi$ on
\[
\om_0(\underline{U}^{\ab})
\]
has eigenvalues\footnote{We mean eigenvalues of any $\coeff'_0$-linear power of $\phi$.} in $w_{<0}(\coeff'_0) \coloneqq \bigcup_{m<0} w_m(\coeff'_0)$.

Let $\underline{\Uu}$ denote the completed universal enveloping algebra of $\underline{U}$ (\ref{prounip_in_T}), and let $\underline{\aA}$ denote its augmentation ideal. Then for each $n \in \Zb_{>0}$, $\underline{\aA}^{n}/\underline{\aA}^{n+1}$ is a quotient object of a tensor power of $U^\m{ab}=\underline{\aA}/\underline{\aA}^2$ in $\Tc$. In particular, since $w_{<0}(\coeff'_0)$ is closed under multiplication, all Frobenius eigenvalues acting on $\om_0(\underline{\aA}^{n}/\underline{\aA}^{n+1})$ and hence on $\om_0(\underline{\aA})$ are in $w_{<0}(\coeff'_0)$

Consequently, the proof of \cite[Theorem 3.1]{BesserColeman} shows that the map $U(\om_0) \to U(\om_0)$ given on points by
\[
g  \mapsto g\inv \phi(g)
\]
is an isomorphism. Specifically, the operator $S$ on $\om_0(\underline{\aA}^{n}/\underline{\aA}^{n+1})$ given by
\[
S(y) = \phi(y) - y
\]
considered in loc. cit. has eigenvalues of the form $\alpha-1$ for $\alpha \in w_{<0}(\coeff'_0)$; since $1 \in w_0(\coeff'_0)$, and the different $w_m(\coeff'_0)$ are disjoint (Definition \ref{defn:theory_of_weights}), the operator $S$ is invertible. In view of this, the proof goes through with no change. 

\subsubsection{}
\label{g2end}

We recall the $\underline{U}-\underline{U^{\gr}}$-bitorsor $\underline{\Spl^W}$ (\ref{groupoid_action_on_bitorsor}). Applying $\om_{\phi}$, we get a $\phi$-equivariant composition action
\[
\circ \colon
U(\om_0) \times \underline{\Spl^W}(\om_0) \times U^{\gr}(\om_0) \to \underline{\Spl^W}(\om_0).
\]

Consequently, the proof of \cite[Corollary 3.2]{BesserColeman} applies without change to show that $\underline{\Spl^W}(\om_0)$ possesses a unique $\coeff'_0$-point $\tannloop^\m{cr}_0$ such that $\phi(\tannloop^\m{cr}_0)=\tannloop^\m{cr}_0$.

\subsubsection{}\label{g2end2}

We must analyze the $\phi$-action on $\underline{\Spl^W_H}(\om_0)$ to show that $\tannloop^\m{cr}_0$ satisfies the conditions of the theorem. Let $E \in \Tc$. Recall there is an isomorphism \ref{rem:gr_commute}(*)
\[
\om_0^{\gr}(E) \simeq \om_0(\gr^W E).
\]

For $g \in \Gc(T)$, note that $g^{\m{ss}}_E = g_{\gr^W E}$ under this isomorphism. By the action of $\Gc$ on $\underline{\Spl^W}(\om_0)$ described in \ref{groupoid_action_on_bitorsor}, the map
\[\tag{**}
\underline{\Spl^W}(\om_0) \times \om_0^{\gr}(E) \to \om_0(E)
\]
is $\Gc$-equivariant. Since the action on $\om_0^{\gr}(E)$ is the usual $\Gc$-action on $\om_0(\gr^W E)$, this upgrades to a diagram
\[
\underline{\Spl^W} \times \gr^W E \to E
\]
in $\Tc$.

Applying $\om_{\phi}$, we get an action
\[
\om_{\phi}(\underline{\Spl^W}) \times \om_{\phi}(\gr^W E) \to \om_{\phi}(E)
\]
in the Tannakian category $\ph\m{Mod}(R)$. In other words, this says that (**) is $\phi$-equivariant. The action of $\phi$ on $\underline{\Spl^W}(\om_0)$ is the one that fixes $\tannloop^\m{cr}_0$, the action on $E$ is the $\phi$ from \ref{thm:arith_frob_path}(*), and the action on $\om_0^{\gr}(E) \simeq \om_0(\gr^W E)$ is the $\phi$-module structure on $\om_{\phi}(\gr^W E)$, which is none other than $\gr^W(\phi)$. Thus for $v \in \om_0^{\gr}(E)$,
\[
\phi(\tannloop^\m{cr}_0(\alpha))
=
(\phi(\tannloop^\m{cr}_0))(\gr^W(\phi)(\alpha))
=
\tannloop^\m{cr}_0(\gr^W(\phi)(\alpha)),
\]
as desired. \qed

\subsection{Hodge and Frobenius together}
\label{sec:hodge_frob_together}

%% **** What is the subscript F? Are the notations for the fields consistent?

Suppose $\Tc$ is a weight-filtered Tannakian category over $\coeff$ with an admissible Hodge-filtered Frobenius-equivariant fiber functor (\ref{defn:hodge_frob},\ref{rem:adm}) $\om_{H,\phi} = (\om_H,\om_{\phi},\eta)$ over $(\coeff',\coeff'_0,\phi)$ on $\Tc$. Let $U \coloneqq U^W$, $L^W$, $P_F$, and $F^0 U \coloneqq U^W_H$ be as in 
\S \ref{WH_groups_diagram}.

%% **** A little bit confusing about whether this is going to give something over k or k'. See also comment about this in the next section.

By Theorem \ref{thm:hodge-filtered_weight-splitting}, there exists a Hodge-filtered weight-splitting
\[
\tannloop^\m{H} \colon \om^{\gr} \to \om
\]
over $\coeff'$.

We refer to $\tannloop^\m{H}$ as an \emph{arithmetic Hodge path}. Arithmetic Hodge paths form a trivial torsor under $F^0 U$. By Theorem \ref{thm:arith_frob_path} there's a unique Frobenius-equivariant weight-splitting
\[
\tannloop^\m{cr}_0 \colon \om_0^{\gr} \to \om_0
\]
which we refer to as the \emph{arithmetic crystalline path}. We let
\[
\tannloop^\m{cr} \colon \om^{\gr} \to \om
\]
denote the path obtained from $\tannloop^\m{cr}_0$ by base-changing from $\coeff'_0$ to $\coeff'$ and transporting along $\eta$.

\begin{sDefinition}\label{defn:period_loops}
	In the situation above, we define the \emph{right-hand unipotent period loop of $\Tc$} by 
	\[
	\uU^R \coloneqq  \tannloop^\m{H} \circ (\tannloop^\m{cr})^{-1} \in U,
	\]
	We define the \emph{left-hand unipotent $\pf$-adic period loop of $\Tc$} by\footnote{Evidently these are inverse to one another; our terminology allows us to avoid picking a favorite.} 
	\[
	\uU^L \coloneqq  {\tannloop^\m{cr}} \circ (\tannloop^\m{H})^{-1} = (\uU^R)^{-1} \in U.
	\]
	
	When we wish to get rid of the dependence on choice of arithmetic Hodge path, we may regard $\uU^R$ as a $k'$-point of the homogeneous space $F^0U \backslash U$ and $\uU^L$ as a $k'$-point of the homogeneous space
	$
	U/F^0U.
	$
\end{sDefinition}

\subsubsection{\yesrem} While the terms `left' and `right' might seem reversed, we will see in Theorems \ref{thm:right_commutes} and \ref{thm:left_commutes} that they correspond to Bloch--Kato logarithms with values in $F^0 \backslash \omega(\pi)$ and $\omega(\pi) / F^0$, respectively.

%%%%%%%%%%%%%%%%%%
\section{Geometry of Localization-Realization Maps for Selmer Schemes, in an axiomatic setting}\label{sec:geom_loc_real}
%%%%%%%%%%%%%%%%%%

In this section, while staying in an abstract setting, we prove that the period loops of Definition \ref{defn:period_loops} have the desired relationship to the non-abelian Bloch--Kato Selmer maps of Chabauty--Kim.

\label{SectionGeom1}

\begin{sProposition}
	\label{prop:reductive_triv_cohom}
	Let $\GG$ be a reductive group over a field $\coeff$, let $\pi$ be a unipotent $\coeff$-group equipped with a $\GG$-action, and let $R$ be an arbitrary $\coeff$-algebra. Then $H^1(\GG_R; \pi_R) = \{*\}$.  
\end{sProposition}

\begin{proof}
	The category of quasi-coherent representations of $\GG_R$ is canonically equivalent to the category of quasi-coherent sheaves on the algebraic stack $B\GG_R$ over $\Spec R$; the equivalence interchanges group cohomology with cohomology of quasi-coherent sheaves. This helps to put the problem in a heavily studied context.\footnote{Our use of algebraic stacks here is frivolous: we could equally argue directly with injective resolutions in the category of quasi-coherent representations \`a la \cite{SGA3I}.}
	In particular, we have compatibility of cohomology with flat base-change by \cite[\href{https://stacks.math.columbia.edu/tag/0765}{Tag 0765}]{stacks-project}: if $V$ is a representation of $\GG$ over $\coeff$ then 
	\[
	H^i(\GG_R; V_R) = H^i(\GG; V) \otimes R = 0
	\]
	for $i >0$. Assuming for an induction on the level of unipotence that $H^1(\GG_R; \pi_R) = 0$ and that
	\[
	0 \to V \to \pi' \to \pi \to 0
	\]
	is an extension of unipotent $\coeff$-groups, we have a transitive action of $\{*\} = H^1(\GG_R; V_R)$ of the fibers of $H^1(\GG_R; \pi'_R) \to H^1(\GG_R; \pi_R) = \{*\}$ and the proposition follows. 
\end{proof}

\begin{ssLemma}\label{lemma:reductive_triv_cohom}
	If $\pi$ is a unipotent group in a weight-filtered Tannakian category $\Tc$ and $P$ a right $\pi$-torsor in $\Tc$, then $\gr^W P$ is a trivial $\gr^W \pi$-torsor.
\end{ssLemma}

\begin{proof}
	Since $\gr^W$ is an exact tensor functor from $\Tc$ to $\Tc^{\m{ss}}$, we may apply it to $P$ to get a $\gr^W \pi$-torsor in $\Tc^{\m{ss}}$. Since $\Tc^{\m{ss}}$ is semisimple, we may show by induction on the unipotence level of $\pi$ as in the proof of Proposition \ref{prop:reductive_triv_cohom} that $H^1(\Tc^{\m{ss}};\gr^W \pi)=0$. Thus $\gr^W P$ is a trivial $\gr^W \pi$-torsor in $\Tc^{\m{ss}}$ and hence also in $\Tc$.
\end{proof}

%%% We actually need this in Lemma \ref{h0.5} in order to deal with Frobenius

% It might be possible to prove an analogue for a $\coeff$-algebra $R$, at least if we assume $\Tc^{\m{ss}}$ is algebraic, using the same cohomology and base-change formula for the stack

\subsection{Prounipotent groups in a filtered Tannakian category}\label{sec:recoll_prounip_T}

Let $\Tc$ be a  strictly filtered Tannakian category over a field $\coeff$. We refer to \S \ref{algebraic_geometry_tann_cat} for basics on pro- and ind- categories and algebraic geometry in a Tannakian category.

\subsubsection{}\label{sec:pro_ind_gr-fin} We say that an object $\phi \colon I \to \Tc$ of $\Pro \Tc$ (resp. $\Ind \Tc$) for $I$ a cofiltered (resp. filtered) indexing category is \emph{graded-finite} if the the composed functor $\gr^W_n \circ \phi$ has a limit (resp. colimit) in $\Tc$ for all $n \in \Zb$. If $M \in \Pro{\Tc}$ (resp. $\in \Ind{\Tc}$), we can define the tensor algebra $T(M) \in \Pro{\Tc}$ (resp. $\in \Ind{\Tc}$) as the projective (resp. inductive) limit of the truncated tensor algebras $T^{\le n}(M) \coloneqq \bigoplus_{i=0}^n M^{\otimes i}$ under the obvious projection (resp. inclusion) maps.

\subsubsection{}\label{fib_func_pro_ind} For a graded-finite object $M$, $\om^{W,\gr}(M)$ is naturally in $\Gr(\LF(S))$.

\subsubsection{} There is a `duality' contravariant auto-equivalence on $\Gr(\LF(S))$ sending $\{M_i\}_{i \in \Zb} \in \Gr(\LF(S))$ to $\{M_{-i}^{\vee}\}_{i \in \Zb}$, and $\om^{W,\gr}(M)$ commutes with duality on graded-finite objects.

\subsubsection{} If $\pi$ is a prounipotent group in $\Tc$ (\ref{prounip_in_T}), we say that $\pi$ is \emph{graded-finite} if $\Lie \pi$ is. We say that $\pi$ is a \emph{unipotent group} if $\Lie \pi \in \Tc$.

\subsubsection{}\label{om_H_of_group_scheme} If $\om_H$ is a Hodge-filtered fiber functor over $\coeff'$ and $\pi$ a prounipotent group in $\Tc$, then we get a filtration $F^i \om(\Lie \pi)$. We denote by $F^i \om(\pi)$ the image of $F^i \om(\Lie \pi)$ under the exponential map $\exp \colon \om(\Lie \pi) \to \om(\pi)$ (which is an isomorphism of schemes over $\coeff'$).

\subsubsection{}\label{om_H_of_group_scheme_i>=0}

In the notation of \ref{om_H_of_group_scheme}, suppose $i \ge 0$. Then $F^i \om(\pi)$ is a subgroup scheme; in this case, we let $\om(\pi)/F^i$ (resp. $F^i \backslash \om(\pi)$) denote the homogeneous space $\om(\pi)/F^i \om(\pi)$ (resp. $F^i \om(\pi) \backslash \om(\pi)$).

Let $I_{\pi}$ denote the kernel of the counit $\Oc(\pi) \to \mathbbm{1}_{\Tc}$. By \cite[Lemme 7.6]{Deligne89}, $F^i\om(\pi)$ is defined by the ideal generated by $F^{1-i} \om(I(\pi))$. If $i=0$, note that this is the ideal generated by $F^1 \om(\Oc(\pi)) = F^1 \Oc(\om(\pi))$.

\begin{sDefinition}
	\label{defn:effective}
	For $\Tc$ as above, an \emph{effective\footnote{Here, `effective' accords with Voevodsky's conventions, by which the \textit{homological} motives functor lands in the triangulated category of \textit{effecti/ve motives}.} prounipotent $\Tc$-group} is a graded-finite prounipotent group $\pi$ in $\Tc$ such that $W_{-1} \Lie \pi =  \Lie \pi$ (equivalently, $W_0 \CoLie \pi = 0$).
	
	\subsubsection{\notreallyrem}\label{rem:effective}
	If $\pi$ is effective, then $W_{-i} (\Lie \pi)^{\otimes i} = (\Lie \pi)^{\otimes i}$ for $i \ge 0$, so \[T^{\le n+1}(\Lie \pi) \to T^{\le n}(\Lie \pi)\] induces an isomorphism on $\gr^W_{-i}$ for $n \ge i$. Therefore, $T(\Lie \pi)$ and hence $\Uc \pi$ and $\Oc(\pi)$ are graded-finite.
	
	More precisely, $W_{-1} T^{\le n}(\Lie \pi) = \oplus_{i=1}^n (\Lie \pi)^{\otimes i}$ and $W_0 T^{\le n}(\Lie \pi) = T^{\le n}(\Lie \pi)$, so \[T(\Lie \pi) = W_0 T(\Lie \pi) \twoheadrightarrow \gr^W_0 T(\Lie \pi) = \mathbbm{1}_{\Tc}\] is the counit. Since the kernel of $T(\Lie \pi) \to \Uc \pi$ is in $W_{-1} T(\Lie \pi)$, the same is true for $\Uc \pi$. By duality, we have $W_{-1} \Oc(\pi) = 0$ and $\gr^W_0 \Oc(\pi) = W_0 \Oc(\pi)$ is the image of the unit $\mathbbm{1}_{\Tc} \xrightarrow{\eta} \Oc(\pi)$.

\end{sDefinition}

\begin{ssProposition}
	\label{prop:effective_gr_point}
	
	Let $\pi$ be an effective unipotent $\Tc$-group, $\omega$ a fiber functor over $\coeff'$, and $R$ an arbitrary $\coeff'$-algebra. Suppose $\Tc$ is weight-filtered, so that $\Gb \coloneqq L^W(\om)$ is reductive. Let $P \in H^1(\Tc;\pi)(R)$. Then $\omgr(P)$ possesses a unique $\GG_R$-fixed $R$-point $\gamma^\m{gr}_R$. 
\end{ssProposition}

\begin{proof}
	
	By Proposition \ref{prop:reductive_triv_cohom}, $\omgr_R(P)$ is trivial when regarded as a $\GG_R$-equivariant $\omgr_R(\pi)$-torsor, so we may assume that $\omgr(P) = \omgr_R(\pi) = \omgr(\pi)_R$ ($\Gb_R$-equivariant $\omgr(\pi)_R$-torsors). By \ref{rem:effective}, $W_0 \Oc(\pi) = \mathrm{Im}(\eta)$, so $W_0 \Oc(\omgr(\pi)_R) = \mathrm{Im}(\omgr(\eta)_R) = R 1_{\Oc(\omgr(\pi)_R)}$.
	
	A $\GG_R$-fixed point corresponds to an $R$-algebra homomorphism from $\Oc(\omgr(\pi)_R)$ to $R=R_0$ respecting the $\GG_R$-action and thus the grading; in particular, it must kill $\Oc(\omgr(\pi)_R)_{>0}$. But $R 1_{\Oc(\omgr(\pi)_R)}$ surjects onto $\Oc(\omgr(\pi)_R)/\Oc(\omgr(\pi)_R)_{>0} = \gr_0^W \Oc(\omgr(\pi)_R)$, so such a homomorphism is determined by $R$-linearity and thus unique. (More precisely, it is the counit $\epsilon$ of $\Oc(\omgr(\pi)_R)$, corresponding to the identity of the group $\omgr(\pi)_R(R)$.)
\end{proof}

\subsubsection{\notreallyrem}\label{rem:unique_Gb_equiv} Proposition \ref{prop:effective_gr_point} implies there is a unique $\Gb_R$-equivariant isomorphism $\omgr(P) \simeq \omgr(\pi)_R$ of $\omgr(\pi)_R$-torsors. By applying this to $R=\coeff'$, we find that the trivialization in Lemma \ref{lemma:reductive_triv_cohom} is unique when $\pi$ is effective.

\begin{ssCorollary}\label{reductive_triv_cohom2}
	
	If $\pi$ is an effective prounipotent $\Tc$-group, then the conclusion of Lemma \ref{lemma:reductive_triv_cohom} (resp. Proposition \ref{prop:effective_gr_point}) holds (in particular, any torsor becomes trivial in $H^1(\GG_R;-)$).
	
\end{ssCorollary}

\begin{proof}
	
	Note that since $W_0 \Lie \pi = \Lie \pi$, the subobjects $W_{-n} \Lie \pi$ are Lie ideals in $\Tc$ for $n \ge 0$. Setting $\pi_n$ to be the unipotent group in $\Tc$ associated to $\Lie \pi / W_{-n-1} \Lie \pi$, we have $\pi = \varprojlim_{n} \pi_n$.
	
	Let $P$ be a right $\pi$-torsor (resp. $\pi_R$-torsor) in $\Tc$ (resp. $\Tc_R$), and let $P_n$ be the pushout of $P$ along $\pi \twoheadrightarrow \pi_n$ (resp. $\pi_R \twoheadrightarrow (\pi_n)_R$). Then $P = \varprojlim_n P_n$ as a scheme in $\Tc$ (resp. $\Tc_R$).
	
	By graded-finiteness, each $\pi_n$ is unipotent, so we may apply \ref{rem:unique_Gb_equiv} (resp. Proposition \ref{prop:effective_gr_point}) to find that $\gr^W P_n$ has a unique trivialization (resp. $\omgr(P_n)$ has a unique $\Gb_R$-fixed $R$-point). The uniqueness implies that these piece together to give a trivialization of $\gr^W P$ (resp. $\Gb_R$-fixed $R$-point of $\omgr(P)$). Then we can assume $\gr^W P = \gr^W \pi$ (resp. $\omgr(P)=\omgr(\pi)_R$) as $\gr^W \pi$-torsors in $\Tc^\m{ss}$ (resp. $\Gb_R$-equivariant $\omgr(\pi)_R$-torsors), and the rest of the proof proceeds as in Proposition \ref{prop:effective_gr_point} and \ref{rem:unique_Gb_equiv}.
\end{proof}

\subsection{Representability of cohomology schemes}\label{sssec:representability_of_coh_schemes}

Let us consider the association $R \mapsto H^1(\Tc;\pi)(R)$ defined on $\coeff$-schemes $R$ in \ref{cohomology_in_T}. Given a $\coeff$-map $R \to R'$, we get a map $(\Fib{\Tc})_{R'} \to (\Fib{\Tc})_{R}$ of $\coeff$-stacks, thus an associated pullback functor $\Tc_{R} \to \Tc_{R'}$, and thus a map $H^1(\Tc;\pi)(R) \to H^1(\Tc;\pi)(R')$. In particular, this association is a functor from $\coeff$-schemes to (pointed) sets. We wish to discuss when this is representable by a $\coeff$-scheme.

For an object $M \in \Tc$ and $i \in \Zb_{\ge 0}$, we have the functor $R \mapsto H^i(\Tc;M)(R)$.

\begin{ssLemma}\label{ssLemma:base_change_H^i}
Suppose $\Tc$ is neutral,\footnote{This assumption should be unnecessary; for example, it would follow by compatibility of cohomology with flat base-change for fpqc-algebraic stacks. However, \cite[\href{https://stacks.math.columbia.edu/tag/0765}{Tag 0765}]{stacks-project} proves it only for (fppf-)algebraic stacks.} and suppose $M$ is an object $\Tc$. Then $H^i(\Tc;M)(R) = H^i(\Tc;M) \otimes_{\coeff} R$. In particular, if $\dim_{\coeff} H^i(\Tc;M)<\infty$, then $R \mapsto H^i(\Tc;M)(R)$ is representable by the vector group associated to the vector space $H^i(\Tc;M)$.
\end{ssLemma}

\begin{proof}
Let $\om$ be a neutral fiber functor, $\om_R$ its base-change to $R$, and recall that
\[
H^i(\Tc;M)(R) = H^i(\pi_1(\Tc,\om_R);\om_R(M)).
\]

The right-hand-side is computed as the cohomology of the complex $C^{\bullet}(\pi_1(\Tc,\om_R);\om_R(M))$, where
\[
C^k(\pi_1(\Tc,\om_R);\om_R(M))
\]
is the set of homomorphisms of $R$-schemes from $\pi_1(\Tc,\om_R)^k$ to $\om_R(M)$. Thus
\begin{eqnarray*}
C^k(\pi_1(\Tc,\om_R);\om_R(M)) &=& \om_R(M) \otimes_R \Oc(\pi_1(\Tc,\om_R)^k)\\
&=& (\om(M) \otimes_{\coeff} R) \otimes_R \Oc(\pi_1(\Tc,\om)^k_R)\\
&=&
\om(M) \otimes_{\coeff} \Oc(\pi_1(\Tc,\om)^k) \otimes_{\coeff} R\\
&=&
C^k(\pi_1(\Tc,\om);\om(M)) \otimes_{\coeff} R.
\end{eqnarray*}

Since $R$ is flat over $\coeff$, tensoring with $R$ commutes with taking cohomology, so  we find that $H^i(\Tc;M)(R) = H^i(\Tc;M) \otimes_{\coeff} R$, and we are done.
\end{proof}

\subsubsection{}\label{functoriality_alg_grp_coh}

We make a remark about functoriality in $R$ of $C^{\bullet}(\pi_1(\Tc,\om_R);\om_R(M))$ and thus $H^k(\pi_1(\Tc,\om_R);\om_R(M))$. If $R'$ is an $R$-algebra (i.e., we have a homomorphism $R \to R'$ of $\coeff$-algebras), then $C^k(\pi_1(\Tc,\om_{R'});\om_{R'}(M))$, a priori the set of maps of $R'$-schemes from $\pi_1(\Tc,\om_{R'})$ to $\om_{R'}(M)$, may be identified with the set of maps of maps of $R$-schemes from $\pi_1(\Tc,\om_{R})$ to $\om_{R'}(M)$. As a result, the map of $R$-schemes with $G_R$-action $\om_{R}(M) \to \om_{R'}(M)$ defines a map of complexes $C^{\bullet}(\pi_1(\Tc,\om_{R});\om_{R}(M)) \to C^{\bullet}(\pi_1(\Tc,\om_{R'});\om_{R'}(M))$ and thus a map on $H^k$.

The same is true for $\pi$ in place of $M$, at least for $k \le 1$.

\begin{ssTheorem}\label{thm:representability}
Let $\Tc$ be a neutral weight-filtered Tannakian category and $\pi$ an effective pro-unipotent $\Tc$-group. Suppose that $\dim_{\coeff} \Ext^1_{\Tc}(\mathbbm{1}_{\Tc},\Gr^W_{-n}\Lie \pi)<\infty$ for all $n \in \Zb_{\ge 1}$. Then the functor on $\coeff$-algebras $R \mapsto H^1(\Tc;\pi)(R)$ is representable by an affine $\coeff$-scheme.
\end{ssTheorem}

\begin{proof}
The proof follows that of \cite[Theorem 2.2.1]{BettsWeightFil2023}, which itself follows that of \cite[Proposition 2]{kim05}.

%The compatibility of the Lie bracket on $\Lie \pi$ with the weight filtration implies that $W_{-n} \Lie \pi$ is a Lie ideal in $\Pro{\Tc}$, and we denote the corresponding normal $\Tc$-subgroup scheme by $W_{-n} \pi$.

We let $\om$ be a neutral fiber functor and set $G \coloneqq \pi_1(\Tc;\om)$. We set $U \coloneqq \om(\pi)$, $U_n \coloneqq \om(\pi_n)$, and $V_n \coloneqq \om(\gr^W_{-n} \pi) = \om(\gr^W_{-n} \Lie{\pi})$, all with $G$-action. The compatibility of the Lie bracket on $\Lie \pi$ with the weight filtration implies that $V_n$ is central in $U_n$ (in particular, abelian).

The proof proceeds by induction on $n$, showing that \[R \mapsto H^1(\Tc;\pi_n)(R) = H^1(G_R;(U_n)_R)\] is representable by a finite-type affine scheme and that $H^0(\Tc_R;(U_n)_R)=0$ for all $n$. For $n=0$, $U_n$ is trivial, so the result is immediate (it is represented by $\Spec{\coeff}$).

Suppose $n \ge 1$, and $R \mapsto H^1(G_R;(U_{n-1})_R)$ is representable by a finite-type affine $\coeff$-scheme. We apply \cite[Proposition 43]{SerreCohGal94} to the short exact sequence $0 \to (V_n)_R \to (U_n)_R \to (U_{n-1})_R \to 0$, noting as in \cite[Footnote 7]{BettsWeightFil2023} that the result still works because $U_n \to U_{n-1}$ and hence $(U_n)_R \to (U_{n-1})_R$ is algebraically split (and hence the map on chains $C^1(G_R;(U_n)_R) \to C^1(G_R;(U_{n-1})_R)$ is surjective). We get a long exact sequence
\begin{equation}\label{eqn:les_nab_coh}
\begin{split}
0 \to H^0(G_R;(V_n)_R) \to H^0(G_R;(U_n)_R) \to H^0(G_R;(U_{n-1})_R) \to \\
H^1(G_R;(V_n)_R)
\to H^1(G_R;(U_n)_R) \to H^1(G_R;(U_{n-1})_R) \xto{\del_R} H^2(G_R;(V_n)_R),
\end{split}
\end{equation}
with the usual interpretation given that $H^1(G_R;(U_n)_R)$ and $H^1(G_R;(U_{n-1})_R)$ are pointed sets rather than groups.

For a map $R \to R'$ of $\coeff$-algebras, we get a commutative diagram 
\[
\xymatrix{
0 \ar[r] & (V_n)_R \ar[r] \ar[d] & (U_n)_R \ar[r] \ar[d] & (U_{n-1})_R \ar[r] \ar[d] & 0\\
0 \ar[r] & (V_n)_{R'} \ar[r] & (U_n)_{R'} \ar[r] & (U_{n-1})_{R'} \ar[r] & 0
}
\]
that is $G_R$-equivariant, hence by \ref{functoriality_alg_grp_coh} induces a map of the corresponding exact sequences (\ref{eqn:les_nab_coh}).

We have $H^0(G;V_n) = \Ext^0_{\Tc}(\mathbbm{1}_{\Tc},V_n) = 0$ by Lemma \ref{lemma:strictness_crit} because $\mathbbm{1}_{\Tc}$ and $V_n$ are concentrated in degrees $0$ and $-n$, respectively, and $n \ge 1$. Thus $V_n^G = 0$ in a scheme-theoretic sense, hence  $(V_n)_R^{G_R} = H^0(G_R;(V_n)_R) = 0$. By the induction hypothesis, $H^0(G_R;(U_{n-1})_R)=0$, so by the long exact sequence, $H^0(G_R;(U_n)_R)=0$.

Since $H^0(G_R;(U_{n-1})_R)=0$,
\[
0 \to H^1(G_R;(V_n)_R)
\to H^1(G_R;(U_n)_R) \xto{\beta'} H^1(G_R;(U_{n-1})_R) \xto{\del_R} H^2(G_R;(V_n)_R)
\]
is exact. Thus
\[
0 \to H^1(G_R;(V_n)_R)
\to H^1(G_R;(U_n)_R) \xto{\beta} \del_R^{-1}(0) \to 0
\]
is exact for every $R$, in the sense that the abelian group $H^1(G_R;(V_n)_R)$ acts simply transitively on the fibers of the surjective map $\beta$.

We note, as in \cite[Corollary 2.2.5]{BettsWeightFil2023}, the functor $R \mapsto H^2(G_R;(V_n)_R)$ is, by Lemma \ref{ssLemma:base_change_H^i}, a subfunctor of a representable functor. More precisely, choosing a basis $\{b_i\}_{i \in I}$ of the $\coeff$-vector space $H^2(G;V_n)$, we have $H^2(G_R;(V_n)_R) \cong \bigoplus_{i \in I} R$ functorially in $R$, which is a subfunctor of the representable functor $\Ab^{H^2(G;V_n)} \colon R \mapsto \prod_{i \in I} R$ (which, as noted in loc.cit., is represented by $\Spec \Sym H^2(G;V_n)^{\vee}$). Thus $\del_R^{-1}(0)$ is the kernel of the composite $H^1(G_R;(U_{n-1})_R) \xto{\del_R} H^2(G_R;(V_n)_R) \hookrightarrow \Ab^{H^2(G;V_n)}$, which is a morphism of representable presheaves and thus a morphism of schemes. Furthermore, $R \mapsto \{0\} \subseteq \prod_{i \in I} R$ is a closed subscheme (defined by the ideal generated by $H^2(G;V_n)^{\vee}$ in $\Sym H^2(G;V_n)^{\vee}$), so $R \mapsto \del_R^{-1}(0)$ is a closed subscheme of $R \mapsto H^1(G_R;(U_{n-1})_R)$, hence representable by a finite-type affine scheme.

Let $A$ represent $R \mapsto \del_R^{-1}(0)$. Then we have an element $\id_A \in \del_A^{-1}(0) = \homo_{\coeff}(A,A)$. Let $\alpha \in H^1(G_A;(U_n)_A)$ map onto $\id_A$. Then $\alpha$ defines by Yoneda's Lemma a map of functors $\alpha \colon \del_R^{-1}(0) \to H^1(G_R;(U_n)_R)$ that is a section of $\beta$. Using the action $\star$ of $H^1(G_R;(V_n)_R)$ on $H^1(G_R;(U_n)_R)$, we define a map
\[
H^1(G_R;(V_n)_R) \times \del_R^{-1}(0) \to H^1(G_R;(U_n)_R)
\]
given by $(x,y) \mapsto x \star \alpha(y)$. Since $\alpha$ is a section of $\beta$, and $H^1(G_R;(V_n)_R)$ acts transitively, this defines an isomorphism $H^1(G_R;(V_n)_R) \times \del_R^{-1}(0) \simeq H^1(G_R;(U_n)_R)$ of functors from $\coeff$-algebras $R$ to sets. Since $R \mapsto H^1(G_R;(V_n)_R)$ is representable, finite-type, and affine by Lemma \ref{ssLemma:base_change_H^i} and the finite-dimensionality of $H^1(G;V_n) \cong \Ext^1_{\Tc}(\mathbbm{1}_{\Tc},\Gr^W_{-n}\Lie \pi)$, and we showed $R \mapsto \del_R^{-1}(0)$ is representable, finite-type, and affine, the same is true of $R \mapsto H^1(G_R;(U_n)_R)$.

Finally, as in the proof of Corollary \ref{reductive_triv_cohom2}, a right $\pi_R$-torsor in $\Tc_R$ is an inverse limit of $\pi_n$-torsors $P_n$, so it follows that $H^1(\Tc;\pi) = \varprojlim_n H^1(\Tc;\pi_n)$. An inverse limit (as presheaves) of affine schemes is again an affine scheme, so we are done.
\end{proof}

\subsubsection{Remark}\label{rem:finite-type_affine}

The proof shows more precisely that $H^1(\Tc;\pi)$ is finite-type when $\pi$ is finite-type and in general an inverse limit of finite-type affine schemes.

\subsubsection{Remark}\label{rem:indep_of_om} Although the proof of Theorem \ref{thm:representability} used a neutral fiber functor $\om$, the result is independent of $\om$. In particular, if we had two fiber functors $\{\om_j\}_{j=1,2}$ over $\coeff$ that were not isomorphic, the resulting $\coeff$-scheme structure on $H^1(\Tc;\pi)$ would be independent of $j$. Once again, we suspect Theorem \ref{thm:representability} to be true without the assumption that $\Tc$ is neutral.

\subsubsection{Remark} Without the assumption $\dim_{\coeff} \Ext^1_{\Tc}(\mathbbm{1}_{\Tc},\Gr^W_{-n}\Lie \pi)<\infty$, the functor $R \mapsto H^1(\Tc;\pi)(R)$ should still be \emph{ind-representable}, as is the case for infinite-dimensional vector spaces.

\subsection{Base-change of filtered fiber functors}\label{base_change_fil_gr}

In order to define the non-abelian Bloch--Kato logarithm as a map of schemes, we must define it for torsors with coefficients in an arbitrary $\coeff'$-scheme $S$. It follows that we must define the Hodge filtration on an arbitrary such torsor, and for this we must be able to base-change a filtered fiber functor from $\Tc$ to $\Tc_S$ for an arbitrary $\coeff$-scheme $S$. This has been done in \cite[Lemma 5.3]{ZieglerGrFil15} when $S$ is a finite field extension of $\coeff'$ and in \cite[Proposition 4.7(i)]{ZieglerFiltGeneralBase} when the Tannakian category in question is neutral.

If $\phi$ is a graded or filtered fiber functor over a field $\coeff'$ with underlying fiber functor $\omega$, and $S$ is a $\coeff$-scheme, we recall from \ref{base_change_fiber_functors} that we have a base-changed fiber functor $\om_{S_{\coeff'}/S} \colon \Tc_S \to \LF(S_{\coeff'})$, and given a map $S \to \Spec{\coeff'}$, we get $\om_{S/S} \colon \Tc_S \to \LF(S)$. We would like to do the same for $\phi$.% (see \ref{base_change_of_filtered} for more motivation).

\subsubsection{}\label{base_change_of_graded}

If $\phi$ is graded and we are given $S \to \Spec{\coeff'}$, we can base-change the natural map $(\Gm)_{\coeff'} \to \pi_1(\Tc,\om)$ of $\coeff'$-schemes to a map $(\Gm)_{S} \to \pi_1(\Tc,\om)_S$. Using the equivalence $\Tc_{S} \simeq \Rep(\pi_1(\Tc,\om)_S)$, we get a grading on $\om_{S/S}(M)$ for every $M \in \Tc_S$, which upgrades $\om_{S/S}$ to a graded fiber functor $\phi_{S/S} \colon \Tc_S \to \GrLF(S)$.

We describe an alternative, stack-theoretic version of this construction that will be necessary in \ref{base_change_of_filtered}. Note that $\GrLF(\coeff') = \LF((B\Gm)_{\coeff'})$, so that $\phi$ is the same as a fiber functor valued in $\LF((B\Gm)_{\coeff'})$, hence by Theorem \ref{theorem:tannakian_reconstr_groupoids} is isomorphic to the pullback functor associated to a map $(B\Gm)_{\coeff'} \xrightarrow{F_{\LF}^{-1}(\phi)} \Fib{\Tc}$ of $\coeff$-stacks.

As in \ref{base_change_fiber_functors} but with $X = (B\Gm)_{\coeff'}$, we can base-change along $S \to \Spec{\coeff}$ to get a map $(B\Gm)_{S_{\coeff'}} \to \Fib{\Tc}_S$ of $S$-stacks and thus a graded fiber functor
\[
\phi_{S_{\coeff'}/S} \colon \Tc_S \to \GrLF(S_{\coeff'})
.\] Given $S \to \Spec{\coeff'}$, we can base-change along $S \to S_{\coeff'}$ to get a map $(B\Gm)_{S} \to \Fib{\Tc}_S$ of $S$-stacks and thus a graded fiber functor
\[
\phi_{S/S} \colon \Tc_S \to \GrLF(S).
\]

\subsubsection{}\label{base_change_of_filtered}

Now suppose $\phi$ is filtered. By \cite[Example 2.13]{wedhorn2024extensionliftinggbundlesstacks}, for a ring $A$, the category $\LF([\Gm \backslash \Ab^1]_A)$ is the same as $\FilLF^{\fin}(A)$ (Definition \ref{filtrations_in_families}). In fact, it holds for a general scheme $S$ in place of $A$ with `finite projective' replaced by `finite locally free' in (a) and (b) of loc.cit. In particular, $\phi$ defines an exact tensor functor from $\Tc$ to $\LF([\Gm \backslash \Ab^1]_{\coeff'})$, whence by Theorem \ref{theorem:tannakian_reconstr_groupoids} a map $[\Gm \backslash \Ab^1]_{\coeff'} \xrightarrow{F_{\LF}^{-1}(\phi)} \Fib{\Tc}$ of $\coeff$-stacks.

As in \ref{base_change_of_graded}, we get a map $[\Gm \backslash \Ab^1]_{S_{\coeff'}} \to \Fib{\Tc}_S$ of $S$-stacks and thus a filtered fiber functor
\[
\phi_{S_{\coeff'}/S} \colon \Tc_S \to \FilLF(S_{\coeff'}).
\] Given $S \to \Spec{\coeff'}$, can base-change along $S \to S_{\coeff'}$ to get a map $[\Gm \backslash \Ab^1]_{S} \to \Fib{\Tc}_S$ of $S$-stacks and thus a filtered fiber functor
\[
\phi_{S/S} \colon \Tc_S \to \FilLF(S).
\]

Note that $\phi_{S/S}$ is compatible with base-change from $\Tc$ to $\Tc_S$, so it is the unique extension described in \cite[Proposition 4.7(i)]{ZieglerFiltGeneralBase} when $S=\Spec{R'}=\Spec{R''}$ is affine.

%%% The last sentence above is the same as applying the adjunction associated with base-change of stacks to $S \to \Fb{\Tc}$ to get a map $S \to \Fib{\Tc_S}=\Fib(\Tc)_S$.

%contained in $R$, we have $\om \in \operatorname{Fib}{\Tc}(\coeff') \to \operatorname{Fib}{\Tc}(R)=(\operatorname{Fib}{\Tc})_R(R)$, and by abuse of notation, we denote by $\om$ its image in the latter.

%%%% By the definition of a vector bundle on a stack, an object $M \in \Tc_S$ provides a functor $\operatorname{Fib}{\Tc}(S') = \operatorname{Fib}{\Tc_S}(S') \to \LF(S')$. Fixing the object $\om \in (\operatorname{Fib}{\Tc})_S(S)$, we get a functor $\om \colon \Tc_S \to \LF('s)$.

\subsection{Setup for Abstract Chabauty--Kim diagram}\label{sec:abstract_setup} For the rest of \S \ref{sec:geom_loc_real}, we let $\Tc$ denote a weight-filtered Tannakian category over $\coeff$, $\om_{H,\phi} = (\om_H,\om_{\phi},\eta)$ over $(\coeff',\coeff'_0,\phi)$ a Hodge-filtered Frobenius-equivariant fiber functor, $G$, $G^{\gr}$, $U \coloneqq U^W$, $L^W$, $P_F$, $F^0 U \coloneqq U^W_H$ as in \S \ref{WH_groups_diagram}, $\tannloop^\m{H}$, $\tannloop^\m{cr}_0$, and $\tannloop^\m{cr}$  as in \S \ref{sec:hodge_frob_together}. We set $\GG \coloneqq L^W = \pi_1(\Tc^{\m{ss}},\om) = \pi_1(\Tc^{\m{ss}},\om^{\gr})$, which is reductive.
%
%We assume for simplicity that $\phi$ acts trivially on $\coeff'$.\footnote{We don't expect this assumption to be necessary. However, in the applications we have in mind, this can always be achieved by replacing Frobenius by a power which fixes the base field.}
% 
We have the usual diagram of pro-algebraic $\coeff'$-groups (\ref{sec:graded_galois_groups_diagram}):
\[
\xymatrix{
	U^{\gr}
	\subseteq
	G^{\gr} \ar[r]_-{\m{ss}^\m{gr}}  & \GG 
	\ar@/_2ex/[l]_-{r} 
	& G
	\supseteq U
	\ar[l]^-{\m{ss}}.
}
\]
We have
\[
\uU^R \coloneqq  \tannloop^\m{H} \circ {\tannloop^\m{cr}}^{-1} \in U(\coeff')
\]
\[
\uU^L \coloneqq  \tannloop^\m{cr} \circ {\tannloop^\m{H}}^{-1} \in U(\coeff')
\]
as in Definition \ref{defn:period_loops}, which is uniquely determined modulo left (resp. right) multiplication by $F^0U$.

\subsection{Frobenius-invariant paths of torsors}
\label{21210A}
Let $\pi$ be an effective prounipotent $\Tc$-group, and let $P$ be a right $\pi$-torsor. The proof of \cite[Corollary 3.2]{BesserColeman} applies without change to show that $\om_0(P)$ possesses a unique $\phi$-fixed $\coeff'_0$-point $\gamma^\m{cr}_0$. The next lemma concerns the interaction between $\gamma_0^\m{cr} \in \om_0(P)(\coeff'_0)$ and the $\GG$-fixed point $\gamma_0^\m{gr} \coloneqq \gamma^{\gr}_{\coeff'_0}$ of Corollary \ref{reductive_triv_cohom2} through the Frobenius equivariant Tannakian path
\[
\tannloop_0^\m{cr} \colon \om_0^\m{gr} \to \om_0
\]
of Theorem \ref{thm:arith_frob_path}.

\begin{ssLemma}
	\label{h0.5}
	In the situation and the notation of \S \ref{21210A}
	\[
	\gamma_0^\m{gr} = (\tannloop_0^\m{cr})\inv (\gamma_0^\m{cr}) \in \omgr_0(P)(\coeff'_0).
	\]
\end{ssLemma}

\begin{proof}
	The associated graded $\gr^W(\pi)$ is an effective prounipotent group object of $\Tc^\m{ss}$, and $\gr^W(P)$ is a $\gr^W(\pi)$-torsor in $\Tc^\m{ss}$.
 %For $M \in \Tc$, the weight filtration in $\Tc$ induces a filtration on the filtered $\phi$-module $\om_{\phi}(M)$, and $\om_{\phi}(\gr^W M)$ is the associated graded in filtered $\phi$-modules. Thus
 The Frobenius automorphism on $\om_0(\gr^W(P)) \simeq \omgr_0(P)$ is the associated graded $\gr^W(\phi)$ of the Frobenius automorphism $\phi$ of $\om_0(P)$.
	
	It follows, on the one hand, that $\omgr_0(P)$ contains a unique $\coeff'_0$-point fixed by $\gr^W(\phi)$; by the commutativity of Diagram \ref{thm:arith_frob_path}(*), $(\tannloop^\m{cr}_0)\inv (\gamma^\m{cr}_0)$ is fixed by $\gr^W(\phi)$, hence it \textit{is} that point. On the other hand, by Corollary \ref{reductive_triv_cohom2}, it follows that $\gr^W(P)$ is trivial; i.e., there's an isomorphism of affine $\Tc$-schemes $\gr^WP \simeq \gr^W \pi$. Under such an isomorphism, the image of the unique Frobenius-fixed point of $\omgr_0(P)$ must again be the unique Frobenius-fixed point of $\omgr_0(\pi)$, namely the group identity $1$. Since $1$ is $\GG$-fixed, so is $(\tannloop^\m{cr}_0)\inv (\gamma^\m{cr}_0)$, and we are done.
\end{proof}

\subsubsection{}
\label{21223Aa}
Thus, if $\pi$ is an effective prounipotent
%\footnote{We switch here from working with prounipotent groups to unipotent groups, which will allow us to apply Proposition \ref{prop:reductive_triv_cohom}.}
$\Tc$-group and $P$ is a $\pi$-torsor, then $\tannloop_0^\m{cr}(\gamma_0^\m{gr})$ is the unique $\phi$-fixed $\coeff'_0$-point of $\om_0(P)$. For $R$ an arbitrary $\coeff'_0$-algebra and $P$ a $\pi_R$-torsor in $\Tc_R$, we \textit{define}
\[
\tag{*}
\gamma^\m{cr}_R \coloneqq \tannloop^\m{cr}(\gamma^\m{gr}_R) \in (\om_0)_{R/R}(P)(R).
\]

\subsection{Hodge paths of torsors}\label{hodge_paths_torsors}

Let $\pi$ be an effective prounipotent $\Tc$-group, $R$ a $\coeff'$-algebra, $P \in H^1(\Tc;\pi)(R)$, and $\gamma^\m{gr}_R$ as in Proposition \ref{prop:effective_gr_point} and Corollary \ref{reductive_triv_cohom2}. For a choice of $\tannloop^\m{H}$, we define
\[
\gamma_R^\m{H} \coloneqq \tannloop^\m{H}(\gamma_R^{\gr}) \in \om(P)(R).
\]

We spend the rest of \S \ref{hodge_paths_torsors} proving:

\begin{ssProposition}\label{prop:two_hodge_choices_differ}
	Any two choices of $\gamma_R^\m{H}$ differ by an element of $F^0 \om(\pi)(R)$.
\end{ssProposition}

In doing so, we will in fact define a subscheme $F^0 \om(P) \subseteq \om(P)$ and show that $\gamma_{\coeff'}^\m{H} \in F^0 \om(P)$.

\subsubsection{}\label{sec:F^0_of_affine_scheme}

Let $R$ be a ring, $A$ a filtered $R$-algebra,
%(i.e., algebra in $\Fil(\mathbf{Mod}_R)$ with the monoidal structure of \ref{defn:filtered_objects})
and $Y \coloneqq \Spec{A}$. In such a case, we define an $R$-subscheme
\[
F^0 Y \coloneqq \Spec{A/(A \cdot F^1 A)} \subseteq Y.
\]

Note that if $Y=\om(\pi)$ for $\pi$ a prounipotent group in $\Tc$, this notation is consistent with \ref{om_H_of_group_scheme_i>=0}.

\begin{ssLemma}\label{lemma:F^0_of_affine_scheme}
	Suppose $\{A_i\}_{i=1,2,3}$ are filtered $R$-algebras, $Y_i \coloneqq \Spec{A_i}$, and
	\[
	f \colon Y_1 \times_R Y_2 \to Y_3
	\]
	is a map for which $f^* \colon A_3 \to A_1 \otimes_R A_2$ is a filtered map. Then $f$ restricts to a map $F^0 f \colon F^0 Y_1 \times_R F^0 Y_2 \to F^0 Y_3$.
\end{ssLemma}

\begin{proof}
	
	Indeed, the image of $F^1 A_3$ lies in $\sum_{p+q=1} F^p A_1 \otimes_R F^q A_2$, which is contained in $F^1 A_1 \otimes_{R} A_2 + A_1 \otimes_{R} F^1 A_2$. The latter generates the ideal of the subscheme $F^0 Y_1 \times_R F^0 Y_2$ in $Y_1 \otimes_R Y_2$, so we are done.
\end{proof}

\subsubsection{}\label{F^0_of_schemes_in_T}

If $S$ is an arbitrary $\coeff$-scheme, then by \ref{base_change_of_filtered}, we can base-change $\om_H$ to a filtered fiber functor
\[
(\om_H)_{S_{\coeff'}/S} \colon \Tc_S \to \FilLF(S_{\coeff'}).
\]

\emph{From now through \ref{it_is_torsor_under_F^0}, if $M$ is an object of or affine scheme in $\Tc_S$, we let $\om(M)$ denote $\om_{S_{\coeff'}/S}(M)$; by abuse of notation, we refer to $\fil^n (\om_H)_{S_{\coeff'}/S}(M)$ as $F^i \om(M)$.}

By \ref{lemma:F^0_of_affine_scheme}, if $S$ is affine, $\{Y_i\}$ are affine schemes in $\Tc_S$, and we have a map $f \colon Y_1 \times Y_2 \to Y_3$ of affine schemes in $\Tc_S$, then we get an associated map of $S_{\coeff'}$-schemes:
\[
F^0 \om(f) \colon F^0 \om(Y_1) \times_{S_{\coeff'}} F^0 \om(Y_2) \to F^0 \om(Y_3).
\]

\subsubsection{}\label{F^0_of_torsor}

In particular, given a torsor action $a_P \colon P \times \pi_S \to P$ in $\Tc_S$, $\om(a_P)$ restricts to a map $F^0 \om(P) \times_{S_{\coeff'}} F^0 \om(\pi_S) \xto{F^0 \om(a_P)} F^0 \om(P)$. Since $F^0 \om(\pi_S) = (F^0 \om(\pi))_S$ is a subgroup scheme, we thus have an action of $F^0 \om(\pi_S)$ on $F^0 \om(P)$ over $S_{\coeff'}$.

\begin{ssLemma}\label{lemma:F^0_is_torsor}
$F^0 \om(P)$ is an $F^0 \om(\pi_S)$-torsor over $S_{\coeff'}$.
\end{ssLemma}

The proof spans \ref{sssec:transitive_suffices}-\ref{it_is_torsor_under_F^0}. First, we show that this lemma implies Proposition \ref{prop:two_hodge_choices_differ}.

\subsubsection{Proof that Lemma \ref{lemma:F^0_is_torsor} $\implies$ Proposition \ref{prop:two_hodge_choices_differ}}\label{proof_of_prop:two_hodge_choices_differ}

We return to notation as in Proposition \ref{prop:two_hodge_choices_differ}, in which $R$ is a $\coeff'$-algebra, $P$ a torsor under $\pi_R$ in $\Tc_R$, and $\om(P)$ (resp. $\om(\pi_R)$) refers to $\om_{R/R}(P)$ (resp. $\om_{R/R}(\pi_R)$).

Corollary \ref{reductive_triv_cohom2} provides a unique $\Gb_R$-trivialization
	\[
	\omgr(P) \simeq \omgr(\pi_R).
	\]
Under this isomorphism, $\gamma^{\gr}$ corresponds to the identity of $\omgr(\pi_R) \simeq \omgr(\pi)_R \simeq \om(\gr^W \pi)_R$, which is necessarily in $F^0\om(\gr^W \pi)$ because it is $\exp(0)$.

Then $\tannloop^\m{H} \colon \omgr(P) \to \om(P)$ respects the Hodge filtration, so in particular, it sends $F^0 \omgr(P)$ to $F^0 \om(P)$. Thus $\gamma_R^\m{H} \in F^0 \om(P)$ for any choice of $\tannloop^\m{H}$.

By \ref{F^0_of_torsor}, $F^0 \om_{R_{\coeff'}/R}(P)$ is a torsor under $F^0 \om_{R_{\coeff'}/R}(\pi_R)$. We can base-change along $R_{\coeff'} \to R$ to see that $F^0 \om(P) = F^0 \om_{R/R}(P)$ is a torsor under $F^0 \om(\pi_R)$. In particular, any two choices of $\gamma_R^\m{H}$ differ by an element of $F^0 \om(\pi_R)$, and we are done.

\subsubsection{Proof of Lemma \ref{lemma:F^0_is_torsor}}\label{sssec:transitive_suffices}

We already know the action $F^0 \om(P) \times_{S_{\coeff'}} F^0 \om(\pi_S) \to F^0 \om(P)$ has trivial stabilizers, as the same is true of $\om(P) \times_{S_{\coeff'}} \om(\pi_S) \to \om(P)$. Therefore, it suffices to show the action is transitive. We devote \ref{inverse_torsor_map}-\ref{it_is_torsor_under_F^0} to showing this transitivity.

\subsubsection{}\label{inverse_torsor_map}

Let $Q$ be a (pseudo-)torsor under a group scheme $H$ over a scheme $X$ with action map $a_Q$. We can define the \emph{difference map} $d_Q \colon Q \times_X Q \to H$ as follows:

For an $X$-scheme $T$ and $a,b \in P(T)$, there is a unique $g \in \pi(T)$ such that $ag = b$. We let the image of $(a,b) \in (P \times_X P)(T)$ be $g \in \pi(T)$.

The difference map is the unique map for which the composition
\[
Q \times_X Q \xrightarrow{\Delta_Q \times \id_Q} (Q \times_X Q) \times_X Q = Q \times_X (Q \times_X Q)
\xrightarrow{\id_Q \times d_Q} Q \times_X H \xrightarrow{a_Q} Q
\]
is projection onto the second factor. 

\subsubsection{}\label{inverse_torsor_map_functoriality}

We mention two functoriality properties of the difference map. First, if $X'$ is an $X$-scheme, and $Q' \coloneqq Q_{X'}$ is the base-changed torsor under $H' \coloneqq H_{X'}$ over $X'$, then $d_{Q'}$ is the base-change of $d_Q$ to $X'$. This follows simply from the description by taking $T$ to be an $X'$-scheme instead of an $X$-scheme. Second, if $Q_i$ is a torsor under $H_i$ for $i=1,2$ and $e_1 \colon H_1 \to H_2$ is an isomorphism of group schemes over $X$ and $e_2 \colon Q_1 \to Q_2$ is an isomorphism that intertwines with $e_1$ and the actions, then the following diagram commutes:
\[
\xymatrix{
Q_1 \times_X Q_1 \ar[d]_-{e_2 \times e_2} \ar[r]^-{d_{Q_1}} & H_1 \ar[d]^-{e_1}\\
Q_2 \times_X Q_2 \ar[r]^-{d_{Q_2}} & H_2
}
\]

\subsubsection{}\label{inverse_torsor_map_in_on_fibers}

%%**** is $R$ over $\coeff$ or $\coeff'$? We can take torsors for the former, but then there's a question about whether $\om$ is taken only for $R$ over $\coeff'$ or whether it outputs stuff over $R \otimes_{\coeff} \coeff'$. And compare with the $\omega_{-/-}$ notation

Let $P$ be as in \ref{F^0_of_torsor}. Then \ref{inverse_torsor_map} for $H=\om(\pi_R)$, $Q=\om(P)$, and $X=S_{\coeff'}$ provides a map of $S_{\coeff'}$-schemes $d_{\om(P)} \colon \om(P) \times_{S_{\coeff'}} \om(P) \to \om(\pi_S)$. We wish to show in \ref{inverse_torsor_map_in_T} that the difference map can be upgraded to a map $P \times P \to \pi_S$ of $\Tc_S$-schemes.

\subsubsection{}\label{upgrading_to_map_in_T}
As described in \cite[5.11]{Deligne89} (in the case $S$ is a field), to upgrade this to a map in $\Tc_S$, it suffices to define such a map for every $\om' \in \Fib(\Tc_S)(S')$ for an $S$-scheme $S'$, compatible with base-change and isomorphisms between fiber functors. In fact, the same is true for general $S$ by the description of affine schemes in $\Tc_S$ as stacks affine over $\Fib(\Tc_S)$ in \ref{alg_geo_stacks}.

\subsubsection{}\label{inverse_torsor_map_in_T}

Given $\om'$, we have a torsor action in $S'$-schemes $\om'(P) \times_{S'} \om'(\pi_S) \to \om'(P)$. We define $(d_P)_{\om'}$ as the inverse torsor map associated to $H=\om'(\pi_R)$, $Q=\om'(P)$, and $X=S'$.

If $\om' \in \Fib(\Tc_S)(S')$ and $S''$ is an $S'$-scheme, let $\om''$ be the image of $\om'$ in $\Fib(\Tc_S)(S'')$. Then the torsor action $\om''(P) \times_{S''} \om''(\pi_S) \to \om''(P)$ is just the base-change along $S'' \to S'$ of $\om'(P) \times_{S'} \om'(\pi_S) \to \om'(P)$. It follows by \ref{inverse_torsor_map_functoriality} that $(d_P)_{\om''}$ is the base-change along $S'' \to S'$ of $(d_P)_{\om'}$.

If $\om'_1,\om'_2 \in \Fib(\Tc_S)(S')$ and $g \colon \om'_1 \to \om'_2$ is an isomorphism, we get isomorphisms $g(\pi_S) \colon \om'_1(\pi_S) \to \om'_2(\pi_S)$ and $g(P) \colon \om'_1(P) \to \om'_2(P)$ that intertwine with the torsor action, hence a commutative diagram
\[
\xymatrix{
\om'_1(P) \times_{S'} \om'_1(P) \ar[d]_-{g(P) \times g(P)} \ar[r]^-{(d_P)_{\om'_1}} & \om'_1(\pi_S) \ar[d]^-{g(\pi_S)}\\
\om'_2(P) \times_{S'} \om'_2(P) \ar[r]^-{(d_P)_{\om'_2}} & \om'_2(\pi_S)
}
\]

We have thus shown the compatibilities required by \ref{upgrading_to_map_in_T}, so we have an inverse torsor map $d_P \colon P \times P \to \pi_S$ in $\Tc_S$.

\subsubsection{}\label{it_is_torsor_under_F^0}

Assuming now $S$ is affine, we may apply $(\om_H)_{S_{\coeff'}/S}$ to $d_P \colon P \times P \to \pi_S$ to get a filtered map
\[
(d_P)_{\om} \colon \om(P) \times \om(P) \to \om(\pi_S) = \om_{S_{\coeff'}/S}(\pi_S)
\]
and thus by Lemma \ref{lemma:F^0_of_affine_scheme}, a map of $S_{\coeff'}$-schemes
\[
F^0 (d_P)_{\om} \colon F^0 \om(P) \times F^0 \om(P) \to F^0 \om(\pi_S)
\]
for which the composition
\begin{equation}\label{eqn:F^0_d_comp}
F^0 \om(P) \times F^0 \om(P) \xrightarrow{\Delta \times \id} F^0 \om(P) \times F^0 \om(P) \times F^0 \om(P)
\xrightarrow{\id \times F^0 (d_P)_{\om}} F^0 \om(P) \times F^0 \om(\pi_S) \xrightarrow{(a_P)_{\om}} F^0 \om(P)
\end{equation}
is projection onto the second factor.

Transitivity now follows: let $S'$ be an $S_{\coeff'}$-scheme, $\om'_H$ the base-change of $(\om_H)_{S_{\coeff'}/S}$ to $S'$, and $p_1,p_2 \in F^0 \om'(P)$. Then base-changing (\ref{eqn:F^0_d_comp}) shows that $F^0 (d_P)_{\om}(p_1,p_2)$ is an element of $F^0 \om'(\pi_S)$ sending $p_1$ to $p_2$, so the action is indeed transitive, and we are done.

\subsection{Abstract Chabauty--Kim diagram: objects}
\label{sec:abstract_CK_diagram}
Let $\pi$ be an effective prounipotent $\Tc$-group. We construct morphisms $\be^R$, $\de$, $\tau^\m{cr}$, $\ev_{\uU^R}/F^0$ of \textit{pointed set}-valued functors\footnote{Under the conditions of Theorem \ref{thm:representability}, they are all representable by affine $\coeff'$-schemes, but we will not need that for the rest of \S \ref{sec:geom_loc_real}.} on $\coeff'$-algebras as in the following diagram
\[
\tag{*}
\xymatrix{
	H^1 (\Tc; \pi)_{\coeff'} 
	\ar[r]^-{\be^R} \ar[d]_\de^\simeq &
	\om(\pi)/F^0 
	\\
	Z^1
	\big(
	U^{\gr}; \omgr(\pi)
	\big)^{\GG}
	\ar[r]_-{\tau^\m{cr}} 
	& 
	Z^1 \big( U; \om(\pi) \big)^{\GG_{\cri}}.
	\ar[u]_-{\ev_{\uU^R}/F^0}
}
\]
%Throughout \ref{sec:abstract_CK_diagram}-\ref{thm:right_commutes} we default to using the right period points $\uU$ and we allow ourselves to drop the decoration $R$ to avoid clouding the notation. 

We define $H^1(\Tc; \pi)$ as in \ref{cohomology_in_T} and $\om(\pi)/F^0$ is as in \ref{om_H_of_group_scheme}. For a $\coeff'$-algebra $R$, we define
\[
Z^1\big(U(\omgr); \omgr(\pi)\big)(R) = 
Z^1\big(U(\omgr)_R; \omgr(\pi)_R\big)
\]
and
\[
Z^1\big(U(\om); \om(\pi)\big)(R) = 
Z^1\big(U(\om)_R; \om(\pi)_R\big).
\]
The $\GG$-fixed subsheaf $Z^1\big(U(\omgr); \omgr(\pi)\big)^\GG$ is defined in the usual way: an $R$-point $c$ of $Z^1\big(U(\omgr); \omgr(\pi)\big)$ is said to be fixed by $\GG$ if it’s fixed by the action of $\GG(R')$ on $Z^1\big(U(\omgr); \omgr(\pi)\big)(R')$ for every $R$-algebra $R'$. We let $\GG_{\cri}$ denote the image of $\Gb \subseteq G^{\gr}$ under the transport-of-structure isomorphism $\tau^{\cri} \colon G^{\gr} \xto{\sim} G$ induced by $\tannloop^{\cri} \colon \omgr \xto{\sim} \om$.

\subsection{Abstract Chabauty--Kim diagram: maps}
\label{sec:abstract_CK_maps}
Having defined the vertices of \ref{sec:abstract_CK_diagram}(*) we now define the arrows. Let $R$ be a $\coeff'$-algebra and let
\[
P \in H^1 (\Tc; \pi)_{\coeff'}(R).
\]
We define
\[
\de P \colon U(\omgr)_R \to \omgr(\pi)_R
\]
on $u \in U(\omgr)(R')$ ($R'$ an arbitrary $R$-algebra) by 
\[
(\de P)(u) = (\gamma^\m{gr}_{R'})\inv \cdot u(\gamma^\m{gr}_{R'});
\]
note that $\de$ is an isomorphism by \cite[Theorem A.4]{CorwinTSV}.
We choose $\gamma \in F^0 \om(P)$ arbitrary and define
\[
\be^R (P) \coloneqq  (\gamma^\m{cr}_R)\inv \cdot \gamma \in \big(\om(\pi)/F^0 \big)(R).
\]
The map $\tau^\m{cr}$ is transport of structure induced by $\tannloop^\m{cr}$. The map $\ev_{\uU^R}/F^0$ is given by evaluation at $\uU^R$ followed by projection $\om(\pi) \twoheadrightarrow \om(\pi)/F^0$.

\begin{sTheorem}
	\label{thm:right_commutes}
	In the situation and the notation of \ref{sec:abstract_setup}, Diagram \ref{sec:abstract_CK_diagram}(*) commutes. 
\end{sTheorem}

\begin{proof}
	Fix $R$ an arbitrary $\coeff'$-algebra and $P \in H^1(\Tc; \pi)(R)$.
	\begin{align*}
\ev_{\uU^R}/F^0(\tau^\m{cr}(\delta(P)))
		&=&
		(\tau^\m{cr}\de P)(\uU^R) & \pmod{F^0}\\
        &=&
		\tannloop^\m{cr} \big(\de P(\tau^{\cri}(\uU^R))\big) & \pmod{F^0}\\
		&=&
		\tannloop^\m{cr} \big(\de P((\tannloop^{\m{cr}})\inv \circ \uU^R \circ \tannloop^\m{cr})\big) & \pmod{F^0}\\
		&=&
		\tannloop^\m{cr} \big(\de P((\tannloop^{\m{cr}})\inv \circ \tannloop^\m{H} \circ {\tannloop^\m{cr}}^{-1} \circ \tannloop^\m{cr})\big) & \pmod{F^0}\\
		&=&
		\tannloop^\m{cr} \big(\de P((\tannloop^{\m{cr}})\inv \circ \tannloop^\m{H})\big) & \pmod{F^0}\\
		&=&
		\tannloop^\m{cr} \big((\gamma^\m{gr}_R)\inv \cdot (\tannloop^{\m{cr}})\inv(\tannloop^\m{H}(\gamma^\m{gr}_R))\big) & \pmod{F^0}\\
		&=& \tannloop^\m{cr} \big((\gamma^\m{gr}_R)\inv \cdot (\tannloop^{\m{cr}})\inv(\ga^\m{H}_R)\big) & \pmod{F^0}\\
		&=& (\gamma^\m{cr}_R)\inv \cdot \ga^\m{H}_R & \pmod{F^0}\\
		&=& \be^R(P),
	\end{align*}
    the last step because $\gamma = \ga^\m{H}_R \in F^0 \om(P)$ by Proposition \ref{prop:two_hodge_choices_differ}.
\end{proof}

\subsubsection{\yesrem}
	The theorem implies that $\ev_{\uU^R} \circ \tau^{\cri} \pmod{F^0}$ is independent of the choice of $\uU^R \in U(\coeff')$.

\subsection{Remark on splittings}
\label{right_drawback}
One downside of this approach is that the splitting given by $\tau^{\cri}$ does not agree with the splitting used in the setting of mixed Tate motives (e.g., in \cite{MTMUE}, \cite{PolGonI}), in which the weight filtration is canonically split by the Hodge filtration.

Alternatively, we may replace $\tau^\m{cr}$ by $\tau^\m{H}$; in fact, we can then circumnavigate $\tau^\m{H}$ altogether. Our goal in \S \ref{left_variant}-\ref{thm:left_commutes} is to record this variant.

\subsection{Left-handed variant}\label{left_variant}
With notation as in \S \ref{sec:abstract_CK_maps}-\ref{sec:abstract_CK_diagram}, we define the \emph{(left) localization-realization map}
\[
\m{LocReal}=\be^L \colon H^1(\Tc; \pi)_{\coeff'} \to 
F^0\backslash \om(\pi)
\]
by choosing $\gamma \in F^0 \om(P)$ and setting
\[
\be^L(P) = (\gamma)\inv \cdot \gamma^\m{cr}_R.
\]
Letting $\tau^\m{H}$ denote transport-of-structure along $\tannloop^\m{H}$, $\Gb_H$ the image of $\Gb \subseteq G^{\gr}$ under $\tau^\m{H} \colon G^{\gr} \xto{\sim} G$, and $\uU^L_{\gr} \coloneqq (\tannloop^\m{H})\inv \circ \tannloop^\m{cr}$, we get a diagram
\[
\tag{*}
\xymatrix{
	H^1(\Tc; \pi)_{\coeff'}
	\ar[r]^-{\de}_-\simeq
	\ar@/_9ex/[drr]_{\be^L}
	&
	Z^1\big( U^{\gr}; \omgr(\pi) \big)^\GG
	\ar[r]^-{\tau^\m{H}}_-\simeq
	\ar[d]_{\ev_{\uU_\m{gr}^L}}
	&
	Z^1\big( U; \om(\pi) \big)^{\GG_H}
	\ar[d]^{\ev_{\uU^L}}
	\\
	& F^0 \backslash \omgr(\pi)
	\ar[r]_-{\tannloop^\m{H}}^-\simeq
	&
	F^0 \backslash \om(\pi).
}
\]

\begin{sTheorem}
	\label{thm:left_commutes}
	Diagram \ref{left_variant}(*) commutes.
\end{sTheorem}

\begin{proof}
	This is similar to \ref{thm:right_commutes}. If $c^\m{gr} \in Z^1\big( U^\m{gr}; \pi^\m{gr} \big)^\GG(R)$, then $c\coloneqq\tau^\m{H}(c^\m{gr})$ is given by the composition indicated in the diagram below left. 
	\[
	\xymatrix{
		U \ar[r]^c \ar[d]_{(\tannloop^\m{H})\inv(\cdot)\tannloop^\m{H}}
		& \om(\pi)
		&\uU^L = \tannloop^\m{cr} \circ (\tannloop^\m{H})\inv
		\ar@{|->}[d]
		& (\gamma^\m{H})\inv \cdot \gamma^\m{cr}_R
		\\
		U^\m{gr} 
		\ar[r]_{c^\m{gr}}
		&
		\omgr(\pi)
		\ar[u]_{\tannloop^\m{H}}
		& \uU_\m{gr}^L
		\ar@{|->}[r]
		& (\gamma^\m{gr}_R)\inv\big( (\tannloop^\m{H})\inv \tannloop^\m{cr} \big)(\gamma^\m{gr}_R)
		\ar@{|->}[u]
	}
	\]
	Suppose now that $c = \de(P)$ comes from a $\pi_R$-torsor $P \in H^1(\Tc; \pi)(R)$. Then, as depicted above right, we have
	\[
	\ev_{\uU^L} \big(\tau^\m{H}(\de P) \big) =
	c(\uU^L) = \tannloop^\m{H}(c^{\gr}(\uU^L_{\gr})) = (\gamma^\m{H}_R)\inv \cdot \gamma^\m{cr}_R \equiv \be^L(P) \pmod{F^0}.
	\]
\end{proof}

\subsubsection{\yesrem}
	As before, the theorem implies that $\ev_{\uU^L} \circ \tau^\m{H} = \tannloop^\m{H} \circ \ev_{\uU^L_{\gr}} \pmod{F^0}$ is independent of the choice of $\uU^L \in U(\coeff')$.

\subsubsection{\yesrem}\label{left_vs_right} The right-handed variant (\S \ref{sec:abstract_CK_maps}-\ref{thm:right_commutes}, given by $\beta^R(P) = (\gamma^{\cri})^{-1} \cdot \gamma^\m{H}$, is the original variant of Kim (\cite[p.105]{kim09}). As pointed out to us by M. L\"udtke, it agrees with the original Bloch--Kato logarithm of \cite{BlochKato} in the abelian case.

On the other hand, the left-handed variant (\S \ref{left_variant}-\ref{thm:left_commutes}), given by $\beta^L(P) = (\gamma^\m{H})^{-1} \cdot \gamma^{\cri}$, is the more recent trend in Chabauty--Kim theory, found e.g., in \cite[Lemma 42]{BesserHeidelberg}, \cite[1.3]{MTMUE}, \cite[Proof of Proposition 3.2]{Sakugawa_nabBK17}, and \cite[Remark 3.1.13,\S 4.2.1]{BettsWeightFil2023}. As pointed out to us by M. L\"udtke, this is preferably because it leads to Coleman integrals from the basepoint to the varying point rather than the other way around.

\begin{sRemark}
	If $M \in \Tc^{ss}$ is semisimple, and $\om^W$ is a weight-filtered fiber functor, then $\om^W(M)$ is canonically split. Moreover, any splitting of $\om^W$ induces the canonical splitting on $\om^W(M)$. In particular, in the setting of paragraph \ref{left_variant}, if $\pi$ is semisimple, then the isomorphism between the graded and ungraded versions of the $p$-adic Albanese variety induced by $\tannloop^\m{H}$ (arrow labeled $\tannloop^\m{H}$ in Diagram \ref{left_variant}(*)) is just the canonical splitting of the weight filtration and thus independent of choice of $\tannloop^\m{H}$.
\end{sRemark}

%%%%%%%%%%%%%%%%%%%%%%%
%%%%%%%%%%%%%%%%%%%%%XXXXXXXXXXXXXXXX
%%%%%%%%%%%%%%%%%%%%%%%%%%

% Material on motivic iterated integrals removed after version 0.55, see folder pPerSelArt_V055_backup230712_last

\section{\texorpdfstring{$p$}{p}-Adic Periods for Motivic Structures and Galois Representations}\label{sec:p-adic_periods_for_motivics_galoisr}

In this section, we apply the abstract considerations of \S \ref{sec:arith_hodge_path}-\ref{sec:geom_loc_real} to the categories of arithmetic interest mentioned in \S \ref{sec:examples_cats}. We discuss in particular our $p$-adic period conjecture and the relation to syntomic regulators. We leave out all discussion of non-abelian cohomology and the Chabauty--Kim method, as those are the subject of \S \ref{sec:ck_application}.

%\commentDC{DC: Why is the rest of Section 5 in Section 5? I thought the idea is to first prove an abstract theorem about a weight-filtered Tannakian category with Frobenius and Hodge fiber functor, THEN it applies equally to the various categories of Section 3. I thought maybe this material should go in a new section (new Section 6? Section 8?) along with motivic iterated integrals and the period conjecture (which I figured was the main contribution of this article other than to Chabauty--Kim). ID: Seems ok to me as is.}

%\commentDC{This new sectioning solves the issue mentioned above}

\subsection{Categories of motivic structures}

\begin{sDefinition}
	\label{category_of_structures}
	Let $\basefield$ be a number field. We define a \emph{category of motivic structures over $\basefield$} (resp. \emph{category of weight-semisimple strongly geometric $p$-adic Galois representations}) to be a full thick Tannakian subcategory $\Tc$ of the category $\mathbf{SM}_\basefield(\Qb)$ (resp. $\Rep_{\Qp}^{\mathrm{sg}}(G_{\basefield})^{\wss}$).
	
	Let $S \subseteq S_f(\basefield)$. We define a \emph{category of lisse motivic structures over $\Oc_{\basefield,S}$} (resp. \emph{category of lisse $p$-adic Galois representations over $\Oc_{\basefield,S}$}) to be a full thick Tannakian subcategory $\Tc$ of the category $\mathbf{SM}_{\basefield,S}(\Qb)$ (Definition \ref{defn:SM_beta}) (resp. $\Rep_{\Qp}^{\mathrm{sf,S}}(G_{\basefield})^{\wss}$).
	
	The thickness property ensures that if $M \in \Tc \subseteq \Cc$ where \[\Cc \in \{\mathbf{SM}_\basefield(\Qb),\mathbf{SM}_{\basefield,S}(\Qb),\Rep_{\Qp}^{\mathrm{sg}}(G_{\basefield})^{\wss},\Rep_{\Qp}^{\mathrm{sf,S}}(G_{\basefield})^{\wss}\},\] then
	\[
	\Ext^1_\Tc(\Qb(0),M) = \Ext^1_{\Cc}(\Qb(0),M).
	\]
	
	%Note that thickness is not an obstacle to the existence of such $\Tc$, because taking the thick envelope of a full Tannakian subcategory does not change the pro-reductive part of $\pi_1(\Tc, \dR)$ and therefore does not change the component group.

\end{sDefinition}

%\begin{sRemark}
%In view of Olsson \cite{OlssonTowards}, Remark 1.13 of Deligne \cite{Deligne89} no longer applies, and the category of systems might accordingly be upgraded somewhat to include not only the crystalline Frobenii as a `` bizarre appendix'', but also comparison isomorphisms \`a la $p$-adic Hodge theory. For our purposes this is unnecessary, and we will not dwell on this foundational aspect here. 
%\end{sRemark}

\begin{ssRemark}
	A category of systems is evidently weight-filtered by \ref{rem:filtr_sub}.
\end{ssRemark}

\begin{ssExample}
	\label{category_of_lisse_structures}
	For $A$ an abelian variety over $\basefield$ with good reduction outside $S$, let $\Tc \coloneqq \mathbf{SM}_{\Oc_{\basefield,S}}(\Qb,A)$ (\ref{sssec:SM_A,S}). Then $\Tc$ is a category of lisse motivic structures over $\Oc_{\basefield,S}$.
\end{ssExample}

\subsubsection{Notation}\label{notation:GZ_UZ}

For a category $\Tc$ of (lisse) motivic structures (resp. Galois representations) and \[\Xi \in \{G,U,G^{\gr},U^{\gr},L^W=\Gb,U^W_H = F^0 U\},\] we set
\[
\Xi(\Tc) \coloneqq \Xi(\om^{\dR}),
\]
which is a pro-algebraic group over $\basefield$ (resp. $\basefield_{\pf}$). For \[\Tc = \mathbf{SM}_\basefield(\Qb),\mathbf{SM}_{\basefield,S}(\Qb),\mathbf{SM}_{\basefield}(\Qb,A),\mathbf{SM}_{\basefield,S}(\Qb,A),\] we refer to $\Xi(\Tc)$, respectively, by \[\Xi(\basefield), \Xi(\Oc_{\basefield,S}), \Xi(A), \Xi(\Oc_{\basefield,S},A).\]

In the latter cases, note that $\Gb(A) = \Gb(A,\Oc_{\basefield,S})$ is independent of $S$.

%\begin{sExample}
%\label{X55}
%More generally, if $\Ac \to Z$ is an abelian scheme with connected $p$-adic monodromy group, then $\Tc = \Rep_{\Qp}^{\mathrm{sg}}(G_{\basefield},\Ac_\basefield)$ is admissible. This is expected to happen, for example, if $\Ac_\basefield$ is an abelian variety of dimension $g$ with Mumford--Tate group $\operatorname{GSp}_{2g}$, by the Mumford--Tate Conjecture. In this latter case, a similar argument as in \ref{example:E} shows $\Tc = \mathbf{SM}_{Z}(\Qb,\Ac_{\basefield})$ is admissible as well.\end{sExample}

\subsection{Period loops for motivic structures and Galois representations}\label{sec:period_loops_mot_gal} Let $\Tc$ be a category of lisse motivic structures (resp. lisse $p$-adic Galois representations) over $\Oc_{\basefield,S}$, $\pf \in S_p(\basefield) \setminus S$, and $U=U(\Tc)$. Using the Hodge-filtered Frobenius equivariant fiber functor of \ref{sec:hodge_frob_spm} (resp. \ref{hodge_frob_rep}) and a choice of arithmetic Hodge path, Definition \ref{defn:period_loops} gives us the right-hand and left-hand \emph{unipotent $\pf$-adic period loop of $\Tc$} 
\[
\uU^R_{\pf} = \uU^R_{\pf}(\Tc) \in U(\basefield_{\pf})
\]
\[
\uU^L_{\pf} = \uU^L_{\pf}(\Tc) \in U(\basefield_{\pf})
\]

Dually, we get \emph{$\pf$-adic period homomorphisms}
\[
\per_{\pf}^R = \per_{\pf}^R(\Tc) \colon \Oc(U)  \to \basefield_{\pf}
\]
\[
\per_{\pf}^L = \per_{\pf}^L(\Tc) \colon \Oc(U)  \to \basefield_{\pf},
\]
whose restrictions to $\Oc(F^0 \backslash U)$ (resp. $\Oc(U/F^0)$) are independent of the choice of arithmetic Hodge path.

%% ***** I think this is fine but just check the notation about the field agrees with defn:period_loops. like make sure k is somehow implicit there. Maybe want to introduce notation like $\om_{k'}$ for a fiber functor \om over k. Then use terminology like $\om^{\dR}$ instead of $\dR$. Also make consistent between $\om^\m{H}$ and $\om^{\dR}$ (former is maybe okay in Sec 4-5, but not in Sec 3 and 6).

\subsubsection{\yesrem}\label{rem:compat_period_loops}
By Remark \ref{rem:compat_hodge_frob}, the versions for motivic structures and Galois representations, respectively, correspond under $\mathrm{real}_{\et,\ell}$.

We now formulate a $p$-adic version of the Kontsevich--Zagier Period Conjecture, generalizing that of \cite{YamashitaBounds10}:

\begin{ssConjecture}\label{conj:period_conjecture}
Let $\Tc$ be a category of lisse motivic structures over $\Oc_{\basefield,S}$ and $\pf \notin S$. Then $\uU^R_{\pf}$ (resp. $\uU^L_{\pf}$) is dense in $F^0 \backslash U$ (resp. $U/F^0$). Equivalently, $\per_{\pf}^R$ (resp. $\per_{\pf}^L$) is injective on $\Oc(F^0 \backslash U)$ (resp. $\Oc(U/F^0)$).
\end{ssConjecture}

\subsubsection{Remark}\label{rem:}

This should be compared with the closely-related conjectures of Ancona--Fr\u{a}\c{t}il\u{a} \cite[Conjectures 8.1,8.3]{AnFr2025algebraic}. Theirs would be more likely candidates for the ``correct'' $p$-adic period conjecture in the pure case. However, if one wants an actual conjecture and not a meta-conjecture whose formulation depends on other conjectures (c.f. \cite[\S 9]{AnFr2025algebraic}), then the conjecture of Ancona-Fr\u{a}\c{t}il\u{a} is largely limited to the setting of pure motives. Indeed, in order to formulate a similar conjecture for our categories of mixed motivic structures, these categories would need good special fibers over finite fields. For a fuller discussion of these issues and for a comparison in the mixed Tate setting, see \cite{dancohen2025andreperiodsmixedtate}.

\subsection{$p$-adic periods and $p$-adic regulators}

\subsubsection{Introduction}

In this subsection, we establish a relationship between our $p$-adic period point $\mathfrak{u}_\pf$ and syntomic regulators. In the case of mixed Tate motives, i.e. $\Tc=\mathbf{MT}(\Oc_{\basefield,S},\Qb)$, we have an embedding 
\[
\bigoplus_{n} K_{2n-1}(\Oc_{\basefield,S})^{(n)} \otimes \Qb(-n)_{\dR} = \bigoplus_{n} \Ext^1_\Tc(\Qb(0),\Qb(n)) \otimes \Qb(-n)_{\dR} \subseteq \Oc(U(\om^{\dR})).
\]

Restricting the period map \[\Oc(U) \xrightarrow{\per_p} \Qp\] to \[K_{2n-1}(\Oc_{\basefield,S})^{(n)} \otimes \Qb(-n)_{\dR}\] for a fixed $n$, we equivalently get a map
\[
K_{2n-1}(\Oc_{\basefield,S})^{(n)} \to \Qp(n)_{\dR}.
\]

For $n \le 0$, the left-hand side is zero, so this is not interesting. But for $n \ge 1$, the right-hand side is precisely the syntomic cohomology
\[
H^1_{\syn}(\Oc_{\pf},n),
\]
and the map is identified with the syntomic regulator.

\subsubsection{Bloch--Kato Exponential Map}\label{Bloch--Kato_exp}

With notation as in \S \ref{sec:period_loops_mot_gal},\footnote{Our notation reflects the case of motivic structures, but Proposition \ref{prop:periods_BK_log} holds in the case of $p$-adic Galois representations as well.} let $M$ be an object of $\Tc$ with $\gr^W_0 M = 0$. Then $1-\phi$ induces an isomorphism on $M_{\cri}$, so we have by \cite[Corollary 3.8.4]{BlochKato} an isomorphism \[\operatorname{D}_{\dR}(M_p)/F^0 \xto{\exp_{\mathrm{BK}}} H^1_f(G_{\pf};M_p),\] whose inverse is denoted $\log_{\mathrm{BK}}$. Note that $\iota_{\pf}$ (\ref{PSMP}) induces an isomorphism $M_{\dR,\pf} \coloneqq M_{\dR} \otimes_{\basefield} \basefield_{\pf} \xto{\sim} \operatorname{D}_{\dR}(M_p)$. We consider the composition
\begin{align*}\label{eqn:BK_composition}
\Ext^1_{\Tc}(\QQ(0),M) = \Ext^1_{\mathbf{SM}_{\basefield,S}(\Qb)}(\QQ(0),M) \xto{\mathrm{real}_{\et}}\\
H^1_{f,S}(G_{\basefield};M_p) \xto{\m{loc}_{\pf}}
H^1_f(G_{\pf};M_p)
\xrightarrow[\sim]{\log_{\mathrm{BK}}} M_{\dR,\pf}/F^0.
\end{align*}

Dually, this gives an element
\begin{eqnarray*}
LM_{\mathrm{BK}} &\in& \homo(\Ext^1_{\Tc}(\QQ(0),M),\Qb) \otimes_{\Qb} M_{\dR,\pf}/F^0\\
&\simeq& \homo(\Ext^1_{\Tc}(\QQ(0),M) \otimes_{\Qb} (M_{\dR}/F^0)^{\vee},K_{\pf}).
\end{eqnarray*}

We would like to relate $LM_{\mathrm{BK}}$ to $\per_{\pf}^L \colon \Oc(U/F^0) \to K_{\pf}$.

\subsubsection{Extensions and Graded Galois Groups}\label{extensions_and_graded}

Let $\om=\om^{\dR}$. For $M$ as in \ref{Bloch--Kato_exp}, we first define a canonical map
\[
a_M^{\gr} \colon \Ext^1_{\Tc}(\QQ(0),M) \otimes_{\Qb} \omgr(M)^{\vee} \to \Oc(U^{\gr}).
\]

Given $c \otimes f \in \Ext^1_{\Tc}(\QQ(0),M) \otimes_{\Qb} \omgr(M)^{\vee}$, we choose an extension
\[
0 \to M \to E \to \QQ(0) \to 0,
\]
representing $c$ and then apply $\omgr$ to get an extension
\[
0 \to \omgr(M) \to \omgr(E) \to \Qb \to 0
\]
of $G^{\gr}(\Tc)$-representations. Since $\Gb$ is reductive, the sequence
\[
0 \to \omgr(M)^{\Gb} \to \omgr(E)^{\Gb} \to \Qb^{\Gb} \to 0
\]
is again exact. But $\omgr(M)^{\Gb}=0$ since $\gr^W_0 M = 0$, so $\omgr(E)^{\Gb} = \Qb^{\Gb} = \Qb$. Therefore, there is a unique $1_E \in \omgr(E)^{\Gb}$ mapping to $1 \in \Qb$, which is the $\gamma^{\gr}$ of Proposition \ref{prop:effective_gr_point} for $M$ effective. Since $U^{\gr}(\Tc)$ acts trivially on $\Qb$, we have for every $u \in U^{\gr}(\Tc)$ that \[-1_E+u(1_E) \in \omgr(M).\] We define $a_M^{\gr}(c \otimes f)$ to be the function \[u \in U^{\gr}(\Tc) \mapsto f(-1_E+u(1_E)).\]

For $M$ effective, there is another description of $a_M^{\gr}$ via the map $\delta$ of \S \ref{sec:abstract_CK_maps}. We note that
\[\Ext^1_{\Tc}(\QQ(0),M)_{\basefield} \overset{\footnotemark}{\simeq} \footnotetext{Given an extension $0 \to M \to E \xto{\rho} \QQ(0) \to 0$, the associated torsor is $\rho^{-1}(1)$.}\, H^1(\Tc;M)_{\basefield} = H^1(G^{\gr};\omgr(M)) \xrightarrow[\delta]{\sim} Z^1(U^{\gr}; \omgr(M))^{\GG}.\]

Given \[c \otimes f \in \Ext^1_{\Tc}(\QQ(0),M) \otimes_{\Qb} \omgr(M)^{\vee} \simeq Z^1(U^{\gr}; \omgr(M))^{\GG} \otimes_{\basefield} \omgr(M)^{\vee},\] the functional $f \colon \omgr(M) \to \basefield$ is an element of $\Oc(\omgr(M))$. We have $a_M^{\gr}(c \otimes f) = c^{\#}(f)$.

\subsubsection{Extensions and Ungraded Galois Groups}\label{extensions_and_ungraded}

Let us now fix a choice of arithmetic Hodge path $\tannloop^\m{H} \colon \omgr \to \om$. This gives us an isomorphism $\tau^\m{H} \colon G^{\gr} \to G$ restricting to an isomorphism $U^{\gr} \to U$ and $1_H \coloneqq \tannloop^\m{H}(1_E) \in \om(E)$. Note that $1_H \in F^0 \om(E)$.

We define
\[
a_M \colon \Ext^1_{\Tc}(\QQ(0),M) \otimes_{\Qb} \om(M)^{\vee} \to \Oc(U)
\]
depending on $\tannloop^\m{H}$ (but independent if $M$ is semisimple), as follows.

The simplest way is to transport everything via the isomorphisms $\omgr(V) \to \om(V)$ and $U^{\gr} \to U$ induced by $\tannloop^\m{H}$. But we may give a description parallel to the one above by letting $a_M(c \otimes f)$ be the function sending $u \in U$ to $f(-1_H+u(1_H))$.

\begin{ssProposition}\label{prop:periods_BK_log}
For any choice of $\tannloop^\m{H}$, the composition
\[
\Ext^1_{\Tc}(\QQ(0),M) \otimes_{\Qb} (\om(M)/F^0)^{\vee}
\hookrightarrow
\Ext^1_{\Tc}(\QQ(0),M) \otimes_{\Qb} \om(M)^{\vee} \xto{a_M}
\Oc(U) \xto{\per_{\pf}^L} K_{\pf}
\]
is
\[-LM_{\mathrm{BK}} \in \homo(\Ext^1_{\Tc}(\QQ(0),M) \otimes_{\Qb} (M_{\dR}/F^0)^{\vee},K_{\pf}).\]
\end{ssProposition}

\begin{proof}
By definition of $a_M$, we have $\per_{\pf}^L(a_M(c \otimes f)) = f(-1_H+u_{\pf}^L(1_H))$. But $u_{\pf}^L(1_H) = \tannloop^{\cri}({\tannloop^\m{H}}^{-1}(1_H)) = \tannloop^{\cri}(1_E)$ is by Lemma \ref{h0.5} the unique $\phi$-invariant element of $\om(E)$ mapping to $1 \in \Qb$, which we denote by $1_{\cri}$. Thus \[\per_{\pf}^L(a_M(c \otimes f)) = f(-1_H+1_{\cri}).\] 

Now $\log_{\m{BK}}(c) = 1_H - 1_{\cri}$, so indeed
\[
\per_{\pf}^L(a_M(c \otimes f)) = -\log_{\m{BK}}(c \otimes f),
\]
and we are done.
\end{proof}

\subsubsection{Syntomic Cohomology and Regulators}\label{sssec:syntomic_coh_reg} We recall the theory of syntomic regulators. A common theme in results about Iwasawa theory and $p$-adic L-functions is the computation of syntomic regulators of specific elements in algebraic $K$-theory or motivic cohomology. See for example \cite{BerDarRot15,DarmonRotger17,LoefSkiZer20,loeffler2021padic}.

While ultimately the connection with $\per_{\pf}$ is just a reformulation of Proposition \ref{prop:periods_BK_log}, we thought it fruitful to elucidate the connection given the importance of syntomic regulators. In Example \ref{example:syn_reg_K_2}, we describe an application of this to computing $\per_{\pf}$ for use in the Chabauty--Kim method.

Let $\Xc \to V \coloneqq \Spec{\Oc_\pf}$ smooth and projective, $X = \Xc_{\basefield_{\pf}}$, and let $\kappa$ be the residue field. For $k,n \ge 0$, we have a quasi-isomorphism (\cite[Proposition 9.9]{BesserRigSyn00I},\cite[A.2]{ErtlNiziol19}):
\begin{equation}\label{eqn:integral_syntomic}
\RGamma^{\rig}_{\syn}(\Xc,n) \simeq \RGamma_f(\Xc,\mathcal{K}(n)) \coloneqq [\RGamma_{\rig}(\Xc_{\kappa}/\kappa)^{\ph=p^n} \to \RGamma_{\dR}(X)/F^n]
\end{equation}
where the first is defined in \cite[\S 6]{BesserRigSyn00I} and the second in \cite[\S 2.1]{NiziolImage97} and \cite[\S 3]{NiziolSmooth01}. We denote the cohomology groups by $H^k_{\syn}(\Xc,n)$.

%% \coloneqq [\RGamma_{\rig}(\Xc_{\kappa}/\kappa) \oplus F^n \RGamma_{\dR}(X) \xrightarrow{(x,y) \mapsto ((p^r-\phi)(x),\m{sp}(y)-x)} \RGamma_{\rig}(\Xc_{\kappa}/\kappa)^{\oplus 2}]\\

There are syntomic Fontaine--Messing ``period morphisms''
\[
\alpha_{k,n} \colon H^{k}_{\syn}(\Xc,n) \to H^k_{\et}(X;\Qp(n)),
\]
known as `$l$' in \cite[\S 4-5]{NiziolSmooth01}, and syntomic Chern classes (\cite[\S 7]{BesserRigSyn00I} or \cite[Theorem 2.1]{NiziolImage97})
\[
c_{n,m}^{\syn} \colon K_m(\Xc) \to H^{2n-m}_{\syn}(\Xc,n)
\]
that agree with the usual \'etale Chern classes (\cite[Proposition 3.4]{NiziolImage97}). Via the Chern character (see e.g. \cite[\S 7.1]{SouleOp85}), one gets maps
\begin{equation}\label{eqn:ch_syn}
K_{2n-k}(\Xc)^{(n)}_{\Qb} \xrightarrow{\reg_{\syn}} H^{k}_{\syn}(\Xc,n).
\end{equation}
compatible with cup products, which we call \emph{syntomic regulators}. The composition (first isomorphism because $\Xc$ is regular)
\[
%\cyc_{\syn} \colon 
\CH^{n}(\Xc,2n-k)_{\Qb} \simeq K_{2n-k}(\Xc)^{(n)}_{\Qb} \xrightarrow{\reg_{\syn}} H^{k}_{\syn}(\Xc,n) \xrightarrow{\rho_{\syn}} H^k_{\et}(X;\Qp(n))
\]
is thus the same as
\[
\CH^n(\Xc,2n-k)_{\Qb} \to \CH^n(X,2n-k)_{\Qb} \xrightarrow{\cyc_{\et}} H^k_{\et}(X;\Qp(n)).
\]

Niziol (\cite[\S 5]{NiziolImage97}) shows that (\ref{eqn:integral_syntomic}) is computed by a spectral sequence with $E^{i,j}_2 = H^i_{f}(G_{\pf};H^j_{\et}(\overline{X};\Qp(n)))$, and this spectral sequence maps to the Hochschild--Serre spectral sequence $E^{i,j}_2 = H^i(G_{\pf};H^j_{\et}(\overline{X};\Qp(n)))$ computing $H^k_{\et}(X;\Qp(n))$ via the inclusion $H^i_{f}(G_{\pf};-) \hookrightarrow H^i(G_{{\pf}};-)$ on $E_2$-terms and $\alpha_{k,n}$ on the limit.

The description (\ref{eqn:integral_syntomic}) induces a map \[H^{k-1}_{\dR}(X)/F^n \xrightarrow{\del} H^{k}_{\syn}(\Xc,n).\] By \cite[Proposition 9.11]{BesserRigSyn00I}, the composition \[H^{k-1}_{\dR}(X)/F^n \xrightarrow{\del} H^{k}_{\syn}(\Xc,n) \xrightarrow{\alpha_{k,n}} H^k_{\et}(X;\Qp(n)) \to H^k_{\et}(\overline{X};\Qp(n))\] is $0$. Letting $H^{k}_{\syn}(\Xc,n)_0 \coloneqq \Ker(H^{k}_{\syn}(\Xc,n) \to H^k_{\et}(\overline{X};\Qp(n)))$, the induced map
\begin{equation}\label{eqn:beta_delta}
H^{k-1}_{\dR}(X)/F^n \xrightarrow{\del} H^{k}_{\syn}(\Xc,n)_0 \to H^1_{f}(G_{\pf};H^{k-1}_{\et}(\overline{X};\Qp(n)))
\end{equation}
is $\exp_{\m{BK}}$.\footnote{Note that \cite{BesserRigSyn00I} cites \cite{NiziolSmooth01} in the proof of Proposition 9.11. The map $\beta$ in \cite[Proposition 2.2]{NiziolSmooth01} is the same as that of \cite[(1.17.2)]{BlochKato}, so the resulting exponential is the same and not, e.g., its negative.}

%Ask Disegni for pointers. This motivation is important for selling this paper. Could put this into intro if desired. ***** Ask Darmon (or Rotger) about applications of these ideas?

%EXPLAIN IDEA: Allows you to take p-adic regulator of something motivic that it not in K-theory because it doesn't satisfy a cocycle condition. A little like defining p-adic regulator on Bloch complex, except this works in general. (QUESTION: Can we use it to do what dJ did for usual Bloch complex but instead for elliptic?)

\subsubsection{Syntomic regulators and the $p$-adic period map}

Now suppose that $X$ is a smooth projective variety over $\basefield$ with smooth projective model $\Xc$ over the localization $\Oc_{\basefield,\pf}$ whose base-change to $V$ we denote $\Xc_{\pf}$.

Recall from \S \ref{sec:sm_mot_coh} that we have a map \[\AJ_{\mot} \colon \CH^{n}_{\mathrm{hom}}(X,2n-k)_{\Qb} \to \Ext^1_{{\mathbf{SM}_\basefield(\Qb))}}(\Qb,h_{\mathbf{SM}}^{k-1}(X)(n)).\] %We set $K_{2n-k}^{\mathrm{hom}}(X)^{(n)}_{\Qb} \coloneqq \che^{-1}(\CH^{n}_{\mathrm{hom}}(X,2n-k)_{\Qb})$. Note that this is the kernel of the map $K_{2n-k}(X)^{(n)}_{\Qb} \xrightarrow{\che_{\et}} H^k_{\et}(X;\Qp(n)) \to H^k_{\et}(\overline{X};\Qp(n))$. We let $K_{2n-k}^{\mathrm{hom}}(\Xc)^{(n)}_{\Qb}$ denote its preimage under the natural map $K_{2n-k}(\Xc)^{(n)}_{\Qb} \to K_{2n-k}(X)^{(n)}_{\Qb}$.
By \cite[Theorems 3.1, 3.2]{NiziolImage97}, the image of
\[
\CH^{n}_{\mathrm{hom}}(\Xc,2n-k)_{\Qb} \to \CH^{n}_{\mathrm{hom}}(X,2n-k)_{\Qb} \xrightarrow{\AJ_{\et}} H^1(G_{\basefield};H_{\et}^{k-1}(\overline{X};\Qp(n)))
\]
lands in $H^1_f(G_{\basefield};H_{\et}^{k-1}(\overline{X};\Qp(n)))$, which implies that the image of
\[
\CH^{n}_{\mathrm{hom}}(\Xc,2n-k)_{\Qb} \to \CH^{n}_{\mathrm{hom}}(X,2n-k)_{\Qb} \xrightarrow{\AJ_{\mot}} \Ext^1_{{\mathbf{SM}_\basefield(\Qb)}}(\Qb,h_{\mathbf{SM}}^{k-1}(X)(n))
\]
lands in $\Ext^1_{\mathbf{SM}_{\Oc_{\basefield,\pf}}(\Qb)}(\Qb,h_{\mathbf{SM}}^{k-1}(X)(n))$.

We may define
\begin{align*}
\CH^{n}_{\mathrm{hom}}(\Xc,2n-k)_{\Qb} \otimes H_{k-1}^{\dR}(X)(n)
\xrightarrow{\AJ_{\mot} \otimes \id}
\Ext^1_{{\mathbf{SM}_{\Oc_{\basefield,\pf}}(\Qb)}}(\Qb,h_{\mathbf{SM}}^{k-1}(X)(n)) \otimes H_{k-1}^{\dR}(X)(n)\\
=
\Ext^1_{{\mathbf{SM}_{\Oc_{\basefield,\pf}}(\Qb)}}(\Qb,h_{\mathbf{SM}}^{k-1}(X)(n)) \otimes \om^{\dR}(h_{\mathbf{SM}}^{k-1}(X)(n))^{\vee}
\xhookrightarrow{a_{h_{\mathbf{SM}}^{k-1}(X)(n)}}
\Oc(U)
\xrightarrow{\per_{\pf}^L}
K_{\pf}
\end{align*}

Thus $\per_{\pf}^L$ defines an associated map
\begin{equation}\label{eqn:epsilon^L}
\epsilon^L \colon CH^{n}_{\mathrm{hom}}(\Xc,2n-k)_{\Qb} \to H^{k-1}_{\dR}(X) \otimes_{\basefield} K_{\pf} = H^{k-1}_{\dR}(X_{\pf}),
\end{equation} and we take the composition
\begin{equation}\label{eqn:integral_syn_reg1}
\begin{split}
	\CH^{n}_{\mathrm{hom}}(\Xc,2n-k)_{\Qb} \xto{\epsilon^L} H^{k-1}_{\dR}(X_{\pf})(n)
	\twoheadrightarrow\\
	H^{k-1}_{\dR}(X_{\pf})(n)/F^0
	\simeq
	H^{k-1}_{\dR}(X_{\pf})/F^n
	\xrightarrow{\delta}
	H^{k}_{\syn}(\Xc_{\pf},n).
\end{split}
\end{equation}

\begin{ssProposition}\label{prop:periods_syntomic_reg}
The map of (\ref{eqn:integral_syn_reg1}) defined using our $\per_{\pf}^L$ is the negative of the syntomic cycle class map
\begin{equation}\label{eqn:integral_syn_reg2}
	\begin{split}
		\CH^{n}_{\mathrm{hom}}(\Xc,2n-k)_{\Qb} \hookrightarrow
		\CH^{n}(\Xc,2n-k)_{\Qb}\\
		\simeq
		K_{2n-k}(\Xc)^{(n)}_{\Qb}
		\to
		K_{2n-k}(\Xc_{\pf})^{(n)}_{\Qb}
		\xrightarrow{\reg_{\syn}}
		H^{k}_{\syn}(\Xc_{\pf},n).
	\end{split}
\end{equation}
\end{ssProposition}

\begin{proof}

Note that the images of both lie in $H^k_{\syn}(\Xc_{\pf},n)_0$. We have an injection
\[
\beta \colon H^k_{\syn}(\Xc_{\pf},n)_0 \simeq H^1_{f}(G_{\pf};H^{k-1}_{\et}(\overline{X};\Qp(n))) \hookrightarrow H^1(G_{\pf};H^{k-1}_{\et}(\overline{X};\Qp(n))),
\]
where the first isomorphism is from the spectral sequence for $H^k_{\syn}(\Xc_{\pf},n)$, noting that $H^i_f$ vanishines for $i > 1$. We compose both (\ref{eqn:integral_syn_reg1}) and (\ref{eqn:integral_syn_reg2}) with $\beta$ and prove the results are negatives of each other.

By compatibility between syntomic and \'etale Chern classes, (\ref{eqn:integral_syn_reg2}) just becomes $\AJ_{\et}$.

By Proposition \ref{prop:periods_BK_log}, the map \[\CH^{n}_{\mathrm{hom}}(\Xc,2n-k)_{\Qb} \xto{\epsilon^L} H^{k-1}_{\dR}(X_{\pf})(n) \to
H^{k-1}_{\dR}(X_{\pf})(n)/F^0\] in (\ref{eqn:integral_syn_reg1}) is \[-\log_{\m{BK}} \circ \mathrm{real}_{\et} \circ \AJ_{\mot} = -\log_{\m{BK}} \circ \AJ_{\et}.\] But $\beta \circ \delta$ is $\exp_{\m{BK}}$ as in (\ref{eqn:beta_delta}), so in total we get $\exp_{\m{BK}} \circ (-\log_{\m{BK}} \circ \AJ_{\et})$, which is $-\AJ_{\et}$ as desired.

\end{proof}

\subsubsection{\yesrem}\label{remark:syntomic_indep_of_hodge_path} Proposition \ref{prop:periods_syntomic_reg} implies that the composition (\ref{eqn:integral_syn_reg1}) is independent of choice of $\tannloop^\m{H}$. Ultimately, this is a result of taking the projection $H^{k-1}_{\dR}(X_{\pf})(n)
	\twoheadrightarrow H^{k-1}_{\dR}(X_{\pf})(n)/F^0$.

\begin{ssExample}\label{example:syn_reg_K_2}
Let $E$ be an elliptic curve over $K$ with good model $\Ec$ over $Z \subseteq \Spec{\Oc_K}$ containing $\pf$. Let $\sigma \in K_2(\Ec)$, defined explicitly so that its image in $K_2(K(E)) \simeq K_2^M(K(E))$ is $f \otimes g$, for $f,g \in K(E)^{\times}$. Let $\om_0,\om_1$ be a basis of differential forms of the second kind on $E$. Let $M = h_1^{\mathbf{SM}}(E)(1) = h^1_{\mathbf{SM}}(E)(2)$, and let $k=n=m=2$. We have $\om_0 \otimes (\om_0 \wedge \om_1),\om_1 \otimes (\om_0 \wedge \om_1) \in M_{\dR}^{\vee} = H^1_{\dR}(E)(-1)$.

Via $\AJ_{\mot}$, we get a class $\sigma_{\mathbf{SM}} \in \Ext^1_{\mathbf{SM}_{Z}(\Qb))}(\Qb,h_{\mathbf{SM}}^{1}(E)(2))$. Fixing a choice of arithmetic Hodge path, we send $\sigma_{\mathbf{SM}} \otimes \om_0 \otimes (\om_0 \wedge \om_1), \sigma_{\mathbf{SM}} \otimes \om_1 \otimes (\om_0 \wedge \om_1)$ via \ref{extensions_and_ungraded} to elements $\sigma_0,\sigma_1 \in \Oc(U(Z))$.

By Proposition \ref{prop:periods_syntomic_reg}, the numbers $\per_{\pf}^L(\sigma_0),\per_{\pf}^L(\sigma_1)$ may be computed in terms of the syntomic regulator $K_2(\Ec) \to H^2_{\syn}(\Ec,2)_0 \simeq M_{\dR}/F^0 = M_{\dR}$.

The ratio $\frac{\per_{\pf}^L(\sigma_0)}{\per_{\pf}^L(\sigma_1)}$ (which is independent of $\sigma$) is a ratio of coefficients of a Chabauty--Kim function, specifically $-\frac{2c_2}{c_1}$ in the notation of \cite[Theorem 1.4]{CorwinMECK}. The first author and M. L\"udtke have already computed such a function for a particular elliptic curve by testing the general form of the theorem of loc.cit. on known $\Zb[\frac{1}{2}]$-points, as described in \cite{CLinitialreport2026}; this thus produces a computation of the syntomic regulator \emph{via the Chabauty--Kim method}.

This syntomic regulator may alternatively be computed using $p$-adic integration via the method of \cite{BesserRigSyn00II,BesserArbSyn21}, and comparing them is the subject of work-in-progress by the first author and M. L\"udtke.

\end{ssExample}

%%%%%%%%%%%%%%%%%%
\section{Applications to Chabauty--Kim}\label{sec:ck_application}
%%%%%%%%%%%%%%%%%

Our goal in this section is to explain how our main result (Theorem \ref{thm:left_commutes}) fits into the Chabauty--Kim method, specifically the computation of the localization map from the global Selmer variety to the local Selmer variety. We focus on categories of motivic structures rather than Galois representations; the latter is dealt with in the parallel work \cite{CorwinTSV}.

\subsection{Fundamental group and path torsors}\label{sec:fund_grp_and_path_torsors}

We explain how to use the category $\mathbf{SM}_{\basefield,S}(\Qb)$ to study the Chabauty--Kim method via motivic structures. %\footnote{The case of Galois representations is discussed in \cite{CorwinTSV}.}
In particular, we construct a motivic version $\pi_1^{\mathbf{SM}}(X,b)$ of the unipotent fundamental group of $X$ in the category $\mathbf{SM}_{\basefield,S}(\Qb)$ and define a unipotent Kummer map
\[
X(\Oc_{\basefield,S}) \to H^1(\mathbf{SM}_{\basefield,S}(\Qb);\pi_1^{\mathbf{SM}}(X,b)).
\]
The key is to define path torsors between any two $a,b \in X(\Oc_{\basefield,S})$. The Betti, de Rham, and $\ell$-adic realizations of these path torsors appear already in \cite[\S 10]{Deligne89}, which, modulo the lack of a $p$-adic comparison in loc.cit., allows one to define a version $\pi_1^{\mathbf{SPM}}(X,b)$ of the unipotent fundamental group of $X$ in the category $\mathbf{SPM}_{\basefield,S}(\Qb)$. The key step is thus to show that this lifts to the level of Voevodsky's geometric motives.

We first describe the explicit approach of \cite[\S 4.2]{GonMPMTM01} (c.f. also \cite[\S 3.12]{DelGon05}), originally due to Beilinson, to constructing an inductive object of $\DM_{\gm}(\basefield,\Qb)$ that realizes $\pi_1^{\mathbf{SPM}}(X,b)$. We then make remarks about how more modern technology allows one to get this construction in a cleaner and more abstract way.

\subsubsection{}\label{sec:path_torsor_SM}

Let $X$ be an arbitrary smooth connected variety over a field $\basefield$ of characteristic $0$. Goncharov (\cite[\S 4.2]{GonMPMTM01}) defines explicit complexes of varieties $P^n(X;a,b)$ associated to $X,a,b$, and we let $\Pc^n_{\mathcal{M}}(X;a,b) \in \DM_{\gm}(\basefield,\Qb)$ denote the object associated to $P^n(X;a,b)[n]$, so that \[H^0_B(\Pc^n_{\mathcal{M}}(X;a,b)) = H^n_B(P^n_{\mathcal{M}}(X;a,b)) =  \Pc^n(X(\Cb);a,b) \coloneqq \Pc(X(\Cb);a,b)/\Ic_{a,b}^{n+1}\] in the notation of \cite[\S 4.1]{GonMPMTM01}, and so that we get a map
\[
p_n \colon \Pc^n_{\mathcal{M}}(X;a,b) \to \Pc^{n-1}_{\mathcal{M}}(X;a,b)
\]
for each $n$ coming from the boundary map $P^n(X;a,b)_{\bullet} \to P^{n-1}(X;a,b)_{\bullet}[-1]$ on p.45 of lot.cit.

For $\basefield$ a number field, we note that
\[
\Pc^n_{\mathbf{SM}}(X;a,b) \coloneqq H^0(\mathrm{real}_{\mathbf{SM}}(\Pc^n_{\mathcal{M}}(X;a,b))) \in \mathbf{SM}_\basefield(\Qb)
\]
by definition of the latter. The inductive limit
\[
\varinjlim_n \Pc^n_{\mathbf{SM}}(X;a,b)^{\vee}
\]
along the duals of $p_n$ is a Hopf algebra in $\operatorname{Ind}{\mathbf{SM}_\basefield(\Qb)}$. We denote its spectrum by $\pi_1^{\mathbf{SM}}(X;a,b)$, which is an affine scheme in $\mathbf{SM}_\basefield(\Qb)$. By Theorems 4.1, 4.3, and 4.4 of \cite{GonMPMTM01}, its Betti, $\ell$-adic, and de Rham components give the usual Betti, $\ell$-adic \'etale, and de Rham unipotent path torsors of $X$. For $a=b$, we denote it by $\pi_1^{\mathbf{SM}}(X,a)$.

\subsubsection{Remark}

There are multiple recent constructions of a path-space torsor between any two points in a model-theoretic or infinity-categorical setup: see ${_\eta B_\xi}$ of \cite[2.2]{DCSRationalMotivicPathSpaces22}, $P_X(x,y)$ of \cite[3.12]{IwanariMotivicRational20}, and a simple modification (c.f. Remark \ref{rem:less_ad_hoc_tangential}(\ref{item:ayoub})) $\pi_1^{\un}(X,x_0,x_1)_{\Qb}$ of $\pi_1^{\mathrm{geo}}(X,x_0,x_1)_{\Qb}$ of \cite[Construction 4.5.1]{AyoubAnabelianPresentationPreprint}. In all cases, these give objects in the associated homotopy category, which in all cases is equivalent to $DM(\basefield,\Qb)$ (also denoted $DA(\basefield,\Qb)$), the triangulated category of motivic complexes over $\basefield$ with rational coefficients, of which $\DM_{\gm}(\basefield,\Qb)$ is the full subcategory of compact objects.

\subsection{Tangential basepoints}\label{sec:tangential_basepoints}

%*** Ideal is to use a more geometric construction for tangential basepoints, maybe using log geometry (sent an email to Martin about this).

We may also include the case of tangential basepoints if $X$ is a curve with a $\basefield$-point, as follows. Letting $X'$ denote the smooth compactification of $X$ and $Y=X' \setminus X$, a tangential basepoint is a pair $b=(\overline{b},v)$ with $\overline{b} \in Y(\basefield)$ and $v \in T_{\overline{b}}(X') \setminus \{0\}$. Work of Deligne (\cite[\S 15]{Deligne89}) produces a pro-algebraic group $\pi_1^{\mathbf{SPM}}(X;a,b)$ in $\mathbf{SPM}_\basefield(\Qb)$ when $a$ or $b$ (or both) is a tangential basepoint. We prove:

\begin{ssTheorem}\label{thm:motivic_torsor_tangential}
If $X(\basefield) \neq \emptyset$, $\pi_1^{\mathbf{SPM}}(X;a,b)$ is in $\mathbf{SM}_\basefield(\Qb)$ (and we denote it by $\pi_1^{\mathbf{SM}}(X;a,b)$).
\end{ssTheorem}

\subsubsection{Proof of Theorem \ref{thm:motivic_torsor_tangential}}
We modify the proof of \cite[Th\'eor\`eme 4.4]{DelGon05}. First, note that as in the proof of \cite[4.11]{DelGon05}, it suffices to treat the case where only one of $a$ or $b$ is tangential (this uses the fact that $X(\basefield)$ is nonempty). Let us assume $b$ is tangential. The same argument also shows it suffices to prove it for a single $a \in X(\basefield)$ for each $b$. Let us suppose that only $b=(\overline{b},v)$ is tangential.

Let $g$ be the genus of $X'$ and $c$ an auxiliary closed point of $X'$ not equal to $a$ or $\overline{b}$. We choose $m \in \Zb$ such that $n \coloneqq \deg(m[c]) > 2g$. Then the divisors $-[\overline{b}]+m[c]$, $-[a]-[\overline{b}]+m[c]$, and $-2[\overline{b}]+m[c]$ are all of degree greater than $2g-2$, so by Riemann--Roch, we have
\begin{align*}
\dim{H^0(X';-[\overline{b}]+m[c])} = \deg(-[\overline{b}]+m[c]) = n-1\\
\dim{H^0(X';-[a]-[\overline{b}]+m[c])} = \deg(-[a]-[\overline{b}]+m[c]) = n-2\\
\dim{H^0(X';-2[\overline{b}]+m[c])} = \deg(-2[\overline{b}]+m[c]) = n-2.
\end{align*}

Since $\basefield$ is of characteristic $0$ and therefore infinite, we may find $u \in H^0(X';-[\overline{b}]+m[c]) \setminus (H^0(X';-[a]-[\overline{b}]+m[c]) \cup H^0(X';-2[\overline{b}]+m[c])) \subseteq \basefield(X)=\basefield(X')$. In other words, $u$ is a uniformizer at $\overline{b}$ and is defined and non-vanishing at $a$.

Viewing $u$ as a map from $X'$ to $\Pb^1$, let $U$ denote preimage of $\Gm \subseteq \Pb^1$ in $X$. Note that $a \in U$.% Let $U' = U \cup \{\overline{b}\}$.

Then $a \in U(\basefield)$ and $b$ is a tangential basepoint of $U$, so we may consider $\pi_1^{\mathbf{SPM}}(U;a,b)$. The map $\pi_1^{\mathbf{SPM}}(U;a,b) \to \pi_1^{\mathbf{SPM}}(X;a,b)$ induced by the inclusion $U \hookrightarrow X$ is surjective, so it suffices to prove that $\pi_1^{\mathbf{SPM}}(U;a,b)$ is in $\mathbf{SM}_\basefield(\Qb)$.

The scheme $\pi_1^{\mathbf{SM}}(U,a)$ is isomorphic to its Lie algebra via the exponential map and therefore may be viewed as a pro-vectorial scheme in $\mathbf{SM}_\basefield(\Qb)$. We consider the twist
\[
\Lie \pi_1^{\mathbf{SM}}(U,a)(-1) = \underline{Hom}(\Qb(1),\Lie \pi_1^{\mathbf{SM}}(U,a))
\]
which is a pro-vectorial scheme in $\mathbf{SM}_\basefield(\Qb)$ but without a Lie bracket. The $\underline{Hom}$ refers to internal Hom of pro-vectorial schemes in $\mathbf{SM}_\basefield(\Qb)$. We similarly refer to $\pi_1^{\mathbf{SM}}(U,a)(-1)$ as the same scheme as $\pi_1^{\mathbf{SM}}(U,a)$ but with twisted $\mathbf{SM}_\basefield(\Qb)$-structure, isomorphic to $\Lie \pi_1^{\mathbf{SM}}(U,a)(-1)$ via the same exponential map.

There is a canonical map $\Qb(1) \to \pi_1^{\mathbf{SPM}}(U,b)$ given by local monodromy, and we compose with the change of basepoint map $(p,\gamma) \mapsto p \gamma p^{-1}$ to get a map
\[
\pi_1^{\mathbf{SPM}}(U;a,b) \times \Qb(1) \to \pi_1^{\mathbf{SPM}}(U;a,b) \times \pi_1^{\mathbf{SPM}}(U,b) \to \pi_1^{\mathbf{SPM}}(X,a).
\]

This produces a map
\[
\pi_1^{\mathbf{SPM}}(U;a,b) \to \underline{Hom}(\Qb(1),\Lie \pi_1^{\mathbf{SM}}(U,a)) \simeq \Lie \pi_1^{\mathbf{SM}}(U,a)(-1) \simeq \pi_1^{\mathbf{SM}}(U,a)(-1).
\]

This map is invariant under the right action of the local monodromy $\Qb(1) \subseteq \pi_1^{\mathbf{SPM}}(U,b)$. The existence of a tangential basepoint $b$ implies that $U$ is affine, so its geometric fundamental group is free, so that \cite[Lemma 4.7]{DelGon05} applies and thus \cite[Lemma 4.5]{DelGon05} shows that
\[
\pi_1^{\mathbf{SPM}}(U;a,\overline{b}) \coloneqq  \pi_1^{\mathbf{SPM}}(U;a,b)/\Qb(1)
\]
embeds into $\pi_1^{\mathbf{SM}}(U,a)(-1)$. This implies that $\pi_1^{\mathbf{SPM}}(U;a,\overline{b})$ is in $\mathbf{SM}_\basefield(\Qb)$.

The map $u \colon U \to \Gm$ induces a map $\pi_1^{\mathbf{SPM}}(U,b) \to \pi_1^{\mathbf{SPM}}(\Gm,(0,u(v)))$ sending the local monodromy $\Qb(1)$ at $b$ isomorphically onto $(\Gm,(0,u(v))) \simeq \Qb(1)$ because $u$ is a uniformizer. It also induces a map $\pi_1^{\mathbf{SPM}}(U;a,b) \twoheadrightarrow \pi_1^{\mathbf{SPM}}(\Gm;u(a),(0,u(v)))$ that is equivariant for the $\Qb(1)$-actions on both, the latter of which makes $\pi_1^{\mathbf{SPM}}(\Gm;u(a),(0,u(v)))$ into a $\Qb(1)$-torsor.

The product map
\[
\pi_1^{\mathbf{SPM}}(U;a,b) \to \pi_1^{\mathbf{SPM}}(\Gm;u(a),(0,u(v))) \times \pi_1^{\mathbf{SPM}}(U;a,\overline{b})
\]
is an isomorphism. Finally, the torsor $\pi_1^{\mathbf{SPM}}(\Gm;u(a),(0,u(v)))$ is the same as $\pi_1^{\mathbf{SPM}}(\Gm;u(a),du(v)) \simeq \pi_1^{\mathbf{SPM}}(\Gm;1,\frac{du(v)}{u(a)})$ by \cite[15.51]{Deligne89}. Since the latter is between elements of $\Gm(\basefield)$, it is in $\mathbf{SM}_\basefield(\Qb)$. This implies that $\pi_1^{\mathbf{SPM}}(U;a,b)$ is as well. \qed

\begin{ssRemark}\label{rem:less_ad_hoc_tangential}
One expects a less ad hoc proof of Theorem \ref{thm:motivic_torsor_tangential} to not require the assumption $X(\basefield) \neq \emptyset$. In fact, there are two approaches:
\begin{enumerate}
	\item\label{item:ayoub} Modify \cite[Construction 4.5.6]{AyoubAnabelianPresentationPreprint} by replacing $\rm{LS}_{\rm{geo}}(X,\Qb)^{\heartsuit,\otimes}$ with the subcategory of unipotent objects (the thick subcategory generated by the unit object), which uses the theory of motivic nearby cycles developed in \cite[\S 3]{AyoubSixOpII}.\footnote{This construction may be phrased in the language of \cite{AyoubHopf14}, as explained to the authors by Ayoub, by replacing the triangulated category $\mathbf{SmDA}$ in 2.4.1 by its thick subcategory $\mathbf{SmDA}^{\un}$ generated by $\mathbbm{1}_{\mathbf{SmDA}}$ to get a functor
		\[
		\phi^{\un}_{\Delta*} \colon \mathbf{SmDA}^{\un}(X,\Qb) \to \mathbf{SmDA}^{\un}(X^2,\Qb),
		\]
		taking $\phi_{\Delta*}^{\un}(\mathbbm{1}_{\mathbf{SmDA}})$ to get a Hopf algebroid over $X$, then taking its fiber over $(a,b)$ using Corollaire 2.40 or \cite[Definition 3.2.18]{AyoubAnabelianPresentationPreprint} in the tangential case.}
	\item\label{item:log_structure} Put an appropriate log structure on the complex $P^n(X;a,b)$ and then use logarithmic motives as in \cite{TurDorvault25Tangential}
%  then apply \cite{HowellThesisLogMotive} or \cite{ShuklinLogMotive2024}.
\end{enumerate}
%The second approach is mentioned in \cite[1.7.4]{dupont2024logarithmicmorphismstangentialbasepoints}, and C. Dupont has informed us that S. Tur-Dorvault will complete this in the near future.
Note also \cite{park2024logmotivicnearbycycles}, which provides a bridge between the two approaches. Using either of these, all of \S \ref{sec:CK_setup} applies even if $X(\basefield) = \emptyset$.
\end{ssRemark}

\subsubsection{}\label{sec:path_torsor_\basefield,S}

Let $w \in S_f(\basefield)$. Let $X$ be a smooth geometrically connected variety over $\basefield$ with smooth compactification $X'$. Suppose $X$ has strong good reduction at $w$, i.e., it has a smooth model $\Xc'$ over $\Oc_w$ in which the closure $\Xc$ of $X$ has complement $\Yc$ \'etale over $\Oc_w$.

Let $b=(\overline{b},v)$ be a tangential basepoint and $\Ic \subseteq \Oc(\Xc')$ the ideal sheaf of functions vanishing at $\overline{b}$ (equivalently, at its closure in $\Xc'$). We say that $b$ is \emph{integral at $w$} if $v$ sends the image of $H^0(\Xc';\Ic/\Ic^2) \to (\Ic/\Ic^2)_{\overline{b}} = T_{\overline{b}}(X')^{\vee}$ to a subset of $\Oc_{w} \subseteq \basefield_{w}$ not contained in $\mf_{w}$.

Let each of $a,b$ be basepoints integral at $w$. Then the $\ell$-adic component of of $\pi_1^{\mathbf{SM}}(X;a,b)$ is unramified (resp. crystalline) at $w \notin S_{\ell}(\basefield)$ (resp. $w \in S_{\ell}(\basefield)$) by \cite[proof of Theorem 2.8]{WildeshausRoP}\footnote{For the case of tangential basepoints, one can argue as in the proof of Theorem \ref{thm:motivic_torsor_tangential} (with the caveat that $\Xc(Z)$ must be nonempty if $a$ and $b$ are both tangential).} (resp. \cite[Theorem 1.11]{OlssonTowards}).

In particular if $X$ has strong good reduction outside $S \subseteq S_f(\basefield)$, and $a,b$ are integral outside $S$, then $\pi_1^{\mathbf{SM}}(X;a,b)$ is an affine scheme in $\mathbf{SM}_{\basefield,S}(\Qb)$.

\subsection{The Chabauty--Kim diagram}\label{sec:CK_setup}

We describe the general setup for Chabauty--Kim. Let $\basefield$ be a number field, $Z=  \Spec{\Oc_{\basefield,S}}$ an open subscheme of $\Spec{\Oc_\basefield}$, and $\Xc \to Z$ a smooth morphism whose generic fiber $X$ is a hyperbolic curve along with smooth proper compactification $\Xc \hookrightarrow \Xc' \to Z$ such that the complement $\Yc = \Xc' \setminus \Xc$ is \'etale over $Z$. We denote by $Y,X',X$ their generic fibers, respectively.

\subsubsection{}\label{SM_fund_grp_of_X}
We fix a $Z$-integral basepoint $b$ and consider the unipotent fundamental group
\[
\pi_1^{\mathbf{SM}}(X,b)
\]
of \S \ref{sec:path_torsor_SM}, which is a pro-algebraic group in $\mathbf{SM}_{\basefield,S}(\Qb)$. Then $\Lie{\pi_1^{\mathbf{SM}}(X,b)}$ is a Lie algebra in $\operatorname{Pro}{\mathbf{SM}_{\basefield,S}(\Qb)}$ whose associated graded for the descending central series filtration is a quotient of the free pronilpotent Lie algebra on $h_1^{\mathbf{SM}}(X) \in \mathbf{SM}_{\basefield,S}(\Qb)$. Letting $J$ be the (generalized) Jacobian of $X$ such that $h_1^{\mathbf{SM}}(X)=h_1^{\mathbf{SM}}(J)$, we find that
\[
\Lie{\pi_1^{\mathbf{SM}}(X,b)} \in \operatorname{Pro}{\mathbf{SM}_{\basefield,S}(\Qb,J)}
\]
and thus that $\pi_1^{\mathbf{SM}}(X,b)$ is a prounipotent group in $\mathbf{SM}_{\basefield,S}(\Qb,J)$.

Given $z \in \Xc(Z)$, the construction of \ref{sec:path_torsor_\basefield,S} gives a torsor $\pi_1^{\mathbf{SM}}(X;z,b)$ under $\pi_1^{\mathbf{SM}}(X,b)$ in $\mathbf{SM}_{\basefield,S}(\Qb,J)$. This defines a map
\[\kappa \colon \Xc(Z) \to H^1(\mathbf{SM}_{\basefield,S}(\Qb);\pi_1^{\mathbf{SM}}(X,b)).\]

\subsubsection{}\label{finite-type_motivic_quotient}
Let $\ppi$ be a finite-type quotient of $\pi_1^{\mathbf{SM}}(X,b)$ in $\mathbf{SM}_{\basefield,S}(\Qb,J)$,\footnote{such as the $n$th quotient $\pi_1^{\mathbf{SM}}(X,b)_n$ for the descending central series} and set $\ppi_{\dR} \coloneqq \om^{\dR}(\pi)$, $\ppi_{\dR,\pf} \coloneqq \om^{\dR,\pf}(\pi) = \ppi_{\dR} \otimes_{\basefield} \basefield_{\pf}$, and $\ppi_p \coloneqq \mathrm{real}_{\et,p}(\ppi)$. Then we have an object
\[
H^1(\mathbf{SM}_{\basefield,S}(\Qb,J);\ppi) = H^1(\mathbf{SM}_{\basefield,S}(\Qb);\ppi).
\]

A priori, it is just a functor from $\Qb$-algebras to pointed sets. Conjecture \ref{conj:FPR_adm} implies that $\Ext^1_{\mathbf{SM}_{\basefield,S}(\Qb,J)}(\mathbbm{1},M)$ is finite-dimensional for any $M \in \mathbf{SM}_{\basefield,S}(\Qb,J)$. Since $\mathbf{SM}_{\basefield,S}(\Qb,J)$ is neutralized by the Betti realization, Theorem \ref{thm:representability} would imply that it is an affine scheme (of finite type by Remark \ref{rem:finite-type_affine}).\footnote{By Remark \ref{rem:indep_of_om}, the $\Qb$-structure doesn't depend on our use of the Betti realization; in the case $\basefield=\Qb$, we get the exact same $\Qb$-structure using the de Rham realization to prove Theorem \ref{thm:representability}.}
\subsubsection{}\label{sssec:CK_diagram}

The projection $\pi_1^{\mathbf{SM}}(X,b) \twoheadrightarrow \ppi$ defines a map \[H^1(\mathbf{SM}_{\basefield,S}(\Qb);\pi_1^{\mathbf{SM}}(X,b)) \to H^1(\mathbf{SM}_{\basefield,S}(\Qb);\ppi),\] which we can compose with $\kappa$ to get $\kappa_{\ppi}$.

We choose a prime $p$ and a place $\pf \in S_p(\basefield)\setminus S$. For $z \in \Xc(\Oc_{\pf})$, we have an \'etale unipotent torsor $\pi_1^{\et,\mathrm{un}}(X;z,b)$ with $G_{\pf}$-action, defining a map $\kappa_{\pf} \colon \Xc(\Oc_{\pf}) \to H^1_f(\basefield_{\pf};\pi_1^{\et,\mathrm{un}}(X,b))$ and thus to $H^1_f(\basefield_{\pf};\ppi_{p})$. This produces the following diagram, whose \'etale incarnation first appeared (at least for $\Pi = \pi_1^{\mathbf{SM}}(X,b)_n$) in \cite[p.120]{kim09}:\footnote{One may also consider a version that takes the product over all $\pf \in S_p(\basefield)$ of the bottom-right two objects as in \cite[\S 3.1 (5)]{DograUnlikely24}.}

\begin{equation}\label{eqn:CK_diagram}
\xymatrix{
	\Xc(Z) \ar[r] \ar[d]^{\kappa_{\ppi}} & \Xc(\Oc_{\pf}) \ar[d]_{\kappa_{\pf,\ppi}} \ar[dr]^{\int_{\pf}} \\
	H^1(\mathbf{SM}_{\basefield,S}(\Qb,J);\ppi) \ar[r]_-{\mathrm{loc}_{\ppi}} & H^1_f(\basefield_{\pf};\ppi_{p}) \ar[r]^-{\sim}_-{\log_{\mathrm{BK}}} &  \operatorname{Res}_{\basefield_{\pf}/\Qp} F^0 \backslash \ppi_{\dR,\pf}.
}
\end{equation}

Note that $\kappa_{\ppi}$ maps into the $\Qb$-points of $H^1(\mathbf{SM}_{\basefield,S}(\Qb,J);\ppi)$, $\kappa_{\pf,\ppi}$ maps into the $\Qp$-points of $H^1_f(\basefield_{\pf};\ppi_{p})$, and $\mathrm{loc}_{\ppi}$ is really a map from the base-change of $H^1(\mathbf{SM}_{\basefield,S}(\Qb,J);\ppi)$ from $\Qb$ to $\Qp$.

\subsubsection{}\label{sec:BK_log}

The map $\log_{\mathrm{BK}}$ is defined as follows. Let $R$ be a $\Qp$-algebra and $P \in H^1_f(K_{\pf};\pi_p)(R) = H^1_f(K_{\pf};(\pi_p)_R)$. We note $P_{\cri} \coloneqq \operatorname{D}_{\cris}(P_p)$ and $P_{\dR} \coloneqq \operatorname{D}_{\dR}(P_p)$ over $R_{\pf,0} \coloneqq R \otimes_{\Qp} K_{\pf,0}$ and $R_{\pf} \coloneqq R \otimes_{\Qp} K_{\pf}$, respectively. By \cite[Lemma 1]{kim09}, there is a unique Frobenius-invariant point $\ga^{\cri} \in P_{\cri}(R_{\pf,0})$. Choosing $\ga^{\m{H}}$ in $F^0 T_{\dR}(R_{\pf})$, we set
\[
\log_{\mathrm{BK}}([P_p]) \coloneqq (\ga^{\m{H}})^{-1} \ga^{\cri} \in F^0 \backslash \pi_{\dR,\pf}(R_{\pf}) = \operatorname{Res}_{\basefield_{\pf}/\Qp} (F^0 \backslash \ppi_{\dR,\pf})(R).
\]

By Lemma \ref{h0.5}, the map $\beta^L$ of \ref{left_variant} is the same as $\log_{\mathrm{BK}} \circ \mathrm{loc}_{\ppi}$.

\subsubsection{}\label{sssec:maps_to_kim}

In \cite[p.120]{kim09}, Kim describes a diagram
\begin{equation}\label{eqn:original_CK_diagram}
\xymatrix{
	\Xc(R) \ar[r] \ar[d] & \Xc(R_v) \ar[d] \ar[r] & Res^{F_v}_{\Qp} (U_n^{DR}/F^0)(\Qp) \\
	H^1(G_T,U_n^{et}) \ar[r] & H^1(G_v,U_n^{et})(\Qp) \ar[ur].
}
\end{equation}

Let $\Pi=\Pi_n$, the $n$th quotient by the descending central series. In our notation, $v=\pf$, $F=\basefield$, $R=\Oc_{\basefield,S}$, $R_v = \Oc_{\pf}$, $T = S \cup S_p(\basefield)$, $U_n^{et} = (\Pi)_p$, and $U_n^{DR} = \Pi_{\dR,\pf}$. Then our diagram (\ref{eqn:CK_diagram}) is nearly the same as (\ref{eqn:original_CK_diagram}) except for two differences
\begin{enumerate}
    \item \label{item:global_different} The object $H^1(\mathbf{SM}_{\basefield,S}(\Qb,J);\ppi)$ is different from $H^1(G_T,U_n^{et})$ and in particular lies over $\Qb$ instead of $\Qp$
    \item \label{item:BK_log_different} The map in Kim's diagram outputs $(\ga^{\cri})^{-1} \ga^{\m{H}}$ instead of our $(\ga^{\m{H}})^{-1} \ga^{\cri}$.
\end{enumerate}

Nonetheless, we claim that our diagram maps to Kim's diagram. For (\ref{item:BK_log_different}), we can apply inversion on $\Pi_{\dR,\pf} = U_n^{DR}$ to switch between the two versions. For (\ref{item:global_different}), we describe a map $\mathrm{r}_{\et,p} \colon H^1(\mathbf{SM}_{\basefield,S}(\Qb,J);\ppi)_{\Qp} \to H^1(G_T,U_n^{et})$ defined using the functor $\mathrm{real}_{\et,\ell}$ (for $\ell=p$) of \ref{sec:fiber_func_PSPM}.

First, note that $\mathrm{real}_{\et,p}$ restricts to
\[
\mathrm{real}_{\et,p} \colon \mathbf{SM}_{\basefield,S}(\Qb,J) \to \Rep_{\Qp}^{\mathrm{sf,S}}(G_{\basefield},J).
\]
Since the functor is $\Qb$-linear, and the target is a Tannakian category over $\Qp$, we can view it as a functor out of $\mathbf{SM}_{\basefield,S}(\Qb,J)_{\Qp}$. In particular, we get a map
\[
H^1(\mathbf{SM}_{\basefield,S}(\Qb,J);\ppi)_{\Qp} \overset{\ref{cohomology_in_T}}{=}
H^1(\mathbf{SM}_{\basefield,S}(\Qb,J)_{\Qp};\ppi_{\Qp})
\xrightarrow{\mathrm{real}_{\et,p}}
H^1(\Rep_{\Qp}^{\mathrm{sf,S}}(G_{\basefield},J);\ppi_p).
\]

Finally, we recall by \cite[]{CorwinTSV} that there is a map $H^1(\Rep_{\Qp}^{\mathrm{sf,S}}(G_{\basefield},J);\ppi_p) \to H^1(G_{\basefield,T};\ppi_p) = H^1(G_T,U_n^{et})$ that on $\Qp$-points sends a torsor in $\Rep_{\Qp}^{\mathrm{sf,S}}(G_{\basefield},J)$ to the corresponding torsor in $\Qp$-schemes with $G_{\basefield,T}$-action. Composing this with the map defined by $\mathrm{real}_{\et,p}$, we get our map $\mathrm{r}_{\et,p}$. Since the motivic path torsor of \S \ref{sec:fund_grp_and_path_torsors}-\ref{sec:tangential_basepoints} has 'etale realization equal to the usual \'etale path torsor, the map $\kappa_{\ppi}$ is compatible with Kim's map.

We thus have a commutative diagram
\begin{equation*}\label{eqn:map_between_CK_diagrams}
\xymatrix{H^1(\mathbf{SM}_{\basefield,S}(\Qb,J);\ppi)_{\Qp} \ar[r]_-{\mathrm{loc}_{\ppi}} \ar[d]^-{\mathrm{r}_{\et,p}} & H^1_f(\basefield_{\pf};\ppi_{p}) \ar[r]^-{\sim}_-{\log_{\mathrm{BK}}} \ar@{=}[d] &  \operatorname{Res}_{\basefield_{\pf}/\Qp} F^0 \backslash \ppi_{\dR,\pf} \ar[d]^-{(\cdot)^{-1}} \\
H^1(G_T,U_n^{et}) \ar[r] & H^1(G_v,U_n^{et}) \ar[r]^-{\sim}
& 
Res^{F_v}_{\Qp} (U_n^{DR}/F^0),
}
\end{equation*}
where the top row is as in (\ref{eqn:CK_diagram}) and the bottom row as in (\ref{eqn:original_CK_diagram}). It is compatible with the Kummer maps and thus fits into a map from (\ref{eqn:CK_diagram}) to (\ref{eqn:original_CK_diagram}).

%Let $G^{\et}_p = G^{\et}_p(\Oc_{K,S},J) \coloneqq \pi_1(\mathbf{SM}_{\basefield,S}(\Qb,J),\om^{\et,p}$. Notice that $G_T$ acts on $\om^{\et,p}(M)=M_{p}$ for every $M \in \mathbf{SM}_{\basefield,S}(\Qb,J)$, so by definition of Tannakian $\pi_1$, we get a map
%\[
%G_T \to G^{\et}_p(\Qp).
%\]

%In order to avoid base-changing $\mathrm{real}_{\et,p}$ to an arbitrary $\Qp$-algebra, we define it in the following way. Notice that $G_T$ acts on $\om^{p}(M)=M_{p}$ for every $M \in \mathbf{SM}_{\basefield,S}(\Qb,J)$

\subsection{Conclusion}\label{ssec:conclusion}

We put ourselves in the situation of \S \ref{sec:CK_setup}. Now that we have formulated our Chabauty--Kim diagram in terms of the category $\mathbf{SM}_{\basefield,S}(\Qb,J)$, let us apply the theory of \S \ref{sec:geom_loc_real}. Specifically, we apply \S \ref{sec:abstract_setup} to the weight-filtered Tannakian category $\Tc = \mathbf{SM}_{\basefield,S}(\Qb,J)$ and Hodge-filtered Frobenius equivariant fiber functor $(\om^{\dR,\pf}_H,\om^{\cri},\mathrm{id})$ over $(\basefield_{\pf},\basefield_{\pf,0},\Frob_{\pf})$ of \ref{sec:THE_CATS_FOR_CK}.

Recall by \S \ref{left_variant}(*) that
\[
H^1(\mathbf{SM}_{\basefield,S}(\Qb,J);\ppi)_{\basefield} \simeq H^1(G;\ppi_{\dR}) \simeq Z^1(U^{\gr};\om^{\dR,\gr}(\ppi))^{\Gb(J)} \overset{\tau^\m{H}}{\simeq} Z^1(U;\om^{\dR}(\ppi))^{\Gb_H(J)}.
\]

\begin{thm*}\label{thm:final_theorem}

The diagram (\ref{eqn:CK_diagram}), which maps to Kim's diagram (\ref{eqn:original_CK_diagram}) by \ref{sssec:maps_to_kim}, may be identified with
\begin{equation}\label{eqn:CK_diagram2}
\xymatrix{
	\Xc(Z) \ar[r] \ar[d]^{\kappa_{\ppi}} & \Xc(\Oc_{\pf})  \ar[d]^{\int_{\pf}} \\
	Z^1(U;\ppi_{\dR})^{\Gb_H(J)} \ar[r]_-{F^0 \backslash \ev_{\uU_{\pf}^L}} & F^0 \backslash \operatorname{Res}_{\basefield_{\pf}/\Qp} (\ppi_{\dR})_{\basefield_{\pf}},
}
\end{equation}
\end{thm*}
\begin{proof}
This is an immediate corollary of Theorem \ref{thm:left_commutes}.
\end{proof}

This provides a foundation for future work in non-abelian Chabauty beyond Quadratic Chabauty. In particular, it puts us in a generalization of the setup of \cite[\S 2.4.2]{PolGonI}. As in that case, one may write $\ev_{\uU_{\pf}^L}$ as the pullback along $\uU_{\pf}^L \colon \Spec{\basefield_{\pf}} \to U$ itself of the universal evaluation map
\[
Z^1(U;\ppi_{\dR})^{\Gb_H(J)} \times_{\basefield} U \to F^0 \backslash \om^{\dR}(\pi) \times_{\basefield} U.
\]

The computation of this image is thus a purely geometric problem over the function field $\basefield(U)$.

\appendix

\section{Mixed Motives}\label{app:mixed_motives}

\label{sec:ab_mixed_motives}

\subsection{}
Abelian categories of mixed motives are the ultimate candidate and main motivation behind the notion of weight-filtered Tannakian category. As they are mostly conjectural, we have relegated them to this short appendix.

\subsection{}
\label{sec:beilinson_mot}
%\commentDC{ID: is Beilinson Motives the terminology from FG? DC: Not just from there. Also in the Handbook of K-theory reference. ID: Fine. In C--D this refers to one particular construction.}
As originally conjectured in \cite[5.10]{Beilinson87} and described in more detail in \cite[\S 4.4.2, \S 5.3, Conjecture 33]{HandbookKTheory}, there should be a category $\MM_\basefield$ of \emph{mixed motives over $\basefield$} for any field $\basefield$. The rationalization $\mathrm{MMot}(\basefield,\Qb) = \MM_\basefield \otimes \Qb$ should be a Tannakian category, and for any embedding $\iota \colon \basefield \hookrightarrow \Cb$, there should be a corresponding rational singular cohomology realization functor $\om^{\iota}$ that is a $\Qb$-rational fiber functor. One similarly expects $\basefield$-rational (resp. $\Ql$-rational) fiber functors $\om^{\dR}$ (resp. $\om^{\ell}=\om^{\et,\ell}$) from de Rham (resp. $\ell$-adic) cohomology, as long as $\basefield$ has characteristic $0$ (resp. different from $\ell$).

\subsection{}
\label{sec:derived_beilinson}
The derived category of bounded complexes of $\mathrm{MMot}(\basefield,\Qb)$ is expected to be naturally isomorphic to $\DM_{\gm}(\basefield,\Qb)$, as implied for example in \cite[\S A.1-A.2]{CisDegTriCat}. This natural isomorphism should further intertwine $\om^{\iota}$, $\om^{\dR}$, and $\om^{\et,\ell}$ of \ref{sec:beilinson_mot} with $\om^{B,v}$, $\om^{\dR}$, and $\om^{\et,\ell}$ of \ref{sec:fiber_func_PSPM} via this isomorphism and $\mathrm{real}_{\mathbf{SPM}}$ of Theorem \ref{thm:real}.

\subsection{}
\label{sec:motives_weight_filtr}

The properties listed in \cite[\S 4.4.2]{HandbookKTheory} imply that the weight filtration makes $\mathrm{MMot}(\basefield,\Qb)$ into a weight-filtered Tannakian category, as follows. 

First, we verify that the filtration is strict. For this, we verify the hypotheses of Lemma \ref{lemma:strictness_crit}. For $M,N \in \Ob(\Tc)$ and $m \neq n$, recall that $\gr_m M$ and $\gr_n N$ are in the $m$th and $n$th graded parts, respectively, of the graded category $\mathrm{Mot}_{\mathrm{num}}(\basefield,\Qb)$, which is the socle (full subcategory of semisimple objects) of $\mathrm{MMot}(\basefield,\Qb)$. Any subquotient of either is semisimple and thus a subquotient in $\mathrm{Mot}_{\mathrm{num}}(\basefield,\Qb)$, hence of weights $m$ and $n$, respectively. In particular, no two nonzero such subquotients can be isomorphic.

%**It is strict because** (something about socle being graded, i.e. no morphisms between motives of different weights? should I make a prop about socle being graded implies strictness?)

Finally, every $\gr^W_n M \in \mathrm{Mot}_{\mathrm{num}}(\basefield,\Qb)$ and is therefore semisimple, so the filtration indeed satisfies Definition \ref{defn:weight_filtr_cat}.

%Choosing $\ell \neq \char{k}$, \cite[\S 5.3, Conjecture 33.5]{HandbookKTheory} implies that there is a $\Ql$-linear fiber functor $\om_{\ell}$ for which 

\subsubsection{}
\label{sec:motives_subcategory}
By \ref{rem:filtr_sub}, any Tannakian subcategory of $\mathrm{MMot}(\basefield,\Qb)$ should also be a weight-filtered Tannakian category. (In particular, we are interested in subcategories of the form $\Mm_{\Pp}$ as described in \cite[\S 3]{GonMEM}.) While this observation may seem trivial, it is important because certain subcategories of the conjectural $\mathrm{MMot}(\basefield,\Qb)$ have been constructed either entirely non-conjecturally or at least more explicitly.

\subsection{}
\label{sec:MTM}
In particular, the category $\mathrm{MT}(\Spec{\basefield},\QQ)$ of \cite[1.1]{DelGon05} is a weight-filtered Tannakian category that is expected to be a Tannakian subcategory of $\mathrm{MMot}(\basefield,\Qb)$. The former follows from the weight filtration described in loc.cit., including the exactness of the associated graded and the fact that the associated graded is a direct sum of objects of the form $\Qb(n)$, which are all simple. This category (and its unramified subcategories) has already found extensive applications to the Chabauty--Kim method, as in \cite{MTMUE,MTMUEII,PolGonI}.

\subsection{}
\label{sec:MEM}
Candidates have been constructed for categories of mixed elliptic motives (\cite{Patashnick13},\cite{Cao18}). Such a category is expected to be the thick Tannakian subcategory of $\mathrm{MMot}(\basefield,\Qb)$ generated by $h_1(E)$ for a fixed elliptic curve $E/\basefield$.

\subsection{}
\label{sec:motives_hodge}
As mentioned in \ref{sec:beilinson_mot}, if $\chara{\basefield}=0$, then algebraic de Rham cohomology provides a fiber functor $\om^{\dR}$ over $\basefield$.
%\commentDC{ID: On which category? DC: On the category $\mathrm{MMot}(k,\Qb)$, which is the topic of this section.}
This functor comes equipped with a Hodge filtration that makes $\om^{\dR}$ into a Hodge-filtered fiber functor as in \ref{defn:hodge_filtr}.
%\comment{ID: Why "conjecturally"? DC: You're probably right - I probably don't need to say conjecturally.}

\subsection{}
\label{sec:motives_frobenius}
The existence of a Frobenius-equivariant fiber functor
%\comment{ID: Do you mean a Hodge-filtered Frobenius-equivariant fiber functor? DC: Sorry, I meant Frobenius-equivariant.}
is a bit more subtle. For this, we must restrict to motives over a number field $\basefield$ that are unramified at some finite place $\pf$ of $\basefield$.
%\commentDC{ID: I'm not sure I agree with this. We do need a $\phi$ and an $F$, but we don't need admissibility. DC: If you don't have crystalline representations, then how can $D_{\cris}$ give you a fiber functor? It doesn't have the right dimensions, so in particular shouldn't be exact, especially if you have a non-crystalline extension of crystalline representations. ID: You certainly don't need to restrict attention to crystalline stuff. I believe the state of the art in terms of the crystalline realization functor is in Deglise--Niziol.}
The unramifiedness corresponds to crystallinity when considering $p$-adic \'etale cohomology for $p$ the characteristic of $\pf$. 

\subsubsection{}
\label{sec:unramified_motives}
Under the given notation, it is most natural to refer to such a subcategory as $\mathrm{MMot}(\Oc_{\basefield,\pf},\Qb)$, where $\Oc_{\basefield,\pf}$ denotes the localization of the ring $\Oo_\basefield$ at the prime ideal $\pf$. Such a category should sit naturally into a square
\[
\xymatrix{\mathrm{MMot}(\Oc_{\basefield,\pf},\Qb) \ar[r] \ar[d] & \mathrm{MMot}(k,\Qb) \ar[d]\\
\DM_{\gm}(\Spec \Oc_{\basefield,\pf}, \QQ) \ar[r] & \DM_{\gm}(\Spec{\basefield}, \QQ),
}
\]
where the right vertical arrow is given by the inclusion functor $\mathrm{MMot}(k,\Qb) \hookrightarrow D^b(\mathrm{MMot}(k,\Qb))$ sending an object $M$ to the complex $M$ concentrated in degree $0$ followed by the expected equivalence of \ref{sec:derived_beilinson}, and $\DM_{\gm}(\Spec \Oc_{\basefield,\pf}, \QQ)$ is as in \cite[Definition 11.1.10]{CisDegTriCat}.

In the special case of $\mathrm{MT}(\Spec{\basefield},\QQ)$, we may define $\mathrm{MT}(\Spec \Oc_{\basefield,\pf},\QQ)$ as the union over $Z \ni \pf$ of $\mathrm{MT}(Z,\QQ)$ (\cite[1.6]{DelGon05}), and we expect that $\mathrm{MT}(\Spec \Oc_{\basefield,\pf},\QQ) = \mathrm{MMot}(\Oc_{\basefield,\pf},\Qb) \cap \mathrm{MT}(\Spec{\basefield},\QQ)$. 

\subsubsection{}
\label{sec:crystalline_realization}
Whether working with the conjectural $\mathrm{MMot}(\Oc_{\basefield,\pf},\Qb)$ or the non-conjectural $\mathrm{MT}(\Spec \Oc_{\basefield,\pf},\QQ)$, we expect a $\basefield_{\pf,0}$-linear crystalline realization fiber functor $\om^{\cris}$ that intertwines via \ref{sec:derived_beilinson} with $R_{\rig}$ (\cite[17.2.23]{CisDegTriCat}). This functor should come equipped with natural isomorphisms $\om^{\cris} \otimes_{\basefield_{\pf,0}} \basefield_{\pf} \simeq \om^{\dR} \otimes_\basefield \basefield_{\pf}$ and $\om^{\cris} \simeq \operatorname{D}_{\cris} \circ \om^{\et,p}$. The latter property in fact can serve as a definition of $\om^{\cris}$, for which one need only assert that $\om^{\et,p}(M)$ is crystalline at $\pf$ for all $M \in \mathrm{MMot}(\Oc_{\basefield,\pf},\Qb)$.

\subsubsection{}
\label{sec:crystalline_realization_frob_equiv}

%\comment{ID: Let's not do that, since the referee will write that all our results are conditional. DC: I agree it should be phrased differently. Like ``Supposing, as widely expected, that...''. This is just a remark and is only one small of example of our main result. And remember, we have a few other categories - it's just the part about mixed motives that's conditional. ID: I changed the wording a bit.}
It is widely expected that $\{\om^{\et,\ell}(M)\}_{\ell}$ form a strongly compatible system of Galois representations and that $\om^{\et,\ell}(\gr^W_n M)$ is pure of weight $n$ for any $M \in \mathrm{MMot}(\basefield,\Qb)$ (\cite[p.82]{TaylorGaloisReps}). If this is the case, then $\om^{\cris}$ is a Frobenius-equivariant fiber functor over $(\basefield_{\pf,0},\Frob_{\pf})$.

\section{Background and Technicalities on Tannakian Categories}\label{appendix:tannakian_categories}

Here, we review the theory of Tannakian categories along with some technicalities on stacks and groupoids that we need mainly in \S \ref{sec:arith_frob_path}-\ref{sec:geom_loc_real}.

\subsection{Basics on Tannakian categories}\label{sec:tann_basic}

Let $\coeff$ be a field. This material can all be found in \cite{SaavedraRivanoBook,DelMil_TC22,DelTann90}.

\subsubsection{Preliminaries}\label{sec:tann_prelim} An \emph{additive category} is a category in which all finite biproducts (simultaneous coproducts and products) exist; such a category is automatically enriched over the monoidal category of abelian groups with tensor product. A \emph{$\coeff$-linear category} is an additive category further enriched over $\coeff$-vector spaces and $\otimes_\coeff$. We assume notions of \emph{symmetric monoidal category} (what is called \emph{ACU $\otimes$-category} in \cite[I.2.4.3]{SaavedraRivanoBook} and \emph{tensor category} in \cite[II.1.1]{DelMil_TC22}), \emph{dualizable object} (\cite[2.1.2]{DelTann90}), and \emph{abelian category}. A symmetric monoidal category is \emph{rigid} (also called \emph{autonomous}) if all objects are dualizable.

The term \emph{tensor category} seems to have many different definitions in the literature, which range from simply meaning symmetric monoidal category (above) to \cite[2.1]{DelTann90}, which requires abelian, rigid, and the condition that $k \xhookrightarrow{\sim} \End(\mathbbm{1})$. Therefore, we use the term only in the context `X tensor category' for $\emptyset \neq X \subseteq \{\mathrm{additive}, \mathrm{abelian}, k-\mathrm{linear}\}$; this means a symmetric monoidal category with structure $X$, such that both are compatible (e.g., the tensor product respects biproducts and/or is $k$-linear).

Thus a \emph{rigid abelian $\coeff$-linear tensor category} is a rigid symmetric monoidal abelian $\coeff$-linear category in which the tensor product is biadditive (respects biproducts) and $\coeff$-bilinear.% Note that \cite[2.1]{DelTann90} reserves the term \emph{tensor category} for this case along with the condition $k \xhookrightarrow{\sim} \End(\mathbbm{1})$.

Any abelian category is in particular an \emph{exact category}; that is, there is a notion of short exact sequence satisfying the conditions listed in \cite[Definition 2.1]{BuehlerExactCats10}.

%A \emph{tensor functor} between linear tensor categories refers to a linear symmetric monoidal functor.

\subsubsection{Fiber functors}\label{sec:fiber_functors} Let $\Tc$ be a symmetric monoidal $\coeff$-linear exact category. 
A \emph{fiber functor} on $\Tc$ with values in a $\coeff$-scheme $S$ is an exact tensor functor $\Tc \to \QCoh(S)$. 
Note that if $\Tc$ is rigid, then any such functor lands in $\LF(S)$ (\S \ref{grfil_in_families}).\footnote{By \cite[\href{https://stacks.math.columbia.edu/tag/0FNU}{Tag 0FNU}(1) $\Leftrightarrow$ (3)]{stacks-project} and \cite[\href{https://stacks.math.columbia.edu/tag/05P2}{Tag 05P2}(1) $\Leftrightarrow$ (4)]{stacks-project}, $\LF(S)=\QCoh(S)_{\rig}$, and by \cite[2.7]{DelTann90}, a symmetric monoidal functor sends dualizable objects to dualizable objects.}. 
We denote the groupoid of exact $\coeff$-linear symmetric monoidal functors and symmetric monoidal natural transformations by $\operatorname{Hom}^{\otimes}_{\mathrm{ex}}$, so that the groupoid of $S$-valued fiber functors is
\[
\operatorname{Hom}^{\otimes}_{\mathrm{ex}}(\Tc,\QCoh(S))
\]
or, when $\Tc$ is rigid,
\[
\operatorname{Hom}^{\otimes}_{\mathrm{ex}}(\Tc,\LF(S)).
\]

If $\om$ is a fiber functor with values in $S$, and $S'$ is an $S$-scheme, we denote by $\om_{S'}$ or $\om_{S'/\coeff}$ the composition of $\om$ with the pullback functor $\QCoh(S) \to \QCoh(S')$. Note that $\om_{S/\coeff}=\om$.

\subsubsection{Definition}\label{defn:tannakian_category} A \emph{Tannakian category $\Tc$ over $\coeff$} is a rigid abelian $\coeff$-linear tensor category possessing a fiber functor over a nonempty $\coeff$-scheme such that the natural map $\coeff \to \End(\mathbbm{1})$ given by $\alpha \mapsto \alpha \id_{\mathbbm{1}}$ is an isomorphism.

For the rest of \S \ref{sec:tann_basic}, we let $\Tc$ be a Tannakian category over $\coeff$.

\subsubsection{Tannakian fundamental groups}

If $\om$ is a fiber functor on $\Tc$ with values in a $\coeff$-scheme $S$, then $\pi_1(\Tc,\om) \coloneqq \underline{\Aut}^{\otimes}(\om)$, the group tensor automorphisms of $\om$, is naturally a flat affine group schemes over $S$, known as the \emph{Galois group of $\Tc$ at $\om$}. Its functor of points is defined, for an $S$-scheme $S'$, by
\[
\pi_1(\Tc,\om)(S') = \underline{\Aut}^{\otimes}(\om_{S'}).
\]

\subsubsection{Tannakian subcategories}

A \emph{Tannakian subcategory} is a full abelian subcategory containing a unit object and closed under tensor products, duals, and subquotients.

In particular, if $\Tc' \subseteq \Tc$ is a Tannakian subcategory and $\om$ is a fiber functor on $\Tc$, then $\pi_1(\Tc,\om) \to \pi_1(\Tc',\om)$ is faithfully flat (\cite[6.6]{Deligne89} references a proof in the neutral case and states that it is true in general; c.f. also \cite[\S 2]{BertolinMotivicGalois08}).

\subsubsection{Algebraic Tannakian Categories}\label{sec:algebraic_tannakian}

Let $\om$ be a fiber functor over $\coeff'$. Following \cite[Remark 3.10]{DelMil_TC22} or \cite[III.3.3.1]{SaavedraRivanoBook}, $\Tc$ is \emph{algebraic} if any of the following equivalent conditions hold:
\begin{itemize}
\item $\pi_1(\Tc,\om)$ is finite-type as an $\coeff'$-scheme (i.e., it is an \emph{algebraic} rather than \emph{pro-algebraic} group)
\item $\Tc$ is tensor-generated by a single object (``has a tensor generator'')
\item $\Tc$ is tensor-generated by a finite collection of objects
\end{itemize}

From the first condition, it is easy to see that a Tannakian subcategory of an algebraic Tannakian category is algebraic.

\subsubsection{Neutral Tannakian categories} We say that $\Tc$ is \emph{neutral} if it possesses a fiber functor $\om$ with values in $\coeff$. In such a case, Tannakian reconstruction states that the natural functor
\[
\Tc \to \Rep(\pi_1(\Tc,\om))
\]
is an equivalence of Tannakian categories. We discuss a generalization to the non-neutral case in \ref{tann_reconstr_groupoid}.

\subsection{Algebraic geometry in a Tannakian category}\label{algebraic_geometry_tann_cat}

\subsubsection{}\label{pro_ind}

We work with objects of $\Pro{\Tc}$ and $\Ind \Tc$.  Each is an abelian $\coeff$-linear tensor category.

\subsubsection{}\label{pro_ind_fiber_functor}

If $\om \colon \Tc \to \LF(S)$ is a fiber functor, we can extend $\omega$ to exact tensor functors $\Pro \Tc \to \QCoh(S)$ (resp. $\Ind \Tc \to \QCoh(S)$) using the fact that $\QCoh(S)$ is complete (resp. cocomplete).

\subsubsection{}\label{scheme_group_object_in_T}

We recall that an \emph{affine scheme in $\Tc$} is a spectrum in the sense of \cite[5.4]{Deligne89} of a commutative ring in $\Ind{\Tc}$. A \emph{group scheme} $\pi$ in $\Tc$ is a group object in affine schemes in $\Tc$; equivalently, it is the spectrum of a commutative Hopf algebra object $\Oc(\pi)$ in $\Ind{\Tc}$. If $\om$ is a fiber functor over a $\coeff$-scheme $S$, then $\om(\pi) \coloneqq \Spec{\om(\Oc(\pi))}$ is a group scheme over $S$. C.f. \ref{alg_geo_stacks} for another definition of (affine) schemes in $\Tc$.

\subsubsection{}\label{pro_ind_duality} Since $(\Pro{\Tc})^{\mathrm{op}}$ is naturally equivalent to $\Ind \Tc^{\mathrm{op}}$, and duality is an equivalence $\Tc \simeq \Tc^{\mathrm{op}}$, duality gives an equivalence 
\[
\Ind{\Tc} \simeq (\Pro{\Tc})^{\m{op}}.
\]

\subsubsection{}\label{prounip_in_T}  If $\pi$ is a group scheme in $\Tc$, we say that $\pi$ is \emph{prounipotent} if $\omega(\Oc(\pi))$ is the Hopf algebra of a prounipotent group for some (equivalently any) fiber functor $\om$ over a field. The dual of $\Oc(\pi)$ under \ref{pro_ind_duality} is a Hopf algebra in $\Pro{\Tc}$ denoted $\Uc \pi$, as it is naturally the (completed) universal envelopping algebra of $\pi$ under any fiber functor. Then the kernel of $\Uc \pi \xrightarrow{\Delta'} \Uc \pi$ is a Lie algebra $\Lie \pi$ in $\Pro{\Tc}$, and $\Uc \pi$ is a quotient of the tensor algebra $T(\Lie \pi)$.

\subsubsection{}\label{torsors_in_tannakian}

If $\pi$ is a group scheme in $\Tc$, a \emph{torsor in $\Tc$ under $\pi$} is a $\Tc$-scheme $P$ with a group action $\pi \times P \to P$ in $\Tc$-schemes for which $\Oc(P)$ is not the zero object of $\Tc$, and $\pi \times P \to P \times P$ is an isomorphism (equivalently, if some fiber functor $\om$ over some field extension of $\coeff$ makes $\om(P)$ into a torsor under $\om(\pi)$). We let $H^1(\Tc;\pi)(\coeff)$ denote the pointed set of isomorphism classes of torsors in $\Tc$ under $\pi$, with $\pi$ as a torsor under itself forming the distinguished point of this set.

\subsection{Tannakian reconstruction in the non-neutral case}
\label{sec:reconstruction_non_neutral}

% \cite[\S 1]{DelTann90}.

Given a Tannakian category $\Tc$ over $\coeff$ and a fiber functor $\om$ over a $k$-scheme $U$, we get a \emph{$\coeff$-groupoid acting on $U$} as follows (\cite[1.11]{DelTann90}): We let $\pr_1,\pr_2$ denote the projections from $U \otimes_{\coeff} U$ to $U$ and set
\[
\Gc = \Gc(\om) \coloneqq \underline{\Isom}^{\otimes}(\pr_2^* \om, \pr_1^* \om).
\]

Then $\Gc$ is a scheme over $U \times_{\coeff} U$ with projections to $U$ denoted $b$ and $s$. For a $k$-scheme $T$, an element $g \in \Gc(T)$ thus gives a pair $b(g),s(g) \in U(T)$ of $T$-points of $U$, and $g \in \Isom^{\otimes}(s(g)^* \om,b(g)^* \om)$. There is thus a composition map
\[
\circ \colon \Gc \times_{{^s}U^b} \Gc \to \Gc,
\]
where, if $g_1,g_2 \in \Gc$, with $s(g_1)=t_1$, $b(g_1)=s(g_2)=t_2$, and $b(g_2)=t_3$, then
\[
g_2 \circ g_1 \colon t_1^* \om \xto{g_1} t_2^* \om \xto{g_2} t_3^* \om.
\]

Then $(\Gc,\circ,s,b)$ forms a transitive (meaning $\Gc$ is fpqc over $U \times_{\coeff} U$) $k$-groupoid acting on $U$ in the sense of \cite[1.6]{DelTann90}. In the notation of \cite[\href{https://stacks.math.columbia.edu/tag/0231}{Tag 0231}]{stacks-project}, $(U,\Gc,s,b,\circ)$ is a groupoid scheme over $\Spec{\coeff}$.

\subsubsection{}\label{rep_of_groupoid}

A \emph{representation} of $\Gc$ is the data of a quasi-coherent sheaf $V$ on $U$ and maps $\rho(g) \colon V_{s(g)} \to V_{b(g)}$ for every $g \in \Gc(T)$ and $k$-scheme $T$, compatible with $\circ$ and change of $T$. We note $\Rep(U \colon \Gc)$ (\cite[1.7]{DelTann90}) the category of representations for which $V$ is a vector bundle on $S$ (assuming $S$ is non-empty).

\subsubsection{}\label{tann_reconstr_groupoid}
By \cite[Theorem 1.12(ii)]{DelTann90}, the functor $\om$ induces an equivalence between $\Tc$ and $\Rep(U \colon \Gc)$.

Similarly, we get an equivalence between $\Ind \Tc$ and $\Ind \Rep(U \colon \Gc)$. In particular, if $Y$ is a (group) scheme over $U$, to show that $Y$ is $\om$ of a corresponding (group) scheme in $\Tc$ (c.f. \cite[5.4]{Deligne89} or \ref{algebraic_geometry_tann_cat}), it suffices to define an action $\rho$ of $\Gc$ on $Y$. Such an action means giving, for each $g \in \Gc(T)$ for a $\coeff$-scheme $T$, a (homo)morphism $\rho(g) \colon Y_{s(g)} \to Y_{b(g)}$, compatible with $\circ$ and change of $T$.

%, we consider not just cohomology sets but cohomology schemes; it follows that we must consider torsors not just over the base ring $\coeff$ but over any $\coeff$-algebra $R$. For this, we must explain what is the base change $\Tc_R$ of $\Tc$ to $R$ and discuss schemes in this category.

\subsection{Background on stacks}\label{sec:stacks_background}

We need some background on stacks (in particular, the gerbes associated to Tannakian categories) for the purposes of \S \ref{app:base_change_tann}, where we discuss the base-change $\Tc_S$ for a $\coeff$-scheme $S$.

\subsubsection{Stacks}\label{sec:stacks}

A \emph{stack} refers to a stack for the fpqc topology, i.e., a $2$-sheaf of groupoids on the category of affine schemes with the fpqc topology. Following e.g. \cite[Notation]{wedhorn2024extensionliftinggbundlesstacks}, a \emph{prestack} is just a $2$-presheaf of groupoids; i.e., a $(2,1)$-functor from the category of affine schemes to the $(2,1)$-category of groupoids. A scheme, via its functor of points, is a stack valued in discrete groupoids (i.e., those with no non-identity morphisms).

\subsubsection{Sheaves on stacks}\label{sec:sheaves_on_stacks}

A \emph{quasi-coherent sheaf} on a prestack $\Xc$ consists of a collection of functors $\Xc(S) \to \QCoh(S)$ for each scheme $S$ that are compatible (in a $2$-categorical rather than strict sense) with morphisms $S_1 \to S_2$. A \emph{vector bundle} on $\Xc$ is a quasi-coherent sheaf for which each functor $\Xc(S) \to \QCoh(S)$ lands in $\LF(S)$.

The category $\QCoh(\Xc)$ (resp. $\LF(\Xc)$) is the category of quasi-coherent sheaves (resp. vector bundles) on $\Xc$ with an obvious notion of morphism; equivalently, it is the $2$-limit of the categories $\QCoh(S)$ (resp. $\LF(S)$) over the opposite of the fibered category associated to $\Xc$.

The category $\QCoh(\Xc)$ is an abelian tensor category. The subcategory $\LF(\Xc)$ is the full subcategory of dualizable objects; it is an additive tensor category but not usually abelian; but it has an exact structure inherited from that of $\QCoh(\Xc)$.

%We recall some terminology and notation for sheaves on stacks from \cite[Notation]{wedhorn}. Let $\Xc$ be a stack for the fpqc topology (i.e., a $2$-sheaf of groupoids on the category of schemes with the fpqc topology). We let $\LF(\Xc)$ (also known as $\operatorname{Vec}(\Xc)$) denote the category of vector bundles on $\Xc$, i.e., sheaves of $\Oc_{\Xc}$-modules that are locally free of finite type. These are precisely the dualizable objects of $\QCoh(\Xc)$.  It is always an exact symmetric monoidal category (the notion of exact sequence is inherited by embedding it into the abelian category of sheaves of modules on $\Xc$).

%%% the category of sheaves is defined in https://stacks.math.columbia.edu/tag/06WB
%% structure sheaf defined in https://stacks.math.columbia.edu/tag/06TU, category of sheaves mentioned in https://stacks.math.columbia.edu/tag/06TN

\subsubsection{Algebraic stacks}\label{sec:algebraic_stacks}

The term \emph{algebraic stack} (\cite[\href{https://stacks.math.columbia.edu/tag/026O}{Tag 026O}]{stacks-project}) used throughout the Stacks Project and by most algebraic geometers is what we call an \emph{fppf-algebraic} stack: it is a stack with representable, separated, and quasi-compact diagonal that has an fppf (equivalently, smooth) cover by a scheme.

We need the notion of \emph{fpqc-algebraic}, which originally appeared in \cite[Example 1.8]{GoerssPresheaves04} and appears with this name in \cite[Definition 2.29]{goerss2008quasicoherentsheavesmodulistack}. It is a stack whose diagonal is representable, separated, and quasi-compact and which has an \emph{fpqc} cover by a scheme.

\subsubsection{Presentations and groupoids}\label{presentations_groupoids}

If $\Xc$ is fpqc-algebraic and $U \to \Xc$ is an fpqc covering by a scheme over a scheme $S$, then there is a groupoid in algebraic spaces over $S$ (in the sense of \S \ref{sec:reconstruction_non_neutral} or \cite[\href{https://stacks.math.columbia.edu/tag/0231}{Tag 0231}]{stacks-project}) $(R \coloneqq U \times_{\Xc} U,\circ,s,b)$ acting on $U$ with $s,b$ fpqc such that $\Xc = [U/R]$ (in the sense of \cite[\href{https://stacks.math.columbia.edu/tag/044Q}{Tag 044Q}]{stacks-project}). Conversely, if $(R,\circ,s,b)$ is a groupoid scheme over $S$ acting on $U$ for which $s$ (equivalently $b$) is fpqc, then the quotient $\Xc \coloneqq [U/R]$ is an fpqc-algebraic stack over $S$ with fpqc covering by $U$. Such an $(R,\circ,s,b)$ is known as a(n \emph{fpqc) presentation} of $\Xc$.

Then $\QCoh(\Xc)$ is the same as the category of representations of $(R,\circ,s,b)$ as in \ref{rep_of_groupoid}, and $\LF(\Xc) \simeq \Rep(X \colon R)$.

\subsubsection{Gerbes and Tannakian categories}\label{sec:gerbes_tann}

Let $\Tc$ be a Tannakian category over $\coeff$. Consider the functor $\Fib{\Tc}$ from $k$-schemes to groupoids sending $S$ to the groupoid
\[
\Fib{\Tc}(S) \coloneqq \operatorname{Hom}^{\otimes}_{\mathrm{ex}}(\Tc,\LF(S))
\]
of fiber functors on $\Tc$ with values in $S$. Then $\operatorname{Fib}{\Tc}$ is an fpqc-algebraic stack (while it is fppf-algebraic iff $\Tc$ is algebraic) over $\coeff$.

More precisely, as stated in \cite{SaavedraRivanoBook} and proved in \cite{DelTann90}, it is an affine gerbe over $\coeff$.

Given an object $M$ of $\Tc$, we get a functor $\Fib{\Tc}(S) = \operatorname{Hom}^{\otimes}_{\mathrm{ex}}(\Tc,\LF(S)) \to \LF(S)$ given by $\om \mapsto \om(M)$. These form an object of $\LF(\Fib{\Tc})$ as in \ref{sec:sheaves_on_stacks}, and, by varying $M \in \Tc$, this pieces together to a functor \[\Tc \to \LF(\Fib{\Tc}).\] Then the main theorem of Tannakian reconstruction says that this functor is an equivalence of $k$-linear tensor categories.

More formally, the associations $\Tc \mapsto \Fib{\Tc}$ and $\Gc \mapsto \LF(\Gc)$ give an equivalence between the $2$-category of Tannakian categories over $\coeff$ and exact tensor functors and the $2$-category of affine gerbes over $\coeff$.

\subsubsection{Gerbes and groupoid schemes}\label{sec:gerbes_groupoids}

If $\Tc$ is neutral with fiber functor $\om$, then $\Fib{\Tc} = BG$, where $G=\pi_1(\Tc,\om)$. Here, $BG \coloneqq [\Spec{\coeff}/G]$ (also denoted $TORS(G)$ in \cite[II]{DelMil_TC22}) is the classifying stack of the group scheme $G$, and $BG(S)$ is the groupoid of $G$-torsors on $S$.

If $\om$ is a fiber functor over an affine $k$-scheme $U$, then we get a transitive $k$-groupoid $(\Gc,\circ,s,b)$ acting on $U$ as in \S \ref{sec:reconstruction_non_neutral}. The object $\om \in \Fib{\Tc}(U)$ is an fpqc cover of $\Fib{\Tc}$ by $U$, and $(\Gc,\circ,s,b)$ is the $k$-groupoid (for $S=\Spec{\coeff}$) as in \ref{presentations_groupoids}.

Conversely, if we start with a transitive $k$-groupoid $(\Gc,\circ,s,b)$ acting on $U$, then $[U/\Gc]$ is a gerbe equivalent to $\Fib{\Rep(U \colon \Gc)}$.

\subsection{Base-change of Tannakian categories and fiber functors}\label{app:base_change_tann}

Since we work with cohomology (Selmer) \emph{schemes} in \S \ref{sec:abstract_CK_diagram}-\ref{sec:abstract_CK_maps}, we must allow ourselves to work in the base-change $\Tc_S$ of $\Tc$ to an arbitrary $k$-scheme $S$.

\subsubsection{}\label{arbitrary_S} 

%%% see https://stacks.math.columbia.edu/tag/04T3 for the groupoid in algebraic spaces. the page after it shows that [X/R] is an algebraic stack.

As $\Fib{\Tc}$ is a gerbe over $\Spec{\coeff}$, we can base-change it to a gerbe $(\operatorname{Fib}{\Tc})_S$ over $S$. Then we define \[\Tc_S \coloneqq \LF((\operatorname{Fib}{\Tc})_S).\] Pullback along the map $(\operatorname{Fib}{\Tc})_S \to \operatorname{Fib}{\Tc}$ gives an exact tensor functor from $\Tc$ to $\Tc_S$; for $M \in \Tc$, we write $M_S$ for its image in $\Tc_S$.

Then $\Tc_S$ is an $\Oc(S)$-linear category, but it is not generally abelian. More specifically, it is a rigid $\Oc(S)$-linear tensor category, and it has an exact structure inherited from $\QCoh(\Fib{\Tc}_S)$ for which $\otimes$ is exact. C.f. \cite[\S 2.4]{ZieglerFiltGeneralBase} and \cite[\S 10.2]{SchaeppiFormalTannaka13}.

\subsubsection{Fiber functors}\label{fiber_functors_maps_of_stacks}

Given an $S$-scheme $S'$ and a map $f \colon S' \to (\operatorname{Fib}{\Tc})_S$ of stacks over $S$, we get a pullback functor $F_{\LF}(f) = f^* \colon \Tc_S = \LF((\operatorname{Fib}{\Tc})_S) \to \LF(S')$. This functor is an exact $\Oc(S)$-linear tensor functor. We refer to such a functor as an \emph{$S'$-valued fiber functor on $\Tc_S$} and denote the groupoid of such functors by $\operatorname{Hom}^{\otimes}_{\mathrm{ex}}(\Tc_S,\LF(S'))$.

\subsubsection{Base change of fiber functors}\label{base_change_fiber_functors}

If $\om=\om_{X/\coeff}$ is a fiber functor over a $\coeff$-scheme $X$, it corresponds to a map $X \to \Fib{\Tc}$ of stacks over $\coeff$.

We can base-change to get a map $X \times_{\coeff} S \to \Fib{\Tc} \times_{\coeff} S = (\Fib{\Tc})_S$ of stacks over $S$, hence an $\Oc(S)$-linear pullback functor \[\om_{S_X / S} \colon \Tc_S \to \LF(S_X),\] where $S_X \coloneqq S \times_{\coeff} X$.

If $X$ starts out as an $S$-scheme, we can pull back along $X \to S_X$ to get an $X$-valued fiber functor $\om_{X/S}$ on $\Tc_S$. Alternatively, we can use the base-change adjunction $\Hom_\coeff(X,\Fib{\Tc}) \simeq \Hom_S(X,(\Fib{\Tc})_S)$, and $\om_{X/S}$ is the object on the right corresponding to $\om$ on the left.

Similarly, if $S$ starts out as a $X$-scheme over $\coeff$, we can compose $S \to X \to \Fib{\Tc}$ to get an $S$-valued fiber functor $\om_S=\om_{S/\coeff}$ on $\Tc$, then apply the previous sentence to get an $S$-valued fiber functor $\om_{S/S}$ on $\Tc_S$. Alternatively, the map $S \to X$ over $\coeff$ induces $S \to S_X$ over $S$, and pulling back $\om_{S_X / S}$ along this gives $\om_{S/S}$.

\subsubsection{}

If $S' \to (\Fib{\Tc})_S$ is fpqc, then by \ref{presentations_groupoids}, we get a groupoid in $S$-schemes $(S',R,s,t,c)$ and an identification $\Tc_S \equiv \operatorname{Rep}(S',R,s,t,c)$. The pullback $\Tc_S = \LF((\Fib{\Tc})_S) \to \LF(S')$ is a fiber functor on $\Tc_S$ corresponding to forgetting the groupoid action.

\subsubsection{}\label{appendix_reference_to_grfil_section_6} In \S \ref{base_change_fil_gr}, we discuss base-change of graded and filtered fiber functors.

\subsubsection{Tannakian reconstruction}

By definition of $\Fib{\Tc}$, for any $\coeff$-scheme $X$, we have
\[
\Fib{\Tc}(X) = \operatorname{Hom}^{\otimes}_{\mathrm{ex}}(\Tc,\LF(X)).
\]

The content in the theory of Tannaka duality is that this is always an affine gerbe when $\Tc$ is a Tannakian category, and that if $\Tc$ is the category of representations of an affine group (resp. groupoid) scheme, then $\Fib{\Tc}$ is representable by the corresponding classifying (resp. quotient) stack.

Since the above equality is of groupoids and not just sets, this says in particular that the automorphism group of a given fiber functor is representable (as $X$ varies) by a group scheme, and furthermore this group scheme is $G$ if $\Tc=\Rep(G)$.

Let us replace $\Tc$ by $\Tc_S$ for a $\coeff$-scheme $S$ and $X$ by an arbitrary prestack $\Xc$ of groupoids on $S$ schemes. We have a priori different groupoids \[(\operatorname{Fib}{\Tc})_S(\Xc) \coloneqq \Hom_{S}(\Xc,(\operatorname{Fib}{\Tc})_S) = \Hom_{\coeff}(\Xc,\operatorname{Fib}{\Tc}) \eqqcolon \operatorname{Fib}{\Tc}(\Xc)\] and 
\[
(\Fib{\Tc_S})(\Xc) \coloneqq \operatorname{Hom}^{\otimes}_{\mathrm{ex}}(\Tc_S,\LF(\Xc))
\]
as well as a natural functor
\[
(\operatorname{Fib}{\Tc})_S(\Xc) \xrightarrow{F_{\LF}} (\Fib{\Tc_S})(\Xc).
\]

\subsubsection{Remark on terminology}

In \cite[A.21,A.23]{wedhorn2024extensionliftinggbundlesstacks}, for $G$ a group scheme over $\coeff$ and $\Tc = \Rep(G)$ (i.e., if $\Tc$ is neutral), an object of $(\operatorname{Fib}{\Tc})_S(\Xc)$ is known as a \emph{$G$-bundle over $\Xc$}, while an object of $\operatorname{Hom}^{\otimes}_{\mathrm{ex}}(\Tc_S,\LF(\Xc))$ is known as a \emph{Tannakian $G$-bundle over $\Xc$}.

\begin{ssTheorem}\label{theorem:tannakian_reconstr_groupoids}
If $S$ is affine, then the functor $F_{\LF}$ is an equivalence of groupoids.
\end{ssTheorem}

\subsubsection{Proof}\label{adams_resolution_proof}

An affine groupoid scheme $(R,\circ,s,b)$ over $S$ acting on $U$ is the same as a \emph{Hopf algebroid} $(A,\Gamma)$, where $A=\Oc(U)$ and $\Gamma=\Oc(R)$, and a representation of $(R,\circ,s,b)$ is the same as a comodule over $(A,\Gamma)$.

We recall (\cite[Definition 3.1]{GoerssHopkins00},\cite[Definition 1.4.3]{HoveyHtpyComods04}) that a Hopf algebroid $(A,\Gamma)$ is an \emph{Adams Hopf algebroid} if $\Gamma$ is flat over $A$ and if $\Gamma$, viewed as a comodule over $(A,\Gamma)$, is a filtered system $\{\Gamma_{\alpha}\}$ of comodules for which $\Gamma_{\alpha}$ is finitely generated and projective over $A$. An fpqc-algebraic stack is an \emph{Adams stack} (\cite[Definition 6.5(2)]{goerss2008quasicoherentsheavesmodulistack},\cite[\S 1.3]{schaeppi2012characterizationcategoriescoherentsheaves}) if the associated Hopf algebroid for some fpqc presentation (as in \ref{presentations_groupoids}) is an Adams Hopf algebroid.

Choose a fiber functor $\om$ of $\Tc$ over an extension $\coeff'/\coeff$. As in \ref{sec:gerbes_groupoids}, we get a transitive $\coeff$-groupoid $(\Gc,\circ,s,b)$ acting on $U=\Spec{\coeff'}$, so in the associated Hopf algebroid, $A=\coeff'$ (while $\Gamma=\Oc(\Gc)$). By \cite[Corollaire 3.9]{DelTann90}, every representation of $\Gc$ is a filtered union of subrepresentations that are locally free of finite rank over $A$. In particular, this is true of $\Gamma$, so $(A,\Gamma)$ is an Adams Hopf algebroid, and $\Fib{\Tc}=[U/\Gc]$ is an Adams stack.

An affine scheme $S=\Spec{R}$ is obviously an Adams stack (it is presented by a Hopf algebroid with $A=\Gamma=R$, so $\Gamma$ is free of finite rank already). It follows by \cite[Theorem 1.8]{SchaeppiIndAbelian2014} that $(\Fib{\Tc})_S = \Fib{\Tc} \times_{\coeff} S$ is also an Adams stack.

The result then follows as in the proof of \cite[Proposition A.30]{wedhorn2024extensionliftinggbundlesstacks}.\footnote{That is, if $\Xc$ is an Adams stack, then $(\operatorname{Fib}{\Tc})_S(\Xc) \xto{F_{\QCoh}} \operatorname{Hom}^{\otimes}(\QCoh((\operatorname{Fib}{\Tc})_S),\QCoh(\Xc))$ induces an equivalence of $(\operatorname{Fib}{\Tc})_S(\Xc)$ with $\operatorname{Hom}^{\otimes}_{\mathrm{tame}}(\QCoh((\operatorname{Fib}{\Tc})_S),\QCoh(\Xc))$ by \cite[1.2.1]{schaeppi2012characterizationcategoriescoherentsheaves}, which is the same as $\operatorname{Hom}^{\otimes}_{\mathrm{cocont}}(\QCoh((\operatorname{Fib}{\Tc})_S),\QCoh(\Xc))$ by \cite[1.3.2]{schaeppi2012characterizationcategoriescoherentsheaves}. 
Applying \cite[Theorem 3.4.2]{SchaeppiConstructingColimits2020} to $\mathscr{A} = \Tc_S \coloneqq \LF(\operatorname{Fib}{\Tc})_S) \subseteq \QCoh(\operatorname{Fib}{\Tc})_S)$ and $\mathscr{D} = \QCoh(\Xc)$, this is equivalent to $\operatorname{Hom}^{\otimes}_{\m{r-ex}}(\Tc_S,\QCoh(\Xc))$ via restriction along $\LF(\operatorname{Fib}{\Tc})_S) \subseteq \QCoh(\operatorname{Fib}{\Tc})_S)$ by \cite[3.4.3]{SchaeppiConstructingColimits2020}. Since tensor functors respect duals (\cite[2.7]{DelTann90}), this is precisely $\operatorname{Hom}^{\otimes}_{\m{ex}}(\Tc_S,\LF(\Xc)) = \Fib{\Tc_S}(\Xc)$, so we are done in this case. The case of general $\Xc$ follows by the same colimit argument as in the second half of the proof of \cite[Proposition A.30]{wedhorn2024extensionliftinggbundlesstacks}.}

% here is a sketch of the proof: It seems that in fact, if BG_0 has the resolution property, S_0 is coherent (e.g., Noetherian), and S is affine, then BG also has the resolution property. By Remark 6.1.2 of https://arxiv.org/pdf/1206.2764, BG_0 has the strong resolution property and therefore the Adams property. An affine scheme obviously has the Adams property, so assuming S is affine, then by 1.8 of https://arxiv.org/pdf/1211.3678 (or argue concretely by tensoring the filtered colimit of dualizable comodules with O(S) over O(S_0)), so does BG = BG_0 \times_{S_0} S. Thus BG satisfies (strong) resolution. As a consequence, your argument for A.30 applies.

%If we assume that $\Tc$ is neutral, then by \cite[Proposition A.30]{wedhorn}, the map $\operatorname{Fib}{\Tc}(\Xc) \to \operatorname{Hom}^{\otimes}_{\mathrm{ex}}(\Tc,\LF(\Xc))$ is an equivalence
%***** Wait, A.28 needs finite type. Is there anywhere that gives what I want in my case? Maybe the thing in stacks project? Maybe Remark 2.20 of https://arxiv.org/pdf/2401.02279?
%** Do we need neutral?

%% https://d-nb.info/1007550295/34 has a nice survey of resolution property
%% https://arxiv.org/pdf/1410.1716 has some good stuff about tensor categories as a foundation for algebraic geometry

\subsubsection{}\label{alg_geo_over_S}

We can mimic \ref{algebraic_geometry_tann_cat} to define affine schemes in $\Tc_S$. The category $\Ind{\Tc_S}$ inherits a symmetric monoidal structure, and a commutative algebra in $\Ind{\Tc_S}$ is just a commutative monoid object. The category $\operatorname{Aff}(\Tc_S)$ of affine schemes in $\Tc_S$ is defined to be the opposite of the category of commutative algebras in $\Ind{\Tc}$; the corresponding object in $\operatorname{Aff}(\Tc_S)$ is called the \emph{spectrum} of the given algebra in $\Ind{\Tc}$.

A group scheme in $\Tc_S$ is just a group object in $\operatorname{Aff}(\Tc_S)$; equivalently, it is the spectrum of a Hopf algebra object in $\Ind{\Tc_S}$.

Following the notation of \ref{arbitrary_S}, if $X$ is an affine scheme in $\Tc$, we write $X_S$ for the affine scheme in $\Tc_S$ given by $\Spec{\Oc(X)_S}$. If $\pi$ is a group scheme in $\Tc$, we similarly write $\pi_S$ for the corresponding group scheme in $\Tc_S$.

\subsubsection{}\label{alg_geo_stacks}

Let $Z \in \operatorname{Aff}(\Tc_S)$ and $\Ac$ its coordinate algebra in $\Ind{\Tc_S}$. By taking its inductive limit, $\Ac$ gives rise to an algebra in $\QCoh(\Fib(\Tc)_S)$. We can take relative Spec (\cite[10.2.1]{OlssonStacksBook16}) to get a stack $ \underline{\Spec}_{\Fib(\Tc)_S} (\Ac)$ with a morphism $ \underline{\Spec}_{\Fib(\Tc)_S} (\Ac) \to \Fib(\Tc)_S$ that is affine by \cite[Theorem 10.2.4]{OlssonStacksBook16}.

%If $S=\Spec{R}$ is affine, then $\Fib(\Tc)_S$ has the resolution property, so $\QCoh(\Fib(\Tc)_S) = \Ind{\Tc_S}$

Thus $\operatorname{Aff}(\Tc_S)$ is equivalent to the full subcategory of the category of affine morphisms of stacks $\pi \colon \mathscr{Y} \to \Fib(\Tc)_S$ for which $\pi$ is locally free. If $S$ is the spectrum of a field, then this latter condition is always satisfied.

%% **** I might want to explain why Ind(Vect) embeds into QCoh. Related to the inclusion after Example 1.7 of https://arxiv.org/pdf/2309.01260

\subsubsection{}\label{cohomology_in_T} If $\pi$ is a group scheme in $\Tc$, then for a $\coeff$-scheme $S$, we denote by $H^1(\Tc;\pi)(S)$ the pointed set of isomorphism classes of right $\pi_S$-torsors in $\Tc_S$; that is, a (right) $\pi_S$-torsor is a nonempty affine scheme $P$ in $\Tc_S$ with a (right) group action map $P \times \pi_S \to P$ in schemes over $\Tc_S$ for which $P \times \pi_S \to P \times P$ is an isomorphism.

If $\om_S$ is a fiber functor over $S$, then we have a bijection
\[
H^1(\Tc;\pi)(S) \simeq H^1(\pi_1(\Tc,\om_S);\om_S(\pi)),
\]
where the latter is defined by algebraic cocycles modulo algebraic coboundaries.

If $M$ is an object of $\Tc$ and $i \in \Zb_{\ge 0}$, we denote by $H^i(\Tc;M)(S)$ the group $\Ext^i_{\Tc_S}(\mathbbm{1},M_S)$. If $\om_S$ is a fiber functor over $S$, then we have a bijection
\[
H^i(\Tc;M)(S) \simeq H^i(\pi_1(\Tc,\om_S);\om_S(M)).
\]

If $R$ is a $\coeff$-algebra, then $H^1(\Tc;\pi)(R)$ (resp. $H^i(\Tc;M)(R)$) denotes the above for $S=\Spec{R}$.

%%%%%%%%%%%%%%%%

%%%%%%%%%%%%%%%%%%%%
%%%%%%%%%%%%%%%%%%%
\bibliography{pPerSel_Refs}%%%%%%

@MISC {SawinMO_Pure_motives_compatible_systems,
    TITLE = {Pure motives and compatible systems of $\ell$-adic representations},
    AUTHOR = {Will Sawin (https://mathoverflow.net/users/18060/will-sawin)},
    HOWPUBLISHED = {MathOverflow},
    NOTE = {URL:https://mathoverflow.net/q/212374 (version: 2015-07-26)},
    EPRINT = {https://mathoverflow.net/q/212374},
    URL = {https://mathoverflow.net/q/212374}
}

@article {UnverCyclotomic,
    AUTHOR = {\"{U}nver, Sinan},
     TITLE = {Cyclotomic {$p$}-adic multi-zeta values},
   JOURNAL = {J. Pure Appl. Algebra},
  FJOURNAL = {Journal of Pure and Applied Algebra},
    VOLUME = {223},
      YEAR = {2019},
    NUMBER = {2},
     PAGES = {489--503},
      ISSN = {0022-4049},
   MRCLASS = {11M32 (14F30)},
  MRNUMBER = {3850552},
MRREVIEWER = {Bin Zhang},
       DOI = {10.1016/j.jpaa.2018.04.002},
       URL = {https://doi.org/10.1016/j.jpaa.2018.04.002},
}

@article {AnFr2025algebraic,
    AUTHOR = {Ancona, Giuseppe and Fr\u{a}\c{t}il\u{a}, Drago\c{s}},
     TITLE = {Algebraic classes in mixed characteristic and {A}ndr\'{e}'s
              {$p$}-adic periods},
   JOURNAL = {J. Inst. Math. Jussieu},
  FJOURNAL = {Journal of the Institute of Mathematics of Jussieu. JIMJ.
              Journal de l'Institut de Math\'{e}matiques de Jussieu},
    VOLUME = {24},
      YEAR = {2025},
    NUMBER = {4},
     PAGES = {1093--1138},
      ISSN = {1474-7480,1475-3030},
   MRCLASS = {14G45 (11G35 14F30)},
  MRNUMBER = {4922232},
       DOI = {10.1017/S147474802400029X},
       URL = {https://doi.org/10.1017/S147474802400029X},
}

@incollection {AndreGaloisMotTrans09,
    AUTHOR = {Andr\'{e}, Yves},
     TITLE = {Galois theory, motives and transcendental numbers},
 BOOKTITLE = {Renormalization and {G}alois theories},
    SERIES = {IRMA Lect. Math. Theor. Phys.},
    VOLUME = {15},
     PAGES = {165--177},
 PUBLISHER = {Eur. Math. Soc., Z\"{u}rich},
      YEAR = {2009},
      ISBN = {978-3-03719-073-9},
   MRCLASS = {11J81 (12F10 14C15 14F40)},
  MRNUMBER = {2588609},
MRREVIEWER = {Tam\'{a}s\ Szamuely},
       DOI = {10.4171/073-1/4},
       URL = {https://doi.org/10.4171/073-1/4},
}

@misc{stacks-project,
    shorthand    = {Stacks},
    author       = {The {Stacks Project Authors}},
    title        = {\textit{Stacks Project}},
    howpublished = {\url{https://stacks.math.columbia.edu}},
    year         = {2018},
  }

@article {AyoubSixOpII,
    AUTHOR = {Ayoub, Joseph},
     TITLE = {Les six op\'{e}rations de {G}rothendieck et le formalisme des
              cycles \'{e}vanescents dans le monde motivique. {II}},
   JOURNAL = {Ast\'{e}risque},
  FJOURNAL = {Ast\'{e}risque},
    NUMBER = {315},
      YEAR = {2007},
     PAGES = {vi+364 pp. (2008)},
      ISSN = {0303-1179},
      ISBN = {978-2-85629-245-7},
   MRCLASS = {14C25 (14F20 14F42 18A40 19E15)},
  MRNUMBER = {2438151},
MRREVIEWER = {Christian Haesemeyer},
}

@article {AyoubHopf14,
    AUTHOR = {Ayoub, Joseph},
     TITLE = {L'alg\`ebre de {H}opf et le groupe de {G}alois motiviques d'un
              corps de caract\'{e}ristique nulle, {II}},
   JOURNAL = {J. Reine Angew. Math.},
  FJOURNAL = {Journal f\"{u}r die Reine und Angewandte Mathematik. [Crelle's
              Journal]},
    VOLUME = {693},
      YEAR = {2014},
     PAGES = {151--226},
      ISSN = {0075-4102},
   MRCLASS = {14F42 (11R32 18D10 18E30)},
  MRNUMBER = {3259032},
MRREVIEWER = {Florence Lecomte},
       DOI = {10.1515/crelle-2012-0090},
       URL = {https://doi.org/10.1515/crelle-2012-0090},
}

@unpublished{AyoubAnabelianPresentationPreprint,
    author = {Ayoub, Joseph},
    title = {Anabelian Presentation of the Motivic {G}alois Group in Characteristic Zero},
    year = {2024},
    note = {Preprint available from https://user.math.uzh.ch/ayoub/PDF-Files/Anabel.pdf}
}

@article {nabsd,
    AUTHOR = {Balakrishnan, Jennifer S. and Dan-Cohen, Ishai and Kim,
              Minhyong and Wewers, Stefan},
     TITLE = {A non-abelian conjecture of {T}ate-{S}hafarevich type for
              hyperbolic curves},
   JOURNAL = {Math. Ann.},
  FJOURNAL = {Mathematische Annalen},
    VOLUME = {372},
      YEAR = {2018},
    NUMBER = {1-2},
     PAGES = {369--428},
      ISSN = {0025-5831},
   MRCLASS = {11D45 (11G50 14F35 14H52)},
  MRNUMBER = {3856816},
MRREVIEWER = {Christopher Frei},
       DOI = {10.1007/s00208-018-1684-x},
       URL = {https://doi.org/10.1007/s00208-018-1684-x},
}

@article {BalakrishnanDograI,
    AUTHOR = {Balakrishnan, Jennifer S. and Dogra, Netan},
     TITLE = {Quadratic {C}habauty and rational points, {I}: {$p$}-adic
              heights},
      NOTE = {With an appendix by J. Steffen M\"{u}ller},
   JOURNAL = {Duke Math. J.},
  FJOURNAL = {Duke Mathematical Journal},
    VOLUME = {167},
      YEAR = {2018},
    NUMBER = {11},
     PAGES = {1981--2038},
      ISSN = {0012-7094},
   MRCLASS = {14G05 (11G50 14G40)},
  MRNUMBER = {3843370},
MRREVIEWER = {Ariyan Javanpeykar},
       DOI = {10.1215/00127094-2018-0013},
       URL = {https://doi.org/10.1215/00127094-2018-0013},
}

@article {CursedCurve,
    AUTHOR = {Balakrishnan, Jennifer and Dogra, Netan and M\"{u}ller, J. Steffen
              and Tuitman, Jan and Vonk, Jan},
     TITLE = {Explicit {C}habauty-{K}im for the split {C}artan modular curve
              of level 13},
   JOURNAL = {Ann. of Math. (2)},
  FJOURNAL = {Annals of Mathematics. Second Series},
    VOLUME = {189},
      YEAR = {2019},
    NUMBER = {3},
     PAGES = {885--944},
      ISSN = {0003-486X},
   MRCLASS = {14G05 (11G18 11G50 11Y50)},
  MRNUMBER = {3961086},
MRREVIEWER = {Christopher Frei},
       DOI = {10.4007/annals.2019.189.3.6},
       URL = {https://doi.org/10.4007/annals.2019.189.3.6},
}

@article {BalakrishnanDograII,
    AUTHOR = {Balakrishnan, Jennifer S. and Dogra, Netan},
     TITLE = {Quadratic {C}habauty and rational points {II}: {G}eneralised
              height functions on {S}elmer varieties},
   JOURNAL = {Int. Math. Res. Not. IMRN},
  FJOURNAL = {International Mathematics Research Notices. IMRN},
      YEAR = {2021},
    NUMBER = {15},
     PAGES = {11923--12008},
      ISSN = {1073-7928},
   MRCLASS = {11G50 (14G40)},
  MRNUMBER = {4294137},
MRREVIEWER = {Christopher Frei},
       DOI = {10.1093/imrn/rnz362},
       URL = {https://doi.org/10.1093/imrn/rnz362},
}

@article {BertolinMotivicGalois08,
    AUTHOR = {Bertolin, Cristiana},
     TITLE = {Motivic {G}alois theory for 1-motives},
   JOURNAL = {Ann. Sci. Math. Qu\'{e}bec},
  FJOURNAL = {Annales des Sciences Math\'{e}matiques du Qu\'{e}bec},
    VOLUME = {32},
      YEAR = {2008},
    NUMBER = {2},
     PAGES = {105--124},
      ISSN = {0707-9109},
   MRCLASS = {14L15 (14C15)},
  MRNUMBER = {2562038},
MRREVIEWER = {Claudio Pedrini},
}

@article {BerDarRot15,
    AUTHOR = {Bertolini, Massimo and Darmon, Henri and Rotger, Victor},
     TITLE = {Beilinson-{F}lach elements and {E}uler systems {II}: the
              {B}irch-{S}winnerton-{D}yer conjecture for
              {H}asse-{W}eil-{A}rtin {$L$}-series},
   JOURNAL = {J. Algebraic Geom.},
  FJOURNAL = {Journal of Algebraic Geometry},
    VOLUME = {24},
      YEAR = {2015},
    NUMBER = {3},
     PAGES = {569--604},
      ISSN = {1056-3911},
   MRCLASS = {11R42 (11G40)},
  MRNUMBER = {3344765},
MRREVIEWER = {Jeanine Van Order},
       DOI = {10.1090/S1056-3911-2015-00675-0},
       URL = {https://doi.org/10.1090/S1056-3911-2015-00675-0},
}

@inproceedings {BesserRigSyn00I,
    AUTHOR = {Besser, Amnon},
     TITLE = {Syntomic regulators and {$p$}-adic integration. {I}. {R}igid
              syntomic regulators},
 BOOKTITLE = {Proceedings of the {C}onference on {$p$}-adic {A}spects of the
              {T}heory of {A}utomorphic {R}epresentations ({J}erusalem,
              1998)},
   JOURNAL = {Israel J. Math.},
  FJOURNAL = {Israel Journal of Mathematics},
    VOLUME = {120},
      YEAR = {2000},
    NUMBER = {part B},
     PAGES = {291--334},
      ISSN = {0021-2172},
   MRCLASS = {14F43 (11G25 14F30 14G22 19F27)},
  MRNUMBER = {1809626},
MRREVIEWER = {Takao Yamazaki},
       DOI = {10.1007/BF02834843},
       URL = {https://doi.org/10.1007/BF02834843},
}

@inproceedings {BesserRigSyn00II,
    AUTHOR = {Besser, Amnon},
     TITLE = {Syntomic regulators and {$p$}-adic integration. {II}. {$K_2$}
              of curves},
 BOOKTITLE = {Proceedings of the {C}onference on {$p$}-adic {A}spects of the
              {T}heory of {A}utomorphic {R}epresentations ({J}erusalem,
              1998)},
   JOURNAL = {Israel J. Math.},
  FJOURNAL = {Israel Journal of Mathematics},
    VOLUME = {120},
      YEAR = {2000},
    NUMBER = {part B},
     PAGES = {335--359},
      ISSN = {0021-2172},
   MRCLASS = {14F43 (11G20 14C35 14G20 19F27)},
  MRNUMBER = {1809627},
MRREVIEWER = {Takao Yamazaki},
       DOI = {10.1007/BF02834844},
       URL = {https://doi.org/10.1007/BF02834844},
}

@article{BesserColeman,
	Author = {Besser, Amnon},
	Coden = {MAANA},
	Doi = {10.1007/s002080100263},
	Fjournal = {Mathematische Annalen},
	Issn = {0025-5831},
	Journal = {Math. Ann.},
	Mrclass = {11S80 (11G25 14F30)},
	Mrnumber = {1883387 (2003d:11176)},
	Mrreviewer = {Bruno Chiarellotto},
	Number = {1},
	Pages = {19--48},
	Title = {Coleman integration using the {T}annakian formalism},
	Url = {http://dx.doi.org/10.1007/s002080100263},
	Volume = {322},
	Year = {2002}}

@incollection{BesserFurusho,
	author = {Besser, Amnon and Furusho, Hidekazu},
	booktitle = {Primes and knots},
	doi = {10.1090/conm/416/07884},
	mrclass = {11G55 (11M41 11S80)},
	mrnumber = {2276133},
	mrreviewer = {Jan Nekov\'{a}\v{r}},
	pages = {9--29},
	publisher = {Amer. Math. Soc., Providence, RI},
	series = {Contemp. Math.},
	title = {The double shuffle relations for {$p$}-adic multiple zeta values},
	url = {https://doi.org/10.1090/conm/416/07884},
	volume = {416},
	year = {2006}}

@incollection{BesserHeidelberg,
	author = {Besser, Amnon},
	booktitle = {The arithmetic of fundamental groups---{PIA} 2010},
	doi = {10.1007/978-3-642-23905-2_1},
	mrclass = {11S80 (11G25 14F30 14G22)},
	mrnumber = {3220512},
	mrreviewer = {Alessandra Bertapelle},
	pages = {3--52},
	publisher = {Springer, Heidelberg},
	series = {Contrib. Math. Comput. Sci.},
	title = {Heidelberg lectures on {C}oleman integration},
	url = {http://dx.doi.org/10.1007/978-3-642-23905-2_1},
	volume = {2},
	year = {2012}}

@incollection {BesserArbSyn21,
    AUTHOR = {Besser, Amnon},
     TITLE = {The syntomic regulator for {$K_2$} of curves with arbitrary
              reduction},
 BOOKTITLE = {Arithmetic {L}-functions and differential geometric methods},
    SERIES = {Progr. Math.},
    VOLUME = {338},
     PAGES = {75--89},
 PUBLISHER = {Birkh\"{a}user/Springer, Cham},
      YEAR = {[2021] \copyright 2021},
   MRCLASS = {19F27 (11G20)},
  MRNUMBER = {4311239},
MRREVIEWER = {Evangelia Gazaki},
       DOI = {10.1007/978-3-030-65203-6\_3},
       URL = {https://doi.org/10.1007/978-3-030-65203-6_3},
}

@article {BettsEtAl,
    AUTHOR = {Best, Alex J. and Betts, L. Alexander and Kumpitsch, Theresa
              and L\"{u}dtke, Martin and McAndrew, Angus W. and Qian, Lie and
              Studnia, Elie and Xu, Yujie},
     TITLE = {Refined {S}elmer equations for the thrice-punctured line in
              depth two},
   JOURNAL = {Math. Comp.},
  FJOURNAL = {Mathematics of Computation},
    VOLUME = {93},
      YEAR = {2024},
    NUMBER = {347},
     PAGES = {1497--1527},
      ISSN = {0025-5718},
   MRCLASS = {14G05 (11G55 11Y50)},
  MRNUMBER = {4709209},
       DOI = {10.1090/mcom/3898},
       URL = {https://doi.org/10.1090/mcom/3898},
}

@article {BettsWeightFil2023,
    AUTHOR = {Betts, L. Alexander},
     TITLE = {Weight filtrations on {S}elmer schemes and the effective
              {C}habauty-{K}im method},
   JOURNAL = {Compos. Math.},
  FJOURNAL = {Compositio Mathematica},
    VOLUME = {159},
      YEAR = {2023},
    NUMBER = {7},
     PAGES = {1531--1605},
      ISSN = {0010-437X},
   MRCLASS = {11G35 (14F35 14G05 14G20)},
  MRNUMBER = {4604872},
MRREVIEWER = {Francesc Bars},
       DOI = {10.1112/S0010437X2300725X},
       URL = {https://doi.org/10.1112/S0010437X2300725X},
}

@article {BettsCorwinLeonhardt,
    AUTHOR = {Betts, L. Alexander and Corwin, David and Leonhardt, Marius},
     TITLE = {Bounds on the {C}habauty-{K}im locus of hyperbolic curves},
   JOURNAL = {Int. Math. Res. Not. IMRN},
  FJOURNAL = {International Mathematics Research Notices. IMRN},
      YEAR = {2024},
    NUMBER = {12},
     PAGES = {9705--9727},
      ISSN = {1073-7928},
   MRCLASS = {11G30 (11G10 14G05 14H30)},
  MRNUMBER = {4761774},
       DOI = {10.1093/imrn/rnae067},
       URL = {https://doi.org/10.1093/imrn/rnae067},
}

@incollection{BlochKato,
	author = {Bloch, Spencer and Kato, Kazuya},
	booktitle = {The {G}rothendieck {F}estschrift, {V}ol.\ {I}},
	mrclass = {11G40 (11G09 14C35 14F30 14G10)},
	mrnumber = {1086888 (92g:11063)},
	mrreviewer = {Ehud de Shalit},
	pages = {333--400},
	publisher = {Birkh\"auser Boston, Boston, MA},
	series = {Progr. Math.},
	title = {{$L$}-functions and {T}amagawa numbers of motives},
	volume = {86},
	year = {1990}}

@article{BrownMTMZ,
	author = {Brown, Francis},
	coden = {ANMAAH},
	doi = {10.4007/annals.2012.175.2.10},
	fjournal = {Annals of Mathematics. Second Series},
	issn = {0003-486X},
	journal = {Ann. of Math. (2)},
	mrclass = {19Exx (11M32 14Fxx)},
	mrnumber = {2993755},
	number = {2},
	pages = {949--976},
	title = {Mixed {T}ate motives over {$\mathbb Z$}},
	url = {http://dx.doi.org/10.4007/annals.2012.175.2.10},
	volume = {175},
	year = {2012}}

@incollection{BrownDecomp,
	Author = {Brown, Francis C. S.},
	Booktitle = {Galois-{T}eichm{\"u}ller theory and arithmetic geometry},
	Mrclass = {11M32 (13B05 16T15)},
	Mrnumber = {3051238},
	Pages = {31--58},
	Publisher = {Math. Soc. Japan, Tokyo},
	Series = {Adv. Stud. Pure Math.},
	Title = {On the decomposition of motivic multiple zeta values},
	Volume = {63},
	Year = {2012}}

@article {BrownNotes17,
    AUTHOR = {Brown, Francis},
     TITLE = {Notes on motivic periods},
   JOURNAL = {Commun. Number Theory Phys.},
  FJOURNAL = {Communications in Number Theory and Physics},
    VOLUME = {11},
      YEAR = {2017},
    NUMBER = {3},
     PAGES = {557--655},
      ISSN = {1931-4523},
   MRCLASS = {81Q30 (11R32)},
  MRNUMBER = {3713352},
       DOI = {10.4310/CNTP.2017.v11.n3.a2},
       URL = {https://doi.org/10.4310/CNTP.2017.v11.n3.a2},
}

@misc{BrownUnit,
      title={Integral points on curves, the unit equation, and motivic periods}, 
      author={Francis Brown},
      year={2017},
      eprint={1704.00555},
    note = {arXiv \href{https://arxiv.org/abs/1704.00555}{1704.00555}},
      archivePrefix={arXiv},
      primaryClass={math.NT},
      url={https://arxiv.org/abs/1704.00555}, 
}

@article {BuehlerExactCats10,
    AUTHOR = {B{\"u}hler, Theo},
     TITLE = {Exact categories},
   JOURNAL = {Expo. Math.},
  FJOURNAL = {Expositiones Mathematicae},
    VOLUME = {28},
      YEAR = {2010},
    NUMBER = {1},
     PAGES = {1--69},
      ISSN = {0723-0869},
   MRCLASS = {18E10 (18-02 18E30)},
  MRNUMBER = {2606234},
MRREVIEWER = {Sunil K. Chebolu},
       DOI = {10.1016/j.exmath.2009.04.004},
       URL = {https://doi.org/10.1016/j.exmath.2009.04.004},
}

@article {Cao18,
    AUTHOR = {Cao, Jin},
     TITLE = {Motives for an elliptic curve},
   JOURNAL = {Math. Ann.},
  FJOURNAL = {Mathematische Annalen},
    VOLUME = {372},
      YEAR = {2018},
    NUMBER = {1-2},
     PAGES = {189--227},
      ISSN = {0025-5831},
   MRCLASS = {14F42 (14C15 14H52 18E30)},
  MRNUMBER = {3856811},
MRREVIEWER = {Jens Hornbostel},
       DOI = {10.1007/s00208-018-1690-z},
       URL = {https://doi.org/10.1007/s00208-018-1690-z},
}

@article {ChatUnv13,
    AUTHOR = {Chatzistamatiou, Andre and \"{U}nver, Sinan},
     TITLE = {On {$p$}-adic periods for mixed {T}ate motives over a number
              field},
   JOURNAL = {Math. Res. Lett.},
  FJOURNAL = {Mathematical Research Letters},
    VOLUME = {20},
      YEAR = {2013},
    NUMBER = {5},
     PAGES = {825--844},
      ISSN = {1073-2780},
   MRCLASS = {11S31 (11F85)},
  MRNUMBER = {3207355},
MRREVIEWER = {Jamshid Derakhshan},
       DOI = {10.4310/MRL.2013.v20.n5.a2},
       URL = {https://doi.org/10.4310/MRL.2013.v20.n5.a2},
}

@article {ChiarBruLaz19,
    AUTHOR = {Chiarellotto, Bruno and Lazda, Christopher and Mazzari,
              Nicola},
     TITLE = {The filtered {O}gus realisation of motives},
   JOURNAL = {J. Algebra},
  FJOURNAL = {Journal of Algebra},
    VOLUME = {527},
      YEAR = {2019},
     PAGES = {348--365},
      ISSN = {0021-8693},
   MRCLASS = {14F42 (11F75 11G09 14F30 14J28)},
  MRNUMBER = {3924438},
MRREVIEWER = {Damian R\"{o}ssler},
       DOI = {10.1016/j.jalgebra.2019.03.006},
       URL = {https://doi.org/10.1016/j.jalgebra.2019.03.006},
}

@book {CisDegTriCat,
    AUTHOR = {Cisinski, Denis-Charles and D\'{e}glise, Fr\'{e}d\'{e}ric},
     TITLE = {Triangulated categories of mixed motives},
    SERIES = {Springer Monographs in Mathematics},
 PUBLISHER = {Springer, Cham},
      YEAR = {[2019] \copyright 2019},
     PAGES = {xlii+406},
      ISBN = {978-3-030-33241-9; 978-3-030-33242-6},
   MRCLASS = {14F42 (14C15 14C35 18G80 19D55)},
  MRNUMBER = {3971240},
MRREVIEWER = {Igor A. Rapinchuk},
       DOI = {10.1007/978-3-030-33242-6},
       URL = {https://doi.org/10.1007/978-3-030-33242-6},
}

@article {PolGonI,
    AUTHOR = {Corwin, David and Dan-Cohen, Ishai},
     TITLE = {The polylog quotient and the {G}oncharov quotient in
              computational {C}habauty--{K}im {T}heory {I}},
   JOURNAL = {Int. J. Number Theory},
  FJOURNAL = {International Journal of Number Theory},
    VOLUME = {16},
      YEAR = {2020},
    NUMBER = {8},
     PAGES = {1859--1905},
      ISSN = {1793-0421},
   MRCLASS = {11G55 (11S80 14F42)},
  MRNUMBER = {4143688},
       DOI = {10.1142/S1793042120500967},
       URL = {https://doi.org/10.1142/S1793042120500967},
}

@misc{CorwinMECK,
      title={Explicit Motivic Mixed Elliptic {C}habauty-{K}im}, 
      author={David Corwin},
      year={2021},
      eprint={2102.08371},
      archivePrefix={arXiv},
      primaryClass={math.NT},
      url={https://arxiv.org/abs/2102.08371},
    note = {arXiv \href{https://arxiv.org/abs/2102.08371}{2102.08371}},
}

@misc{CorwinTSV,
  author = {Corwin, David},
  title = {{{T}annakian {S}elmer Varieties}},
  howpublished = "\url{https://www.math.bgu.ac.il/~corwind/files/research/TSV.pdf}",
  year = {2022}, 
}

@unpublished{CLinitialreport2026,
author = {Corwin, David and L\"udtke, Martin},
title = {{Computations in Mixed Elliptic Motivic Chabauty--Kim: Initial Report}},
  howpublished = "\url{https://martinluedtke.github.io/files/corwin_luedtke_initial_report_2026.pdf}",
note ={\url{https://martinluedtke.github.io/files/corwin_luedtke_initial_report_2026.pdf}},
year = {2026},
}

@article {DegNiz18,
    AUTHOR = {D\'{e}glise, Fr\'{e}d\'{e}ric and Nizio\polishL{}, Wies\polishL{}awa},
     TITLE = {On {$p$}-adic absolute {H}odge cohomology and syntomic
              coefficients. {I}},
   JOURNAL = {Comment. Math. Helv.},
  FJOURNAL = {Commentarii Mathematici Helvetici. A Journal of the Swiss
              Mathematical Society},
    VOLUME = {93},
      YEAR = {2018},
    NUMBER = {1},
     PAGES = {71--131},
      ISSN = {0010-2571},
   MRCLASS = {14F30 (11G25 14C30 14F42 14G20)},
  MRNUMBER = {3777126},
MRREVIEWER = {Adolfo Quir\'{o}s},
       DOI = {10.4171/CMH/430},
       URL = {https://doi.org/10.4171/CMH/430},
}

@article {CKTwo,
    AUTHOR = {Dan-Cohen, Ishai and Wewers, Stefan},
     TITLE = {Explicit {C}habauty-{K}im theory for the thrice punctured line
              in depth 2},
   JOURNAL = {Proc. Lond. Math. Soc. (3)},
  FJOURNAL = {Proceedings of the London Mathematical Society. Third Series},
    VOLUME = {110},
      YEAR = {2015},
    NUMBER = {1},
     PAGES = {133--171},
      ISSN = {0024-6115,1460-244X},
   MRCLASS = {11G55 (11D45 14F42)},
  MRNUMBER = {3299602},
MRREVIEWER = {Piotr\ Kraso\'{n}},
       DOI = {10.1112/plms/pdu034},
       URL = {https://doi.org/10.1112/plms/pdu034},
}

@article {MTMUE,
    AUTHOR = {Dan-Cohen, Ishai and Wewers, Stefan},
     TITLE = {Mixed {T}ate motives and the unit equation},
   JOURNAL = {Int. Math. Res. Not. IMRN},
  FJOURNAL = {International Mathematics Research Notices. IMRN},
      YEAR = {2016},
    NUMBER = {17},
     PAGES = {5291--5354},
      ISSN = {1073-7928},
   MRCLASS = {14F42 (11D45 11G55 14H30)},
  MRNUMBER = {3556439},
MRREVIEWER = {Piotr Kraso\'n},
       DOI = {10.1093/imrn/rnv239},
       URL = {http://dx.doi.org/10.1093/imrn/rnv239},
}

@article {MTMUEII,
    AUTHOR = {Dan-Cohen, Ishai},
     TITLE = {Mixed {T}ate motives and the unit equation {II}},
   JOURNAL = {Algebra Number Theory},
  FJOURNAL = {Algebra \& Number Theory},
    VOLUME = {14},
      YEAR = {2020},
    NUMBER = {5},
     PAGES = {1175--1237},
      ISSN = {1937-0652},
   MRCLASS = {14F42 (11D45 11G55 14F30 14F35 14G05 14Q05)},
  MRNUMBER = {4129385},
       DOI = {10.2140/ant.2020.14.1175},
       URL = {https://doi.org/10.2140/ant.2020.14.1175},
}

@article {PolGonII,
    AUTHOR = {Dan-Cohen, Ishai and Corwin, David},
     TITLE = {The polylog quotient and the {G}oncharov quotient in
              computational {C}habauty--{K}im theory {II}},
   JOURNAL = {Trans. Amer. Math. Soc.},
  FJOURNAL = {Transactions of the American Mathematical Society},
    VOLUME = {373},
      YEAR = {2020},
    NUMBER = {10},
     PAGES = {6835--6861},
      ISSN = {0002-9947},
   MRCLASS = {14G05 (11G55 14F42)},
  MRNUMBER = {4155193},
       DOI = {10.1090/tran/7964},
       URL = {https://doi.org/10.1090/tran/7964},
}

@article {DCSRationalMotivicPathSpaces22,
    AUTHOR = {Dan-Cohen, Ishai and Schlank, Tomer},
     TITLE = {Rational motivic path spaces and {K}im's relative unipotent
              section conjecture},
   JOURNAL = {Rend. Semin. Mat. Univ. Padova},
  FJOURNAL = {Rendiconti del Seminario Matematico della Universit\`a di
              Padova},
    VOLUME = {148},
      YEAR = {2022},
     PAGES = {117--172},
      ISSN = {0041-8994,2240-2926},
   MRCLASS = {14G05 (14C15 14F42 55P62)},
  MRNUMBER = {4542375},
MRREVIEWER = {Shusuke\ Otabe},
       DOI = {10.4171/rsmup/97},
       URL = {https://doi.org/10.4171/rsmup/97},
}

@misc{dancohen2025andreperiodsmixedtate,
      title={On {A}ndr\'e periods of mixed {T}ate motives}, 
      author={Ishai Dan-Cohen},
      year={2025},
      eprint={2502.17404},
      archivePrefix={arXiv},
      primaryClass={math.AG},
      url={https://arxiv.org/abs/2502.17404}, 
    note = {arXiv \href{https://arxiv.org/abs/2502.17404}{2502.17404}},
}

@article {DarmonRotger17,
    AUTHOR = {Darmon, Henri and Rotger, Victor},
     TITLE = {Diagonal cycles and {E}uler systems {II}: {T}he {B}irch and
              {S}winnerton-{D}yer conjecture for {H}asse-{W}eil-{A}rtin
              {$L$}-functions},
   JOURNAL = {J. Amer. Math. Soc.},
  FJOURNAL = {Journal of the American Mathematical Society},
    VOLUME = {30},
      YEAR = {2017},
    NUMBER = {3},
     PAGES = {601--672},
      ISSN = {0894-0347},
   MRCLASS = {11G05 (11G40)},
  MRNUMBER = {3630084},
MRREVIEWER = {Rolf Berndt},
       DOI = {10.1090/jams/861},
       URL = {https://doi.org/10.1090/jams/861},
}

@article {DeligneHodgeII,
    AUTHOR = {Deligne, Pierre},
     TITLE = {Th\'{e}orie de {H}odge. {II}},
   JOURNAL = {Inst. Hautes \'{E}tudes Sci. Publ. Math.},
  FJOURNAL = {Institut des Hautes \'{E}tudes Scientifiques. Publications
              Math\'{e}matiques},
    NUMBER = {40},
      YEAR = {1971},
     PAGES = {5--57},
      ISSN = {0073-8301},
   MRCLASS = {14C30 (14F15)},
  MRNUMBER = {498551},
MRREVIEWER = {J. H. M. Steenbrink},
       URL = {http://www.numdam.org/item?id=PMIHES_1971__40__5_0},
}

@article {WeilII,
    AUTHOR = {Deligne, Pierre},
     TITLE = {La conjecture de {W}eil. {II}},
   JOURNAL = {Inst. Hautes \'{E}tudes Sci. Publ. Math.},
  FJOURNAL = {Institut des Hautes \'{E}tudes Scientifiques. Publications
              Math\'{e}matiques},
    NUMBER = {52},
      YEAR = {1980},
     PAGES = {137--252},
      ISSN = {0073-8301},
   MRCLASS = {14G13 (10H10)},
  MRNUMBER = {601520},
MRREVIEWER = {Spencer J. Bloch},
       URL = {http://www.numdam.org/item?id=PMIHES_1980__52__137_0},
}

@incollection{Deligne89,
	address = {New York},
	author = {Deligne, Pierre},
	booktitle = {Galois groups over {Q} ({B}erkeley, {CA}, 1987)},
	mrclass = {14G25 (11G35 11M06 11R70 14F35 19E99 19F27)},
	mrnumber = {1012168 (90m:14016)},
	mrreviewer = {James Milne},
	pages = {79--297},
	publisher = {Springer},
	series = {Math. Sci. Res. Inst. Publ.},
	title = {Le groupe fondamental de la droite projective moins trois points},
	volume = {16},
	year = {1989}}

@incollection {DelTann90,
    AUTHOR = {Deligne, P.},
     TITLE = {Cat\'{e}gories tannakiennes},
 BOOKTITLE = {The {G}rothendieck {F}estschrift, {V}ol. {II}},
    SERIES = {Progr. Math.},
    VOLUME = {87},
     PAGES = {111--195},
 PUBLISHER = {Birkh\"{a}user Boston, Boston, MA},
      YEAR = {1990},
   MRCLASS = {14A99 (12H05 18A99)},
  MRNUMBER = {1106898},
MRREVIEWER = {James Milne},
}

@article{DelGon05,
	Author = {Deligne, P. and Goncharov, A. B.},
	Coden = {ASENAH},
	Doi = {10.1016/j.ansens.2004.11.001},
	Fjournal = {Annales Scientifiques de l'\'Ecole Normale Sup\'erieure. Quatri\`eme S\'erie},
	Issn = {0012-9593},
	Journal = {Ann. Sci. \'Ecole Norm. Sup. (4)},
	Mrclass = {11G55 (14F42 14G10 19F27)},
	Mrnumber = {2136480 (2006b:11066)},
	Mrreviewer = {Tam{\'a}s Szamuely},
	Number = {1},
	Pages = {1--56},
	Title = {Groupes fondamentaux motiviques de {T}ate mixte},
	Url = {http://dx.doi.org/10.1016/j.ansens.2004.11.001},
	Volume = {38},
	Year = {2005}}

@misc{DelMil_TC22,
  author = {Deligne, Pierre and Milne, JS},
  title = {{Tannakian Categories}},
  howpublished = "\url{https://www.jmilne.org/math/xnotes/tc2022.pdf}",
  year = {2022}, 
  note = "[Online; accessed 5-Aug-2025]"
}

@article {DograUnlikely24,
    AUTHOR = {Dogra, Netan},
     TITLE = {Unlikely intersections and the {C}habauty-{K}im method over
              number fields},
   JOURNAL = {Math. Ann.},
  FJOURNAL = {Mathematische Annalen},
    VOLUME = {389},
      YEAR = {2024},
    NUMBER = {1},
     PAGES = {1--62},
      ISSN = {0025-5831,1432-1807},
   MRCLASS = {11G30 (11G10 14G05)},
  MRNUMBER = {4735940},
MRREVIEWER = {Tristan\ Phillips},
       DOI = {10.1007/s00208-023-02638-2},
       URL = {https://doi.org/10.1007/s00208-023-02638-2},
}

@article {ErtlNiziol19,
    AUTHOR = {Ertl, Veronika and Nizio\l , Wies\l awa},
     TITLE = {Syntomic cohomology and {$p$}-adic motivic cohomology},
   JOURNAL = {Algebr. Geom.},
  FJOURNAL = {Algebraic Geometry},
    VOLUME = {6},
      YEAR = {2019},
    NUMBER = {1},
     PAGES = {100--131},
      ISSN = {2313-1691},
   MRCLASS = {14F42 (11G25 14F20)},
  MRNUMBER = {3904801},
MRREVIEWER = {Matthias Wendt},
       DOI = {10.1109/jas.2018.7511189},
       URL = {https://doi.org/10.1109/jas.2018.7511189},
}

@misc{eskandari2024tannakianfundamentalgroupsblended,
      title={Tannakian fundamental groups of blended extensions}, 
      author={Payman Eskandari},
      year={2024},
      eprint={2407.01379},
      archivePrefix={arXiv},
      primaryClass={math.AG},
      url={https://arxiv.org/abs/2407.01379},
    note = {arXiv \href{https://arxiv.org/abs/2407.01379}{2407.01379}},
}

@incollection {FPR91,
    AUTHOR = {Fontaine, Jean-Marc and Perrin-Riou, Bernadette},
     TITLE = {Autour des conjectures de {B}loch et {K}ato: cohomologie
              galoisienne et valeurs de fonctions {$L$}},
 BOOKTITLE = {Motives ({S}eattle, {WA}, 1991)},
    SERIES = {Proc. Sympos. Pure Math.},
    VOLUME = {55},
     PAGES = {599--706},
 PUBLISHER = {Amer. Math. Soc., Providence, RI},
      YEAR = {1994},
   MRCLASS = {11F85 (11F33 11F67 11G09 11G40 14G10 19F27)},
  MRNUMBER = {1265546},
MRREVIEWER = {Alexey A. Panchishkin},
}

@incollection {Fontaine94,
    AUTHOR = {Fontaine, Jean-Marc},
     TITLE = {Repr\'{e}sentations {$l$}-adiques potentiellement semi-stables},
      NOTE = {P\'{e}riodes $p$-adiques (Bures-sur-Yvette, 1988)},
   JOURNAL = {Ast\'{e}risque},
  FJOURNAL = {Ast\'{e}risque},
    NUMBER = {223},
      YEAR = {1994},
     PAGES = {321--347},
      ISSN = {0303-1179},
   MRCLASS = {14F20 (11G25 14F30 14G20)},
  MRNUMBER = {1293977},
MRREVIEWER = {Adolfo Quir\'{o}s},
}

@misc{fontaineouyang,
title = {Galois Representations},
author = {Jean-Marc Fontaine and Yi Ouyang},
year = {2022},
url = {http://staff.ustc.edu.cn/~yiouyang/galoisrep.pdf},
note = {http://staff.ustc.edu.cn/~yiouyang/galoisrep.pdf}
}

@article{FurushoI,
	Author = {Furusho, Hidekazu},
	Coden = {INVMBH},
	Doi = {10.1007/s00222-003-0320-9},
	Fjournal = {Inventiones Mathematicae},
	Issn = {0020-9910},
	Journal = {Invent. Math.},
	Mrclass = {11G55 (11M41 11S80)},
	Mrnumber = {2031428 (2005b:11095)},
	Mrreviewer = {Jan Nekov{\'a}{\v{r}}},
	Number = {2},
	Pages = {253--286},
	Title = {{$p$}-adic multiple zeta values. {I}. {$p$}-adic multiple polylogarithms and the {$p$}-adic {KZ} equation},
	Url = {http://dx.doi.org/10.1007/s00222-003-0320-9},
	Volume = {155},
	Year = {2004}}

@article{FurushoJafari,
	author = {Furusho, Hidekazu and Jafari, Amir},
	doi = {10.1112/S0010437X0600265X},
	fjournal = {Compositio Mathematica},
	issn = {0010-437X},
	journal = {Compos. Math.},
	mrclass = {11G55 (11M41 11S80)},
	mrnumber = {2360311},
	mrreviewer = {Jan Nekov\'{a}\v{r}},
	number = {5},
	pages = {1089--1107},
	title = {Regularization and generalized double shuffle relations for {$p$}-adic multiple zeta values},
	url = {https://doi.org/10.1112/S0010437X0600265X},
	volume = {143},
	year = {2007}}

@article{FurushoII,
	Author = {Furusho, Hidekazu},
	Coden = {AJMAAN},
	Doi = {10.1353/ajm.2007.0024},
	Fjournal = {American Journal of Mathematics},
	Issn = {0002-9327},
	Journal = {Amer. J. Math.},
	Mrclass = {11G55 (11M41 11S80)},
	Mrnumber = {2343385 (2009b:11111)},
	Mrreviewer = {Pierre A. Lochak},
	Number = {4},
	Pages = {1105--1144},
	Title = {{$p$}-adic multiple zeta values. {II}. {T}annakian interpretations},
	Url = {http://dx.doi.org/10.1353/ajm.2007.0024},
	Volume = {129},
	Year = {2007}}

@article{FurushoPenagonHexagon,
	author = {Furusho, Hidekazu},
	doi = {10.4007/annals.2010.171.545},
	fjournal = {Annals of Mathematics. Second Series},
	issn = {0003-486X},
	journal = {Ann. of Math. (2)},
	mrclass = {16T20 (17B01 20E18)},
	mrnumber = {2630048},
	mrreviewer = {Pierre A. Lochak},
	number = {1},
	pages = {545--556},
	title = {Pentagon and hexagon equations},
	url = {https://doi.org/10.4007/annals.2010.171.545},
	volume = {171},
	year = {2010}}

@online{FresanGil,
	author = {Jos\'e Ignacio Burgos Gil and Javier Fres\'an},
	title = {Multiple zeta values: from numbers to motives},
  howpublished = "\url{http://javier.fresan.perso.math.cnrs.fr/mzv.pdf}",
	url = {http://javier.fresan.perso.math.cnrs.fr/mzv.pdf},
note = {To appear in Clay Mathematics Proceedings}
}

@incollection {GoerssHopkins00,
    AUTHOR = {Goerss, Paul G. and Hopkins, Michael J.},
     TITLE = {Andr\'{e}-{Q}uillen (co)-homology for simplicial algebras over
              simplicial operads},
 BOOKTITLE = {Une d\'{e}gustation topologique [{T}opological morsels]: homotopy
              theory in the {S}wiss {A}lps ({A}rolla, 1999)},
    SERIES = {Contemp. Math.},
    VOLUME = {265},
     PAGES = {41--85},
 PUBLISHER = {Amer. Math. Soc., Providence, RI},
      YEAR = {2000},
   MRCLASS = {18D50 (13D03 55P43 55U35)},
  MRNUMBER = {1803952},
MRREVIEWER = {Haynes R. Miller},
       DOI = {10.1090/conm/265/04243},
       URL = {https://doi.org/10.1090/conm/265/04243},
}

@incollection {GoerssPresheaves04,
    AUTHOR = {Goerss, Paul G.},
     TITLE = {({P}re-)sheaves of ring spectra over the moduli stack of
              formal group laws},
 BOOKTITLE = {Axiomatic, enriched and motivic homotopy theory},
    SERIES = {NATO Sci. Ser. II Math. Phys. Chem.},
    VOLUME = {131},
     PAGES = {101--131},
 PUBLISHER = {Kluwer Acad. Publ., Dordrecht},
      YEAR = {2004},
   MRCLASS = {55N22 (55P42)},
  MRNUMBER = {2061853},
MRREVIEWER = {Andrey Yu. Lazarev},
       DOI = {10.1007/978-94-007-0948-5\_4},
       URL = {https://doi.org/10.1007/978-94-007-0948-5_4},
}

@misc{goerss2008quasicoherentsheavesmodulistack,
      title={Quasi-coherent sheaves on the moduli stack of formal groups}, 
      author={Paul G. Goerss},
      year={2008},
      eprint={0802.0996},
      archivePrefix={arXiv},
      primaryClass={math.AT},
      url={https://arxiv.org/abs/0802.0996},
    note = {arXiv \href{https://arxiv.org/abs/0802.0996}{0802.0996}},
}

@incollection{GonMot,
	address = {Providence, RI},
	author = {Goncharov, Alexander B.},
	booktitle = {Motives ({S}eattle, {WA}, 1991)},
	mrclass = {19F15 (11G09 11R42 11R70 19E20)},
	mrnumber = {1265551 (94m:19003)},
	mrreviewer = {Richard M. Hain},
	pages = {43--96},
	publisher = {Amer. Math. Soc.},
	series = {Proc. Sympos. Pure Math.},
	title = {Polylogarithms and motivic {G}alois groups},
	volume = {55},
	year = {1994}}

@inproceedings{GonICM,
	address = {Basel},
	author = {Goncharov, Alexander B.},
	booktitle = {Proceedings of the {I}nternational {C}ongress of {M}athematicians, {V}ol.\ 1, 2 ({Z}\"urich, 1994)},
	mrclass = {19F27 (11R42 11R70 19E20)},
	mrnumber = {1403938 (97h:19010)},
	mrreviewer = {J. Browkin},
	pages = {374--387},
	publisher = {Birkh\"auser},
	title = {Polylogarithms in arithmetic and geometry},
	year = {1995}}

@incollection {GonMEM,
    AUTHOR = {Goncharov, Alexander},
     TITLE = {Mixed elliptic motives},
 BOOKTITLE = {Galois representations in arithmetic algebraic geometry
              ({D}urham, 1996)},
    SERIES = {London Math. Soc. Lecture Note Ser.},
    VOLUME = {254},
     PAGES = {147--221},
 PUBLISHER = {Cambridge Univ. Press, Cambridge},
      YEAR = {1998},
   MRCLASS = {11G40 (11G09 11G55 14C35 14G10 19F27)},
  MRNUMBER = {1696477},
MRREVIEWER = {Jan Nekov\'{a}\v{r}},
       DOI = {10.1017/CBO9780511662010.005},
       URL = {https://doi-org.libproxy.berkeley.edu/10.1017/CBO9780511662010.005},
}

@misc{GonMPMTM01,
      title={Multiple polylogarithms and mixed {T}ate motives}, 
      author={A. B. Goncharov},
      year={2001},
      eprint={math/0103059},
      archivePrefix={arXiv},
      primaryClass={math.AG},
      url={https://arxiv.org/abs/math/0103059},
    note = {arXiv \href{https://arxiv.org/abs/math/0103059}{0103059}},
}

@article{GonGal,
	author = {Goncharov, Alexander B.},
	coden = {DUMJAO},
	doi = {10.1215/S0012-7094-04-12822-2},
	fjournal = {Duke Mathematical Journal},
	issn = {0012-7094},
	journal = {Duke Math. J.},
	mrclass = {11G55 (11G09 14C30 16W30 19E15 20F34)},
	mrnumber = {2140264 (2007b:11094)},
	mrreviewer = {Matilde Marcolli},
	number = {2},
	pages = {209--284},
	title = {Galois symmetries of fundamental groupoids and noncommutative geometry},
	url = {http://dx.doi.org/10.1215/S0012-7094-04-12822-2},
	volume = {128},
	year = {2005}}

@article {EGAIV3,
    AUTHOR = {Grothendieck, A.},
     TITLE = {\'{E}l\'{e}ments de g\'{e}om\'{e}trie alg\'{e}brique. {IV}.
              \'{E}tude locale des sch\'{e}mas et des morphismes de
              sch\'{e}mas. {III}},
   JOURNAL = {Inst. Hautes \'{E}tudes Sci. Publ. Math.},
  FJOURNAL = {Institut des Hautes \'{E}tudes Scientifiques. Publications
              Math\'{e}matiques},
    NUMBER = {28},
      YEAR = {1966},
     PAGES = {255},
      ISSN = {0073-8301,1618-1913},
   MRCLASS = {14.55},
  MRNUMBER = {217086},
MRREVIEWER = {J.\ P.\ Murre},
       URL = {http://www.numdam.org/item?id=PMIHES_1966__28__255_0},
}

@article {HainWeighted03,
    AUTHOR = {Hain, Richard and Matsumoto, Makoto},
     TITLE = {Weighted completion of {G}alois groups and {G}alois actions on
              the fundamental group of {$\Bbb P^1-\{0,1,\infty\}$}},
   JOURNAL = {Compositio Math.},
  FJOURNAL = {Compositio Mathematica},
    VOLUME = {139},
      YEAR = {2003},
    NUMBER = {2},
     PAGES = {119--167},
      ISSN = {0010-437X},
   MRCLASS = {14G32 (11G55 11R34 14F35 14F42 14H30)},
  MRNUMBER = {2025807},
MRREVIEWER = {Romyar T. Sharifi},
       DOI = {10.1023/B:COMP.0000005077.42732.93},
       URL = {https://doi.org/10.1023/B:COMP.0000005077.42732.93},
}

@incollection {Hain_HodgeDeRham_Chapter2016,
    AUTHOR = {Hain, Richard},
     TITLE = {The {H}odge--de {R}ham theory of modular groups},
 BOOKTITLE = {Recent advances in {H}odge theory},
    SERIES = {London Math. Soc. Lecture Note Ser.},
    VOLUME = {427},
     PAGES = {422--514},
 PUBLISHER = {Cambridge Univ. Press, Cambridge},
      YEAR = {2016},
   MRCLASS = {14D07 (14G35 14H52 20F40 32S35 58A14)},
  MRNUMBER = {3409885},
MRREVIEWER = {Javier A. Fern\'{a}ndez},
}

@article {HainGoldman20,
    AUTHOR = {Hain, Richard},
     TITLE = {Hodge theory of the {G}oldman bracket},
   JOURNAL = {Geom. Topol.},
  FJOURNAL = {Geometry \& Topology},
    VOLUME = {24},
      YEAR = {2020},
    NUMBER = {4},
     PAGES = {1841--1906},
      ISSN = {1465-3060},
   MRCLASS = {58A12 (14C30 17B62 57K20)},
  MRNUMBER = {4173923},
MRREVIEWER = {Yuli B. Rudyak},
       DOI = {10.2140/gt.2020.24.1841},
       URL = {https://doi.org/10.2140/gt.2020.24.1841},
}

@article {HainMatsumoto,
    AUTHOR = {Hain, Richard and Matsumoto, Makoto},
     TITLE = {Universal mixed elliptic motives},
   JOURNAL = {J. Inst. Math. Jussieu},
  FJOURNAL = {Journal of the Institute of Mathematics of Jussieu. JIMJ.
              Journal de l'Institut de Math\'{e}matiques de Jussieu},
    VOLUME = {19},
      YEAR = {2020},
    NUMBER = {3},
     PAGES = {663--766},
      ISSN = {1474-7480},
   MRCLASS = {14F42 (11G55 14H52 19E20)},
  MRNUMBER = {4094705},
       DOI = {10.1017/s1474748018000130},
       URL = {https://doi.org/10.1017/s1474748018000130},
}

@incollection {HoveyHtpyComods04,
    AUTHOR = {Hovey, Mark},
     TITLE = {Homotopy theory of comodules over a {H}opf algebroid},
 BOOKTITLE = {Homotopy theory: relations with algebraic geometry, group
              cohomology, and algebraic {$K$}-theory},
    SERIES = {Contemp. Math.},
    VOLUME = {346},
     PAGES = {261--304},
 PUBLISHER = {Amer. Math. Soc., Providence, RI},
      YEAR = {2004},
   MRCLASS = {18G35 (18E30 55N20 55U35)},
  MRNUMBER = {2066503},
MRREVIEWER = {Paul G. Goerss},
       DOI = {10.1090/conm/346/06291},
       URL = {https://doi.org/10.1090/conm/346/06291},
}

@article {HuberMMRealizDerivCat,
    AUTHOR = {Huber, Annette},
     TITLE = {Realization of {V}oevodsky's motives},
   JOURNAL = {J. Algebraic Geom.},
  FJOURNAL = {Journal of Algebraic Geometry},
    VOLUME = {9},
      YEAR = {2000},
    NUMBER = {4},
     PAGES = {755--799},
      ISSN = {1056-3911},
   MRCLASS = {14F42 (14F30 18E30)},
  MRNUMBER = {1775312},
MRREVIEWER = {Yuichiro Takeda},
}

@book {HuberMSBook17,
    AUTHOR = {Huber, Annette and M\"{u}ller-Stach, Stefan},
     TITLE = {Periods and {N}ori motives},
    SERIES = {Ergebnisse der Mathematik und ihrer Grenzgebiete. 3. Folge. A
              Series of Modern Surveys in Mathematics [Results in
              Mathematics and Related Areas. 3rd Series. A Series of Modern
              Surveys in Mathematics]},
    VOLUME = {65},
      NOTE = {With contributions by Benjamin Friedrich and Jonas von
              Wangenheim},
 PUBLISHER = {Springer, Cham},
      YEAR = {2017},
     PAGES = {xxiii+372},
      ISBN = {978-3-319-50925-9; 978-3-319-50926-6},
   MRCLASS = {14F42 (11G05 14C15 14C30 19E15 32G20)},
  MRNUMBER = {3618276},
MRREVIEWER = {Alberto Collino},
}

@article {IwanariMotivicRational20,
    AUTHOR = {Iwanari, Isamu},
     TITLE = {Motivic rational homotopy type},
   JOURNAL = {High. Struct.},
  FJOURNAL = {Higher Structures},
    VOLUME = {4},
      YEAR = {2020},
    NUMBER = {2},
     PAGES = {57--133},
      ISSN = {2209-0606},
   MRCLASS = {14F42 (14F08 19E15 55P62)},
  MRNUMBER = {4133164},
MRREVIEWER = {David\ A.\ Blanc},
}

@book {JannsenBook90,
    AUTHOR = {Jannsen, Uwe},
     TITLE = {Mixed motives and algebraic {$K$}-theory},
    SERIES = {Lecture Notes in Mathematics},
    VOLUME = {1400},
      NOTE = {With appendices by S. Bloch and C. Schoen},
 PUBLISHER = {Springer-Verlag, Berlin},
      YEAR = {1990},
     PAGES = {xiv+246},
      ISBN = {3-540-52260-3},
   MRCLASS = {14C35 (11G09 11G35 14G25 19E20)},
  MRNUMBER = {1043451},
MRREVIEWER = {James Milne},
       DOI = {10.1007/BFb0085080},
       URL = {https://doi.org/10.1007/BFb0085080},
}

@article {JarossayDynamical,
    AUTHOR = {Jarossay, David},
     TITLE = {Pro-unipotent harmonic actions and dynamical properties of
              {$p$}-adic cyclotomic multiple zeta values},
   JOURNAL = {Algebra Number Theory},
  FJOURNAL = {Algebra \& Number Theory},
    VOLUME = {14},
      YEAR = {2020},
    NUMBER = {7},
     PAGES = {1711--1746},
      ISSN = {1937-0652,1944-7833},
   MRCLASS = {11M32 (11G55 11S80 14G32)},
  MRNUMBER = {4150248},
MRREVIEWER = {Nils\ Matthes},
       DOI = {10.2140/ant.2020.14.1711},
       URL = {https://doi.org/10.2140/ant.2020.14.1711},
}

@misc{JarLilSaeWeiZeh24,
      title={Polylogarithmic motivic {C}habauty-{K}im for $\mathbb{P}^1 \setminus \{ 0,1,\infty \}$: the geometric step via resultants}, 
      author={David Jarossay and David T. -B. G. Lilienfeldt and Francesco Maria Saettone and Ariel Weiss and Sa'ar Zehavi},
      year={2024},
      eprint={2408.07400},
      archivePrefix={arXiv},
      primaryClass={math.NT},
      url={https://arxiv.org/abs/2408.07400}, 
      note = {arXiv \href{https://arxiv.org/abs/2408.07400}{2408.07400}},
}

@article{kim05,
	Author = {Kim, Minhyong},
	Coden = {INVMBH},
	Doi = {10.1007/s00222-004-0433-9},
	Fjournal = {Inventiones Mathematicae},
	Issn = {0020-9910},
	Journal = {Invent. Math.},
	Mrclass = {11G30 (11G55 14F30 14F42)},
	Mrnumber = {2181717 (2006k:11119)},
	Mrreviewer = {Tam{\'a}s Szamuely},
	Number = {3},
	Pages = {629--656},
	Title = {The motivic fundamental group of {$\mathbb P^1\setminus \{0,1,\infty\}$} and the theorem of {S}iegel},
	Url = {http://dx.doi.org/10.1007/s00222-004-0433-9},
	Volume = {161},
	Year = {2005}}

@article{kim09,
	Author = {Kim, Minhyong},
	Coden = {KRMPBV},
	Doi = {10.2977/prims/1234361156},
	Fjournal = {Kyoto University. Research Institute for Mathematical Sciences. Publications},
	Issn = {0034-5318},
	Journal = {Publ. Res. Inst. Math. Sci.},
	Mrclass = {14G05 (11G30 11G40 14F35)},
	Mrnumber = {2512779 (2010k:14029)},
	Mrreviewer = {Ramdorai Sujatha},
	Number = {1},
	Pages = {89--133},
	Title = {The unipotent {A}lbanese map and {S}elmer varieties for curves},
	Url = {http://dx.doi.org/10.2977/prims/1234361156},
	Volume = {45},
	Year = {2009}}

@incollection {KimEffective,
    AUTHOR = {Kim, Minhyong},
     TITLE = {Remark on fundamental groups and effective {D}iophantine
              methods for hyperbolic curves},
 BOOKTITLE = {Number theory, analysis and geometry},
     PAGES = {355--368},
 PUBLISHER = {Springer, New York},
      YEAR = {2012},
   MRCLASS = {11G30 (14H25)},
  MRNUMBER = {2867924},
MRREVIEWER = {Jan Nekov\'{a}\v{r}},
       DOI = {10.1007/978-1-4614-1260-1\_16},
       URL = {https://doi.org/10.1007/978-1-4614-1260-1_16},
}

@incollection{LevineBlochRevisited,
     author = {Levine, Marc},
     title = {Bloch's higher {C}how groups revisited},
     booktitle = {$K$-theory - Strasbourg, 1992},
     series = {Ast\'erisque},
     publisher = {Soci\'et\'e math\'ematique de France},
     number = {226},
     year = {1994},
     pages = {235-320},
     zbl = {0817.19004},
     language = {en},
     url = {http://www.numdam.org/item/AST_1994__226__235_0}
}

@incollection {LoefSkiZer20,
    AUTHOR = {Loeffler, David and Skinner, Christopher and Zerbes, Sarah
              Livia},
     TITLE = {Syntomic regulators of {A}sai-{F}lach classes},
 BOOKTITLE = {Development of {I}wasawa theory---the centennial of {K}.
              {I}wasawa's birth},
    SERIES = {Adv. Stud. Pure Math.},
    VOLUME = {86},
     PAGES = {595--638},
 PUBLISHER = {Math. Soc. Japan, Tokyo},
      YEAR = {[2020] \copyright 2020},
   MRCLASS = {11F41 (11F67 11F80 19F27)},
  MRNUMBER = {4385092},
MRREVIEWER = {Andrew James Graham},
}

@misc{loeffler2021padic,
      title={On $p$-adic regulators for $\m{GSp}(4) \times \m{GL}(2)$ and $\m{GSp}(4) \times \m{GL}(2) \times \m{GL}(2)$}, 
      author={David Loeffler and Sarah Livia Zerbes},
      year={2021},
      eprint={2011.15098},
    note = {arXiv \href{https://arxiv.org/abs/2011.15098}{2011.15098}},
      archivePrefix={arXiv},
      primaryClass={math.NT}
}

@book {LNMC2006,
    AUTHOR = {Mazza, Carlo and Voevodsky, Vladimir and Weibel, Charles},
     TITLE = {Lecture notes on motivic cohomology},
    SERIES = {Clay Mathematics Monographs},
    VOLUME = {2},
 PUBLISHER = {American Mathematical Society, Providence, RI; Clay
              Mathematics Institute, Cambridge, MA},
      YEAR = {2006},
     PAGES = {xiv+216},
      ISBN = {978-0-8218-3847-1; 0-8218-3847-4},
   MRCLASS = {14F42 (19E15)},
  MRNUMBER = {2242284},
MRREVIEWER = {Thomas Geisser},
}

@unpublished {MoonenIntroMT,
AUTHOR = {Moonen, Ben},
Title = {An Introduction to {M}umford-{T}ate {G}roups},
Date = {2004},
Note = {https://www.math.ru.nl/~bmoonen/Lecturenotes/MTGps.pdf},
}

@article {NiziolSmooth01,
    AUTHOR = {Nizio\l , Wies\l awa},
     TITLE = {Cohomology of crystalline smooth sheaves},
   JOURNAL = {Compositio Math.},
  FJOURNAL = {Compositio Mathematica},
    VOLUME = {129},
      YEAR = {2001},
    NUMBER = {2},
     PAGES = {123--147},
      ISSN = {0010-437X},
   MRCLASS = {14F30},
  MRNUMBER = {1863299},
MRREVIEWER = {Jean-Yves \'{E}tesse},
       DOI = {10.1023/A:1014577828651},
       URL = {https://doi.org/10.1023/A:1014577828651},
}

@article {NizUnique20,
    AUTHOR = {Nizio\l , Wies\l awa},
     TITLE = {On uniqueness of {$p$}-adic period morphisms, {II}},
   JOURNAL = {Compos. Math.},
  FJOURNAL = {Compositio Mathematica},
    VOLUME = {156},
      YEAR = {2020},
    NUMBER = {9},
     PAGES = {1915--1964},
      ISSN = {0010-437X},
   MRCLASS = {14F30 (11F80 14F40 14F42)},
  MRNUMBER = {4172627},
MRREVIEWER = {Lance Edward Miller},
       DOI = {10.1112/s0010437x20007344},
       URL = {https://doi.org/10.1112/s0010437x20007344},
}

@article {NiziolImage97,
    AUTHOR = {Niziol, Wieslawa},
     TITLE = {On the image of {$p$}-adic regulators},
   JOURNAL = {Invent. Math.},
  FJOURNAL = {Inventiones Mathematicae},
    VOLUME = {127},
      YEAR = {1997},
    NUMBER = {2},
     PAGES = {375--400},
      ISSN = {0020-9910},
   MRCLASS = {14F30 (14C35 14G20 19F27)},
  MRNUMBER = {1427624},
MRREVIEWER = {Jan Nekov\'{a}\v{r}},
       DOI = {10.1007/s002220050125},
       URL = {https://doi.org/10.1007/s002220050125},
}

@article{OlssonTowards,
	author = {Olsson, Martin C.},
	coden = {MAMCAU},
	doi = {10.1090/S0065-9266-2010-00625-2},
	fjournal = {Memoirs of the American Mathematical Society},
	isbn = {978-0-8218-5240-8},
	issn = {0065-9266},
	journal = {Mem. Amer. Math. Soc.},
	mrclass = {14G17 (11F85 11G25 11S15 14C30)},
	mrnumber = {2791384},
	number = {990},
	pages = {vi+157},
	title = {Towards non-abelian {$p$}-adic {H}odge theory in the good reduction case},
	url = {http://dx.doi.org/10.1090/S0065-9266-2010-00625-2},
	volume = {210},
	year = {2011}}

@book {OlssonStacksBook16,
    AUTHOR = {Olsson, Martin},
     TITLE = {Algebraic spaces and stacks},
    SERIES = {American Mathematical Society Colloquium Publications},
    VOLUME = {62},
 PUBLISHER = {American Mathematical Society, Providence, RI},
      YEAR = {2016},
     PAGES = {xi+298},
      ISBN = {978-1-4704-2798-6},
   MRCLASS = {14D23 (14D22)},
  MRNUMBER = {3495343},
MRREVIEWER = {Stefan Schr\"{o}er},
       DOI = {10.1090/coll/062},
       URL = {https://doi.org/10.1090/coll/062},
}

@misc{park2024logmotivicnearbycycles,
      title={Log motivic nearby cycles}, 
      author={Doosung Park},
      year={2024},
      eprint={2405.14083},
      archivePrefix={arXiv},
      primaryClass={math.AG},
      url={https://arxiv.org/abs/2405.14083}, 
    note = {arXiv \href{https://arxiv.org/abs/2405.14083}{2405.14083}},
}

@article {Patashnick13,
    AUTHOR = {Patashnick, Owen},
     TITLE = {A candidate for the abelian category of mixed elliptic
              motives},
   JOURNAL = {J. K-Theory},
  FJOURNAL = {Journal of K-Theory. K-Theory and its Applications in Algebra,
              Geometry, Analysis \& Topology},
    VOLUME = {12},
      YEAR = {2013},
    NUMBER = {3},
     PAGES = {569--600},
      ISSN = {1865-2433},
   MRCLASS = {14C25 (11G55 14C15 14F42 55P62)},
  MRNUMBER = {3165187},
MRREVIEWER = {Annette Huber},
       DOI = {10.1017/is013007012jkt237},
       URL = {http://dx.doi.org/10.1017/is013007012jkt237},
}

@book {SaavedraRivanoBook,
    AUTHOR = {Saavedra Rivano, Neantro},
     TITLE = {Cat\'{e}gories {T}annakiennes},
    SERIES = {Lecture Notes in Mathematics, Vol. 265},
 PUBLISHER = {Springer-Verlag, Berlin-New York},
      YEAR = {1972},
     PAGES = {ii+418},
   MRCLASS = {14L15 (18D10 20G05)},
  MRNUMBER = {0338002},
MRREVIEWER = {P. Abellanas},
}

@article {Sakugawa_nabBK17,
    AUTHOR = {Sakugawa, Kenji},
     TITLE = {On a non-abelian generalization of the {B}loch-{K}ato
              exponential map},
   JOURNAL = {Math. J. Okayama Univ.},
  FJOURNAL = {Mathematical Journal of Okayama University},
    VOLUME = {59},
      YEAR = {2017},
    NUMBER = {[2016 on cover]},
     PAGES = {41--70},
      ISSN = {0030-1566},
   MRCLASS = {11R34 (14C30 14F30 19F27)},
  MRNUMBER = {3643426},
}

@article {SchaeppiFormalTannaka13,
    AUTHOR = {Sch\"{a}ppi, Daniel},
     TITLE = {The formal theory of {T}annaka duality},
   JOURNAL = {Ast\'{e}risque},
  FJOURNAL = {Ast\'{e}risque},
    NUMBER = {357},
      YEAR = {2013},
     PAGES = {viii+140},
      ISSN = {0303-1179},
      ISBN = {978-2-85629-773-5},
   MRCLASS = {18D05 (18D10 18D30 18E10)},
  MRNUMBER = {3185459},
MRREVIEWER = {Behrang Noohi},
}

@article {SchaeppiIndAbelian2014,
    AUTHOR = {Sch\"{a}ppi, Daniel},
     TITLE = {Ind-abelian categories and quasi-coherent sheaves},
   JOURNAL = {Math. Proc. Cambridge Philos. Soc.},
  FJOURNAL = {Mathematical Proceedings of the Cambridge Philosophical
              Society},
    VOLUME = {157},
      YEAR = {2014},
    NUMBER = {3},
     PAGES = {391--423},
      ISSN = {0305-0041},
   MRCLASS = {18E10 (14F05)},
  MRNUMBER = {3286515},
MRREVIEWER = {Sandra Mantovani},
       DOI = {10.1017/S0305004114000401},
       URL = {https://doi.org/10.1017/S0305004114000401},
}

@article {SchaeppiConstructingColimits2020,
    AUTHOR = {Sch\"{a}ppi, Daniel},
     TITLE = {Constructing colimits by gluing vector bundles},
   JOURNAL = {Adv. Math.},
  FJOURNAL = {Advances in Mathematics},
    VOLUME = {375},
      YEAR = {2020},
     PAGES = {107394, 85},
      ISSN = {0001-8708},
   MRCLASS = {14A20 (16T05 18E10 18F20 18M15)},
  MRNUMBER = {4170218},
MRREVIEWER = {Jon Eivind Vatne},
       DOI = {10.1016/j.aim.2020.107394},
       URL = {https://doi.org/10.1016/j.aim.2020.107394},
}

@misc{schapira2013derivedcategoryfilteredobjects,
      title={Derived category of filtered objects}, 
      author={Pierre Schapira and Jean-Pierre Schneiders},
      year={2013},
      eprint={1306.1359},
      archivePrefix={arXiv},
      primaryClass={math.AG},
      url={https://arxiv.org/abs/1306.1359}, 
    note = {arXiv \href{https://arxiv.org/abs/1306.1359}{1306.1359}},
}

@article {SchneidersQA99,
    AUTHOR = {Schneiders, Jean-Pierre},
     TITLE = {Quasi-abelian categories and sheaves},
   JOURNAL = {M\'{e}m. Soc. Math. Fr. (N.S.)},
  FJOURNAL = {M\'{e}moires de la Soci\'{e}t\'{e} Math\'{e}matique de France. Nouvelle S\'{e}rie},
      YEAR = {1999},
    NUMBER = {76},
     PAGES = {vi+134},
      ISSN = {0249-633X},
   MRCLASS = {18G35 (18E30 18F20)},
  MRNUMBER = {1779315},
MRREVIEWER = {R. H. Street},
}

@misc{schaeppi2012characterizationcategoriescoherentsheaves,
      title={A characterization of categories of coherent sheaves of certain algebraic stacks}, 
      author={Daniel Schäppi},
      year={2012},
      eprint={1206.2764},
      archivePrefix={arXiv},
      primaryClass={math.AG},
      url={https://arxiv.org/abs/1206.2764},
      note = {arXiv \href{https://arxiv.org/abs/1206.2764}{1206.2764}},
}

@book {SerreCohGal94,
    AUTHOR = {Serre, Jean-Pierre},
     TITLE = {Cohomologie galoisienne},
    SERIES = {Lecture Notes in Mathematics},
    VOLUME = {5},
   EDITION = {Fifth},
 PUBLISHER = {Springer-Verlag, Berlin},
      YEAR = {1994},
     PAGES = {x+181},
      ISBN = {3-540-58002-6},
   MRCLASS = {12G05 (11R34)},
  MRNUMBER = {1324577},
       DOI = {10.1007/BFb0108758},
       URL = {https://doi.org/10.1007/BFb0108758},
}

@book {SerreOpenImageI,
    AUTHOR = {Serre, Jean-Pierre},
     TITLE = {Abelian {$l$}-adic representations and elliptic curves},
    SERIES = {Research Notes in Mathematics},
    VOLUME = {7},
      NOTE = {With the collaboration of Willem Kuyk and John Labute,
              Revised reprint of the 1968 original},
 PUBLISHER = {A K Peters, Ltd., Wellesley, MA},
      YEAR = {1998},
     PAGES = {199},
      ISBN = {1-56881-077-6},
   MRCLASS = {11G05 (11F33 11F80)},
  MRNUMBER = {1484415},
}

@article {SouleOp85,
    AUTHOR = {Soul\'{e}, Christophe},
     TITLE = {Op\'{e}rations en {$K$}-th\'{e}orie alg\'{e}brique},
   JOURNAL = {Canad. J. Math.},
  FJOURNAL = {Canadian Journal of Mathematics. Journal Canadien de
              Math\'{e}matiques},
    VOLUME = {37},
      YEAR = {1985},
    NUMBER = {3},
     PAGES = {488--550},
      ISSN = {0008-414X},
   MRCLASS = {18F25 (14C35 14F15 19E08)},
  MRNUMBER = {787114},
MRREVIEWER = {Ulf Rehmann},
       DOI = {10.4153/CJM-1985-029-x},
       URL = {https://doi.org/10.4153/CJM-1985-029-x},
}

@article {TaylorGaloisReps,
    AUTHOR = {Taylor, Richard},
     TITLE = {Galois representations},
   JOURNAL = {Ann. Fac. Sci. Toulouse Math. (6)},
  FJOURNAL = {Annales de la Facult\'{e} des Sciences de Toulouse. Math\'{e}matiques.
              S\'{e}rie 6},
    VOLUME = {13},
      YEAR = {2004},
    NUMBER = {1},
     PAGES = {73--119},
      ISSN = {0240-2963},
   MRCLASS = {11F80 (11F70 11S37)},
  MRNUMBER = {2060030},
MRREVIEWER = {Thomas A. Weston},
       URL = {http://afst.cedram.org/item?id=AFST_2004_6_13_1_73_0},
}

@incollection {TsujiSurvey02,
    AUTHOR = {Tsuji, Takeshi},
     TITLE = {Semi-stable conjecture of {F}ontaine-{J}annsen: a survey},
      NOTE = {Cohomologies $p$-adiques et applications arithm\'{e}tiques, II},
   JOURNAL = {Ast\'{e}risque},
  FJOURNAL = {Ast\'{e}risque},
    NUMBER = {279},
      YEAR = {2002},
     PAGES = {323--370},
      ISSN = {0303-1179},
   MRCLASS = {14F30 (14F40 14F43)},
  MRNUMBER = {1922833},
MRREVIEWER = {Elmar Grosse-Kl\"{o}nne},
}

@misc{TurDorvault25Tangential,
      title={The motivic fundamental groupoid at tangential basepoints}, 
      author={Sofian Tur-Dorvault},
      year={2025},
      eprint={2510.18151},
      archivePrefix={arXiv},
      primaryClass={math.AG},
      url={https://arxiv.org/abs/2510.18151},
    note = {arXiv \href{https://arxiv.org/abs/2510.18151}{2510.18151}},
}

@incollection {VoevodskyTriCa00,
    AUTHOR = {Voevodsky, Vladimir},
     TITLE = {Triangulated categories of motives over a field},
 BOOKTITLE = {Cycles, transfers, and motivic homology theories},
    SERIES = {Ann. of Math. Stud.},
    VOLUME = {143},
     PAGES = {188--238},
 PUBLISHER = {Princeton Univ. Press, Princeton, NJ},
      YEAR = {2000},
   MRCLASS = {14F42 (14C25)},
  MRNUMBER = {1764202},
}

@misc{wedhorn2024extensionliftinggbundlesstacks,
      title={Extension and lifting of G-bundles on stacks}, 
      author={Torsten Wedhorn},
      year={2024},
      eprint={2311.05151},
      archivePrefix={arXiv},
      primaryClass={math.AG},
      url={https://arxiv.org/abs/2311.05151}, 
note = {arXiv \href{https://arxiv.org/abs/2311.05151}{2311.05151}},
}

@book {WildeshausRoP,
    AUTHOR = {Wildeshaus, J\"{o}rg},
     TITLE = {Realizations of polylogarithms},
    SERIES = {Lecture Notes in Mathematics},
    VOLUME = {1650},
 PUBLISHER = {Springer-Verlag, Berlin},
      YEAR = {1997},
     PAGES = {xii+343},
      ISBN = {3-540-62460-0},
   MRCLASS = {11G18 (11G09 11G40 14D07 14G35 19F27)},
  MRNUMBER = {1482233},
MRREVIEWER = {Jan Nekov\'{a}\v{r}},
       DOI = {10.1007/BFb0093051},
       URL = {https://doi.org/10.1007/BFb0093051},
}

@article {YamashitaBounds10,
    AUTHOR = {Yamashita, Go},
     TITLE = {Bounds for the dimensions of {$p$}-adic multiple {$L$}-value
              spaces},
   JOURNAL = {Doc. Math.},
  FJOURNAL = {Documenta Mathematica},
      YEAR = {2010},
    NUMBER = {Extra vol.: Andrei A. Suslin sixtieth birthday},
     PAGES = {687--723},
      ISSN = {1431-0635},
   MRCLASS = {11G55 (11M32 11R42 14F30 14F42)},
  MRNUMBER = {2804269},
MRREVIEWER = {Jan Nekov\'a\v r},
}

@article {ZieglerGrFil15,
    AUTHOR = {Ziegler, Paul},
     TITLE = {Graded and filtered fiber functors on {T}annakian categories},
   JOURNAL = {J. Inst. Math. Jussieu},
  FJOURNAL = {Journal of the Institute of Mathematics of Jussieu. JIMJ.
              Journal de l'Institut de Math\'{e}matiques de Jussieu},
    VOLUME = {14},
      YEAR = {2015},
    NUMBER = {1},
     PAGES = {87--130},
      ISSN = {1474-7480},
   MRCLASS = {18D10},
  MRNUMBER = {3284480},
MRREVIEWER = {Sorin D\u{a}sc\u{a}lescu},
       DOI = {10.1017/S1474748013000376},
       URL = {https://doi.org/10.1017/S1474748013000376},
}

@article{ZieglerFiltGeneralBase,
  author       = {Ziegler, Paul},
  title        = {Filtered Fiber Functors Over a General Base},
  journal      = {Transformation Groups},
issn = {1083-4362},
  year         = {2024},
  note         = {Published online August 8 2024},
  doi          = {10.1007/s00031-024-09875-y},  
  url          = {https://link.springer.com/article/10.1007/s00031-024-09875-y},
}

@book {HandbookKTheory,
     TITLE = {Handbook of {$K$}-theory. {V}ol. 1, 2},
    EDITOR = {Friedlander, Eric M. and Grayson, Daniel R.},
 PUBLISHER = {Springer-Verlag, Berlin},
      YEAR = {2005},
     PAGES = {Vol. 1: xiv+535 pp.; Vol. 2: pp. i--x and 537--1163},
      ISBN = {978-3-540-23019-9; 3-540-23019-X},
   MRCLASS = {19-06},
  MRNUMBER = {2182598},
       DOI = {10.1007/3-540-27855-9},
       URL = {https://doi-org.libproxy.berkeley.edu/10.1007/3-540-27855-9},
}

@book {SGA3I,
     TITLE = {Sch\'{e}mas en groupes ({SGA} 3). {T}ome {I}. {P}ropri\'{e}t\'{e}s
              g\'{e}n\'{e}rales des sch\'{e}mas en groupes},
    SERIES = {Documents Math\'{e}matiques (Paris) [Mathematical Documents
              (Paris)]},
    VOLUME = {7},
    EDITOR = {Gille, Philippe and Polo, Patrick},
      NOTE = {S\'{e}minaire de G\'{e}om\'{e}trie Alg\'{e}brique du Bois Marie 1962--64.
              [Algebraic Geometry Seminar of Bois Marie 1962--64],
              A seminar directed by M. Demazure and A. Grothendieck with the
              collaboration of M. Artin, J.-E. Bertin, P. Gabriel, M.
              Raynaud and J-P. Serre,
              Revised and annotated edition of the 1970 French original},
 PUBLISHER = {Soci\'{e}t\'{e} Math\'{e}matique de France, Paris},
      YEAR = {2011},
     PAGES = {xxviii+610},
      ISBN = {978-2-85629-323-2},
   MRCLASS = {14L15},
  MRNUMBER = {2867621},
}

@Inbook{Beilinson87,
author="Beilinson, A. A.",
editor="Manin, Yuri I.",
title="Height pairing between algebraic cycles",
bookTitle="K-Theory, Arithmetic and Geometry: Seminar, Moscow University, 1984--1986",
year="1987",
publisher="Springer Berlin Heidelberg",
address="Berlin, Heidelberg",
pages="1--26",
isbn="978-3-540-48016-7",
doi="10.1007/BFb0078364",
url="https://doi.org/10.1007/BFb0078364"
}
%%%%%%%%%%%%%%%%%%
\bibliographystyle{alphanum}%%%
%%%%%%%%%%%%%%%%%
%%%%%%%%%%%%%%%%%

\vfill

\Small

\noindent
\textsc{David Corwin} 

\noindent 
\textsc{Department of mathematics}

\noindent
\textsc{Ben-Gurion University of the Negev}

\noindent
\textsc{Be'er Sheva, Israel}

\noindent
\textsc{Email address:} \url{corwin_david@yahoo.com }

\vspace{.5cm}

\noindent
\textsc{Ishai Dan-Cohen} 

\noindent 
\textsc{Department of mathematics}

\noindent
\textsc{Ben-Gurion University of the Negev}

\noindent
\textsc{Be'er Sheva, Israel}

\noindent
\textsc{Email address:} \url{ishaida@bgu.ac.il}

\end{document}